%% file: 1DynkinExtensions-1_V1.tex
\documentclass[a4paper,11pt,reqno]{amsart}

\usepackage{graphicx}
\usepackage{amsmath}
\usepackage{amssymb}
\usepackage{amsfonts}
\usepackage{epsfig}
\usepackage{epigraph}
\usepackage[hang]{caption}
\usepackage{hhline}
\usepackage{comment}
\usepackage{color}

\usepackage{mathrsfs}

\usepackage{graphics}
\usepackage{float}

\usepackage{mathtools}
\usepackage{booktabs}

\graphicspath{{Linkages/}
              {Conj_Classes/}
              {Loctets/}
              {A_classes/}}

\def\colorcom#1#2#3{
 \def#1##1 {
  \begin{color}{#3}
      \bf[\,#2: ##1\,]
  \end{color}}
} \colorcom\RS{RS}{red}

 \DeclarePairedDelimiter\norm{\lVert}{\rVert}%

\makeatletter

\newtheorem{theorem}{Theorem}[section]

\@addtoreset{equation}{section}

\def\rank{{\rm rank}}

\def\botG{\stackrel{}{\Gamma}}

\newtheorem{corollary}[theorem]{Corollary}
\newtheorem{remark}[theorem]{Remark}
\newtheorem{proposition}[theorem]{Proposition}
\newtheorem{lemma}[theorem]{Lemma}

\def\PerfProof{{\it Proof.\ }}

\makeindex

\begin{document}

\title{Classification of linkage systems}
         \author{Rafael Stekolshchik}
         %%\thanks{{\large email: rs2@biu.013.net.il}}

\date{}

%% \mathfrak changed to \mathsf

\begin{abstract}
   In $1972$, R. Carter introduced admissible diagrams to classify conjugacy classes in a finite Weyl group $W$.
   For any two non-orthogonal roots $\alpha$ and $\beta$ corresponding to vertices of admissible diagram,
   we draw the {\it dotted} (resp. {\it solid}) edge $\{\alpha, \beta\}$ if $(\alpha, \beta) > 0$ (resp. $(\alpha, \beta) < 0$).
   The diagram with properties of admissible diagram and possibly containing dotted edges are said to be {\it Carter diagrams}.
   For any Carter diagram ${\Gamma}$, we introduce the {\it partial Cartan matrix} $B_{\Gamma}$,
   which is analogous to the Cartan matrix associated with a Dynkin diagram.
   A {\it linkage diagram} is obtained from $\Gamma$ by adding an extra root $\gamma$, together with its bonds,
   so that the resulting subset of roots is linearly independent.
   With every linkage diagram we associate the {\it linkage label vector} $\gamma^{\nabla}$,
   similar to the \lq\lq{numerical labels}\rq\rq~ introduced by Dynkin
   for the study of irreducible linear representations of the semisimple Lie algebras.
   The linkage diagrams connected under the action of {\it dual partial Weyl group} $W^{\vee}_{S}$ (associated with $B_{\Gamma}$)
   constitute the set $\mathscr{L}(\Gamma)$ dubbed the {\it linkage system}.
   For any simply-laced Carter diagram $\Gamma$, the linkage system $\mathscr{L}(\Gamma)$ is explicitly constructed.
   To obtain linkage diagrams $\theta^{\nabla} \in \mathscr{L}(\Gamma)$, we use an easily verifiable criterion:
       $\mathscr{B}^{\vee}_{\Gamma}(\theta^{\nabla}) < 2$,
   where $\mathscr{B}^{\vee}_{\Gamma}$ is the quadratic form associated with $B^{-1}_{\Gamma}$.

   A Dynkin diagram $\Gamma'$ such that $\rank(\Gamma') = \rank(\Gamma) + 1$ and
   any $\Gamma$-associated root subset $S$ lies in the root system $\varPhi(\Gamma')$,
   will be called the {\it Dynkin extension of the Carter diagram $\Gamma$}
   and denoted by $\Gamma <_D \Gamma'$.
   The linkage system $\mathscr{L}(\Gamma)$ is the union of $\Gamma_i$-components $\mathscr{L}_{\Gamma_i}(\Gamma)$
   taken for all Dynkin extensions of $\Gamma <_D \Gamma_i$.
   The subset $\varPhi(S)$ of roots of $\varPhi(\Gamma')$
   linearly dependent on roots of $S$ is said to be a {\it partial root system}.
   The size of $\mathscr{L}_{\Gamma'}(\Gamma)$ is estimated as follows:
   $\mid \mathscr{L}_{\Gamma'}(\Gamma) \mid  \leq  \mid\varPhi(\Gamma')\mid - \mid\varPhi(S)\mid$.

   Carter diagrams of the same type and the same index are said to be {\it covalent}.
   For any pair $\{\Gamma, \widetilde\Gamma\}$ of covalent Carter diagrams, where $\Gamma$ is the Dynkin diagram,
   we explicitly construct the invertible linear map $M : \mathcal{P} \longrightarrow \mathcal{R}$, where $\mathcal{R}$ (resp. $\mathcal{P}$)
   is the root system (resp. partial root system) corresponding to $\Gamma$ (resp. $\widetilde\Gamma$).
   This allows us to derive some results for the Carter diagrams from analogous results for the Dynkin diagrams.
   In particular, we have $\mid \mathscr{L}(\widetilde\Gamma) \mid = \mid \mathscr{L}(\Gamma) \mid$.

   For the simple-laced Dynkin diagram $\Gamma$, every $A$- $D$- or $E$-component of the linkage system $\mathscr{L}(\Gamma)$
   coincides with a weight system of one of fundamental representations of the simple Lie algebra $\mathfrak{g}$
   associated with $\Gamma$.

   The $8$-cell {\lq\lq}spindle-like{\rq\rq} subsystems in $\mathscr{L}(\Gamma)$,
   called {\it loctets}, play the essential role in describing the linkage systems.
\end{abstract}

\maketitle

\tableofcontents

\newpage
~\\
~\\

\setlength{\epigraphwidth}{100mm}

\input 1Preview.tex
\input 2IntroMain.tex
\input 3transition.tex
\input 4linkageRoot.tex
\input 5loctets.tex
\input 6enumerat.tex

\input 7weights.tex
\input 8D_partialDlak.tex
\input 8E_partialElak.tex
\input 8A_partialAl.tex
\begin{appendix}
%%\vspace{15mm}
\input AB_treesMatrix.tex
\input C_linkageDiagr.tex

\input D_linkageSyst.tex
\end{appendix}
\newpage
\listoffigures
\newpage
\input biblio-2.tex

\printindex

\end{document}

%% file: 1Preview.tex
 \section{\sc\bf Preview of the basic notions}

  \subsection{Admissible and Carter diagrams, solid and dotted edges}

  \subsubsection{The primary root system}
    \label{sec_fulcrum}
    \index{primary root system}
    \index{root system ! - primary}
   We consider a root system $\varPhi$, the finite Weyl group $W$ acting on $\varPhi$,
   the diagram Dynkin $\Gamma_{D}$ associated with $\varPhi$ and the corresponding Cartan matrix {\bf B}.
   The quadratic form $\mathscr{B}$ is associated with the Cartan matrix {\bf B} and $\mathscr{B}$
   determines the inner product $(\cdot,\cdot )$  on the linear space $V$ spanned
   by simple roots of $\varPhi$. Inside the root system $\varPhi$ we will consider different
   root subsets and root subsystems, so we call $\varPhi$ the {\it primary root system}.
   There exist root subsets which can be embedded into different primary root systems. For example,
   $\varPhi(A_7) \subset \varPhi(A_8), \varPhi(D_8), \varPhi(E_8)$.
   In this paper, besides the quadruple of objects $\{\varPhi, W, \Gamma_{D}, {\bf B}\}$ we consider also
   a certain its generalization, the  quadruple $\{\mathcal{P}, W^{\vee}, \Gamma, B_{\Gamma}\}$,
   where $\mathcal{P}$ (resp. $W^{\vee}$, resp. $B_{\Gamma}$) is the {\it partial root system},
   (resp. partial Weyl group, partial Cartan matrix) and $\Gamma$ is the Carter diagram.
   These notions will be explained shortly. The generalized quadruple is our object of interest.
   And all this is studied in the area which is the primary root system $\varPhi$.

   \subsubsection{The admissible diagram}
  \label{sec_adm_diagr}
  \index{root system}

 In $1972$, R.~Carter introduced {\it admissible diagrams} to classify conjugacy classes in a finite Weyl group $W$.
 These diagrams are also used to characterize elements of the Weyl group, see definition in \S\ref{sec_Carter}.
 Each element $w \in W$ can be expressed in the form
 \begin{equation}
   \label{any_roots_0}
    w  = s_{\tau_1} s_{\tau_2} \dots s_{\tau_k}, \text{ where } \tau_i \in \varPhi,
 \end{equation}
 and $s_{\tau_i} \in W$ are reflections
 corresponding to not necessarily simple roots $\tau_i \in \varPhi$.
 %% Here, $\varPhi$ is a certain root system, and $W$ is the Weyl group associated with $\varPhi$.
 Carter proved that $k$ in the decomposition \eqref{any_roots_0} is
 the smallest if and only if the subset of roots $\{\tau_1, \tau_2, \dots, \tau_k\}$ is linearly independent; such a
 decomposition is said to be {\it reduced}.
  We denote by $l_C(w)$ the smallest value $l$ in any expression like
 \eqref{any_roots_0}, we have $l_C(s_{\tau_1} s_{\tau_2} \dots s_{\tau_l}) = l$.
 The set $S = \{\tau_1,\dots,\tau_l\}$ consists of linearly
 independent and \underline{not necessarily simple} roots, see \cite[Lemma $3$]{Ca72}.
 We always have $l_C(w) \leq l(w)$, where $l(w)$ is the smallest value
 $k$ in any expression like \eqref{any_roots_0} such that all roots $\alpha_i$ are \underline{simple}.

 \begin{lemma}{\cite[Lemma 3]{Ca72}}
  \label{lem_lin_indep}
   Let $\alpha_1, \alpha_2, \dots, \alpha_k \in \varPhi$.
   Then $s_{\alpha_1} s_{\alpha_2} \dots s_{\alpha_k}$ is reduced
   if and only if $\alpha_1, \alpha_2, \dots, \alpha_k$ are linearly
   independent. \qed
 \end{lemma}

 We associate with the element $w \in W$ (and with conjugacy class $C(w)$ containing $w$) the diagram
 $\Gamma$ with nodes corresponding to roots in \eqref{any_roots_0}, nodes $\alpha$ and $\beta$  are joined
 by  $n_{{\alpha}{\beta}}\cdot n_{{\beta}{\alpha}}$ bonds, where
\begin{equation}
   \label{n_sr}
       n_{{\alpha}{\beta}} = 2\frac{(\alpha,\beta)}{(\alpha,\alpha)}, \quad
       n_{{\beta}{\alpha}} = 2\frac{(\beta,\alpha)}{(\beta,\beta)}.
 \end{equation}
Number $n_{{\alpha}{\beta}}$ and $n_{{\beta}{\alpha}}$  are integer, and $n_{{\alpha}{\beta}} \cdot n_{{\beta}{\alpha}} \in \{0, 1, 2, 3\}$,
see \cite[Ch.6, \S1, ${n^o}$3]{Bo02}.

\index{admissible diagram}
\index{diagram! - admissible diagram}
 A diagram $\Gamma$ is said to be {\it admissible}, see \cite[p. 7]{Ca72}, if
\begin{equation}
 \label{eq_def_adm}
 \begin{split}
 & (a)  \text{ The nodes of } \Gamma \text{ correspond to a set of linearly independent roots in } \varPhi. \\
 & (b) \text{
  If a subdiagram of } \Gamma \text{ is a cycle, then it contains an even number of nodes. }
 \end{split}
\end{equation}

 Surprisingly,
 the admissible diagrams can contain cycles though the extended Dynkin
 diagram $\widetilde{A}_l$ cannot be a part of any admissible diagram, see Lemma \ref{lem_must_dotted}.
 It turned out that the cycles in admissible diagrams essentially differ from the cycle $\widetilde{A}_l$.
 Namely, in such a cycle, there exist necessarily two pairs of roots: A pair with a
 positive inner product and a pair with a negative inner product. This does not happen for $\widetilde {A}_l$.

 \subsubsection{Solid and dotted edges}
  \label{sec_solid_dotted}
 \index{dotted ! - edge}
 \index{solid ! - edge}
 This observation motivated us to distinguish such pairs of roots:
 Let us draw the {\it dotted} (resp. {\it solid}) edge $\{\alpha, \beta\}$ if
 $(\alpha, \beta) > 0$ (resp. $(\alpha, \beta) < 0$), see Fig. \ref{E6a1_exam_linkages}.
 The diagrams with properties of admissible diagrams and containing
 dotted edges are said to be {\it Carter diagrams}. Up to dotted edges,
 the classification of Carter diagrams coincides with the
 classification of admissible diagrams.

  Recall that $(\alpha, \beta) > 0$
 (resp. $(\alpha, \beta) < 0$) means that the angle between
 roots $\alpha$ and $\beta$ is acute (resp. obtuse).
  Recall that, for the Dynkin diagrams, all angles between simple roots
  are obtuse, all edges are solid, and no special designation is necessary.
  A {\it solid edge} indicates an obtuse angle between roots exactly as for simple roots
  in the case of Dynkin diagrams. A {\it dotted edge} indicates an acute angle between the
  roots considered.

 \subsubsection{The $\Gamma$-associated root subset}
   \label{sec_associated}
   \index{$S$, $\Gamma$-associated root subset}
   \index{$\Gamma$-associated root subset}
   \index{primary root system}
   \index{root system ! - primary}

  For any Carter diagram $\Gamma$ and a primary root system $\varPhi$,
  let $\Gamma'$ be a Dynkin diagram such that $\rank(\Gamma') = \rank(\Gamma) + 1$,
  $\varPhi(\Gamma') \subseteq \varPhi$, where $\varPhi(\Gamma')$ is the root system associated with $\Gamma'$.
  Consider a subset $S \subset \varPhi(\Gamma')$ of linearly independent roots with following properties:
  Vertices of $\Gamma$ are in one-to-one correspondence with roots of $S$;
  solid (resp. dotted) edges of $\Gamma$ are in one-to-one correspondence with pairs
  $\{\alpha, \beta\}$, where $\alpha, \beta \in S$, such that $(\alpha, \beta) < 0$
  (resp. $(\alpha, \beta) > 0$. The root subset $S$ is said to be {\it $\Gamma$-associated}.
  %% This extension is said to be {\it Dynkin extension}, see \S\ref{sec_Dynkin_ext0}.

  \begin{remark}
  {\rm (i) Emphasize that whenever we say that $S$ is the $\Gamma$-associated root subset,
 we mean that elements of $S$ are taken from $\varPhi(\Gamma') \subseteq \varPhi$ , i.e.,
 they are taken from the primary root system $\varPhi$, see \S\ref{sec_fulcrum}.

 (ii) Not necessarily for every Carter diagram $\Gamma$ there exists a set $S$
 with a given properties. For example, there is no $E_6$-associated root subset in the primary root subset $\varPhi = \varPhi(D_n)$
 for any integer $n$, see Lemma \ref{lem_E6_not_in_Dn}.

 }
  \end{remark}

 \subsubsection{The Carter diagrams}
  \label{sec_Carter}

 \index{Carter diagram}
 \index{diagram! - Carter diagram}
 \index{$\Gamma$-associated set of roots}
 \index{$\Gamma$-associated element}
 \index{$S$-associated element}
 \index{$S$ (set of linearly independent roots)}
 \index{$S_{\alpha}$ (subset $\{\alpha_i \mid i = 1,\dots,k\}$)}
 \index{$S_{\beta}$ (subset $\{\beta_j \mid j = 1,\dots,h\}$)}
 \index{Weyl group}
  Any admissible diagram $\Gamma$ is said to be a {\it Carter diagram}
  if any edge connecting a pair of roots $\{\alpha, \beta\}$ with
  inner product $(\alpha, \beta) > 0$ (resp. $(\alpha, \beta) < 0$)
  is drawn as dotted (resp. solid) edge. There is no edge for the inner product $(\alpha, \beta) = 0$. Let
\begin{equation}
   \label{eq_alpha_bet}
     S = \{\alpha_1, \alpha_2, \dots, \alpha_k, \beta_1, \beta_2, \dots, \beta_h \}
 \end{equation}
 be any $\Gamma$-associated set (of not necessarily simple roots),
 where roots of the set $S_{\alpha} := \{\alpha_i \mid i = 1,\dots,k\}$ (resp. $S_{\beta} := \{\beta_j \mid j = 1,\dots,h\}$) are
 mutually orthogonal.
 According to (\ref{eq_def_adm}$(a)$), there exists the set \eqref{eq_alpha_bet} of linearly independent roots.
 Thanks to (\ref{eq_def_adm}$(b)$), such a partition into the sum of two mutually orthogonal sets $S_{\alpha}$
 and $S_{\beta}$ is possible. The partition \eqref{eq_alpha_bet} is said to be {\it bicolored partition}.
 Let $w = w_1 w_2$ be the decomposition of $w$ into the product of two involutions.
 By \cite[Lemma 5]{Ca72} each of $w_1$ and $w_2$ can
 be expressed as a product of reflections corresponding to mutually
 orthogonal roots as follows:
 \begin{equation}
   \label{two_invol}
      w = w_1{w}_2, \quad \text{ where } \quad
      w_1 = s_{\alpha_1} s_{\alpha_2} \dots s_{\alpha_k}, \quad
      w_2 = s_{\beta_1} s_{\beta_2} \dots s_{\beta_h}.
 \end{equation}
 Since $S$ is linearly independent, the decomposition
 \eqref{two_invol} is reduced, see Lemma \ref{lem_lin_indep}, and $k + h = l_C(w)$.
 The decomposition \eqref{two_invol}
 is said to be a {\it bicolored decomposition}.
 The subset of roots corresponding to $w_1$ (resp. $w_2$) is said to be {\it $\alpha$-set} (resp. {\it $\beta$-set}):
 \index{bicolored decomposition}
\begin{equation}
   \label{two_sets}
       \alpha\text{-set} = \{ \alpha_1, \alpha_2, \dots, \alpha_k \}, \quad
       \beta\text{-set} =  \{ \beta_1, \beta_2, \dots, \beta_h \}.
 \end{equation}

 \index{$w = w_1{w}_2$, bicolored decomposition of $w$; $w_1$, $w_2$ - involutions}
 \index{$\alpha$-set, subset of roots corresponding to $w_1$}
 \index{$\beta$-set, subset of roots corresponding to $w_2$}
 \index{$\alpha$-label, coordinate out of $\alpha$-set}
 \index{$\beta$-label, coordinate out of $\beta$-set}
 \index{label ! - $\alpha$-label}
 \index{label ! - $\beta$-label}
 \index{$k$, number of coordinates in $\alpha$-set (= number of $\alpha$\text{-labels})}
 \index{$h$, number of coordinates in $\beta$-set (= number of $\beta$\text{-labels})}
 \index{$l$, number of vertices in the Carter diagram $\Gamma$, $l = k + h$}

 Up to the difference between solid and dotted edges,
 the Carter diagram $\Gamma$ is the admissible diagram.
 Any coordinate of a given vector of linkage labels from an $\alpha$-set (resp. a $\beta$-set)
 we call an {\it $\alpha$-label} (resp. a {\it $\beta$-label}).
 Let $L \subset V$ be the linear subspace spanned by root subsets \eqref{eq_alpha_bet}.
 We denote by $\Pi_w$ the corresponding root basis, whose elements
 are not necessarily simple roots:
 \index{$\Pi_w$, root subset associated with the bicolored decomposition of $w$}
 \begin{equation}
   \label{root_subset_L}
       \Pi_w = \{ \alpha_1, \alpha_2, \dots, \alpha_k, \beta_1, \beta_2, \dots, \beta_h \}.
 \end{equation}

\subsubsection {Three classes of Carter diagrams}
   \label{sec_3_class}
 \index{$W$, Weyl group}
 \index{Weyl group ! - $W$}
 \index{class of diagrams ! - $\mathsf{C4}$}
 \index{class of diagrams ! - $\mathsf{DE4}$}
 \index{${\bf B}$, Cartan matrix associated with a Dynkin diagram}
 \index{Cartan matrix ! - ${\bf B}$}
 \index{label ! - linkage label vector}
 \index{class of diagrams ! - simply-laced connected Carter diagrams}
 \index{$\Gamma$, Carter diagram}
 \index{uniqueness of the associated conjugacy class}
 \index{uniqueness theorem}

  We consider three classes of simply-laced connected Carter diagrams:
  Denote the class of diagrams containing $4$-cycle $D_4(a_1)$ by $\mathsf{C4}$;
  the class of diagrams without cycles and containing $D_4$ (as a subdiagram) by $\mathsf{DE4}$,
  i.e., $\mathsf{DE4}$ is the class consisting of Dynkin diagrams  $E_6$, $E_7$, $E_8$ and $D_l$ for $l \geq 4$;
  the class of Dynkin diagrams $A_l$ by $\mathsf{A}$.
  Conjugate elements in the Weyl group $W$ are associated with the same Carter diagram $\Gamma$.
  For any $\Gamma \in \mathsf{C4} \coprod \mathsf{DE4}$, the {\it uniqueness of the associated conjugacy class} takes place,
  see \cite[Theorem 6.5]{St10}.
  %% There exist Carter diagrams which do not determine a
  %% single conjugacy class in $W$, \cite[Lemma 27]{Ca72}.
  %% For $A_l$, the uniqueness theorem does not hold, see Remark \ref{two_class_Al}.
\begin{remark}
  \label{two_class_Al}
{\rm
    Note that the uniqueness theorem does not hold for $A_l$.
    There are two conjugacy classes of type $A_5$ in $W(E_7)$ and in $W(D_6)$,
    two conjugacy classes of type $A_7$ in $W(E_8)$, see \cite[p. 31]{Ca72} and \cite[Lemma 27]{Ca72}.
    Moreover, for the Carter diagram $A_3$, there exist three conjugacy classes in $W(D_4)$.
   } \qed
\end{remark}

%% \clearpage
\subsection{Linkage diagrams and linkage systems}

 \subsubsection{The partial Cartan matrix $B_{\Gamma}$}
    \label{sec_partial_B}
    \index{$B_{\Gamma}$, partial Cartan matrix}
    \index{Cartan matrix ! - partial Cartan matrix $B_{\Gamma}$}

  Similarly to the Cartan matrix associated with Dynkin diagrams, we determine the Cartan matrix
  for each Carter diagram $\Gamma$ as follows

 \begin{equation}
   \label{canon_dec_2}
   B_{\Gamma} :=
      \left (
        \begin{array}{cccccc}
         (\alpha_1, \alpha_1) & \dots & (\alpha_i, \alpha_k) &
         (\alpha_1, \beta_1)  & \dots &  (\alpha_1, \beta_h) \\
         \dots                & \dots & \dots &
         \dots                & \dots & \dots \\
         (\alpha_k, \alpha_1) & \dots & (\alpha_k, \alpha_k) &
         (\alpha_k, \beta_1)  & \dots &  (\alpha_k, \beta_h) \\
         (\beta_1, \alpha_1)  & \dots & (\beta_1, \alpha_k) &
         (\beta_1, \beta_1)   & \dots & (\beta_1, \beta_h) \\
         \dots                & \dots & \dots &
         \dots                & \dots & \dots \\
         (\beta_h, \alpha_1)  & \dots & (\beta_h, \alpha_k) &
         (\beta_h, \beta_1)   & \dots & (\beta_h, \beta_h) \\
        \end{array}
      \right ),
 \end{equation}
 where $S = \{\alpha_1,\dots,\alpha_k, \beta_1, \dots, \beta_h\}$. Here, subsets of roots
 $\alpha$-set and $\beta$-set, see \eqref{two_sets}, match
 the bicolored partition \eqref{eq_alpha_bet} and bicolored decomposition of a certain element $w \in W$
 corresponding to $\Gamma$. We also write $S = \{ \tau_1,\dots,\tau_{k+h} \}$ if the
 bicolored partition like \eqref{eq_alpha_bet} does not matter.

 \index{quadratic form ! - $\mathscr{B}_{\Gamma}$, associated with the partial Cartan matrix $B_{\Gamma}$}
 \index{$\mathscr{B}_{\Gamma}$, quadratic form associated with the partial Cartan matrix $B_{\Gamma}$}
 \index{$S = \{\tau_1,\dots,\tau_l\}$, $\Gamma$-associated root subset}
 \index{Cartan matrix ! - inverse $B^{-1}_{\Gamma}$ of the partial Cartan matrix}
 \index{$(\cdot, \cdot)$, symmetric bilinear form associated with ${\bf B}$}
 \index{symmetric bilinear form $(\cdot, \cdot)$}
 \index{$(\cdot, \cdot)_{\Gamma}$, symmetric bilinear form associated with $B_{\Gamma}$}
 \index{symmetric bilinear form $(\cdot, \cdot)_{\Gamma}$}
 \index{$S$-associated subspace}
 We call the matrix $B_{\Gamma}$ a {\it partial Cartan matrix} corresponding to
 the Carter diagram $\Gamma$.
 The partial Cartan matrix $B_{\Gamma}$ is well-defined since $(\tau_i, \tau_j)$ in \eqref{canon_dec_2}
 do not depend on the choice of $\Gamma$-associated root subset.
 If $\Gamma$ is a Dynkin diagram,
 the partial Cartan matrix $B_{\Gamma}$ is the Cartan matrix
 associated with $\Gamma$.
 The matrix $B_{\Gamma}$ is positive definite (Proposition \ref{restr_forms_coincide}).
 The symmetric bilinear form associated with
 the partial Cartan matrix $B_{\Gamma}$ is denoted by $(\cdot, \cdot)_{\Gamma}$ and the
 corresponding quadratic form is denoted by $\mathscr{B}_{\Gamma}$. Let $L \subset V$ be the
 subspace spanned by the root subset $S$, the subspace $L$ is said to be {\it S-associated subspace}.
 For $S = \{\tau_1,\dots,\tau_l\}$, we write $L =[\tau_1,\dots,\tau_l]$.

\begin{remark}{\rm
 \label{rem_values}
 In all cases where we consider the partial Cartan matrix $B_{\Gamma}$ for simply-laced Carter diagrams,
 we assume that $(\tau_i, \tau_j)$ takes values in $\{-1, 0,  1 \}$, not in $\{-\frac{1}{2}, 0, \frac{1}{2}\}$.
 In other words, we consider the bilinear form $(\cdot, \cdot)$ obtained by doubling the usual bilinear form
 related to the Cartan matrix. Then linkage labels introduced in \S\ref{sec_linkage_diagr} are integer.
 }
\end{remark}

%%%% \subsubsection{Which came first, the Root System or the Dynkin Diagram?}

\subsubsection{The linkage diagrams}
  \label{sec_linkage_diagr_0}
  \index{Dynkin labels}
  \index{label ! - Dynkin labels}
  \index{label ! - linkage label vector}
  \index{$\gamma$, linkage root}
  \index{linkage root}
  \index{class of diagrams ! - $\Gamma$, Carter diagrams}
  \index{$\Gamma$, Carter diagram}
  \index{connection diagram}

  We consider a class of diagrams, called {\it linkage diagrams}, which
  constitute a subclass of the class of connection diagrams, see \S\ref{sec_connection},
  and generalize the class of Carter diagrams, see \S\ref{sec_Carter}.
  Any linkage diagram is obtained from a Carter diagram $\Gamma$ by adding one extra root $\gamma$,
  with its bonds, so that the roots corresponding to vertices of $\Gamma$ together with $\gamma$ form some
  linearly independent root subset. The extra root $\gamma$ added to the Carter diagram $\Gamma$ is said to be
  a {\it linkage root}, see \S\ref{sec_linkage}.  Any linkage diagram constructed in this way may
  also be a Carter diagram but this is not necessarily so.  The following inclusions hold:

 \begin{equation*}
   \fbox{$\begin{array}{c}
      \text{ Dynkin } \\
      \text{ diagrams of CCl\footnotemark[1] }
   \end{array}$} \quad\subset\quad
   \fbox{$\begin{array}{c}
      \text{ Carter } \\
      \text{ diagrams }
   \end{array}$} \quad\subset\quad
   \fbox{$\begin{array}{c}
      \text{ Linkage } \\
      \text{ diagrams }
   \end{array}$} \quad\subset\quad
   \fbox{$\begin{array}{c}
      \text{ Connection } \\
      \text{ diagrams }
   \end{array}$}
 \end{equation*}

 \footnotetext[1]{The Dynkin diagrams in this article appear in two ways:
 (1) associated with some irreducible root system (customary use);
 (2) representing some conjugacy class (CCl), i.e, the Carter diagram
 which looked like a Dynkin diagram.
 In a few cases Dynkin diagrams represent two (and even three!) conjugacy classes, see Remark \ref{two_class_Al}.}

  With every linkage diagram we associate the vector of {\it linkage labels}.
  The linkage labels are similar to the Dynkin labels,
  see \cite{Sl81}, which are the {\lq\lq}numerical labels{\rq\rq} introduced by Dynkin in \cite{Dy50}
  for the study of irreducible linear representations of the semisimple Lie algebras, \cite{GOV90}, \cite{KOV95}, \cite{Ch84}.

\begin{figure}[h]
\centering
\includegraphics[scale=0.7]{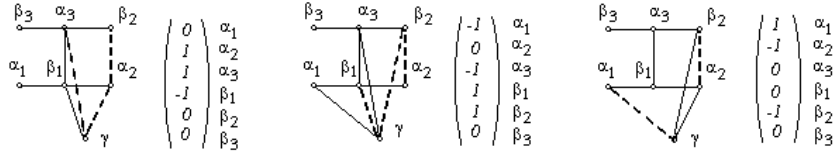}
 \caption{\hspace{3mm}Examples of linkage diagrams and vectors of linkage labels for $E_6(a_1)$.}
%%%%%% The label must come after caption
\label{E6a1_exam_linkages}
\end{figure}

  For any simply-laced Carter diagram $\Gamma$, a linkage label
  takes one of three values $\{-1, 0, 1 \}$. There is one-to-one correspondence between
  linkage diagrams obtained from the given simply-laced Carter diagram $\Gamma$ and vectors of linkage labels
  taking coordinates from the set $\{-1, 0, 1 \}$.
  For this reason, we use terms \underline{linkage labels} and  \underline{linkage diagrams} as
  synonyms. Some linkage diagrams and their linkage labels for the Carter diagram $E_6(a_1)$
  are depicted in Fig. \ref{E6a1_exam_linkages}.
 \index{linkage diagrams ! - examples for $E_6(a_1)$}
  The linkage diagrams are the main characters of this paper. {\bf We give a complete description of linkage diagrams constructed
  for every simply-laced Carter diagram.}
  By abuse of notation, this description is essentially based on the answer to the following question:
 ~\\

  {\it What linkage roots can be added to a given irreducible linearly independent root subset,
  so that the resulting set would also be an irreducible linearly independent root subset?}
 ~\\
 ~\\
 It turns out that the answer to this question is very simple within the framework of the
 quadratic form associated with the partial Cartan matrix $B_{\Gamma}$, see \S\ref{sec_criterion}.

\subsubsection{The linkage systems}
 The group $W^{\vee}_S$, named the {\it dual partial Weyl group}, acts
 on linkage diagrams, see \S\ref{sec_dual_partial}.
 Under the action of $W^{\vee}_S$ the set of
 linkage diagrams (= vectors of linkage labels) constitute the
 diagram called the {\it linkage system} similarly to the weight system
 in the representation theory of semisimple Lie algebras, \cite[p. 30]{Sl81}.
 We denote the linkage system associated with the Carter diagram $\Gamma$ by $\mathscr{L}(\Gamma)$.

\subsection{Dynkin extensions, partial root systems and linkage system components}
   \label{sec_D_ext_1}

\subsubsection{$A$-, $D$-, $E$-types of Dynkin and Carter diagram}
  The Dynkin diagram $A_l$, where $l \ge 1$ (resp. $D_l$, where $l \ge 4$; resp. $E_l$, where $l = 6,7,8$) is said to be the
  {\it Dynkin diagram of $A$-type}
  (resp. {\it $D$-type}, resp. {\it $E$-type}). The Carter diagram $A_l$, where $l \ge 1$ (resp. $D_l$, $D_l(a_k)$,
  where $l \ge 4$, $1 \leq k \leq \big [ \tfrac{l-2}{2} \big ]$; resp. $E_l$, $E_l(a_k)$, where $l = 6,7,8$, $k = 1,2,3,4$)
  is said to be the {\it Carter diagram of $A$-type} (resp. {\it $D$-type}, resp. {\it $E$-type}).

\subsubsection{Partial root system}
   For a $\Gamma$-associated root subset $S$, we denote by $\varPhi(S)$ the subset of roots of $\varPhi$
   \underline{linearly dependent} on roots of $S$. The subset of roots $\varPhi(S)$ is said to be
   a {\it partial root system}.
   There are examples of two non-conjugate $\Gamma$-associated root subsets
   $S_1$ and $S_2$, see \S\ref{sec_example_4cycles}.
   This is a reason why I prefer denote the partial root system by $\varPhi(S)$, not by $\varPhi(\Gamma)$.

\subsubsection{Covalent Carter diagrams}
  \label{sec_coval_list}
 Carter diagrams of the same type and the same index are said to be {\it covalent} Carter diagrams.
 The corresponding root systems (usual or partial) are said to be {\it covalent} root systems.
 For example, the $E$-type covalent diagrams of index $6$ are $\{E_6, E_6(a_1), E_6(a_2)\}$ and
 the $E$-type covalent diagrams of index $7$ are $\{E_7, E_7(a_1), E_7(a_2), E_7(a_3), E_7(a_4))\}$,
 the $D$-type covalent diagrams of index $l$ are
          $\{D_l, D_l(a_1), D_l(a_2), \dots, D_l(a_{\big [ \tfrac{l-2}{2} \big ]}\}$.
~\\

 We are interested in the following covalent pairs of diagrams $\{\widetilde\Gamma, \Gamma\}$, where $\Gamma$ is a Dynkin diagram
 and $\widetilde\Gamma$ is a Carter diagram covalent to $\Gamma$:
   \begin{equation}
     \label{eq_pairs_root_syst_0}
    %% \label{eq_covalent_pairs}
    %% \footnotesize
     \begin{split}
      & \{D_l(a_k), D_l\} \text{ for } l \geq 4 \text{ and } 1 \leq k \leq \big [\tfrac{l-2}{2} \big ], \\
      & \{E_6(a_k), E_6\} \text{ for }  k = 1,2, \\
      & \{E_7(a_k), E_7\} \text{ for }  k = 1,2,3,4.
     \end{split}
   \end{equation}

   Let $\{\widetilde\Gamma, \Gamma\}$ be a covalent pair of Carter diagrams out of list \eqref{eq_pairs_root_syst_0}
   and $\widetilde{S}$ (resp. $S$) be $\widetilde\Gamma$-associated (resp. $\Gamma$-associated) subset,
   $\mathcal{P}$ (resp. $\mathcal{R}$) be the partial root system (resp. root system) spanned by $\widetilde{S}$ (resp. $S$).

   In \S\ref{sec_transition}, we show that there exists the invertible linear map $M : \mathcal{P} \longmapsto \mathcal{R}$,
   Moreover, $\widetilde\tau$ is a root of $\mathcal{P}$ if and only if $\widetilde\tau$ is a root of $\mathcal{R}$
   (Theorem \ref{prop_bijection_root_syst}). For any covalent pair \eqref{eq_pairs_root_syst_0}, we have
   $\mid \mathcal{R} \mid = \mid \mathcal{P} \mid.$ (Corollary \ref{corol_biject}).

  \begin{table}[H]
  %%\footnotesize
  \Small
  \centering
  \renewcommand{\arraystretch}{1.5}
  \begin{tabular} {|c|c|c|c|}
  \hline
      $\Gamma$ &  Number of  &
      \multicolumn{2}{c|}{$E$-components}   \cr
      \cline{3-4}
        &  components   &
      $\mathscr{B}^{\vee}_{\Gamma}(\gamma^{\nabla}$) &  Number of \cr
        &   &  &  linkage diagrams \\
    \hline
       $E_6$, $E_6(a_1)$, $E_6(a_2)$ &   $2$  &  $\frac{4}{3}$  & $2\times27 = 54$  \\
    \hline
       $E_7$, $E_7(a_1)$, $E_7(a_2)$, $E_7(a_3)$, $E_7(a_4)$ &   $1$  &   $\frac{3}{2}$  & $56$ \cr
   \hline    %1
  \end{tabular}
  ~\\ \vspace{2mm}
  \centering
  \begin{tabular} {|c|c|c|c|c|}
  \hline
     $\Gamma$ &  Number of  &
      \multicolumn{3}{c|}{Number of linkage diagrams,}  \cr
       &  components   & \multicolumn{3}{c|}{and $p = \mathscr{B}^{\vee}_{\Gamma}(\gamma^{\nabla})$}    \cr
      \cline{3-5}
         &    & $D$-components   & $E$-components &
          \hspace{1.8mm} In all \hspace{1.8mm} \cr
      \cline{3-4}
          &     & $p = 1$ & $p = \frac{l}{4}$ & \\
    \hline
       %% $D_4$, ~$D_4(a_1)$ &   $3$   & $3 \times 8 = 24$ & -\footnotemark[1]  &  $24$ \\
       $D_4$, ~$D_4(a_1)$ &   $3$   & $3 \times 8 = 24$ & - &  $24$ \\
    \hline
       $D_5$, ~$D_5(a_1)$ &   $3$   &  $10$ & $2\times16 = 32$   & $42$ \\
    \hline
       $D_6$, $D_6(a_1)$, $D_6(a_2)$  &   $3$ & $12$  & $2\times32 = 64$  &    $76$ \\
    \hline
       $D_7$, $D_7(a_1)$, $D_7(a_2)$ &   $3$   & $14$ & $2\times64 = 128$  &  $142$ \\
    \hline
       $D_l$, ~$D_l(a_k)$, $l > 7$ &   $1$  & $2l$ & -  &  $2l$ \\
   \hline    %1
\end{tabular}
  ~\\ \vspace{2mm}
  \centering
\begin{tabular} {|c|c|c|c|c|}
  \hline  %1
    $\Gamma$
     & \multicolumn{4}{c|}{Number of linkage diagrams, and the number $p = \mathscr{B}^{\vee}_{\Gamma}(\gamma^{\nabla})$}  \cr
        \cline{2-5}
    &  $A$-component        & $D$-component         &  $E$-component         &
     \hspace{6.0mm}In all \hspace{6.0mm}  \cr
     %% In all \cr
      \cline{2-4}
          &    $p = \frac{l}{l+1}$ & $p = 2\frac{l-1}{l+1}$    &  $p = 3\frac{l-2}{l+1}$   &      \\
  \hline
        $A_3$     &  $2 \times 4 = 8$ & $6$  &  -     & $14$ \\
  \hline
        $A_4$     &  $2 \times 5 = 10$ & $2 \times 10 = 20$  &  -     & $30$ \\ %% $10 + 20 = 30$    \\
  \hline
        $A_5$     &  $2 \times 6 = 12$ & $2 \times 15 = 30$  & $20$  & $62$ \\ %% $12 + 30 + 20 = 62$ \\
  \hline
        $A_6$     &  $2 \times 7 = 14$ & $2 \times 21 = 42$  & $2 \times 35 = 70$  & $126$ \\ %% $14 + 42 + 70 = 126$   \\
  \hline
        $A_7$     &  $2 \times 8 = 16$ & $2 \times 28 = 56$  & $2 \times 56 = 112$
                        & $184$ \\ %% $\begin{array}{c} 16 + 56 + 112 = 184 \end{array}$  \\
  \hline
        $A_8$     &  $2 \times 9 = 18$ & $2 \times 36 = 72$  & -  & $90$ \\ %% $\begin{array}{c}  $18 + 72  = 90$  \end{array}$ \\
  \hline
        $A_l$, $l > 8$  &  $2\times(l+1)$  & $2 \times \frac{l(l+1)}{2} = l(l+1)$   &  - &
         $(l+1)(l+2)$   \\
   \hline
  \end{tabular}
  \vspace{2mm}
  \caption{\hspace{1mm}The $ADE$ components of the linkage system $\mathscr{L}(\Gamma)$,
  the number of linkage diagrams in every component and values $\mathscr{B}^{\vee}_{\Gamma}(\gamma^{\nabla})$.}
  \label{tab_val_Bmin1}
  \end{table}

\subsubsection{Dynkin extensions}
  \label{sec_Dynkin_ext0}

   We show in Corollary \ref{cor_exist_Dynkin_ext} that for a Carter diagram  $\Gamma$ (rank($\Gamma$) $< 8$) and any
   $\Gamma$-associated root subset $S$ \underline{there exists the Dynkin diagram $\Gamma'$} such that $\rank(\Gamma') = \rank(\Gamma) + 1$
   and $\varPhi(S) \subset \varPhi(\Gamma')$, where $\varPhi(S)$ is the partial root system
   and $\varPhi(\Gamma')$ is the root system associated with $\Gamma'$.
   This pair $\{\Gamma, \Gamma'\}$ is said to be the {\it Dynkin extension  of the Carter diagram $\Gamma$} and
   is denoted by $\Gamma <_D \Gamma'$.

   %% Consider a Carter diagram  $\Gamma$ and a certain $\Gamma$-associated root subset $S$.
   %% Suppose $\Gamma'$ is the Dynkin diagram such that $\rank(\Gamma') = \rank(\Gamma) + 1$ and
   %% $S \subset \varPhi(\Gamma')$, where  $\varPhi(\Gamma')$ is the root system associated with $\Gamma'$.
   %% Such an extension is said to be the {\it Dynkin extension  of the Carter diagram $\Gamma$} and denote it by $\Gamma <_D \Gamma'$.
   %% It is shown that the {\bf Dynkin extension $\Gamma <_D \Gamma'$ exists}
   %% for any simply-laced Carter diagram $\Gamma$ (Corollary \ref{cor_exist_Dynkin_ext}).

  For any pair $\{\widetilde\Gamma, \Gamma\}$ out of list \eqref{eq_pairs_root_syst_0},
  we have $\mid \mathscr{L}(\widetilde\Gamma) \mid ~=~ \mid \mathscr{L}(\Gamma) \mid$ (Corollary \ref{corol_Estim_2}).
%%  {\bf Corollary} (Corollary \ref{cor_exist_Dynkin_ext})
%%  {\it For any Carter diagram out of the list of \S\ref{sec_coval_list},
%%    there exists the Dynkin extension.}

\subsubsection{Root stratum}
  \label{sec_Dynkin_ext}
 \index{regular extension $\Gamma < \widetilde{\Gamma}$}
 \index{$\Gamma < \widetilde{\Gamma}$, regular extension}
 \index{$\Gamma \stackrel{\alpha}{<} \widetilde{\Gamma}$, regular extension}

  Consider the connection diagram $\widetilde{\Gamma}$, see \S\ref{sec_connection},
  obtained from a certain Carter diagram $\Gamma$ by adding only one vertex $\alpha$,
  where $\alpha$ is connected to $\Gamma$ at $v$ points, where $v=1,2$ or $3$. If $\widetilde{\Gamma}$ is also a Carter diagram, this
  extension $\Gamma$ to $\widetilde{\Gamma}$ is said to be a {\it regular extension} and is denoted by
 \begin{equation*}
  \Gamma < \widetilde{\Gamma} \quad \text{ or } \quad  \Gamma \stackrel{\alpha}{<} \widetilde{\Gamma}.
 \end{equation*}
   The Dynkin extension is not necessary regular.  For example, $D_n(a_k) <_D D_{n+1}$
   is a Dynkin extension but $D_n(a_k)$ is not a subdiagram of $D_{n+1}$.
 ~\\

   The subset of roots $\varPhi(\Gamma')\backslash\varPhi(S)$ is said to be the {\it root stratum};
   it consists of roots of $\varPhi(\Gamma')$ \underline{linearly independent} of roots of $S$.

\subsubsection{Linkage system components}
  \label{sec_link_syst_comp}
 \index{ $\varPhi(S)$, partial root system}
 \index{root system ! - partial}
 \index{partial root system}
 \index{conjugation ! - of bases of the root system}
 \index{conjugation ! - non-conjugate $\Gamma$-associated root subsets $S_1$ and $S_2$}
 \index{root stratum}
 \index{root system}
  \index{linkage system component}
  \index{stratum size}
  Let $\mathscr{L}_{\Gamma'}(\Gamma) := \{\gamma^{\nabla} \mid \gamma \in \varPhi(\Gamma')\backslash\varPhi(S)\}$.
  Here, $\gamma^{\nabla}$ is the vector of linkage labels associated with $\gamma$.
  The set $\mathscr{L}_{\Gamma'}(\Gamma)$ is said to be the {\it linkage system component}
 (or, {\it $\Gamma'$-component} of $\mathscr{L}(\Gamma)$).
  The linkage system $\mathscr{L}(\Gamma)$ is the union of the sets $\mathscr{L}_{\Gamma_i}(\Gamma)$
  taken for all Dynkin extensions of $\Gamma$:
  \index{linkage diagrams ! - in $\mathscr{L}_{\Gamma'}(\Gamma)$}
 \begin{equation*}
     \mathscr{L}(\Gamma) \quad = \quad  \bigcup_{\Gamma <_D \Gamma_i} \mathscr{L}_{\Gamma_i}(\Gamma).
 \end{equation*}

 There are several Dynkin extensions and, accordingly,
 several linkage system components for the given Carter diagram $\Gamma$.
 The number of roots in the root stratum $\varPhi(\Gamma')\backslash\varPhi(S)$
 is said to be the {\it stratum size}. We have
 \begin{equation*}
   \mid \varPhi(\Gamma')\backslash\varPhi(S) \mid  ~=~  \mid\varPhi(\Gamma')\mid - \mid\varPhi(S)\mid
     \text{ and }
   \mid \mathscr{L}_{\Gamma'}(\Gamma) \mid  \quad \leq \quad  \mid\varPhi(\Gamma')\mid - \mid\varPhi(S)\mid,
 \end{equation*}
 see \S\ref{sec_part_root_and_comp}. From the latter inequality we get the following proposition:
 ~\\
 ~\\
  {\bf Proposition} (Proposition \ref{prop_Estim})
     {\it For any Carter diagram $\Gamma$, we have $\mid \mathscr{L}(\Gamma) \mid \leq \mathscr{E}$, where
     $\mathscr{E}$ is given by Table \ref{tab_val_Estim}.}
 ~\\

  It is checked in \S\ref{sec_stratif_D}, \S\ref{sec_stratif_E} and \S\ref{sec_stratif_A} that
  there exists at least $\mathscr{E}$ linkage diagrams in the linkage system $\mathscr{L}(\Gamma)$.
  Together with Proposition \ref{prop_Estim}, we get the following
  ~\\
  ~\\
 {\bf Corollary} (Corollary \ref{corol_Estim})
     {\it For any Carter diagram $\Gamma$, we have $\mid \mathscr{L}(\Gamma) \mid = \mathscr{E}$,
     see Table \ref{tab_val_Estim}.}
~\\

   For each Carter diagrams $\Gamma \in  \mathsf{C4} \coprod \mathsf{DE4} \coprod \mathsf{A}$,
   the values $\mathscr{B}^{\vee}_{\Gamma}(\gamma^{\nabla})$, the number of components, and number of
   linkage diagrams are presented in Table \ref{tab_val_Bmin1}, see Theorem \ref{th_full_descr}.
 ~\\

 Every component of $\mathscr{L}(\Gamma)$
 is determined by the type of Dynkin extension $\Gamma <_D \Gamma'$. There are $A$-, $D$-, $E$-types of
 Dynkin extensions. The corresponding components of the linkage systems are called $A$-, $D$-, $E$-components.
 For description of components of the linkage system for each Carter diagram $\Gamma$, see
 Table \ref{tab_val_Bmin1} and Theorem \ref{th_full_descr}.
 ~\\

%% file: 2IntroMain.tex
\section{\sc\bf Introduction and main results}

\subsection{The partial Cartan matrix, inverse quadratic form and linkages}
 \label{sec_part_Cartan}

\subsubsection{Linkages and linkage diagrams}
  \label{sec_linkage}
  \index{$\varPhi$, root system}
  \index{bicolored decomposition}
  \index{$w = w_1{w}_2$, bicolored decomposition of $w$; $w_1$, $w_2$ - involutions}
  \index{$k$, number of coordinates in $\alpha$-set (= number of $\alpha$\text{-labels})}
  \index{$h$, number of coordinates in $\beta$-set (= number of $\beta$\text{-labels})}
  \index{$l$, number of vertices in the Carter diagram $\Gamma$, $l = k + h$}

 Let $w = w_1w_2$ be the bicolored decomposition of some element
 $w \in W$, where $w_1$, $w_2$ are two involutions
 associated, respectively, with an $\alpha$-set $\{ \alpha_1, \dots, \alpha_k \}$ and a $\beta$-set
 $\{ \beta_1, \dots, \beta_h \}$  of roots from the root system $\varPhi$, see
 \eqref{two_invol}, \eqref{two_sets}, and let $\Gamma$ be the Carter
 diagram associated with this bicolored decomposition.
 We consider the {\it extension} of the root basis $\Pi_w$ by means of the root $\gamma \in \varPhi$,
 so that the set of roots
 \index{dotted ! - edge}
 \index{solid ! - edge}
 \begin{equation}
   \label{alpha_beta}
    \Pi_w(\gamma) = \{ \alpha_1, \dots, \alpha_k, \beta_1, \dots, \beta_h, \gamma \}
 \end{equation}
 is linearly independent. Let us multiply $w$ on the right by the reflection
 $s_{\gamma}$ corresponding to $\gamma$ and consider the diagram
 $\Gamma' = \Gamma \cup \gamma$ together with new edges.
By \S\ref{sec_solid_dotted}, these edges are \underline{solid} (resp. \underline{dotted}) for $(\gamma, \tau) = -1$ (resp. $(\gamma, \tau) = 1$),
where $\tau \in \Pi_w$, see \eqref{root_subset_L}.
~\\
~\\
\begin{minipage}{10.6cm}
 \index{$\alpha$-set, subset of roots corresponding to $w_1$}
 \index{$\beta$-set, subset of roots corresponding to $w_2$}
 \index{label ! - linkage label vector}
 \index{$\alpha$-endpoints}
 \index{$\beta$-endpoints}
 \index{linkage label vector}

 \hspace{5mm}The diagram $\Gamma'$ is said to be
 a {\it linkage diagram} and the root $\gamma$  is said to be a {\it linkage root}.
 The roots $\tau$ corresponding to the new edges ($(\gamma, \tau) \neq 0$) are said to be
 {\it endpoints} of the linkage diagram; endpoints lying in an $\alpha$-set (resp. a $\beta$-set)
 are said to be {\it $\alpha$-endpoints} (resp. {\it $\beta$-endpoints}).
 Consider vector $\gamma^{\nabla}$ defined by \eqref{dual_gamma}.
 This vector is said to be the {\it linkage label vector}.
 There is, clearly, the one-to-one correspondence between linkage label
 \end{minipage}
\begin{minipage}{5.4cm}
  \quad
 \begin{equation}
   \label{dual_gamma}
   \gamma^{\nabla} :=
   \left (
    \begin{array}{c}
    (\gamma, \alpha_1) \\
    \dots,     \\
    (\gamma, \alpha_k)  \\
    (\gamma, \beta_1)  \\
    \dots,  \\
    (\gamma, \beta_h) \\
    \end{array}
   \right )
 \end{equation}
\end{minipage}
~\\
     vectors $\gamma^{\nabla}$
 (with labels $\gamma^{\nabla}_i \in \{-1, 0, 1\}$) and simply-laced linkage diagrams
 (i.e., linkage diagrams  such that $(\gamma, \tau) \in \{-1, 0, 1\}$).

\subsubsection{The projection of the linkage root}
  \label{sec_proj_linkage}

 \index{$L$, linear space ! - spanned by the $\Gamma$-associated root subset}
 \index{$\gamma$, linkage root}
 \index{linkage root}
 \index{linkage label vectors ! - $\tau_i^{\nabla}$}

 Let $L \subset V$ be the linear space spanned by the $\Gamma$-associated root subset,
 \underline{$\gamma_L$ be the projection} of the linkage root $\gamma$ on $L$.
 For any  $\tau_i \in S$, we define vectors $\tau_i^{\nabla}$ as follows:
 \index{$L^{\nabla}$, space of linkage labels}
 \index{linkage label vectors ! - $L^{\nabla}$ spanned by vectors $\tau_i^{\nabla}, \tau_i \in \Pi_w$}
 \index{quadratic form ! - inverse quadratic form $\mathscr{B}^{\vee}_{\Gamma}$}
 \index{$\mathscr{B}^{\vee}_{\Gamma}$, inverse quadratic form}
   \begin{equation*}
       \tau_i^{\nabla} := B_{\Gamma}{\tau_i},
   \end{equation*}
 where $B_{\Gamma}$ is the partial Cartan matrix,
 for details see \S\ref{sec_dual_partial}.
 Consider the linear space $L^{\nabla}$ spanned by vectors $\tau_i^{\nabla}$, where $\tau \in \Pi_w$,
 see \S\ref{sec_linkage}.
 For any linkage root $\gamma \in \varPhi$, we denote by $\gamma^{\nabla} \in L^{\nabla}$ the vector
 obtained as follows:
 %%
 %% The partial Cartan matrix $B_{\Gamma}$ maps $\gamma_L \in L$ to $\gamma^{\nabla} \in L^{\nabla}$ as follows:
 \begin{equation*}
   \label{B_pairing}
      \gamma^{\nabla} = B_{\Gamma} \gamma_L, \quad B_{\Gamma}^{-1}\gamma^{\nabla} =  \gamma_L,
 \end{equation*}
 see Proposition \ref{prop_link_diagr_conn}.
 The linear space $L^{\nabla}$ is said to be {\it the space of linkage labels}.
 The vector $\gamma^{\nabla}$ is said to be the {\it linkage label vector}.
 In the space $L^{\nabla}$, we take the {\it inverse quadratic form}
 $\mathscr{B}^{\vee}_{\Gamma}$ associated with the inverse matrix $B_{\Gamma}^{-1}$.
 {\bf The quadratic form $\mathscr{B}^{\vee}_{\Gamma}$  provides
 an easily verifiable criterion that the vector $u$ is the linkage label vector
 for a certain linkage root $\gamma \not\in L$ (i.e., $u = B_{\Gamma}\gamma_L$ for a certain $\gamma \not\in L$).
 This criterion (Theorem \ref{th_B_less_2}) is the following inequality:}
  \begin{equation*}
      \mathscr{B}^{\vee}_{\Gamma}(u) < 2.
  \end{equation*}

  \subsection{The linkage systems and loctets}
  \index{dual partial Weyl group $W^{\vee}_S$}
  \index{Weyl group ! - dual partial $W^{\vee}_S$}

  \subsubsection{The starlike numbering of vertices}
    \label{sec_starlike}
  \label{two_predef_classes}
  \index{numbering near the branch point}
  \index{starlike numbering}
  \index{branch point}

\begin{figure}[H]
\centering
\includegraphics[scale=0.9]{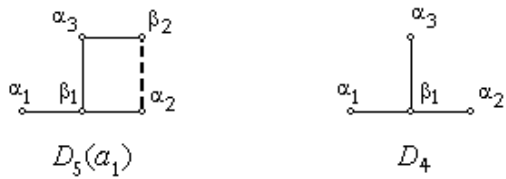}
\vspace{3mm} \caption{\hspace{3mm}The starlike numbering of vertices adjacent to the branch point}
%%%%%% The label must come after caption
\label{D5a1_D4_pattern}
\end{figure}

   Throughout this article, we use a special numbering of vertices  adjacent to the branch point.
   Such a numbering is said to be {\it starlike}. Let $\Gamma$ be the Carter diagram such that
   $\Gamma \in \mathsf{C4} \coprod \mathsf{DE4}$ and $\Gamma \ne D_4(a_1)$. Then $\Gamma$ contains $D_4$ as subdiagram.
   We use the numbering of the vertices shown in Fig \ref{D5a1_D4_pattern}. Namely, the branch point is denoted by $\beta_1$ and the neighboring vertices
   are denoted by $\{\alpha_1, \alpha_2, \alpha_3\}$.
   The starlike numbering is presented, for example, in  Figs. \ref{E7a1_gamma_ij_3} -- \ref{E7a3_gamma_ij_4}.
   This numbering is very important for unification construction of loctets.

\subsubsection{Loctets}

  The group $W^{\vee}_S$, named the dual partial Weyl group, see \S\ref{sec_dual_partial}, acts
  on vectors of linkage labels as follows:
 \begin{equation*}
      (w\gamma)^{\nabla} = w^{*}\gamma^{\nabla} \text{ for any } \gamma^{\nabla} \in L^{\nabla}, \quad w^{*} \in W^{\vee}_S,
 \end{equation*}
 see Proposition \ref{prop_link_diagr_conn}. Under the action of $W^{\vee}_S$ the set of
 linkage diagrams (= linkage label vectors) constitute the
 diagram called the {\it linkage system} similarly to the weight system
 in the representation theory of semisimple Lie algebras, \cite[p. 30]{Sl81}.
 We denote the linkage system associated with the Carter diagram $\Gamma$ by $\mathscr{L}(\Gamma)$.

\begin{figure}[H]
\centering
\includegraphics[scale=1.3]{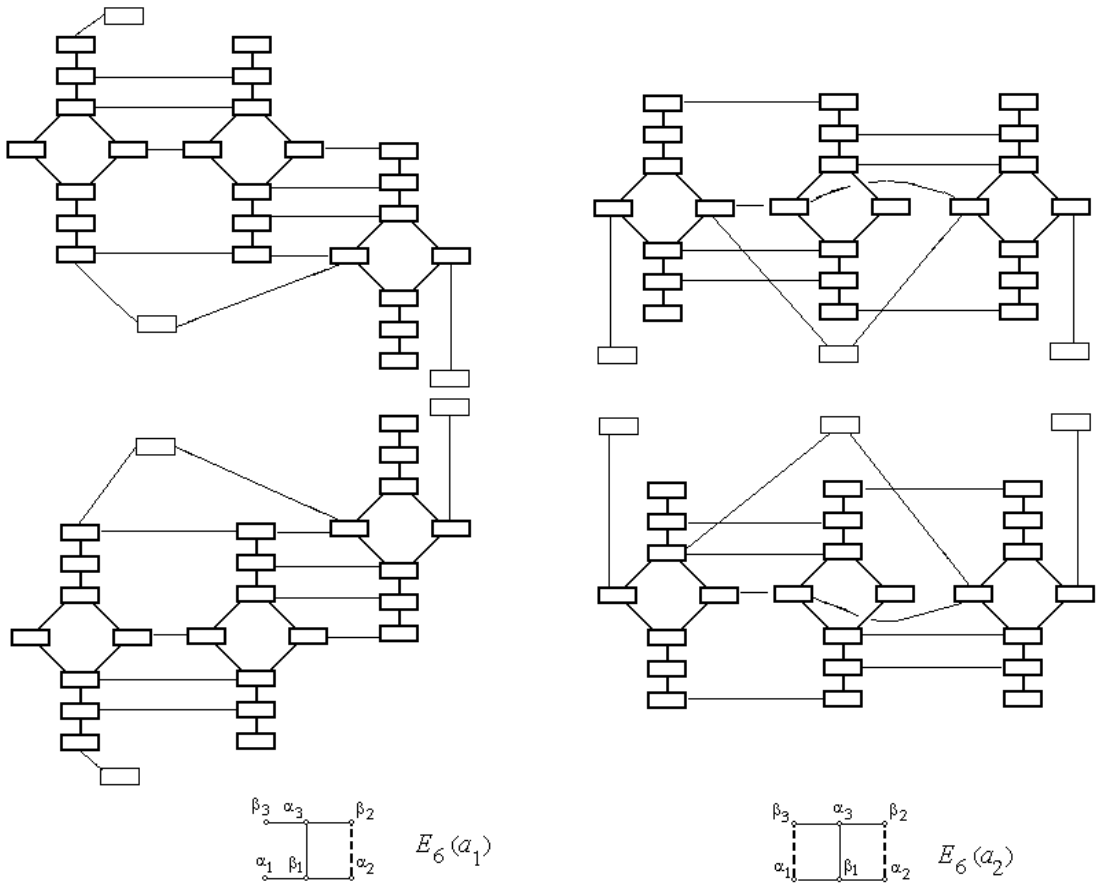}
\caption[\hspace{3mm}Linkage systems $\mathscr{L}(E_6(a_1))$ and
$\mathscr{L}(E_6(a_2))$.
 The $8$-cell bold subdiagrams are loctets]{\hspace{3mm}Linkage systems
 $\mathscr{L}(E_6(a_1))$ and $\mathscr{L}(E_6(a_2))$,
 see Fig. \ref{E6a1_linkages} and Fig. \ref{E6a2_linkages}. The $8$-cell {\bf bold} subdiagrams are {\bf loctets},
 see Fig. \ref{fig_loctets}}
%%%%%% The label must come after caption
\label{linkage_systems_example}
\end{figure}

 The linkage systems $\mathscr{L}(E_6(a_1))$ and $\mathscr{L}(E_6(a_2))$ are depicted in Fig. \ref{linkage_systems_example}.
 The linkage systems for all Carter diagrams are presented in Figs. \ref{D4a1_linkages}--\ref{E7a4_linkages},
  \ref{Dk_al_wind_rose},   \ref{27_weight_diagr_E6__2comp}(top),  \ref{E7pure_linkage_system}, \ref{D5pure_loctets},
  \ref{D6pure_loctets}, \ref{D7a1_linkages_cdef}--\ref{D7a1_D7a2_D7pu_loctets_comp3}.
 Every linkage diagram containing at least one non-zero $\alpha$-label $\alpha_1$, $\alpha_2$ or $\alpha_3$ (see \S\ref{sec_Carter})
 belongs to a certain $8$-cell "spindle-like" linkage subsystem called
 \underline{loctet} (= linkage octet). {\bf The loctets are the main construction blocks for every linkage system.}

 Let  $a_i = (\alpha_i, \gamma)$ (resp. $b_i = (\beta_i, \gamma)$), where $i = 1,2,3$, be coordinates of the linkage label vector $\gamma^{\nabla}$.
 If all $a_i = 0$ (resp. all $b_i = 0$), the linkage diagram $\gamma^{\nabla}$
 is said to be a {\it $\beta$-unicolored ({\rm resp.} $\alpha$-unicolored) linkage diagram}.
 {\bf Every linkage system is the union of several loctets and several $\beta$-unicolored linkage diagrams,}
 see \S\ref{sec_enum_loctets}.
 There are exactly $6$ loctets in the linkage system $\mathscr{L}(E_6(a_k))$ for $k = 1, 2$,
 see Figs. \ref{E6a1_linkages}, \ref{E6a2_linkages}, and in the linkage system
 $\mathscr{L}(E_7(a_k))$ for $k = 1, 2, 3, 4$, see Figs. \ref{E7a1_linkages} -- \ref{E7a4_linkages}.
 For $E_6(a_1)$ and $E_6(a_2)$,  the linkage system is the union of $6$ loctets and $6$ $\beta$-unicolored linkage diagrams,
 altogether $54 = 2 \times 27 = 6 \times 8 + 6$ linkage diagrams, see Fig. \ref{linkage_systems_example}.

%% We use the knowledge of the structure of linkage systems in our proof of the Carter theorem
%% on decomposition of every element in the Weyl group into the product of two involutions, see \cite{St11}.
%%~\\

 \index{$\alpha$-label, coordinate out of $\alpha$-set}
 \index{$\beta$-label, coordinate out of $\beta$-set}
 \index{label ! - $\alpha$-label}
 \index{label ! - $\beta$-label}
 \index{loctet ! - (= linkage octet)}
 \index{linkage diagrams ! - $\alpha$-unicolored}
 \index{linkage diagrams ! - $\beta$-unicolored}
 \index{$\gamma$, linkage root}
 \index{linkage root}
 \index{$W$, Weyl group}
 \index{Weyl group ! - $W$}

\subsection{The Carter diagrams and connection diagrams}
 \label{sect_diagram}
 \index{bicolored decomposition}
 \index{class of diagrams ! - connection diagrams}
 \index{class of diagrams ! - simply-laced connected Carter diagrams}
 \index{connection diagram}

  In this paper we consider several types of diagrams.
  The {\it Carter diagrams} are obtained by a small deformation of {\it admissible diagrams}
  introduced by R.~Carter \cite{Ca72}, they
  describes a bicolored decomposition of some element $w \in W$, see \S\ref{sec_Carter}.
  The {\it connection diagrams} generalize the Carter diagram; each connection diagram describes a decomposition
  (not necessary bicolored) of a certain element $w \in W$,
  and this diagram is supplied with a certain order of reflections, see \S\ref{sec_connection}.
  In both cases all reflections are associated with roots which are not necessary simple.
  The {\it linkage diagram} is a particular case of the connection diagram obtained from a certain Carter diagram
  by adding one extra vertex with its bonds, see \S\ref{sec_linkage}.
  The linkage diagrams are the focus of this paper.

\subsubsection{The Dynkin diagrams}
  \label{sec_Dynkin}
  \index{$E$, linear space spanned by all roots of $\varPhi$}
  \index{$\Pi$, set of all simple roots in $\varPhi$}
  \index{${\bf B}$, Cartan matrix associated with a Dynkin diagram}
  \index{Cartan matrix ! - ${\bf B}$}
  \index{$(\cdot, \cdot)$, symmetric bilinear form associated with ${\bf B}$}
  \index{symmetric bilinear form $(\cdot, \cdot)$}
  \index{quadratic form ! - $\mathscr{B}$, quadratic Tits form associated with the Cartan matrix ${\bf B}$}
  \index{$\mathscr{B}$, quadratic Tits form associated with the Cartan matrix ${\bf B}$}
  \index{root system}

  Let $\varPhi$ be a certain classic root system such that all roots of the same length.
  Then the corresponding Dynkin diagram $\Gamma_D$ is simply-laced.
  Let $\Pi$ the set of all simple roots in $\varPhi$, $V$ be the linear space spanned by all roots,
  $W$ the finite Weyl group acting on $\varPhi$.
  Let ${\bf B}$ be the corresponding Cartan matrix, $(\cdot, \cdot)$ the corresponding symmetric bilinear form
  (= inner product on $V$) and $\mathscr{B}$ the quadratic Tits form associated with ${\bf B}$,  \cite[Ch. 2]{St08}.
  We suppose that each diagonal element of ${\bf B}$ is equal to $2$, see Remark \ref{rem_classic}. The following relation
  is the well-known property connecting roots and the quadratic Tits forms\footnotemark[1]:
  \begin{equation}
     \label{eq_Kac}
      \mathscr{B}(\alpha) = 2 \Longleftrightarrow \alpha \in \varPhi.
  \end{equation}
  For any two non-orthogonal simple roots $\alpha \neq \beta$, we have
  \begin{equation}
     \label{eq_Kac_2}
      (\alpha, \beta) = \Arrowvert \alpha \Arrowvert \Arrowvert \beta \Arrowvert
      \cos(\widehat{\alpha, \beta}) = \sqrt{2} \cdot \sqrt{2} (-\frac{1}{2}) = -1.
  \end{equation}

  \footnotetext[1]{In order for the values of the linkage labels (see \S\ref{sec_linkage})
  be integer, as in \eqref{eq_Kac_2}, we choose the diagonal elements equal $2$.
  Often, the diagonal elements are chosen equal $1$;
  then \eqref{eq_Kac} looks as follows: $\mathscr{B}(\alpha) = 1 \Longleftrightarrow \alpha \in \varPhi$,  see \cite{Kac80}.}

\subsubsection{The connection diagrams}
  \label{sec_connection}
  \index{class of diagrams ! - connection diagrams}
  \index{$l$, number of vertices in the Carter diagram $\Gamma$, $l = k + h$}
  \index{connection diagram}

  Let $\Gamma$ be the diagram characterizing
  connections between roots of a certain set $S$ of linearly
  independent and not necessarily simple roots, $o$
  be the order of reflections in the decomposition \eqref{any_roots_0}.
  The pair $(\Gamma, o)$ is said to be a {\it connection diagram}.
  Here, $\Gamma$ is the diagram describing
  connections between roots as it is described by
  the Dynkin diagrams or by the Carter diagrams, and $o$ is the
  order of elements in the (not necessarily bicolored) decomposition \eqref{any_roots_0},
  see \S\ref{seq_sim_diagram}, \S\ref{sec_basic_lem}.

%% \clearpage
\subsection{The main results}

\subsubsection{Partial Weyl group and dual partial Weyl group}
 \index{$W_S$, partial Weyl group}
 \index{Weyl group ! - partial $W_S$}
 \index{$W^{\vee}_S$, dual partial Weyl group}
 \index{Weyl group ! - dual partial $W^{\vee}_S$}
 \index{$B_{\Gamma}$, partial Cartan matrix}

 %%\index{$B_{\Gamma}$, partial Cartan matrix}
 %%\index{Cartan matrix ! - partial Cartan matrix $B_{\Gamma}$}
 \index{$\gamma^{\nabla}$, linkage diagram}
 \index{$S = \{\tau_1,\dots,\tau_l\}$, $\Gamma$-associated root subset}
 \index{dual reflection $s^{*}_{\tau_i}$}
 \index{$s^{*}_{\tau_i}$, dual reflection}
 \index{linkage label vector}

 Let $S = \{\tau_1, \dots, \tau_l \}$ be a certain $\Gamma$-associated root subset,
 where $\Gamma$ is a certain Carter diagram. %% roots $\tau_i$ are not necessarily simple.
 We introduce the partial Weyl group $W_S$ generated by reflections
 $\{s_{\tau_1}, \dots, s_{\tau_l} \}$, and the dual partial Weyl group $W^{\vee}_S$
 generated by dual reflections  $\{s^*_{\tau_1}, \dots, s^*_{\tau_l} \}$, see \S\ref{sec_dual_partial}.
 The linkage diagrams $\gamma^{\nabla}$ and $(w\gamma)^{\nabla}$ are related as follows:
 $(w\gamma)^{\nabla} = w^{*}\gamma^{\nabla}$, where $w^{*} \in W^{\vee}_S$ (Proposition \ref{prop_link_diagr_conn}).
 The quadratic form $\mathscr{B}_{\Gamma}$ takes a certain constant value for all elements
 $w\gamma$, where $w$ runs over  $W_S$ (Proposition \ref{prop_unique_val}).

\subsubsection{Bijection of covalent root systems}
  \label{sec_Dynkin_ext_2}

 \index{regular extension $\Gamma < \widetilde{\Gamma}$}
 \index{$\Gamma < \widetilde{\Gamma}$, regular extension}
 \index{$\Gamma \stackrel{\alpha}{<} \widetilde{\Gamma}$, regular extension}

   Let $\Gamma <_D \Gamma'$ be a certain Dynkin extension, see \S\ref{sec_Dynkin_ext0}, $S$ be a $\Gamma$-associated root subset.
   Recall that the partial root system $\varPhi(S)$ is the subset of roots of $\varPhi(\Gamma')$
   which are \underline{linearly dependent} on roots of $S$, and the root stratum $\varPhi(\Gamma')\backslash\varPhi(S)$
   consists of roots of $\varPhi(\Gamma')$ which are \underline{linearly independent} of roots of $S$,
   see \S\ref{sec_Dynkin_ext}.
   %% For the case of the the root system, any two bases are conjugate, see \cite[Ch.6. $\S1$, ${n^{\rm o}}5$, Theorem 2]{Bo02}.
   %% For the case of the partial root system, there are not conjugate two $\Gamma$-associated root subsets
   %% $S_1$ and $S_2$, see the example in \cite[$\S$B.1.2]{St10}.
   %% This is a reason why I prefer denote the partial root system by $\varPhi(S)$, not by $\varPhi(\Gamma)$.

  \index{ $\varPhi(S)$, partial root system}
  \index{root system ! - partial}
  \index{partial root system}
  \index{conjugation ! - of bases of the root system}
  \index{conjugation ! - non-conjugate $\Gamma$-associated root subsets $S_1$ and $S_2$}
  \index{root stratum}
  \index{root system}
  \index{linkage system component}
  \index{stratum size}

  Let $\mathscr{L}_{\Gamma'}(\Gamma) := \{\gamma^{\nabla} \mid \gamma \in \varPhi(\Gamma')\backslash\varPhi(S)\}$.
  The set $\mathscr{L}_{\Gamma'}(\Gamma)$ is said to be the {\it linkage system component}
 (or, {\it $\Gamma'$-component} of $\mathscr{L}(\Gamma)$).
  The linkage system $\mathscr{L}(\Gamma)$ is the union of the sets $\mathscr{L}_{\Gamma_i}(\Gamma)$
  taken for all Dynkin extensions of $\Gamma$:
  \index{linkage diagrams ! - in $\mathscr{L}_{\Gamma'}(\Gamma)$}
 \begin{equation*}
     \mathscr{L}(\Gamma) \quad = \quad  \bigcup_{\Gamma <_D \Gamma_i} \mathscr{L}_{\Gamma_i}(\Gamma).
 \end{equation*}
 ~\\
 There are several Dynkin extensions and, accordingly,
 several linkage system components for the given Carter diagram $\Gamma$.
 The number of roots in the root stratum $\varPhi(\Gamma')\backslash\varPhi(S)$
 is said to be the {\it stratum size}. We have
 \begin{equation*}
   \mid \varPhi(\Gamma')\backslash\varPhi(S) \mid  ~=~  \mid\varPhi(\Gamma')\mid - \mid\varPhi(S)\mid
     \text{ and }
   \mid \mathscr{L}_{\Gamma'}(\Gamma) \mid  \quad \leq \quad  \mid\varPhi(\Gamma')\mid - \mid\varPhi(S)\mid.
 \end{equation*}
 From the latter inequality we get the following proposition:
 ~\\
 ~\\
  {\bf Proposition} (Proposition \ref{prop_Estim})
     {\it For any Carter diagram $\Gamma$, we have $\mid \mathscr{L}(\Gamma) \mid \leq \mathscr{E}$, where
    the estimate $\mathscr{E}$ is given by Table \ref{tab_val_Estim}.}
 ~\\

  It is checked in \S\ref{sec_stratif_D}, \S\ref{sec_stratif_E} and \S\ref{sec_stratif_A} that
  there exists at least $\mathscr{E}$ linkage diagrams in the linkage system $\mathscr{L}(\Gamma)$.
  Together with Proposition \ref{prop_Estim}, we get the following
  ~\\
  ~\\
 {\bf Corollary} (Corollary \ref{corol_Estim})
     {\it For any Carter diagram $\Gamma$, we have $\mid \mathscr{L}(\Gamma) \mid = \mathscr{E}$,
     see Table \ref{tab_val_Estim}.}
 ~\\
 ~\\
  By Proposition \ref{prop_Estim} we get also the following
 ~\\
 ~\\
 {\bf Corollary} (Corollary \ref{corol_Estim_2})
   {\it For the Carter diagram $E_i(a_k)$ and $D_i(a_k)$,  we have }
 \begin{equation*}
   \Small
     \begin{split}
      & \mid \mathscr{L}(E_6(a_k)) \mid = \mid \mathscr{L}(E_6) \mid, \text{ where } k = 1,2, \\
      & \mid \mathscr{L}(E_7(a_k)) \mid = \mid \mathscr{L}(E_7) \mid, \text{ where } k = 1,2,3,4, \\
      & \mid \mathscr{L}(D_l(a_k)) \mid = \mid \mathscr{L}(D_l) \mid, \text{ where }
        l \geq 4 \text{ and } 1 \leq k \leq \big [ \tfrac{l-2}{2} \big ].
     \end{split}
   \end{equation*}
 ~\\

   For each Carter diagrams $\Gamma \in  \mathsf{C4} \coprod \mathsf{DE4} \coprod \mathsf{A}$,
   the values $\mathscr{B}^{\vee}_{\Gamma}(\gamma^{\nabla})$, the number of components, and number of
   linkage diagrams are presented in Table \ref{tab_val_Bmin1}, see Theorem \ref{th_full_descr}.
 ~\\

   The essential part of \S\ref{sec_transition} is the bijection of the partial root system and the root system:
 ~\\
 ~\\
   {\bf Theorem} (Theorem \ref{prop_bijection_root_syst})
  {\it Consider a covalent pair of Carter diagrams $\{\widetilde\Gamma, \Gamma\}$, see \S\ref{sec_coval_list}.
  %% out of the list \eqref{eq_covalent_pairs}.
  Let $\mathcal{R}$ (resp. $\mathcal{P}$) be the root system (resp. partial root system) corresponding to
  $\Gamma$ (resp. $\widetilde\Gamma$).
 %%  \begin{equation}
 %%   \label{eq_pairs_root_syst_0}
 %%   \footnotesize
 %%    \begin{split}
 %%     & \{D_l(a_k), D_l\} \text{ for } l \geq 4 \text{ and } 1 \leq k \leq \big [\tfrac{l-2}{2} \big ], \\
 %%     & \{E_6(a_k), E_6\} \text{ for }  k = 1,2, \\
 %%     & \{E_7(a_k), E_7\} \text{ for }  k = 1,2,3,4.
 %%    \end{split}
 %%  \end{equation}
   There exists the map $M : \mathcal{P} \longmapsto \mathcal{R}$
   given by Table \ref{tab_part_root_syst} is the invertible linear map.}
 ~\\
 ~\\
  {\bf Corollary} (Corollary \ref{corol_biject})
 {\rm (i)} {\it For any pair \eqref{eq_covalent_pairs},
 the root system $\varPhi(\Gamma')$ contains the root subsystem $\mathcal{R}$ if and only if
 $\varPhi(\Gamma')$ contains the partial root subsystem $\mathcal{P}$.}

 {\rm (ii)} {\it For any covalent pair \eqref{eq_covalent_pairs}, we have } $ \mid \mathcal{R} \mid = \mid \mathcal{P} \mid.$
~\\
~\\
  {\bf Corollary} (Corollary \ref{cor_exist_Dynkin_ext})
     \label{cor_exist_Dynkin_ext}
  {\it For any Carter diagram out of the list of \S\ref{sec_coval_list}, there exists the Dynkin extension.}

 \subsubsection{Whether a vector is a linkage root?}
    \label{sec_criterion}

 Let $S = \{\tau_1, \dots, \tau_l \}$ be a certain $\Gamma$-associated root subset,
 where roots $\tau_i$ are not necessarily simple.
 Let $L \subset E$ be the linear subspace spanned by all roots of $S$,
 let $\gamma_L$ be the projection of the root $\gamma$ on $L$.
 The main result of \S\ref{sect_quadr_form} is the following theorem which verifies whether or not
 a given vector is a linkage root:
  ~\\
  ~\\
 {\bf Theorem} (Theorem \ref{th_B_less_2}).
     {\it A root $\theta \in \varPhi$ is a linkage root (i.e., $\theta$ is linearly independent
     of roots of $L$)
     %% A vector $u^{\nabla} \in L^{\nabla}$ is a linkage label vector
     %% corresponding to a certain root $\gamma \in \varPhi$, where $\gamma \not\in L$
     %%(i.e., $u^{\nabla} = B_{\Gamma}\gamma_L$),
     if and only if
  \index{label ! - linkage label vector}
     \begin{equation}
       \label{eq_p_less_2_loc}
         %% \mathscr{B}^{\vee}_{\Gamma}(u^{\nabla}) < 2.
          \mathscr{B}^{\vee}_{\Gamma}(\theta^{\nabla}) < 2,
     \end{equation}
   where $\theta^{\nabla} = B_\Gamma\theta_L$ and $\theta_L$ is the projection of $\theta$ on $L$.
   }

 \subsubsection{Three loctet types}
    \index{$\alpha$-label, coordinate out of $\alpha$-set}
    \index{label ! - $\alpha$-label}
    \index{loctet ! - (= linkage octet)}
    \index{loctet ! - of type $L_{ij}$}
    \index{$L_{ij}$, loctet of type $(ij) \in \{(12), (13), (23)\}$}
    \index{$\gamma^{\nabla}_{ij}(n)$, linkage diagram}

\begin{figure}[h]
\centering
\includegraphics[scale=1.3]{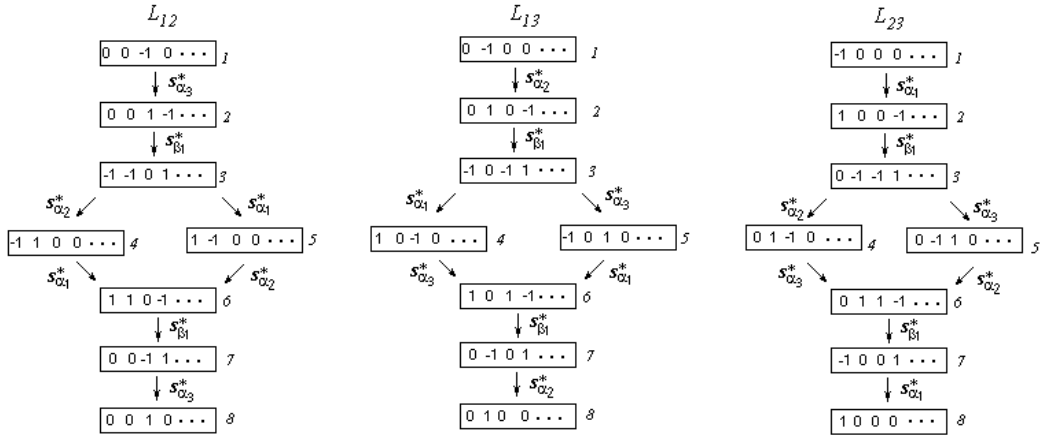}
 \caption{\hspace{3mm}The loctet types $L_{12}$, $L_{13}$ and $L_{23}$.}
%%%%%% The label must come after caption
\label{fig_loctets}
\end{figure}

 In \S\ref{sec_enum_loctets}, by means of inequality \eqref{eq_p_less_2_loc},
 we obtain a complete description of linkage diagrams for all linkage systems.
 We introduce the $8$-cell linkage subsystem called
 {\it loctet (= linkage octet)} as depicted in Fig. \ref{fig_loctets}.
 Recall that we adhere to a starlike numbering of \S\ref{sec_starlike}, so every linkage label vector in Fig. \ref{fig_loctets}
 looks as $\{\alpha_1, \alpha_2, \alpha_3, \beta_1, \dots\}$.
 Here, in Fig. \ref{fig_loctets}, we consider the case where $\alpha$-set contains only $3$ coordinates,
 but in the common case the structure of lockets does not change.
 Loctets are the main construction blocks in the structure of the linkage systems.
  Consider roots $\gamma^{\nabla}_{ij}(n)$ depicted in Fig. \ref{fig_loctets},
  where $\{ij\}$ is associated with type $L_{ij}$ and $n$ is the order number of the linkage diagrams
  in the vertical numbering in Fig. \ref{fig_loctets}.
  We call the octuple of linkages depicted in every connected component
  in Fig. \ref{fig_loctets} the {\it loctet } of type $L_{12}$
  (resp. $L_{13}$, resp. $L_{23}$).

 \index{loctet ! - structure of loctets}
 \index{linkage diagrams ! - $\beta$-unicolored}

 \subsubsection{The structure of loctets}
    \label{sec_struct_loctets}
  ~\\
  {\bf Corollary} (on the structure of loctets and linkage diagrams (Corollary \ref{corol_loctet})).
  {\it
    {\rm (i)} Any linkage diagram containing a non-zero $\alpha$-label $\alpha_1$, $\alpha_2$ or $\alpha_3$
    belongs to one of the loctets of the linkage system.

    {\rm (ii)} Any linkage diagram of the loctet uniquely determines the whole loctet.

    {\rm (iii)} If two loctets have one common linkage diagram, they coincide.

    {\rm (iv)} Every linkage diagram from the linkage system either belongs to one of the loctets
    or is $\beta$-unicolored.
   }

\index{class of diagrams ! - $\mathsf{C4}$}
\index{class of diagrams ! - $\mathsf{DE4}$}
 %%%

   \begin{remark} { \rm %% (i)
     \label{rem_l_less_8}
    The Carter diagrams $E_8(a_k)$ for $k = 1,\dots,4$ and $E_8$  do not represent any conjugacy classes in  $W(D_n)$
    and what is why:

     (a) Among Carter diagrams of $\mathsf{DE4}$ class, see \S\ref{sec_3_class},
     only Carter diagrams $D_l$ (for some $l$) represent conjugacy classes in $W(D_n)$;
     the $E_l$-associated root subsets, where $i=6,7,8$, can not be mapped into
     $W(D_n)$, see Lemma \ref{lem_E6_not_in_Dn}.

     (b) Among Carter diagrams of $\mathsf{C4}$ class,
     only Carter diagrams $D_l(a_k)$ (for some $l, k$) represent conjugacy classes in  $W(D_n)$;
     the $E_l(a_k)$-associated root subsets, where $l=6,7,8$, can not be mapped into  $W(D_n)$,
     see Lemma \ref{lem_E6ak_not_in_Dn}.
     %% For this issue see also \cite[p. 13]{Ca72}.

    Hence, there is no a linearly independent $9$-element root subset containing
    any $E_8(a_k)$-associated or $E_8$-associated  root subset, i.e., \underline{there are no linkages} for any Carter diagram of type
    $E_8(a_k)$ or $E_8$. For this reason, among simply-laced Carter diagrams of $E$-type it suffices
    to consider only diagrams with a number of vertices $l < 8$. }\qed
   \end{remark}

   %% (ii) Let $\Gamma$ be the union of two Carter diagrams $\Gamma_1$ and $\Gamma_2$, $S_1$ (resp. $S_2$)
   %% a $\Gamma_1$-associated (resp. $\Gamma_2$-associated) root subset.
   %% We have the union $S = S_1 \coprod S_2$, and the direct sums $L = L_1 \oplus L_2$ and
   %% $L^{\nabla} = L_1^{\nabla} \oplus L_2^{\nabla}$, see \S\ref{sec_proj_linkage}.
   %% The linkage diagram $\gamma^{\nabla} \in L^{\nabla}$ splits into two summands
   %% $\gamma^{\nabla} = \gamma^{\nabla}_1 + \gamma^{\nabla}_2$,
   %% where $\gamma^{\nabla}_i \in L_i^{\nabla}$. Accordingly, the partial Cartan matrix $B_{\Gamma}$,
   %% and its inverse matrix  $B_{\Gamma}^{-1}$ admit a decomposition into a direct sum,
   %% and $W_S^{\vee} = W_{S_1}^{\vee} \times W_{S_2}^{\vee}$. Then
   %% $\mathscr{L}(\Gamma) = \mathscr{L}(\Gamma_1) \times  \mathscr{L}(\Gamma_2)$. } \qed
  \index{loctet ! - $8$th linkage diagram}
  \index{linkage diagrams ! - $\beta$-unicolored}
  \index{$\gamma^{\nabla}_{ij}(8)$, linkage diagram}
  \subsubsection{On tables and diagrams}
  In \S\ref{sec_calc_gamma_8}, the calculation technique for loctet diagrams  $\gamma^{\nabla}_{ij}(8)$ (= $8$th linkage diagram of
  the loctet) is explained. According to Corollary \ref{corol_loctet}, the whole loctet is uniquely determined
  from $\gamma^{\nabla}_{ij}(8)$. From Tables \ref{sol_inequal_1}--\ref{sol_inequal_6} one can recover the calculation of
  $\gamma^{\nabla}_{ij}(8)$ for Carter diagrams $\Gamma \in \mathsf{C4} \coprod \mathsf{DE4}$.
  In \S\ref{sec_calc_homog},
  the calculation technique for $\beta$-unicolored linkage diagrams is similarly explained.
  From Tables \ref{homog_inequal_1}--\ref{homog_inequal_3} one can recover the calculation of $\beta$-unicolored linkage diagrams
  for Carter diagrams $\Gamma \in \mathsf{C4} \coprod \mathsf{DE4}$.
  In \S\ref{sec_diagr_per_comp}, the loctets per component
  for all linkage systems are listed, see Table \ref{tab_seed_linkages_6}.
  In \S\ref{sec_inv_matr}, the partial Cartan matrix $B_{\Gamma}$, and the matrix $B^{-1}_{\Gamma}$,
  for Carter diagrams $D_l(a_k)$, $D_l$, where $l \leq 8$, and for $E_l(a_k)$, $E_l$, where $l \leq 7$,
  are listed, see Tables \ref{tab_partial Cartan_1}--\ref{tab_partial Cartan_3}.
  In \S\ref{sect_linkage_diagr}, all linkage systems are depicted,
  see Figs. \ref{D4a1_linkages}--\ref{D4_loctets}, \ref{Dk_al_linkages}, \ref{Dk_al_wind_rose}.
  The description of all linkage systems is presented in Theorem \ref{th_full_descr}, see \S\ref{sec_struct_sizes}.

\subsubsection{Linkage systems for Dynkin diagrams and weight systems}
 \label{sec_linkage_syst}
  The linkage system and the weight system for $E_6$ coincide,
  see Fig. \ref{27_weight_diagr_E6__2comp}(top) and Fig. \ref{27_weight_diagr_E6__2comp}(bottom).
  It becomes obvious after recognizing loctets in both diagrams. The comparative figure containing
  both the linkage systems and the weight systems together with all their loctets can be seen in Fig. \ref{27_weights}.

  Similarly, the linkage system and the weight system for $E_7$ coincide,
  see Fig. \ref{E7pure_linkage_system} and Fig. \ref{56_weight_diagr_E7}.

 \index{weight system ! - $(D_l, \overline{\omega}_1)$}
 \index{weight system ! - $(A_l, \overline{\omega}_1)$}
 \index{conjugation ! - non-conjugate conjugacy classes}
 \index{$\mathfrak{g}$, simple Lie algebra}

 We observe that sizes of components in linkage systems for Carter diagrams $E_6$,
 $E_7$ and $D_l$ (and for their covalent $E_6(a_i)$, $E_7(a_i)$ and $D_l(a_i)$, see \S\ref{sec_Dynkin_ext_2},
 Corollary \ref{corol_Estim_2}) are, respectively, $27$, $56$ and $2l$ which coincide with the dimensions of
 the smallest fundamental representations of simple Lie algebras $E_6$,  $E_7$, and $D_l$, respectively.
 Similarly, sizes of components of the Carter diagram $A_l$, where $l \geq 8$,
 are, respectively, $l+1$, $\frac{l(l+1)}{2}$: The $A$-component (resp. $D$-component)
 coincides with the weight system $(A_l, \overline{\omega}_1$) (resp. $(A_l, \overline{\omega}_2$)),
 see \cite[Figs. 1 and 5]{PSV98}.
  These facts get {\it a priori} reasoning in the following theorem:
~\\
~\\
 {\bf Theorem} (Theorem \ref{th_linkages_n_weights}){
 {\it Let $\Gamma$ be a simply-laced Dynkin diagram,  $\mathfrak{g}$ the simple Lie algebra associated with $\Gamma$.
  Every $A$-, $D$- or $E$-component of the linkage system $\mathscr{L}(\Gamma)$ coincides with
  a weight system $\mathscr{W}(\mathfrak{g})$ of one of
  fundamental representations of $\mathfrak{g}$.}

 %% Note the elements of the linkage system $\mathscr{L}(\Gamma)$ and elements of the weight system $\mathscr{W}(\mathfrak{g})$
 %%  are of different nature.

%%{\Small
\begin{remark}
  \label{rem_weight_system}
  %%  \footnotesize
  {\rm
  (i) The weight system corresponding to the highest weight $\overline{\omega}_1$ for type $D_l$
     is taken from \cite[Fig. 4]{PSV98}, see Fig. \ref{Dl_weight_system}.  The number of weights is
     $2l$.  By Theorem \ref{th_linkages_n_weights},  the linkage system has the same diagram. Compare with
     the linkage systems $\mathscr{L}(D_l(a_k))$, see Fig. \ref{Dk_al_linkages}, \ref{Dk_al_wind_rose},
     \ref{D5pure_loctets}, \ref{D6pure_loctets}.

\begin{figure}[H]
\centering
\includegraphics[scale=0.7]{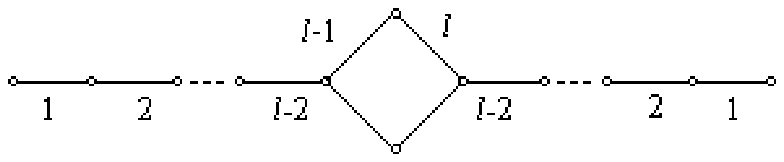}
\vspace{3mm}
\caption{\hspace{3mm}The weight system $(D_l, \overline{\omega}_1)$.}
%%%%%% The label must come after caption
\label{Dl_weight_system}
\end{figure}

 \index{class of diagrams ! - $\mathsf{C4}$}
 \index{$\gamma^{\nabla}$, linkage diagram}

   Note that the linkage system $\mathscr{L}(D_l)$, see Fig. \ref{Dlpu_linkages},
   (= weight system $(D_l, \overline{\omega}_1)$ in Fig. \ref{Dl_weight_system}),
   has exactly the same shape as the Carter diagram $D_l(a_k)$ despite the
   fact that the vertices in these diagrams are of different nature.

  (ii)  There are two conjugacy classes in $W(D_n)$ associated with the Carter diagram $A_l$, see Remark \ref{two_class_Al}.
 Nevertheless, the linkage system for $A_l$ can be constructed. By Theorem \ref{th_linkages_n_weights}
 this linkage system coincides with the weight system for $A_l$. Figures of the weight system
 $(A_l, \overline{\omega}_1)$ can be found in \cite[Fig. 1]{PSV98}. The shape of this weight
 system coincides with the Dynkin diagram $A_{l+1}$.
  } \qed
\end{remark}

 \subsubsection{The structure and sizes of $ADE$ components of linkage systems}
   \label{sec_struct_sizes}
\begin{theorem}[The aggregate description of the structure and sizes of $ADE$ components]
  \label{th_full_descr}
     For each Carter diagrams $\Gamma \in  \mathsf{C4} \coprod \mathsf{DE4} \coprod \mathsf{A}$,
     the number of linkage system components, the number of linkage diagrams in every component,
     the invariant characterizing every component (values $\mathscr{B}^{\vee}_{\Gamma}(\gamma^{\nabla})$)
     are presented in Table \ref{tab_val_Bmin1}.
 \end{theorem}

  \PerfProof The number of linkage diagrams is obtained from
  the results of stratification of root systems in \S\ref{sec_transition} and
  enumeration of loctets and
  $\beta$-unicolored linkage diagrams for every Carter diagram from Table \ref{tab_val_Bmin1}, see \S\ref{sec_enum_loctets}.
  The number of components is obtained from the shape of linkage systems.
  For the Carter diagrams of $\mathsf{C4}$ class, see Figs. \ref{D5a1_linkages}--\ref{E7a4_linkages},
  \ref{D7a1_linkages_cdef}--\ref{D7a2_linkages_comp2}, \ref{D7a1_D7a2_D7pu_loctets_comp3}.
  For the Carter diagrams from $\mathsf{DE4}$ class,
  see Figs. \ref{27_weight_diagr_E6__2comp}(top),  \ref{E7pure_linkage_system}, \ref{D5pure_loctets},
  \ref{D6pure_loctets},  \ref{D7pu_loctets_comp1}--\ref{D7a1_D7a2_D7pu_loctets_comp3}.

  For $D_l(a_k)$ and $D_l$, where $l \geq 8$, the statement is proved in \S\ref{sec_D_lg8}.
  For $D_l(a_k)$ and $D_l$, where $l < 8$, the statement is proved in \S\ref{sec_plis_Dl_Dlak}, see also \S\ref{sect_linkage_diagr}.

  For Carter diagrams $A_l$, the technique of loctets does not work.
  The number of linkage diagrams and number of components are obtained in \S\ref{sec_plis_Al} for $l < 8$
  and in \S\ref{sec_Al_lg8} for $l \geq 8$.

  For the rational number $p$ characterizing every component, see \S\ref{sec_rat_p}.
  \qed

 \index{weight system ! - for $(D_4, \overline{\omega}_1)$, $(D_4, \overline{\omega}_3)$, $(D_4, \overline{\omega}_4)$}
\begin{remark}[On linkage systems $\mathscr{L}(D_4)$ and $\mathscr{L}(D_4(a_1))$]
  \label{rem_DE_except}
{\rm
   (i) For $D_4$ and $D_4(a_1)$, all $3$ components are $D$-components,
   and for each component $\mathscr{B}^{\vee}_{\Gamma}(\gamma^{\nabla}) = 1$.
   For $D_4$, see see Fig. \ref{D4_loctets}; for $D_4(a_1)$,  see Fig. \ref{D4a1_linkages}.

   (ii) To obtain each component of $D_4$
   one can take only $4$ first coordinates for any linkage diagrams in Fig. \ref{fig_loctets},
   see Fig. \ref{D4_loctets}.
   Every component is exactly the \underline{loctet}. These components coincide with $3$
   weight systems of $3$ fundamental representations of semisimple Lie algebra $D_4$:
    $(D_4, \overline{\omega}_1)$, $(D_4, \overline{\omega}_3)$, $(D_4, \overline{\omega}_4)$,
    see \cite[Fig. 10]{PSV98}.
}
\end{remark}

\begin{figure}[H]
\centering
\includegraphics[scale=1.0]{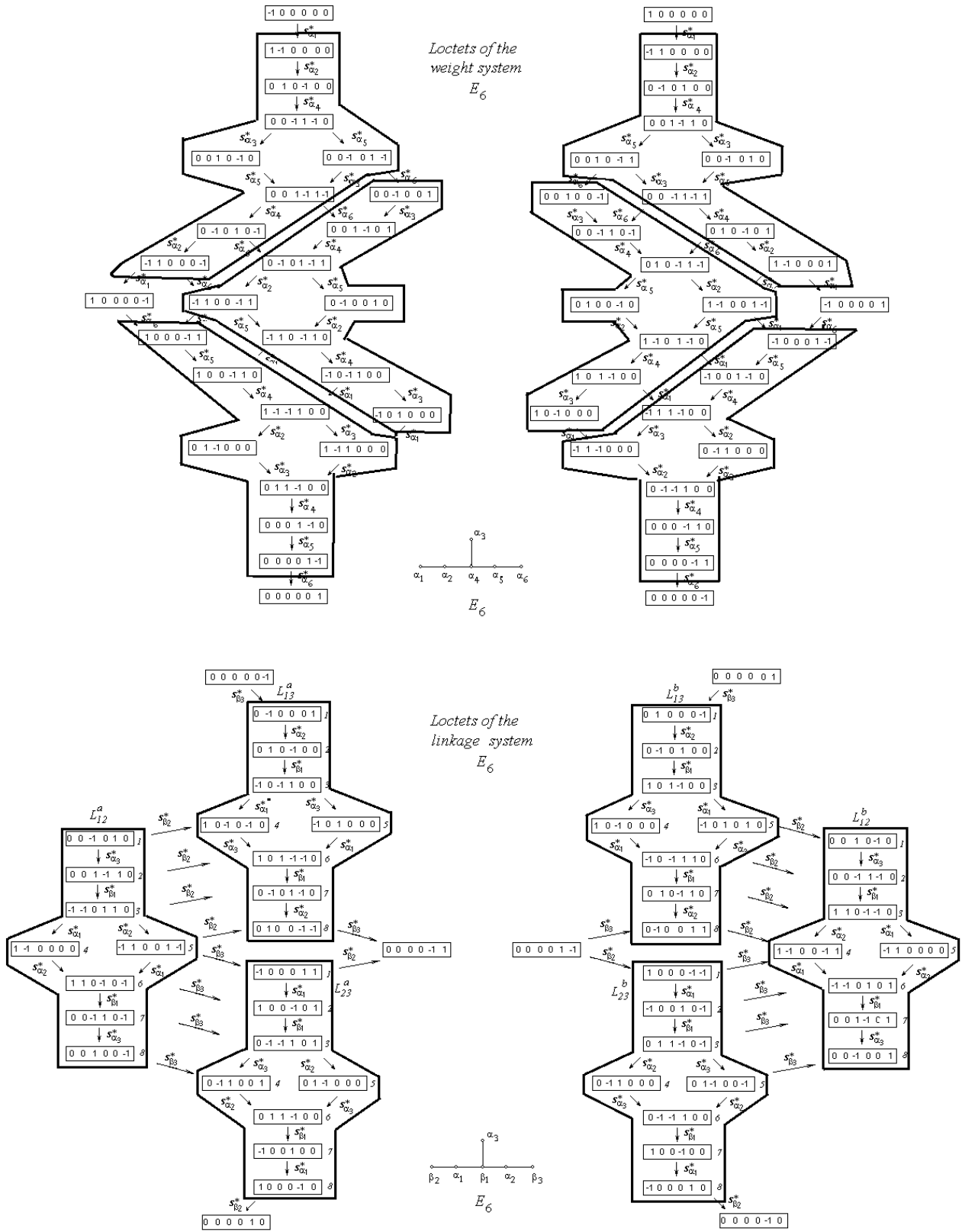}
%%\vspace{3mm}
\caption{\hspace{3mm}Loctets in the weight system and in the linkage system $E_6$}
%%%%%% The label must come after caption
\label{27_weights}
 \index{weight system ! - for $E_6$}
\end{figure}

%%%The rational number $p = \mathscr{B}^{\vee}_{\Gamma}(u^{\nabla})$
\subsubsection{The invariant of the linkage system component}
  \label{sec_rat_p}
  \index{Dynkin extension ! - $\Gamma <_D \Gamma'$}
  \index{$\Gamma <_D \Gamma'$, Dynkin extension}
  \index{linkage system component}

  The linkage system component $\mathscr{L}_{\Gamma'}(\Gamma)$,
  where $\Gamma <_D \Gamma'$ is a Dynkin extension,
  is the $W^{\vee}_S$-orbit on the set of linkage diagrams, see \S\ref{sec_partial_groups}.
  {\bf The rational number $p = \mathscr{B}^{\vee}_{\Gamma}(u^{\nabla})$
  is the invariant characterizing this linkage system component.}
  We call this component {\it the extension of the Carter diagram $\Gamma$ by p}.
  We denote this extension by $\{ \Gamma, p \}$.

  We will see that the {\bf Carter diagrams $D_5(a_1), D_6(a_1), D_6(a_2), D_7(a_1), D_7(a_2), D_5, D_6, D_7$
  have two Dynkin extensions corresponding to the $D$-component and the $E$-component of the linkage system component $\mathscr{L}(\Gamma)$.}
  For example, $D_5(a_1)$ has extension $\{ D_5(a_1), 1 \}$ containing $10$ linkages,
  and extension $\{ D_5(a_1), \frac{5}{4} \}$ containing $32$ linkages,
  see Fig. \ref{D5a1_linkages}, Table \ref{tab_val_Bmin1}.

  We will see also that {\bf Carter diagrams $A_5$, $A_6$, $A_7$ have three Dynkin extensions corresponding to $A$-, $D$- and $E$-component
  of the linkage system $\mathscr{L}(\Gamma)$.} For example, $A_7$ has extension $\{A_7, \frac{1}{8}\}$
  containing $16$ linkages, $\{A_7, \frac{3}{2}\}$ containing $56$ linkages and $\{A_7, \frac{15}{8}\}$
  containing $112$ linkages, see Fig. \ref{A7_to_E8}, Fig. \ref{A5_to_A6D6_pass}(bottom), Table \ref{tab_val_Bmin1}.

\subsection{On connection with the Carter theorem}
  \label{sec_Carter_theorem}

 Denote the set of elements $w \in W$, each of
 which corresponds to an admissible diagram, by $W_0$. The
 existence of an admissible diagram for the element $w$ means that
 $w$ can be decomposed into the product of two involutions as
 follows:
 \begin{equation}
   \label{two_invol_0}
      w = w_1{w}_2, \quad \text{ where } \quad
      w_1 = s_{\alpha_1} s_{\alpha_2} \dots s_{\alpha_k}, \quad
      w_2 = s_{\beta_1} s_{\beta_2} \dots s_{\beta_h},
 \end{equation}
 the roots $\{\alpha_i\mid i = 1,\dots,k\}$ being mutually
 orthogonal, and the roots $\{\beta_j\mid j = 1,\dots,h\}$ being also
 mutually orthogonal.

 Thus, $W_0$ is the subset of elements $w \in W$ that can be
 decomposed as \eqref{two_invol_0}. It turns out that $W_0 = W$. This
 fact is one of the main results of \cite[Theorem C]{Ca72}. I call
 this result {\bf the Carter theorem}.  The proof of the Carter
 theorem is based on the classification of conjugacy classes.
 (By \cite[p. G-21]{Ca70}, N.~Burgoyne carried out a check of classification for
 $E_7$ and $E_8$ with a computer aid).

 \index{Coxeter groups}
 \index{Weyl group}
 I would like to quote Carter  from \cite[p. 4]{Ca00}:

{\it \lq\lq One remarkable feature of the theory of Coxeter groups
and Iwahori-Hecke algebras is the number of key properties which at
present have no uniform proof and can only be proved in a
case-by-case manner. We mention three examples of such properties.
In the first place every element of a Weyl group is a product of two
involutions. This is a key property in Carter's description of the
conjugacy classes in the Weyl group \cite{Ca72}. Secondly, every
element of a finite Coxeter group can be transformed into an element
of minimal length in its conjugacy class by a sequence of
conjugations by simple reflections such that the length does not
increase at any stage. This property is basic to the Geck-Pfeiffer
approach to the conjugacy classes of Coxeter groups. Thirdly, there
are basic properties of Lusztig's a-function which at present have
only case-by-case proofs. The a-function is an important invariant
of irreducible characters of Coxeter groups. One cannot be satisfied
with the theory of Coxeter groups until such case-by-case proofs are
replaced by uniform proofs of a conceptual nature.\rq\rq}

 In \cite{St11}, using the classification of linkage system, we give the proof of the Carter theorem,
 which does not use the classification of conjugacy classes.

%% file: 3transition.tex
~\\
\section{\sc\bf Theorem on covalent Carter diagrams}
  \label{sec_transition}

  \subsection{The size of linkage system $\mathscr{L}(\Gamma)$}
  The following proposition together with corollaries constitute the main result of this section:
  \begin{proposition}
     \label{prop_Estim}
     For the Carter diagram $\Gamma$, we have $\mid \mathscr{L}(\Gamma) \mid \leq \mathscr{E}$, where
     $\mathscr{E}$ is given by Table \ref{tab_val_Estim}.
  \end{proposition}
   In sections \S\ref{sec_stratif_D}-\S\ref{sec_stratif_A}, we check that for each Carter diagram $\Gamma$,
  there exists at least $\mathscr{E}$ linkage diagrams in the linkage system $\mathscr{L}(\Gamma)$.
  Together with Proposition \ref{prop_Estim}, we get the following
  \begin{corollary}
    \label{corol_Estim}
     For the Carter diagram $\Gamma$, we have $\mid \mathscr{L}(\Gamma) \mid = \mathscr{E}$
     (where $\mathscr{E}$ is given by Table \ref{tab_val_Estim}).
  \end{corollary}
 ~\\
  \begin{table}[H]
  \Small
  \centering
  \renewcommand{\arraystretch}{1.5}
  \begin{tabular} {||c|c||c|c||}
  \hline \hline
      The Carter          &   $\mathscr{E}$ & The Carter          & $\mathscr{E}$ \\
        diagram $\Gamma$  &                 & diagram $\Gamma$    & \\
    \hline \hline
       $E_6$, ~$E_6(a_1)$, ~$E_6(a_2)$ &  $54$  & $A_3$    & $14$  \\
    \hline
       $E_7$, $E_7(a_1)$, $E_7(a_2)$, $E_7(a_3)$, $E_7(a_4)$ & $56$ & $A_4$    & $30$  \\
   \hline
       $D_4$, ~$D_4(a_1)$ &  $24$ & $A_5$    & $62$ \\
    \hline
       $D_5$, ~$D_5(a_1)$ &  $42$ & $A_6$    & $126$ \\
    \hline
       $D_6$, ~$D_6(a_1)$, ~$D_6(a_2)$  &  $76$ & $A_7$    & $184$\\
    \hline
       $D_7$, ~$D_7(a_1)$, ~$D_7(a_2)$ & $142$ & $A_8$    & $90$ \\
    \hline
       $D_l$, ~$D_l(a_k)$, $l > 7$ & $2l$ & $A_l$, $l > 8$  & $(l+1)(l+2)$\\
   \hline  \hline
  \end{tabular}
  \vspace{2mm}
  \caption{\hspace{3mm}$\mid \mathscr{L}(\Gamma) \mid \leq \mathscr{E}$}
  \label{tab_val_Estim}
  \end{table}
 ~\\

   By Proposition \ref{prop_Estim} and calculations of this section we get
  \begin{corollary}
    \label{corol_Estim_2}
     For the Carter diagram $E_i(a_k)$ and $D_i(a_k)$,  we have
   \begin{equation*}
     \begin{split}
      & \mid \mathscr{L}(E_6(a_k)) \mid ~=~ \mid \mathscr{L}(E_6) \mid ~=~ 54, \text{ where } k = 1,2, \\
      & \mid \mathscr{L}(E_7(a_k)) \mid ~=~ \mid \mathscr{L}(E_7) \mid ~=~ 56, \text{ where } k = 1,2,3,4, \\
      & \mid \mathscr{L}(D_l(a_k)) \mid ~=~ \mid \mathscr{L}(D_l) \mid,
        \text{ where } l \geq 4 \text{ and } k = 1 \leq k \leq \big [ \tfrac{l-2}{2} \big ].
     \end{split}
   \end{equation*}
  \end{corollary}

   For each Carter diagram $\Gamma \in  \mathsf{C4} \coprod \mathsf{DE4} \coprod \mathsf{A}$,
   the values $\mathscr{B}^{\vee}_{\Gamma}(\gamma^{\nabla})$, the number of components, and number of
   linkage diagrams are presented in Table \ref{tab_val_Bmin1}, see Theorem \ref{th_full_descr}.

  \subsection{Transition between a partial root system $\mathcal{P}$ and a root system $\mathcal{R}$}
    \label{sec_part_root_and_comp}

 \begin{theorem}[On bijection of covalent root systems]
  \label{prop_bijection_root_syst}
  {\rm(i)}
  {\it Consider a covalent pair of Carter diagrams $\{\widetilde\Gamma, \Gamma\}$ out of the following list:
   \begin{equation}
    \label{eq_covalent_pairs}
    %% \footnotesize
     \begin{split}
      & \{D_l(a_k), D_l\} \text{ for } l \geq 4 \text{ and } 1 \leq k \leq \big [\tfrac{l-2}{2} \big ], \\
      & \{E_6(a_k), E_6\} \text{ for }  k = 1,2, \\
      & \{E_7(a_k), E_7\} \text{ for }  k = 1,2,3,4.
     \end{split}
   \end{equation}
  Let $\mathcal{R}$ (resp. $\mathcal{P}$) be the root system (resp. partial root system) corresponding to
  $\Gamma$ (resp. $\widetilde\Gamma$).
   There exists the map $M : \mathcal{P} \longmapsto \mathcal{R}$
   given by Table \ref{tab_part_root_syst} is the invertible linear map.}

 {\rm(ii)} Let $\widetilde\tau$ be a root in $\varPhi$, where $\varPhi$ is the primary root system.
  Then $\widetilde\tau$ is a root of the partial root system $\mathcal{P}$
  if and only if $\widetilde\tau$ is a root of the  root system $\mathcal{R}$.
 \end{theorem}

  \PerfProof
 (i)  For all pairs \{$\mathcal{P}$, $\mathcal{R}$\} from Table \ref{tab_part_root_syst},
 it is easy to check by \eqref{eq_Kac} that vectors images under the map $M$ are roots.
 For example, for case (2), the value $\mathscr{B}(\beta_{k+1})$ is as follows:
 \begin{equation*}
   \Small
   \begin{split}
    \mathscr{B}(\beta_{k+1}) =
             & \mathscr{B}(\tau_1)  + 4\sum\limits_{i=2}^k\mathscr{B}(\tau_i) + \mathscr{B}(\tau_{k+1}) +
                \mathscr{B}(\overline\tau_{k+1}) + \\
             & 4(\tau_1, \tau_2) + 8\sum\limits_{i=2}^{k-1}(\tau_{k-1}, \tau_k) + 4(\tau_k, \tau_{k+1}) + 4(\tau_k, \overline\tau_{k+1}) = \\
             & 2 + 8(k-1) + 2 + 2 - 4  - 8(k-2) - 4 - 4 = 2 + 8(k-1) - 8(k-1) = 2.
   \end{split}
 \end{equation*}
  Thus, $\beta_{k+1}$ is the root.  It only remains to check some inner products.
~\\

  \begin{table}[H]
  %%\footnotesize
  \centering
  \renewcommand{\arraystretch}{1.3}
  \begin{tabular} {|c|c|c|c|}
    \hline
    \hline
      & Partial root system & Root system  & Linear maps  \\
      &     $\mathcal{P}$    &  $\mathcal{R}$  &
            $M : \mathcal{P} \longmapsto \mathcal{R}$ and  $M^{-1} : \mathcal{R} \longmapsto \mathcal{P}$ \\
     &  &  & {\Small ($M$ (resp. $M^{-1}$) acts as the identity on all} \\
     &  &  & {\Small  not mentioned $\alpha_i$ and $\beta_j$)} \\
    \hline
    \hline
     1 &
     $\begin{array}{c} \addlinespace
      \includegraphics[scale=0.9]{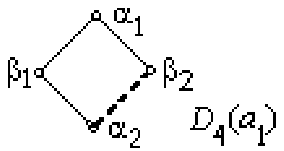}  \end{array}$  &
     $\begin{array}{c}  \addlinespace
       \includegraphics[scale=0.9]{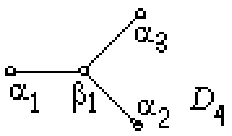}  \end{array}$
           & \footnotesize
         $\begin{array}{cl}
         \addlinespace
         M & \hspace{-3mm}:=
         \begin{cases}
           \alpha_3  = -(\alpha_1 + \beta_1  + \beta_2) %%% \footnotemark[1]
         \end{cases}
         \\  \addlinespace
         M^{-1} &  \hspace{-3mm}:=
          \begin{cases}
           \beta_2 = -(\alpha_1 + \alpha_3 + \beta_1)
          \end{cases}
       \end{array}$ \\
    \hline
    2 &
     $\begin{array}{c} \\ \addlinespace
      \includegraphics[scale=0.9]{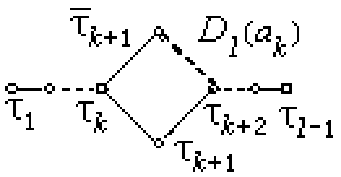}  \end{array}$  &
     $\begin{array}{c} \\ \addlinespace
       \includegraphics[scale=0.9]{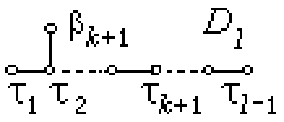}  \end{array}$
        & \footnotesize
     $\begin{array}{cl}
       \addlinespace
        M & \hspace{-3mm}:=
          \begin{cases}
             \beta_{k+1} =
              -(\tau_1 + 2\sum\limits_{i=2}^k\tau_i +  \tau_{k+1} + \overline{\tau}_{k+1}), \\
               \qquad\qquad\qquad\qquad\qquad\qquad\qquad k \geq 2, \\
              \beta_2 = -(\tau_1 + \tau_2 + \overline{\tau}_2).
          \end{cases} \\ \addlinespace
        M^{-1} & \hspace{-3mm}:=
          \begin{cases}
             \overline{\tau}_{k+1} =
              -(\tau_1 + 2\sum\limits_{i=2}^k\tau_i +  \tau_{k+1} + \beta_{k+1}), \\
               \qquad\qquad\qquad\qquad\qquad\qquad\qquad k \geq 2, \\
              \overline{\tau}_2 = -(\tau_1 + \tau_2 + \beta_2).
          \end{cases} \\ \addlinespace
     \end{array}$ \\
    \hline
    \hline
     3 &  $\begin{array}{c} \addlinespace
      \includegraphics[scale=0.8]{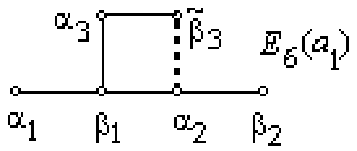}  \end{array}$  &
     $\begin{array}{c} \addlinespace
       \includegraphics[scale=0.8]{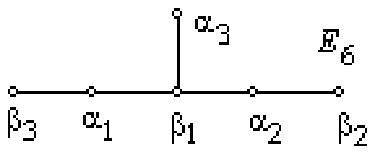}  \end{array}$
            & \footnotesize
         $\begin{array}{cl}
           \addlinespace
         M & \hspace{-3mm}:=
         \begin{cases}
           \beta_3 = -(\alpha_1 + \beta_1 + \alpha_3 + \widetilde\beta_3)
         \end{cases}
         \\  \addlinespace
         M^{-1} &  \hspace{-3mm}:=
          \begin{cases}
            \widetilde\beta_3 = -(\beta_3 + \alpha_1 + \beta_1 + \alpha_3)
          \end{cases}
       \end{array}$ \\
    \hline
     4 &  $\begin{array}{c} \\ \addlinespace
      \includegraphics[scale=0.8]{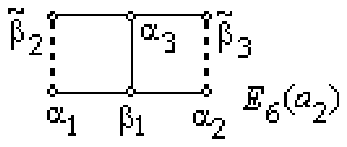}  \end{array}$  &
     $\begin{array}{c} \\ \addlinespace
       \includegraphics[scale=0.8]{E6_root_syst.eps}  \end{array}$
        & \footnotesize
         $\begin{array}{cl}
         \addlinespace
         M & \hspace{-3mm}:=
         \begin{cases}
           \beta_2 = -(\alpha_2 + \beta_1 + \alpha_3 + \widetilde\beta_2) \\
           \beta_3 = -(\alpha_1 + \beta_1 + \alpha_3 + \widetilde\beta_3)
         \end{cases}
         \\  \addlinespace
         M^{-1} &  \hspace{-3mm}:=
          \begin{cases}
            \widetilde\beta_2 = -(\beta_2 + \alpha_2 + \beta_1 + \alpha_3) \\
            \widetilde\beta_3 = -(\beta_3 + \alpha_1 + \beta_1 + \alpha_3)
          \end{cases}
         \end{array}$ \\
    \hline
    \hline
     5 &  $\begin{array}{c}
      \includegraphics[scale=0.8]{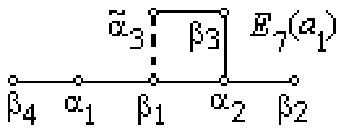}  \end{array}$  &
     $\begin{array}{c}
       \includegraphics[scale=1.0]{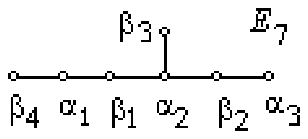}  \end{array}$
        & \footnotesize
         $\begin{array}{cl}
         M & \hspace{-3mm}:=
         \begin{cases}
           \alpha_3 = -(\beta_2 + \alpha_2 + \beta_3 + \widetilde\alpha_3)
         \end{cases}
         \\  \addlinespace
         M^{-1} &  \hspace{-3mm}:=
          \begin{cases}
            \widetilde\alpha_3 = -(\alpha_3 + \beta_2 + \alpha_2 + \beta_3)
          \end{cases}
         \end{array}$
         \\
      \hline
     6 &   $\begin{array}{c}
      \includegraphics[scale=0.9]{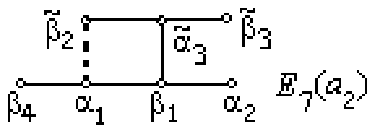}  \end{array}$ &
     $\begin{array}{c}
       \includegraphics[scale=1.0]{E7_root_syst_3.eps}  \end{array}$
        & \footnotesize
       $\begin{array}{cl}
       \addlinespace
        M & \hspace{-3mm}:=
         \begin{cases}
         & \hspace{-2mm}\beta_2 = -( \alpha_2 + \beta_1 + \widetilde\alpha_3 + \widetilde\beta_2), \\
         & \hspace{-2mm}\beta_3 = 2\widetilde\alpha_3 + \beta_1 + \widetilde\beta_3 + \widetilde\beta_2,  \\
         & \hspace{-2mm} \alpha_3 = -\widetilde\beta_3.
          \end{cases}
          \\  \addlinespace
         M^{-1} & \hspace{-3mm}:=
         \begin{cases}
         & \hspace{-2mm}\widetilde\alpha_3 = \alpha_3 + \beta_2 + \alpha_2 + \beta_3, \\
         & \hspace{-2mm}\widetilde\beta_2 =  \alpha_3 + 2\beta_2 + 2\alpha_2 + \beta_1 + \beta_3, \\
         & \hspace{-2mm}\widetilde\beta_3 = -\alpha_3.
          \end{cases}
         \end{array}$
          \\ %% \addlinespace
      \hline
     7 &  $\begin{array}{c} \\
      \includegraphics[scale=0.9]{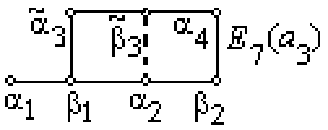}  \end{array}$  &
     $\begin{array}{c} \\
       \includegraphics[scale=1.0]{E7_root_syst_3.eps}  \end{array}$
         & \footnotesize
        $\begin{array}{cl}
        \addlinespace
        M& \hspace{-3mm}:=
         \begin{cases}
         & \hspace{-2mm}\beta_4 = \alpha_4 + \widetilde\beta_3 - \alpha_2 - \beta_1 - \alpha_1, \\
         & \hspace{-2mm}\alpha_3 = \widetilde\alpha_3 + \widetilde\beta_3 - \alpha_2 - \beta_2, \\
         & \hspace{-2mm}\beta_3 = -\widetilde\beta_3.
          \end{cases}
          \\ \addlinespace
         M^{-1} & \hspace{-3mm}:=
         \begin{cases}
         & \hspace{-2mm}\alpha_4 = \beta_4 + \alpha_1  + \beta_1 + \alpha_2 + \beta_3, \\
         & \hspace{-2mm}\widetilde\alpha_3 = \alpha_3 + \beta_2 + \alpha_2 + \beta_3,  \\
         & \hspace{-2mm}\widetilde\beta_3 = -\beta_3.
          \end{cases}
         \end{array}$
        \\
      \hline
     8 &  $\begin{array}{c} \\ \addlinespace
      \includegraphics[scale=0.9]{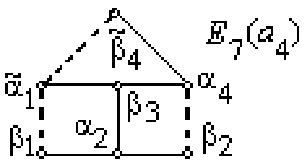} \end{array}$   &
     $\begin{array}{c} \\ \addlinespace
       \includegraphics[scale=1.0]{E7_root_syst_3.eps}   \end{array}$
         & \footnotesize
        $\begin{array}{cl} \addlinespace
        M & \hspace{-3mm}:=
         \begin{cases}
         & \hspace{-2mm}\beta_4 = -(2\alpha_4 + \beta_3 + \widetilde\beta_4 - \beta_2), \\
         & \hspace{-2mm}\alpha_3 = -(\widetilde\alpha_1 + \beta_3 + \alpha_2 + \beta_2), \\
         & \hspace{-2mm}\alpha_1 = -(\beta_1 + \alpha_2 + \beta_2 - \alpha_4 - \widetilde\beta_4).
          \end{cases}
          \\ \addlinespace
         M^{-1} & \hspace{-3mm}:=
         \begin{cases}
         & \hspace{-2mm}\widetilde\alpha_1 =  -(\alpha_3 + \beta_2 + \alpha_2 + \beta_3), \\
         & \hspace{-2mm}\alpha_4 =  -(\beta_4 + \alpha_1 + \beta_1 + \alpha_2 + \beta_3), \\
         & \hspace{-2mm}\widetilde\beta_4 =  \beta_4 + 2\alpha_1 + 2\beta_1 + 2\alpha_2 + \beta_2 + \beta_3.
          \end{cases}
          \\ \addlinespace
         \end{array}$
        \\
      \hline
      \hline
\end{tabular}
  \vspace{2mm}
  \caption{\hspace{3mm}The bijective linear maps between root systems and partial root systems.}
   %% Here $M$ (resp. $M^{-1}$) acts as identity on all not mentioned $\alpha_i$ and $\beta_j$}
  \label{tab_part_root_syst}
  \end{table}

  \underline{$(1)$ {\it Pair } \{$D_4(a_1)$, $D_4$\}.} Let us check relations for $\alpha_3$:
 \begin{equation*}
  \Small
   \begin{array}{cl}
    (\alpha_3, \alpha_1) ~= & -(\mathscr{B}(\alpha_1) + (\beta_1, \alpha_1) + (\beta_2, \alpha_1)) = -(2 - 2) = 0, \\
    (\alpha_3, \alpha_2) ~= & -(\alpha_1, \alpha_2) + (\beta_1, \alpha_2) + (\beta_2, \alpha_2) = -1 + 1 = 0, \\
    (\alpha_3, \beta_1) ~= & -(\alpha_1, \beta_1) + (\beta_1, \beta_1) + (\beta_2, \beta_1) = 1 - 2 = -1. \\
   \end{array}
 \end{equation*}

  \underline{$(2)$ {\it Pair } \{$D_l(a_k)$, $D_l$\}.}
  The root $\beta_{k+1}$ is connected with $S = \{\tau_1, \dots, \tau_l \}$ at point $\tau_2$:

  \begin{equation*}
   \Small
   \begin{split}
   \text{ For } & k \geq 2, \\
      & \beta_{k+1} \perp \tau_1, \text{ since } (\tau_1, \tau_1 + 2\tau_2) = 0, \\
      & (\beta_{k+1}, \tau_2) = -((\tau_2, \tau_1) + 2\mathscr{B}(\tau_2) + 2(\tau_2,\tau_3)) = -(-1 + 4 -2) = -1, \\
      & \beta_{k+1} \perp \tau_i \text{ for } 3 \leq i \leq k-1, \text{ since } (\tau_i, \tau_{i-1} + \tau_i + \tau_{i+1}) = 0, \\
      & \beta_{k+1} \perp \tau_k, \text{ since } (\tau_k, 2\tau_{k-1} + 2\tau_k + \tau_{k+1} +  \overline\tau_{k+1}) = 0, \\
      & \beta_{k+1} \perp \tau_{k+1}, \text{ since } (\tau_{k+1}, 2\tau_k + \tau_{k+1}) = 0, \\
      & \beta_{k+1} \perp \tau_{k+2}, \text{ since } (\tau_{k+2}, \tau_{k+1} +  \overline\tau_{k+1}) = 0, \\
      & \beta_{k+1} \perp \tau_i \text{ for }  i > k+2. \\
   \text{ For } & k = 1, \\
      & \beta_2 \perp \tau_1, \text{ since } (\tau_1, \tau_1 + \tau_2 + \overline\tau_2) = 2 - 1 - 1 = 0, \\
      & (\beta_2, \tau_2) = -((\tau_2, \tau_1) + \mathscr{B}(\tau_2) + (\tau_2,\overline\tau_2)) = -(-1 + 2)  = -1, \\
      & \beta_2 \perp \tau_3, \text{ since } (\tau_2, \tau_3) + (\overline\tau_2, \tau_3) = -1 + 1  = 0, \\
      & \beta_2 \perp \tau_i \text{ for }  i > 3. \\
   \end{split}
 \end{equation*}
 By Remark \ref{rem_max_root} the root $\beta_{k+1}$ is uniquely determined, so the root
 $\overline\tau_{k+1}$ is also uniquely determined by $\beta_{k+1}$.

 \underline{$(3)$ {\it Pair } \{$E_6(a_1)$, $E_6$\}.} For the root $\beta_3 = -(\alpha_1 + \beta_1 + \alpha_3 + \widetilde\beta_3)$, we have
 $\beta_3 \perp  \beta_1, \alpha_2, \alpha_3, \beta_2$ and
 \begin{equation*}
     \Small \begin{array}{c}
   (\beta_3, \alpha_1) = -(\mathscr{B}(\alpha_1) +(\alpha_1, \beta_1)) = -1.
    \end{array}
 \end{equation*}

 \underline{$(4)$ {\it Pair } \{$E_6(a_2)$, $E_6$\}.} For roots $\beta_2 = -(\alpha_2 + \beta_1 + \alpha_3 + \widetilde\beta_2)$ and
 $\beta_3 = -(\alpha_1 + \beta_1 + \alpha_3 + \widetilde\beta_3)$, the following relations hold:
  \begin{equation*}
    \Small
     \begin{split}
     & \beta_2 \perp \alpha_1, \alpha_3, \beta_1, \quad \beta_3 \perp \alpha_2, \alpha_3, \beta_1, \quad
       (\beta_2, \alpha_2) = -1, \quad (\beta_3, \alpha_1) = -1,  \\
     & (\beta_2, \beta_3) = \mathscr{B}(\beta_1) + \mathscr{B}(\alpha_3) +
        (\beta_1, \alpha_1 + 2\alpha_3 + \alpha_2) +
        (\widetilde\beta_3, \alpha_2 + \alpha_3) + (\widetilde\beta_2, \alpha_1 + \alpha_3) =  2 + 2 - 4 + 0 + 0 = 0.
     \end{split}
 \end{equation*}

 \underline{$(5)$ {\it Pair } \{$E_7(a_1)$, $E_7$\}.} Here, $\alpha_3 = -(\beta_2 + \alpha_2 + \beta_3 + \widetilde\alpha_3)$.
 We have $\alpha_3 \perp \alpha_1, \beta_3, \beta_1, \alpha_2, \beta_4$ and
 \begin{equation*}
    \Small
  \begin{array}{c}
     (\alpha_3, \beta_2) = -((\mathscr{B}(\beta_2) + (\beta_2, \alpha_2)) = -1.
         \end{array}
 \end{equation*}

\underline{$(6)$ {\it Pair } \{$E_7(a_2)$, $E_7$\}.}
 For roots $\beta_2 = -( \alpha_2 + \beta_1 + \widetilde\alpha_3 + \widetilde\beta_2)$ and
 $\beta_3 = 2\widetilde\alpha_3 + \beta_1 + \widetilde\beta_3 + \widetilde\beta_2$, we have
 $\beta_2 \perp \beta_1, \alpha_1, \beta_4$ and $\beta_3 \perp \beta_1,\alpha_1, \beta_4$.
 Additionally, the following relations hold:
 \begin{equation*}
  \Small
   \begin{split}
      (\beta_2, \alpha_3) = &(\widetilde\beta_3, \widetilde\alpha_3) = -1, \\
      (\beta_2, \alpha_2) = &-((\alpha_2, \alpha_2) + (\alpha_2, \beta_1)) = -(2 - 1) = -1, \\
      (\beta_2, \beta_3) = &-(\alpha_2 + \beta_1 + \widetilde\alpha_3 + \widetilde\beta_2,
                   2\widetilde\alpha_3 + \beta_1 + \widetilde\beta_3 + \widetilde\beta_2) = \\
                      &-\big ((\alpha_2, \beta_1) + (\beta_1, \beta_1) + 2(\beta_1, \widetilde\alpha_3) +
                       2(\widetilde\alpha_3, \widetilde\alpha_3) + (\alpha_3, \beta_1) + (\alpha_3, \widetilde\beta_3)  \\
                      & + (\alpha_3, \widetilde\beta_2) + (\widetilde\beta_2, \widetilde\beta_2) + 2(\alpha_3, \widetilde\beta_2) \big ) =
                     1 - 2 + 2 - 4 + 1 + 1 + 1  - 2 + 2 = 0. \\
      (\beta_3, \alpha_3) = &-(\widetilde\beta_3,\widetilde\beta_3) - 2(\widetilde\alpha_3, \widetilde\beta_3) = 0. \\
   \end{split}
 \end{equation*}

 \underline{$(7)$ {\it Pair } \{$E_7(a_3)$, $E_7$\}.}
 Here, $\beta_4 \perp \beta_3, \alpha_2, \alpha_3, \beta_1, \beta_2$ and $(\beta_4, \alpha_1) = -1$,
 For example,
 \begin{equation*}
  \Small
   \begin{split}
     & (\beta_4, \beta_1) = (\alpha_4 + \widetilde\beta_3 - \alpha_2 - \beta_1 - \alpha_1, \beta_1) =
        -(\alpha_2 + \beta_1 + \alpha_1, \beta_1) =  0, \\
     & (\beta_4, \beta_2) = (\alpha_4 + \widetilde\beta_3 - \alpha_2 - \beta_1 - \alpha_1, \beta_2) =
        (\alpha_4 - \alpha_2, \beta_2) =  0, \\
     & (\beta_4, \alpha_2) = (\alpha_4 + \widetilde\beta_3 - \alpha_2 - \beta_1 - \alpha_1, \alpha_2) =
        (-\beta_1 -\alpha_2 + \widetilde\beta_3, \alpha_2) =  0. \\
   \end{split}
 \end{equation*}
 Further, $\alpha_3 \perp \beta_3, \alpha_2, \beta_1, \alpha_1$ and $(\alpha_3, \beta_2) = -1$.
 For example,
 \begin{equation*}
  \Small
   \begin{split}
     & (\alpha_3, \alpha_2) = (\widetilde\alpha_3 + \widetilde\beta_3 - \alpha_2 - \beta_2, \alpha_2) =
           -(\alpha_2, \alpha_2) -(\beta_2, \alpha_2) + (\widetilde\beta_3, \alpha_2) = 0, \\
     & (\alpha_3, \beta_1) = (\widetilde\alpha_3 + \widetilde\beta_3 - \alpha_2 - \beta_2, \beta_1) =
          (\beta_1, \widetilde\alpha_3) - (\beta_1, \alpha_2) = 0.
   \end{split}
 \end{equation*}
 Similarly, $\beta_3 \perp \alpha_2, \alpha_1, \beta_1$ and $(\beta_3, \alpha_2) = -1$.

 \underline{$(8)$ {\it Pair } \{$E_7(a_4)$, $E_7$\}.}
 We have $\beta_4 \perp  ~\beta_1, \alpha_2, \beta_3, \beta_2$.
 The remaining cases for $\beta_4$:
\begin{equation*}
  \Small
   \begin{split}
     \qquad
      (\beta_4, \alpha_1) = & ~(2\alpha_4, -\alpha_4 + \beta_2 - \widetilde\beta_4) + (\beta_3, \alpha_2 - \alpha_4) +
        (\widetilde\beta_4, -\widetilde\beta_4 -\alpha_4) -(\beta_2, \beta_2 + \alpha_2 - \alpha_4) = \\
         & 0 + 0 -1 + 0 = -1,  \\
      (\beta_4, \alpha_3) = & ~(2\alpha_4, \beta_3 + \beta_2) + (\beta_3, \beta_3 + \widetilde\alpha_1 + \alpha_2) +
          (\widetilde\beta_4, \widetilde\alpha_1) + (-\beta_2, \beta_2 + \alpha_2) = 0 + 0 + 1 -1 = 0. \\
  \end{split}
 \end{equation*}
 Further, $\alpha_3 \perp ~\alpha_2, \beta_1, \beta_3$ and remaining cases for $\alpha_3$ are
 as follows:
 \begin{equation*}
  \Small
   \begin{split}
      (\alpha_3, \beta_2) = & ~-(\beta_2, \beta_2) - (\beta_2,  \widetilde\alpha_1) = -2 +1 = -1, \\
      (\alpha_3, \alpha_1) = & ~(\widetilde\alpha_1, \beta_1 - \widetilde\beta_4) + (\beta_3, \alpha_2 - \alpha_4)
         + (\alpha_2, \beta_1 + \beta_2 + \alpha_2) + (\beta_2, \beta_2 + \alpha_2 - \alpha_4) = 0. \\
   \end{split}
 \end{equation*}
~\\
 Finally, $\alpha_1 \perp ~\alpha_2, \beta_2, \beta_3$ and
       $(\alpha_1, \beta_1) = ~-(\beta_1, \beta_1) -(\beta_1, \alpha_2) = -1.$
~\\

  \index{$S = \{\tau_1,\dots,\tau_l\}$, $\Gamma$-associated root subset}
  (ii) Let $\widetilde{S} = \{\widetilde\tau_1,\dots,\widetilde\tau_l\}$ (resp. $S = \{\tau_1,\dots,\tau_l\}$)
  be the $\widetilde\Gamma$-associated (resp. $\Gamma$-associated) root subset
  generating the partial root system $\mathcal{P}$ (resp. root system $\mathcal{R}$).
  Let $M$ be matrix given by Table \ref{tab_part_root_syst}, $M$ acts as follows:
\begin{equation}
  \label{eq_matr_N_1}
    M\widetilde\tau_i = \tau_i.
\end{equation}
   Let $\widetilde\tau$ be a certain root in $\mathcal{P}$.  Then
 \begin{equation*}
     \widetilde\tau = \sum\limits_{i = 0}^l t_i\widetilde\tau_i,
 \end{equation*}
  where $t_i$ are some integer coefficients.
  Consider the vector $\tau := M\widetilde\tau$. We have
 \begin{equation*}
     \tau = M\widetilde\tau = \sum\limits_{i = 0}^l t_i{M}\widetilde\tau_i = \sum\limits_{i = 0}^l t_i\tau_i.
 \end{equation*}
  Thus, $\tau$ belongs to the integer lattice spanned by $S$.  Since $M^{-1}$ is the integer matrix then
  the vector $\widetilde\tau = M^{-1}\tau$ also belongs to the integer lattice spanned by $S$.
  In addition, $\widetilde\tau \in \mathcal{P} \subset \varPhi$, i.e., $\widetilde\tau$ is the root in $\varPhi$.
  Therefore, $\widetilde\tau$ is the root in the root system $\mathcal{R}$.
  The converse is similar.

  Thus, linear transformations $M : \mathcal{P} \longmapsto \mathcal{R}$ and $M^{-1} : \mathcal{R} \longmapsto \mathcal{P}$
  map roots of $\mathcal{P}$ onto roots of $\mathcal{R}$ and realize the bijection of $\mathcal{P}$ and $\mathcal{R}$.
 \qed
~\\

The matrix $M$ maps is called the {\it transition matrix} from $\widetilde{S}$ onto $S$.

 \begin{corollary}
  \label{corol_biject}
 {\rm (i)} For any pair \eqref{eq_covalent_pairs},
 the root system $\varPhi(\Gamma')$ contains the root subsystem $\mathcal{R}$ if and only if
 $\varPhi(\Gamma')$ contains the partial root subsystem $\mathcal{P}$.

 {\rm (ii)} For any pair \eqref{eq_covalent_pairs},
 \begin{equation}
    \label{eq_number_of_roots}
     \mid \mathcal{R} \mid = \mid \mathcal{P} \mid.
 \end{equation}
 \qed
 \end{corollary}

  \begin{corollary}
     \label{cor_exist_Dynkin_ext}
      For any Carter diagram out of the list \eqref{eq_covalent_pairs},
      there exists the Dynkin extension.
  \end{corollary}
  \PerfProof For Dynkin diagrams $\Gamma = A_l, D_l, E_l$, the Dynkin extensions exist since
  $\varPhi(E_l) \subset \varPhi(E_{l+1})$ for $l < 8$, $\varPhi(A_l) \subset \varPhi(A_{l+1})$
  and $\varPhi(D_l) \subset \varPhi(D_{l+1})$.
  For non-Dynkin Carter diagrams $\widetilde\Gamma = E_l(a_k), D_l(a_k)$,
  the Dynkin extensions exist by Corollary \ref{corol_biject}(i):
  Since $\varPhi(\Gamma) \subset \varPhi(\Gamma')$ then $\varPhi(\widetilde\Gamma) \subset \varPhi(\Gamma')$.
  \qed

\begin{figure}[H]
\centering
\includegraphics[scale=0.7]{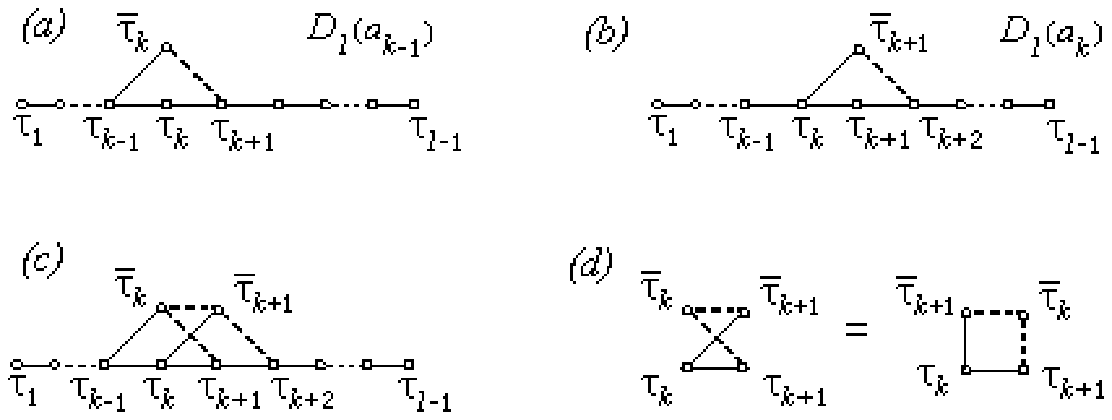}
\caption{ \hspace{3mm}Relation between $\overline\tau_k$ and
$\overline\tau_{k+1}$}
%%%%%% The label must come after caption
\label{Dlak_tk_tk1}
\end{figure}

 \begin{corollary}
  {\rm (i)}
   Let $S = \{\tau_1, \dots, \tau_{l-1}\}$ be an $A_{l-1}$-associated subset of roots in $\varPhi(D_l)$.
   There is the only root $\overline\tau_k$ such $S \cup \overline\tau_k$ is the $D_l(a_{k-1})$-associated subset,
   see Fig. \ref{Dlak_tk_tk1}$(a)$.

  {\rm (ii)} If $\overline\tau_k$ (resp. $\overline\tau_{k+1}$) is a root such that the subset
  $S \cup \overline\tau_k$ (resp. $S \cup \overline\tau_{k+1}$) constitutes the $D_l(a_{k-1})$-associated
  (resp. $D_l(a_k)$-associated) subset, then
 \begin{equation}
    \label{eq_rel_tk_tk1}
         \tau_k + \tau_{k+1} + \overline\tau_{k+1} - \overline\tau_k = 0,
 \end{equation}
  see Fig. \ref{Dlak_tk_tk1}$(d)$.
 \end{corollary}
\PerfProof
 \rm (i) By Remark \ref{rem_max_root}, the maximal root $\mu_{max}$ is uniquely determined.
 The vector $\beta_k$ defined in Table \ref{tab_part_root_syst}(case $(2)$) coincides with
 $\mu_{max}$, i.e., there exists only one root $\beta_k$ connected with $\tau_2$ and orthogonal
 to roots $\{\tau_1,\tau_3, \dots, \tau_l\}$. Thus, the vector
 \begin{equation}
  \label{eq_iniq_def_2}
     \overline\tau_k = -(\tau_1 + 2\sum\limits_{i=2}^{k-1}\tau_i + \tau_k + \beta_k)
 \end{equation}
  is also uniquely determined:
 \begin{equation}
  \label{eq_iniq_def_2a}
     \beta_k = \beta_{k+1}.
 \end{equation}

 \index{dotted ! - edge}
 \rm (ii)  By \eqref{eq_iniq_def_2a}, we have
 \begin{equation*}
   %%\label{eq_iniq_def_3}
    -(\tau_1 + 2\sum\limits_{i=2}^{k-1}\tau_i + \tau_k + \overline\tau_k) =
    -(\tau_1 + 2\sum\limits_{i=2}^k\tau_i + \tau_{k+1} + \overline\tau_{k+1}).
 \end{equation*}
 In other words,
 \begin{equation*}
    \overline\tau_k = \tau_k + \tau_{k+1} + \overline\tau_{k+1},
 \end{equation*}
 i.e., \eqref{eq_rel_tk_tk1} holds. Further,
 \begin{equation*}
    (\overline\tau_k, \overline\tau_{k+1})  = (\tau_k, \overline\tau_{k+1}) +  (\overline\tau_{k+1},\overline\tau_{k+1}) =
        -1 + 2 = 1.
 \end{equation*}
 Therefore, $\overline\tau_k$ and $\overline\tau_{k+1}$ are connected by a dotted edge,
 see Fig. \ref{Dlak_tk_tk1}$(d)$. \qed
~\\

\subsection{Relation of linkage diagrams lying in $\mathscr{L}(\widetilde\Gamma)$ and in $\mathscr{L}(\Gamma)$}
  \label{sec_rel_coval_2}
%% The $\{i,j\}$th element of
%% $B_{\widetilde\Gamma}$ is $(\widetilde\tau_i, \widetilde\tau_j)$. In
%% the transition from the basis $\widetilde{S}$ to the basis $S$ the
%% matrix $B_{\widetilde\Gamma}$ is converted to the form
%% ${}^t{M}B_{\widetilde\Gamma}M$, and its $\{i,j\}$th element is
%% \begin{equation}
%%   \label{eq_matr_M_2}
%%    [{}^t{M}B_{\widetilde\Gamma}M]_{ij} = (M\widetilde\tau_i, M\widetilde\tau_j) = (\tau_i, \tau_j) = [B_{\Gamma}]_{ij},
%% \end{equation}
%% In the other words
%% \begin{equation}
%%   \label{eq_matr_M_3}
%%     {}^t{M}B_{\widetilde\Gamma}M = B_{\Gamma}.
%% \end{equation}
%% Then
%% \begin{equation}
%%   \label{eq_matr_M_4}
%%       B_{\widetilde\Gamma} = {}^t{M}^{-1}B_{\Gamma}M^{-1}, \qquad B_{\widetilde\Gamma}^{-1} = {M}B_{\Gamma}^{-1}({}^tM).
%% \end{equation}

Denote by $\gamma^{\nabla}_{\Gamma}$ and $\gamma^{\nabla}_{\widetilde\Gamma}$ the following vectors:
\begin{equation}
  \label{eq_def_gamma_theta_nabla}
   \gamma^{\nabla}_{\Gamma} =
   \left (
   \begin{array}{c}
       (\gamma, \tau_1) \\
       \dots \\
       (\gamma, \tau_l)
   \end{array}
     \right ), \qquad
      \gamma^{\nabla}_{\widetilde\Gamma} =
   \left (
   \begin{array}{c}
       (\gamma, \widetilde\tau_1) \\
       \dots \\
       (\gamma, \widetilde\tau_l)
   \end{array}
     \right ).
\end{equation}

\begin{lemma}
  \label{lem_contragred_map}
 Let $\{\widetilde\Gamma,  \Gamma\}$ be the covalent pair of Carter diagrams from \eqref{eq_covalent_pairs}.

{\rm(i)} The vector $\gamma^{\nabla}_{\widetilde\Gamma}$ is the
linkage label vector in the linkage system
$\mathscr{L}(\widetilde\Gamma)$ if and only if
$\gamma^{\nabla}_{\Gamma}$ is the linkage label vector in the
linkage system $\mathscr{L}(\Gamma)$.

{\rm(ii)} Let $\gamma^{\nabla}_{\widetilde\Gamma}$,
$\theta^{\nabla}_{\widetilde\Gamma} \in
\mathscr{L}(\widetilde\Gamma)$ and
 $\gamma^{\nabla}_{\Gamma}$, $\theta^{\nabla}_{\Gamma} \in \mathscr{L}(\Gamma)$.
 Then $\gamma^{\nabla}_{\widetilde\Gamma}$, $\theta^{\nabla}_{\widetilde\Gamma}$ coincide if and only if
 $\gamma^{\nabla}_{\Gamma}$, $\theta^{\nabla}_{\Gamma}$ coincide.
\end{lemma}

\PerfProof
  Let $L$ be the subspace spanned by roots of $\widetilde{S} = \{\widetilde\tau_1, \dots \widetilde\tau_l\}$,
  we will write $L = [\widetilde\tau_1, \dots \widetilde\tau_l]$. By Theorem \ref{prop_bijection_root_syst}
  we have
\begin{equation}
 \label{eq_L_two_forms}
   L = [\widetilde\tau_1, \dots \widetilde\tau_l] = [\tau_1, \dots \tau_l].
\end{equation}

  (i) The fact that $\gamma^{\nabla}_{\widetilde\Gamma}$
  is the linkage label vector in the linkage system $\mathscr{L}(\widetilde\Gamma)$
  means that $\gamma$ is the linkage root linearly independent of roots of $\widetilde{S}$, see \S\ref{sec_linkage}.
  In turn, by \eqref{eq_L_two_forms}, this means,
  that $\gamma$ is also linearly independent of roots of $S$ and $\gamma^{\nabla}_{\Gamma}$
  is the linkage label vector in the linkage system $\mathscr{L}(\Gamma)$.

  (ii) By \eqref{eq_labels_proj} and \eqref{eq_def_gamma_theta_nabla}
 \begin{equation}
   \label{eq_gamma_theta}
   \begin{split}
    & \gamma^{\nabla}_{\widetilde\Gamma} = \mathscr{B}_{\widetilde\Gamma}\gamma_L, \qquad
      \theta^{\nabla}_{\widetilde\Gamma} = \mathscr{B}_{\widetilde\Gamma}\theta_L, \\
   %% & \\
    & \gamma^{\nabla}_{\Gamma} = \mathscr{B}_{\Gamma}\gamma_L, \qquad
      \theta^{\nabla}_{\Gamma} = \mathscr{B}_{\Gamma}\theta_L. \\
   \end{split}
 \end{equation}
  where $\gamma_L$ (resp. $\theta_L$) is the projection of $\gamma$ (resp. $\theta$) on $L$.
  Note that projections of $\gamma$ and $\theta$ on $L$ do not depend on choice of the diagram $\widetilde\Gamma$.
  By \eqref{eq_gamma_theta} from $\gamma^{\nabla}_{\widetilde\Gamma} = \theta^{\nabla}_{\widetilde\Gamma}$ we get
  $\mathscr{B}_{\widetilde\Gamma}\gamma_L = \mathscr{B}_{\widetilde\Gamma}\theta_L$.
  Since $\mathscr{B}_{\widetilde\Gamma}$ is positive definite (Proposition \ref{restr_forms_coincide}),
  we have $\gamma _L = \theta_L$. Again, by \eqref{eq_gamma_theta} we have $\gamma^{\nabla}_{\Gamma} = \theta^{\nabla}_{\Gamma}$.
  \qed

 \begin{corollary}
    \label{corol_Estim_2}
     For any covalent pair of Carter diagrams $\{\widetilde\Gamma, \Gamma\}$ out of \eqref{eq_covalent_pairs} we have
   \begin{equation*}
      \mid \mathscr{L}(\widetilde\Gamma) \mid ~=~ \mid \mathscr{L}(\Gamma) \mid.
   \end{equation*}
  \end{corollary}

   This follows directly from Lemma \ref{lem_contragred_map}.
  \qed

  \subsection{Linkage system components}
    \label{sec_linkage_syst_comp}
  \index{Dynkin extension ! - $\Gamma <_D \Gamma'$}
  \index{$\Gamma <_D \Gamma'$, Dynkin extension}
  \index{linkage system component}
  \index{stratum size}
  \index{ $\varPhi(S)$, partial root system}
  \index{root system ! - partial}
  \index{partial root system}
  \index{root system}

   Consider a Carter diagram  $\Gamma$.
   Let $\Gamma <_D \Gamma'$ be a Dynkin extension, $S$ a certain $\Gamma$-associated root subset,
   $S \subset \varPhi(\Gamma')$, where  $\varPhi(\Gamma')$ is the root system associated with $\Gamma'$.
 ~\\

   Recall that the root subset $\varPhi(S)$ consisting of roots of $\varPhi$
   linearly dependent on roots of $S$ is said to be a partial root system.
   The subset of roots $\varPhi(\Gamma')\backslash\varPhi(S)$ is said to be the root stratum, see \S\ref{sec_Dynkin_ext}.
   The number of roots in the root stratum $\varPhi(\Gamma')\backslash\varPhi(S)$
   is said to be the {\it stratum size}:
 \begin{equation}
  \label{eq_regext2}
 \mid \varPhi(\Gamma')\backslash\varPhi(S) \mid \quad = \quad  \mid\varPhi(\Gamma')\mid - \mid\varPhi(S)\mid,
 \end{equation}
  see \S\ref{sec_link_syst_comp}. We define the set of linkage diagrams $\mathscr{L}_{\Gamma'}(\Gamma)$ as follows:
  \index{linkage diagrams ! - in $\mathscr{L}_{\Gamma'}(\Gamma)$}
 \begin{equation}
  \label{eq_regext3}
   \mathscr{L}_{\Gamma'}(\Gamma) = \{\gamma^{\nabla} \mid \gamma \in \varPhi(\Gamma')\backslash\varPhi(S)\}.
 \end{equation}
 The set $\mathscr{L}_{\Gamma'}(\Gamma)$ is said to be the {\it linkage system component}
 (or, {\it $\Gamma'$-component} of $\mathscr{L}(\Gamma)$).
 There are several Dynkin extensions and, accordingly,
 several linkage system components for the given Carter diagram $\Gamma$:
 \begin{equation*}
     \Gamma <_D \Gamma_1, \dots, \Gamma <_D \Gamma_n.
 \end{equation*}
  The linkage system $\mathscr{L}(\Gamma)$ is the union of the sets $\mathscr{L}_{\Gamma_i}(\Gamma)$
  taken for all Dynkin extensions of $\Gamma$:
 \begin{equation*}
     \mathscr{L}(\Gamma) \quad = \quad  \bigcup_{\Gamma <_D \Gamma_i} \mathscr{L}_{\Gamma_i}(\Gamma).
 \end{equation*}
 For example, the Dynkin extensions for $\Gamma = A_7$ are as follows:
 \begin{equation*}
     A_7 <_D A_8, \quad A_7 <_D  D_8, \quad A_7 <_D E_8,
 \end{equation*}
 the corresponding linkage system components are
 $\mathscr{L}_{A_8}(A_7)$, $\mathscr{L}_{D_8}(A_7)$, $\mathscr{L}_{E_8}(A_7)$ and
 the linkage system $\mathscr{L}(A_7)$ is the union of linkage system components of
 $A$-, $D$-, and $E$-type:
 \begin{equation*}
     \mathscr{L}(A_7) ~=~  \mathscr{L}_{A_8}(A_7) \cup
                                      \mathscr{L}_{D_8}(A_7) \cup  \mathscr{L}_{E_8}(A_7).
 \end{equation*}

 For every root $\gamma \in \varPhi(\Gamma')\backslash\varPhi(S)$,
 there exist the linkage diagram $\gamma^{\nabla} \in \mathscr{L}_{\Gamma'}(\Gamma)$.
 There exist pairs of roots $\gamma_1, \gamma_2 \in \varPhi(\Gamma')\backslash\varPhi(S)$
 such that $\gamma_1^{\nabla} = \gamma_2^{\nabla}$, see Lemma \ref{lem_coinc_linkage}.
 Thus, from \eqref{eq_regext2} we deduce that:
 \begin{equation}
   \label{eq_regext4}
     \mid \mathscr{L}_{\Gamma'}(\Gamma) \mid  \quad \leq \quad  \mid\varPhi(\Gamma')\mid - \mid\varPhi(S)\mid.
 \end{equation}

%% file: 4linkageRoot.tex
~\\
\section{\sc\bf Theorem on a linkage root}
  \label{sect_quadr_form}

 \index{${\bf B}$, Cartan matrix associated with a Dynkin diagram}
 \index{Cartan matrix ! - ${\bf B}$}
 \index{$L$, linear space ! - spanned by the $\Gamma$-associated root subset}
 \index{primary root system}
 \index{root system ! - primary}

Let $\Gamma$ be the Carter diagram corresponding
to the bicolored decomposition of $w$ given as in eq. \eqref{two_invol},
and $\Pi_w = \{\tau_1,\dots,\tau_l\}$ the corresponding root basis, see \eqref{root_subset_L}.
Recall that $L$ is the $\Pi_w$-associated subspace:
\begin{equation}
  \label{space_L}
  L = [\tau_1,\dots,\tau_l],
\end{equation}
see \S\ref{sec_proj_linkage}.

 Let ${\bf B}$ be the Cartan matrix corresponding to the primary
 root system $\varPhi$, see \S\ref{sec_fulcrum}.
 %% Let $\Gamma <_D \Gamma'$ be any Dynkin extension and ${\bf B} := {\bf B}(\Gamma')$
 %% be the Cartan matrix corresponding the Dynkin diagram $\Gamma'$.
 %% see \S\ref{sec_Dynkin_ext0}.

 \begin{proposition}
   \label{restr_forms_coincide}
 {\rm(i)}  The restriction of the bilinear form associated with the Cartan matrix
  ${\bf B}$ on the subspace $L$ coincides with the bilinear form associated with the
  partial  Cartan matrix $B_{\Gamma}$, i.e., for any pair of vectors $v, u \in L$, we have
 \begin{equation}
   \label{restr_q}
       (v, u)_{\botG} = (v, u), \text{ and }
       \mathscr{B}_{\Gamma}(v) = \mathscr{B}(v).
 \end{equation}

 {\rm(ii)} For every Carter diagram, the matrix $B_{\botG}$ is positive definite.

 \end{proposition}
\PerfProof
 (i) From \eqref{canon_dec_2} we deduce:
 \begin{equation*}
     (v, u)_{\botG} =   (\sum\limits_i{t_i{\tau_i}}, \sum\limits_j{q_j{\tau_j}})_{\botG} =
      \sum\limits_{i,j}t_i{q}_j(\tau_i, \tau_j)_{\botG} = \sum\limits_{i,j}t_i{q}_j(\tau_i, \tau_j) =
      (v, u).
 \end{equation*}

 (ii) This follows from (i).
\qed

 \index{Cartan matrix ! - generalized}
\begin{remark}[The classical case]
   \label{rem_classic}
 \rm{
  Recall that the $n\times{n}$ matrix $K = (k_{ij})$, where $1 \leq i,j \leq n$, such that
 \begin{equation*}
   \begin{split}
     (C1) & \quad k_{ii} = 2 \text{ for } i = 1,\dots, n, \\
     (C2) & \quad -k_{ij} \in \mathbb{Z} = \{0, 1, 2, \dots \} \text{ for } i \neq j, \\
     (C3) & \quad  k_{ij} = 0 \text{ implies } k_{ji} = 0 \text{ for } i, j = 1, \dots, n
    \end{split}
 \end{equation*}
  is called a {\it generalized Cartan matrix}, \cite{Kac80}, \cite[\S 2.1]{St08}.
  For the Carter diagram $\Gamma$, which is not a Dynkin diagram,
  the condition (C2) fails: The elements $k_{ij}$ associated with dotted edges are positive.

 If the Carter diagram does not contain any cycle,
 then the Carter diagram is the Dynkin diagram, the corresponding conjugacy class
 is the conjugacy class of the Coxeter element, and
 the partial Cartan matrix is the classical Cartan matrix, which is
 the particular case of a generalized Cartan matrix. \qed
 }
\end{remark}

\subsection{Linear dependence and maximal roots}
  \label{sec_linear_dep}
Let $S = \{\tau_1,\dots,\tau_l\}$ be a $\Gamma$-associated subset,
The matrix $B_{\Gamma}$ is well-defined, since
$(\tau_i, \tau_j)$ is $0$ (resp. $-1$, resp. $1$) if the corresponding connection $\{\tau_i, \tau_j\}$ does not exist
(resp. solid edge, resp. dotted edge).
Let $\gamma$ be a root linearly dependent on
$S$ as follows:
\begin{equation}
  \label{eq_lin_depend1}
   \gamma = t_1\tau_1 + \dots + t_l\tau_l.
\end{equation}
Then we have
 \begin{equation}
  \label{eq_lin_depend2}
  \left (
    \begin{array}{c}
      (\gamma, \tau_1) \\
      \dots \\
      (\gamma, \tau_l) \\
    \end{array}
  \right ) =  B_{\Gamma}
  \left (
    \begin{array}{c}
      t_1 \\
      \dots\\
      t_l \\
    \end{array}
   \right ) = B_{\Gamma}\gamma,  \text{ and }
     \left (
    \begin{array}{c}
      t_1 \\
      \dots \\
      t_l \\
    \end{array}
   \right ) = {B}^{-1}_{\Gamma}
    \left (
    \begin{array}{c}
      (\gamma, \tau_1) \\
      \dots \\
      (\gamma, \tau_l) \\
    \end{array}
  \right ).
 \end{equation}
Replacing $\tau_i$ with $-\tau_i$,
we get the coefficient $-t_i$ instead of $t_i$ in the decomposition \eqref{eq_lin_depend1}.

 \index{$\gamma^{\nabla}$, linkage diagram}
\begin{remark}
  \label{rem_max_root}
 \index{dotted ! - edge}
 \index{solid ! - edge}
 \index{quadratic form ! - $\mathscr{B}_{\Gamma}$, associated with the partial Cartan matrix $B_{\Gamma}$}
{\rm
   Let the vector $\gamma$ be \underline{linearly dependent on roots} of $S = \{\tau_1,\dots,\tau_l\}$,
   and let $\gamma$ be connected with only one $\tau_i \in S$.
   We have two frequently occurring cases:
   ~\\

   (i) Suppose $\gamma$ is connected to the same point as the maximal (or minimal) root in the root system $S$.
   In other words,
   the orthogonality relations $(\gamma, \tau_i)$ in eq. \eqref{eq_lin_depend2} coincide with orthogonality relations
   for the maximal (resp. minimal) root while the edge connecting with $\gamma$
   is dotted (resp. solid). Since equation \eqref{eq_lin_depend2} has a unique solution
\begin{equation}
  \label{eq_uniq_sol}
   \gamma = B_{\Gamma}^{-1}\gamma^{\nabla},
\end{equation}
   we deduce that $\gamma$ coincides with the maximal (resp. minimal) root.
   ~\\

   (ii) Consider the necessary condition that $\gamma$ is a root.  We have
   \begin{equation}
     \label{eq_vect_inner_prod}
       \gamma^{\nabla} := \left (
    \begin{array}{c}
      (\gamma, \tau_1) \\
      \dots \\
      (\gamma, \tau_{i}) \\
      \dots \\
      (\gamma, \tau_{l}) \\
    \end{array}
  \right ) =
           \left (
    \begin{array}{c}
      0 \\
      \dots \\
      \pm{1} \\
      \dots \\
      0 \\
    \end{array}
  \right ), %% \qquad \gamma = B_{\Gamma}^{-1}\gamma^{\nabla},
   \end{equation}
  where the sign $+$ (resp. $-$) corresponds to the dotted (resp. solid)
  edge connecting $\gamma$ with $\tau_i$.
  Let $\mathscr{B}_{\Gamma}$ be the quadratic form associated with the partial Cartan matrix $B_{\Gamma}$.
  By \eqref{eq_uniq_sol} and \eqref{eq_vect_inner_prod} the value of $\mathscr{B}_{\Gamma}$ on the root $\gamma$ is as follows
\begin{equation}
  \label{eq_vect_inner_prod_2}
  \mathscr{B}_{\Gamma}(\gamma) = \langle B_{\Gamma}\gamma, \gamma \rangle =
  \langle \gamma^{\nabla}, B_{\Gamma}^{-1}\gamma^{\nabla}  \rangle = b^{\vee}_{i,i},
\end{equation}
 where $b^{\vee}_{i,i}$ is the $i$th diagonal element of $B_{\Gamma}^{-1}$. If $\gamma$ is a root, then
 $\mathscr{B}(\gamma) = 2$, and the necessary condition that $\gamma$ ia a root, which is connected with only one $\tau_i \in S$,
 is the following simple equality:
\begin{equation}
   \label{eq_necessry_cond}
    b^{\vee}_{i,i} = 2.
\end{equation}
}
\end{remark}

 \subsection{The inverse quadratic form $\mathscr{B}^{\vee}_{\Gamma}$}
   \label{sec_inverse_qvadr}
 \subsubsection{Linkage roots}
  \label{sec_linkage_diagr}
 \index{linkage label vector}
 \index{Dynkin labels}
 \index{label ! - Dynkin labels}
 \index{linkage label vector (= linkage diagram)}
 \index{linkage diagrams ! - (= linkage label vectors)}
 \index{label ! - linkage label vector}
 \index{class of diagrams ! - connection diagrams}
 \index{$\gamma^{\nabla}$, linkage diagram}
 \index{linkage root}
 \index{connection diagram}
 \index{quadratic form ! - inverse quadratic form $\mathscr{B}^{\vee}_{\Gamma}$}
 \index{bicolored decomposition}
 \index{$\Pi_w$, root subset associated with the bicolored decomposition of $w$}
 \index{$L$, linear space ! $= [\tau_1,\dots,\tau_l]$}
 \index{$L$, linear space ! $\Pi_w$-associated subspace}
 \index{$\mathscr{B}^{\vee}_{\Gamma}$, inverse quadratic form}
 \index{$\Pi_w$-associated subspace}

  Recall that the linkage diagram is obtained from a Carter diagram $\Gamma$ by adding one extra root $\gamma$,
  with its bonds, so that the roots corresponding to vertices of $\Gamma$ together with $\gamma$ form
  a linearly independent root subset. This extra root $\gamma$ is said to be the {\it linkage root}.
  Any linkage diagram constructed in this way might also be a Carter diagram but this is not necessarily so.
  With every linkage diagram we associate the linkage label vector,
  we use terms \underline{linkage label vectors} and  \underline{linkage diagrams} as
  synonyms, see \S\ref{sec_linkage_diagr_0}.
  The linkage labels are similar to the Dynkin labels, see \S\ref{sec_dominant_weight}.
  Two linkage diagrams and their linkage label vectors for the Carter diagram $E_6(a_1)$
  are depicted in Fig. \ref{E6a1_exam_linkages}.

 We would like to get an answer to the following question:

\begin{equation}
 \begin{array}{l}
  \textit{What linkage roots can be added to the irreducible linearly independent root subset? \qquad~~\quad~}
  \end{array}
\end{equation}
 ~\\
 It turns out that the answer to this question is very simple within the framework of the quadratic
 form associated with the partial Cartan matrix.  Let $\gamma$ be a linkage root for $\Gamma$, and let
  \begin{equation}
     \label{eq_labels}
       \gamma^{\nabla} := \left (
    \begin{array}{c}
      (\gamma, \tau_1) \\
      \dots \\
      (\gamma, \tau_l) \\
    \end{array}
  \right )
 \end{equation}
 be the linkage label vector.
 We denote the space spanned by $L$ and $\gamma$ by $L(\gamma)$, and write
\begin{equation}
  \label{L_extended_gamma}
   L =  [\tau_1,\dots,\tau_l], \qquad
   L(\gamma) =  [\tau_1,\dots,\tau_l,\gamma].
\end{equation}
 Let $\mathscr{B}_{\Gamma}$ be the quadratic form associated with the partial Cartan matrix $B_{\Gamma}$.
 Then
\begin{equation}
  \label{eq_vect_inner_prod_2}
  \mathscr{B}_{\Gamma}(\gamma) = \langle B_{\Gamma}\gamma, \gamma \rangle =
  \langle \gamma^{\nabla}, B_{\Gamma}^{-1}\gamma^{\nabla} \rangle.
\end{equation}
Since $B_{\Gamma}$ is positive definite, the eigenvalues of $B_{\Gamma}$
are positive. Hence the eigenvalues of $B^{-1}_{\Gamma}$ are also
positive, and the matrix $B^{-1}_{\Gamma}$ is positive definite.
We call the quadratic form $\mathscr{B}^{\vee}_{\Gamma}$ corresponding to the matrix
$B^{-1}_{\Gamma}$  the {\it inverse quadratic form}. The form
 $\mathscr{B}^{\vee}_{\Gamma}$ is positive definite.

\subsubsection{The projection of the linkage root}
   \label{sec_proj_ext_vec}
 \index{projection of the linkage root}
 \index{normal extending vector $\mu$}
 \index{$\mu$, normal extending vector}
 \index{$(\cdot, \cdot)_{\Gamma}$, symmetric bilinear form associated with $B_{\Gamma}$}
 \index{symmetric bilinear form $(\cdot, \cdot)_{\Gamma}$}
Let $L^{\perp}$ be the orthogonal complement of $L$ to $L(\gamma)$
in the sense of the symmetric bilinear form $(\cdot, \cdot)$
associated with the primary root system $\varPhi$:
\begin{equation}
   L(\gamma) =  L \oplus L^{\perp}.
\end{equation}
Let $\gamma_L$ be the projection of the linkage root $\gamma$. For
any root $\theta \in L(\gamma)$ such that $\theta \not\in L$, we
have $L(\theta) = L(\gamma)$, and $\theta$ is uniquely decomposed
into the following sum:
\begin{equation}
 \label{eq_theta_decomp}
   \theta =  \theta_L + \mu, \quad \text{ where } \quad \theta_L \in L, \quad \mu \in L^{\perp}.
\end{equation}
Given any vector $\theta$ by decomposition \eqref{eq_theta_decomp}, we introduce also the
{\it conjugate vector} $\overline{\theta}$ as follows:

\begin{equation}
 \label{eq_theta_conjug}
   \theta =  \theta_L + \mu, \quad  \overline{\theta} = \theta_L - \mu, \quad \text{ where } \quad \theta_L \in L, \quad \mu \in L^{\perp}.
\end{equation}

 \begin{figure}[h]
\centering
\includegraphics[scale=0.4]{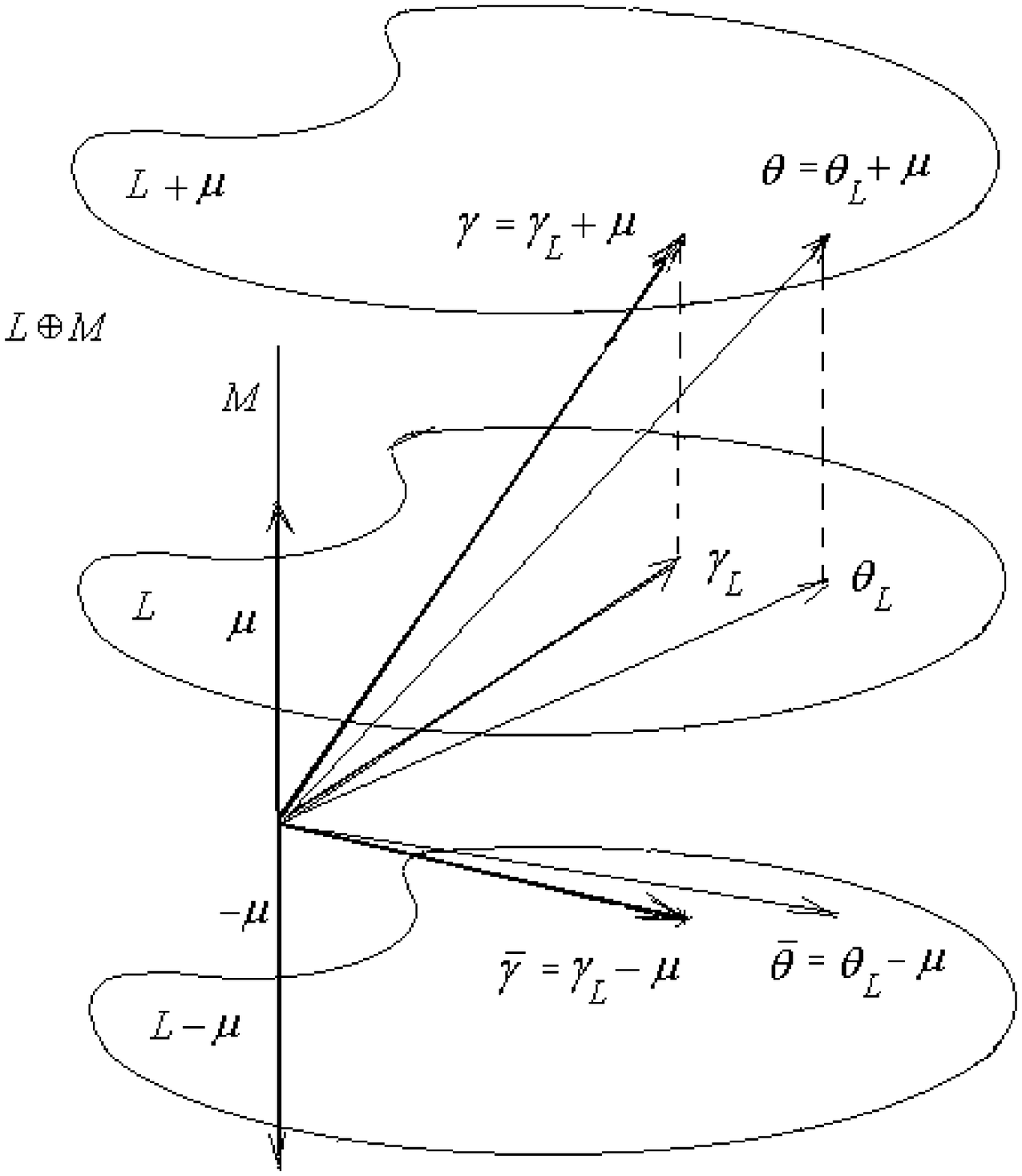}
 \caption{\hspace{3mm}The roots $\gamma = \gamma_L + \mu$ and $\overline{\gamma} = \gamma_L - \mu$.}
%%%%%% The label must come after caption
\label{L_and_mu}
\end{figure}
~\\
 The vector $\mu$ is said to be the {\it normal extending vector}, see Fig. \ref{L_and_mu}.
 The following proposition collects a number of properties of
 the projection $\gamma_L$ of the linkage root $\gamma$, the linkage label vector $\gamma^{\nabla}$, and
 the normal extending vector $\mu$:

 \index{label ! - linkage label vector}
 \index{$\gamma^{\nabla}$, linkage diagram}
 \index{quadratic form ! - inverse quadratic form $\mathscr{B}^{\vee}_{\Gamma}$}
 \index{$\mathscr{B}^{\vee}_{\Gamma}$, inverse quadratic form}
 \index{linkage label vector}
\begin{proposition}
  \label{prop_unique_val}
  \rm{(i)} The linkage label vector $\gamma^{\nabla}$ and the projection $\gamma_L$ are related as follows:
 \begin{equation}
  \label{eq_labels_proj}
   \gamma^{\nabla} = \gamma_L^{\nabla} = B_{\Gamma}\gamma_L.
 \end{equation}

  \rm{(ii)} The component $\mu$ is, up to sign, a fixed vector for any root $\theta \in L(\gamma)\backslash L$.

  \rm{(iii)} The value $\mathscr{B}_{\Gamma}(\gamma_L)$ is constant for any root $\theta \in L(\gamma)\backslash L$.

  \rm{(iv)} The vector $\theta = \theta_L + \mu$ is a root in $L(\gamma)\backslash L$ if and only if
   $\overline{\theta} = \theta_L - \mu$ is a root in  $L(\gamma)\backslash L$.

  \rm{(v)} If $\delta$ is a root in $L(\gamma)\backslash L$ such that $\delta^{\nabla} = \gamma^{\nabla}$,
   then $\delta = \gamma_L + \mu$ or $\gamma_L - \mu$, i.e. $\delta = \gamma$ or $\delta = \overline\gamma$.

  \rm{(vi)} For any root $\theta \in L(\gamma)\backslash L$, we have
\begin{equation}
    \label{length_dual_1}
     \mathscr{B}_{\Gamma}^{\vee}(\theta^{\nabla}) = \mathscr{B}_{\Gamma}(\theta_L),
\end{equation}
~\\
     and, $\mathscr{B}^{\vee}_{\Gamma}(\theta^{\nabla})$ is a constant for all roots  $\theta \in L(\gamma)\backslash L$.

 \end{proposition}

\PerfProof
 (i) Because $(\tau_i, \mu) = 0$ for any $\tau_i \in S$, we have by \eqref{eq_labels} and \eqref{eq_theta_decomp}
 \begin{equation*}
    \begin{split}
    \gamma^{\nabla} = &
      \left (
      \begin{array}{c}
         (\gamma, \tau_1) \\
         \dots \\
         (\gamma, \tau_l) \\
      \end{array}
      \right ) =
      \left (
      \begin{array}{c}
         (\gamma_L + \mu, \tau_1) \\
         \dots \\
         (\gamma_L + \mu, \tau_l) \\
      \end{array}
      \right ) =
      \left (
      \begin{array}{c}
         (\gamma_L, \tau_1) \\
         \dots \\
         (\gamma_L, \tau_l) \\
      \end{array}
      \right ) =  \gamma_L^{\nabla}.
     \end{split}
  \end{equation*}

  By \eqref{eq_lin_depend2}, since $\gamma_L \in L$, it follows that $\gamma_L^{\nabla} = B_{\Gamma}\gamma_L$. Therefore,
  $\gamma^{\nabla} = \gamma_L^{\nabla} = B_{\Gamma}\gamma_L$.
 ~\\

 (ii) Any root $\theta \in L(\gamma)$ such that $\theta \not\in L$ looks as $\tau \pm \gamma$
 for some $\tau \in L$. If $\gamma = \gamma_L + \mu$, where $\gamma_L \in L$ then
 $\theta = \tau  \pm\gamma_L  \pm\mu = \theta_L \pm \mu$, where $\theta_L = \tau \pm\gamma_L \in L$.
 ~\\

 (iii) Let $\mathscr{B}$ be the quadratic form associated with the primary root system $\varPhi$, see \S\ref{sec_fulcrum}.
  By \eqref{eq_theta_decomp} we have $\theta_L \perp \mu$, and
  $\mathscr{B}(\theta) = \mathscr{B}(\theta_L) + \mathscr{B}(\mu)$.
  Here, $\mathscr{B}(\theta) = 2$ since $\theta$ is the root, and by (ii),
  we have $\mathscr{B}(\theta_L)$ is constant.
  By \eqref{restr_q}, we have $\mathscr{B}_{\Gamma}(\theta_L) = \mathscr{B}(\theta_L)$,
  i.e.,  $\mathscr{B}_{\Gamma}(\theta_L)$ is also constant for all $\theta \in L(\gamma)$.
~\\

 (iv) Let $\theta$ be a root, i.e., $\mathscr{B}(\theta) = \mathscr{B}(\theta_L) + \mathscr{B}(\mu) = 2$.
  Then for $\overline{\theta}$, we have
  $\mathscr{B}(\overline{\theta}) = \mathscr{B}(\theta_L) + \mathscr{B}(-\mu) = 2$, and $\overline{\theta}$ is a root as well.
 ~\\

 (v) By  $\delta^{\nabla} = \gamma^{\nabla}$ and \eqref{eq_labels_proj}, we have $B_{\Gamma}\delta_L = B_{\Gamma}\gamma_L$,
  i.e., $\delta_L = \gamma_L$. By heading (ii), we have $\delta = \gamma_L + \mu$ or $\gamma_L - \mu$.
~\\

 (vi) By heading (i), since $\theta_L \in L$, we have $\theta_L^{\nabla} = B_{\Gamma}\theta_L$.  Thus,
 \begin{equation}
      \mathscr{B}^{\vee}_{\Gamma}(\theta^{\nabla}) =
      \langle B_{\Gamma}^{-1}\theta^{\nabla}, \theta^{\nabla} \rangle =
     \langle {\theta_L},  B_{\Gamma}\theta_L \rangle = \mathscr{B}_{\Gamma}(\theta_L). \\
 \end{equation}
~\\
\qed

\subsection{The criterion of a linkage root}
   \index{linkage root}
   \index{linkage root criterion}

 In this section we give a criterion that a given vector is a linkage root.
 %% see \S\ref{sec_inverse_qvadr}.

 \index{$L$, linear space ! - spanned by the $\Gamma$-associated root subset}
\begin{theorem}[Criterion of a linkage root]
   \label{th_B_less_2}
     {\rm(i)} Let $\theta^{\nabla}$ be the linkage label vector
     corresponding to a certain root $\theta \in L(\gamma)$, i.e., $\theta^{\nabla} = B_{\Gamma}\theta_L$.
     The root $\theta$ is a linkage root,
     (i.e., $\theta$ is linearly independent of roots of $L$) if and only if
     \begin{equation}
       \label{eq_p_less_2}
          \mathscr{B}^{\vee}_{\Gamma}(\theta^{\nabla}) < 2.
     \end{equation}

    {\rm(ii)} Let $\theta \in L(\gamma)$ be a root connected only
     with one $\tau_i \in S$. The root $\theta$ is a linkage root if and only if
     \begin{equation}
       \label{eq_bii_less_2}
          b^{\vee}_{i,i} < 2,
     \end{equation}
     where $b^{\vee}_{i,i}$ is the $i$th diagonal element of $B^{-1}_{\Gamma}$.
\end{theorem}

\PerfProof
    (i) Let $\mathscr{B}^{\vee}_{\Gamma}(\theta^{\nabla}) = 2$.
    By Proposition \ref{prop_unique_val},(vi) we have  $\mathscr{B}_{\Gamma}(\theta_L) = 2$,
    and  by Proposition \ref{restr_forms_coincide} we have also $\mathscr{B}(\theta_L) = 2$.
    Since $\theta \in L(\gamma)$ is a root, then
    $\mathscr{B}(\theta) = 2$ and $\mathscr{B}(\mu) = \mathscr{B}(\theta) - \mathscr{B}(\theta_L) = 0$.
    Therefore, $\mu = 0$,  and by \eqref{eq_theta_decomp} $\theta$ coincides with its projections on $L$:
    $\theta = \theta_L$.
    Thus $\theta$ is linearly depends on vectors of $L$.

    Conversely, let $\mathscr{B}^{\vee}_{\Gamma}(\theta^{\nabla}) < 2$, i.e., $\mathscr{B}_{\Gamma}(\theta_L) < 2$.
    As above, we have  $\mathscr{B}(\theta_L) < 2$
    and $\mathscr{B}(\mu) = \mathscr{B}(\theta) - \mathscr{B}(\theta_L) \neq 0$,
    i.e., $\mu \neq 0$ and $\theta$ is linearly independent of roots of $L$.

    (ii) We have
     \begin{equation*}
        \theta^{\nabla} \quad = \quad
                \left (
                 \begin{array}{c}
                    (\theta, \tau_1) \\
                    \dots \\
                    (\theta, \tau_i) \\
                    \dots \\
                    (\theta, \tau_l) \\
                 \end{array}
                  \right )
                \quad = \quad
                 \left (
                 \begin{array}{c}
                    0 \\
                    \dots \\
                    \pm{1} \\
                    \dots \\
                    0 \\
                 \end{array}
                  \right ),
     \end{equation*}
  and $\mathscr{B}^{\vee}_{\Gamma}(\theta^{\nabla}) = b^{\vee}_{i,i}$.
  Thus,  statement (ii) follows from (i). 
  \qed
\begin{remark}{\rm
  There are cases where a given root subset $L$ can be extended to two different subsets $L(\gamma) \neq L(\delta)$
  such that $\mathscr{B}^{\vee}_{\Gamma}(\gamma^{\nabla}) \neq \mathscr{B}^{\vee}_{\Gamma}(\delta^{\nabla})$.
  For example, there are two extensions for $D_5$ with different values of $\mathscr{B}^{\vee}_{\Gamma}(\theta^{\nabla})$:
  \begin{equation}
    \mathscr{B}^{\vee}_{\Gamma}(\gamma^{\nabla}) = \frac{5}{4} \quad \text{ for } \quad  D_5 \stackrel{\gamma}{<} E_6,
    \qquad \qquad
    \mathscr{B}^{\vee}_{\Gamma}(\delta^{\nabla}) = 1 \quad \text{ for } \quad  D_5 \stackrel{\delta}{<} D_6,
  \end{equation}
  see Fig. \ref{two_extensions}. The matrix $B_\Gamma^{-1}$ for $\Gamma = D_5$ is given in Table \ref{tab_partial Cartan_3}.
  For more examples of extensions with different values of
  $\mathscr{B}^{\vee}_{\Gamma}(\theta^{\nabla})$, see Table \ref{tab_val_Bmin1}.
\begin{figure}[h]
\centering
\includegraphics[scale=0.7]{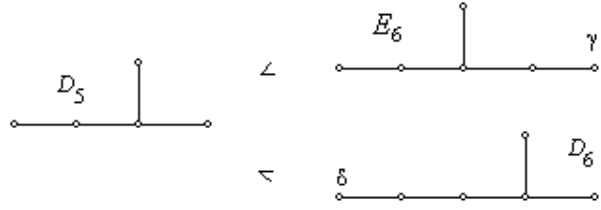}
 \caption{\hspace{3mm}Two extensions: $D_5 \stackrel{\gamma}{<} E_6$ and $D_5 \stackrel{\delta}{<} D_6$}
%%%%%% The label must come after caption
\label{two_extensions}
\end{figure}
}
\end{remark}

%% file: 5loctets.tex
~\\
\section{\sc\bf Loctets}
  \label{sect_loctets}

 \subsection{The dual partial Weyl group associated with a root subset}
   \label{sec_dual_partial}
 \index{$\Pi_w$, root subset associated with the bicolored decomposition of $w$}
 \index{bicolored decomposition}
 \index{$L$, linear space ! $= [\tau_1,\dots,\tau_l]$}
 \index{$L$, linear space ! $\Pi_w$-associated subspace}
 \index{linkage label vectors ! - $\tau_i^{\nabla}$}
 \index{$\Pi_w$-associated subspace}

 As in \S\ref{sec_inverse_qvadr}, let $\Gamma$ be the Carter diagram corresponding
to the bicolored decomposition of $w$ given as in eq.
\eqref{two_invol}, and $\Pi_w = \{\tau_1,\dots,\tau_l\}$ is the
corresponding root basis, see \eqref{root_subset_L}. Let $L$ be the
$\Pi_w$-associated subspace:
\begin{equation*}
  L = [\tau_1,\dots,\tau_l].
\end{equation*}
 For any  $\tau_i \in \Pi_w$, we define vectors $\tau_i^{\nabla}$ as follows:
   \begin{equation}
      \label{eq_tau_vee}
        \tau_i^{\nabla} := B_{\Gamma}{\tau_i}.
   \end{equation}
 Let ${L}^{\nabla}$ be the linear space spanned by vectors $\tau_i^{\nabla}$, where $\tau_i \in \Pi_w$,
 see \S\ref{sec_linkage}. The linear space ${L}^{\nabla}$ is said to be {\it the space of linkage labels}.
   Eq. \eqref{eq_tau_vee} is consistent with \eqref{dual_gamma}.
   The map $\tau_i \longrightarrow \tau_i^{\nabla}$ given by \eqref{eq_tau_vee}
   is expanded to the linear mapping $L \longrightarrow {L}^{\nabla}$, and
   \begin{equation}
      \label{eq_gamma_vee}
       u^{\nabla} = B_{\Gamma}{u} =
       \left (
       \begin{array}{c}
         (u, \tau_1) \\
          \dots \\
         (u, \tau_l) \\
       \end{array}
       \right ) \text{ for any } u \in L.
   \end{equation}
   ~\\
   Consider the restriction of the reflection $s_{\tau_i}$ on the subspace $L$.
   For any $v \in L$, by Proposition \ref{restr_forms_coincide} we have:
 \index{dual reflection $s^{*}_{\tau_i}$}
 \index{$s^{*}_{\tau_i}$, dual reflection}
   \begin{equation}
     \label{refl_tau}
       s_{\tau_i}v = v  - 2\frac{(\tau_i, v)}{(\tau_i, \tau_i)}\tau_i =
                      v  - (\tau_i, v)_{\botG}\tau_i =
                      v  - \langle B_{\Gamma}\tau_i, v \rangle \tau_i =
                      v - \langle \tau_i^\nabla, v \rangle \tau_i.
   \end{equation}
   We define the {\it dual reflection} $s^{*}_{\tau_i}$ acting on a vector $u^{\nabla} \in L^{\nabla}$ as follows:
   \begin{equation}
     \label{refl_tau_dual}
       s^{*}_{\tau_i}u^{\nabla} := u^{\nabla} - \langle u^{\nabla}, \tau_i \rangle \tau_i^{\nabla}.
   \end{equation}
   Let $W_S$ (resp. $W_S^{\vee}$) be the group generated by reflections
   $\{s_{\tau_i} \mid \tau_i \in S \}$ (resp. $\{s^{*}_{\tau_i} \mid \tau_i \in S \})$,
   where $S = \{\tau_1,\dots,\tau_l\}$.

   \begin{proposition}
    \label{prop_contagr_1}
     {\rm(i)} For any $\tau_i \in \Pi_w$, we have
   \begin{equation}
      \label{refl_transp_2}
        s^{*}_{\tau_i} = {}^t{s}_{\tau_i} = {}^t{s}_{\tau_i}^{-1}.
   \end{equation}
     {\rm(ii)} The mapping
   \begin{equation}
        w \rightarrow {}^t{w}^{-1}
   \end{equation}
     determines an isomorphism of $W_S$ onto $W_S^{\vee}$.
   \end{proposition}
 \index{$W_S$, partial Weyl group}
 \index{Weyl group ! - partial $W_S$}
 \index{$W^{\vee}_S$, dual partial Weyl group}
 \index{Weyl group ! - dual partial $W^{\vee}_S$}

   \PerfProof
  (i) By \eqref{refl_tau} and \eqref{refl_tau_dual}, for any $v \in L, u^{\nabla} \in L^{\nabla}$ we have:
   \begin{equation}
     \label{refl_transpon}
     \begin{split}
      & \langle  s^{*}_{\tau_i}u^{\nabla}, v \rangle =
       \langle  u^{\nabla} - \langle u^{\nabla}, \tau_i \rangle \tau_i^{\nabla}, v \rangle =
        \langle u^{\nabla}, v \rangle - \langle u^{\nabla}, \tau_i \rangle \langle v, \tau_i^{\nabla} \rangle, \\
      &   \langle  u^{\nabla}, s_{\tau_i}v \rangle =
        \langle  u^{\nabla} , v - \langle \tau_i^{\nabla}, v \rangle \tau_i \rangle =
        \langle u^{\nabla}, v \rangle - \langle \tau_i^{\nabla}, v \rangle \langle u^{\nabla}, \tau_i \rangle.
     \end{split}
   \end{equation}
   Thus,
   \begin{equation}
      \langle  s^{*}_{\tau_i}u^{\nabla}, v \rangle = \langle  u^{\nabla}, s_{\tau_i}v \rangle,
      \quad \text{ for any } v \in L, u^{\nabla} \in L^{\nabla},
   \end{equation}
   and \eqref{refl_transp_2} holds.
  \qed
  ~\\

 \subsubsection{Groups $W_S$ and $W_S^{\vee}$}
   \label{sec_partial_groups}

  One should note that $W_S$ (and, therefore, $W_S^{\vee}$) is not necessarily Weyl group,
  since the roots $\tau_i \in \Pi_w$ are not necessarily simple: They do not constitute a root subsystem.
  We call $W_S$ the {\it partial Weyl group}, and $W_S^{\vee}$ the {\it dual partial Weyl group
  associated with the root subset $S$}, or, for short, the {\it dual partial Weyl group}.
  By \eqref{eq_tau_vee}, \eqref{eq_gamma_vee} and \eqref{refl_tau_dual}, for any $u^{\nabla} \in L^{\nabla}$, we have
   \begin{equation}
     \label{refl_tau_dual_2}
    \begin{split}
       &  s^{*}_{\tau_i}u^{\nabla} = u^{\nabla} - \langle u^{\nabla}, \tau_i \rangle \tau^{\nabla}_i =
             u^{\nabla} - u^{\nabla}_{\tau_i} (\tau^{\nabla}_i) =
             u^{\nabla} - u^{\nabla}_{\tau_i} (B_{\Gamma}\tau_i),
          \\
       & (s^{*}_{\tau_i}u^{\nabla})_{\tau_k} = u^{\nabla}_{\tau_k} - u^{\nabla}_{\tau_i}(B_{\Gamma}\tau_i)_{\tau_k}
                                  = u^{\nabla}_{\tau_k} - u^{\nabla}_{\tau_i}(\tau_i, \tau_k).
    \end{split}
   \end{equation}
  Then, by \eqref{refl_tau_dual_2} we get
 \index{Weyl group ! - $W$}
 \index{$W$, Weyl group}
 \index{$\varPhi$, root system}
 \index{$E$, linear space spanned by all roots of $\varPhi$}
 \index{$\Pi$, set of all simple roots in $\varPhi$}
 \index{label ! - linkage label vector}
 \index{dotted ! - edge}
 \index{solid ! - edge}
 \index{$(\cdot, \cdot)$, symmetric bilinear form associated with ${\bf B}$}
 \index{symmetric bilinear form $(\cdot, \cdot)$}
 \index{root system}
    \begin{equation}
       \label{dual_refl}
       (s^*_{\tau_i}u^{\nabla})_{\tau_k} =
        \begin{cases}
            -u^{\nabla}_{\tau_i}, & \text{ if } k = i, \\
            u^{\nabla}_{\tau_k} + u^{\nabla}_{\tau_i}, & \text{ if } \{\tau_k, \tau_i \}
            \text{ is a {\it solid} edge, i.e., } (\tau_k, \tau_i) = -1, \\
            u^{\nabla}_{\tau_k} - u^{\nabla}_{\tau_i}, & \text{ if }  \{\tau_k, \tau_i \}
            \text{ is a {\it dotted} edge, i.e., } (\tau_k, \tau_i) = 1, \\
            u^{\nabla}_{\tau_i}, & \text{ if } \tau_k \text{ and } \tau_i
             \text{ are not connected, i.e., } (\tau_k, \tau_i) = 0.
        \end{cases}
    \end{equation}

    Now, we will show that values $(s^*_{\tau_i}u^{\nabla})_{\tau_k}$ belong to the set $\{-1, 0, 1\}$,
    i.e., $s^*_{\tau_i}$ acts in the linkage system $\mathscr{L}({\Gamma})$, and
    \begin{equation}
       W^{\vee}(S) : \mathscr{L}({\Gamma}) \longrightarrow \mathscr{L}({\Gamma}).
    \end{equation}
    Let $\{\tau_k, \tau_i\}$ be a solid edge. If $u^{\nabla}_{\tau_k} = u^{\nabla}_{\tau_i} = 1$
    (resp. $u^{\nabla}_{\tau_k} = u^{\nabla}_{\tau_i} = -1$) then roots $\{-u, \tau_k, \tau_i\}$
    (resp. $\{u, \tau_k, \tau_i\}$) constitute the root system corresponding to the
    extended Dynkin diagram $\widetilde{A}_2$, which is impossible,
    see Lemma \ref{lem_must_dotted}.
    For remaining pairs $\{u^{\nabla}_{\tau_k}, u^{\nabla}_{\tau_i}\}$, we have
    $-1 \leq u^{\nabla}_{\tau_k} + u^{\nabla}_{\tau_i} \leq 1$. Now, let $\{\tau_k, \tau_i\}$ be a dotted edge.
    If $u^{\nabla}_{\tau_k} = 1$ and $u^{\nabla}_{\tau_i} = -1$ (resp. $u^{\nabla}_{\tau_k} = -1$ and $u^{\nabla}_{\tau_i} = 1$)
    then roots $\{u, -\tau_k, \tau_i\}$ (resp. $\{u, \tau_k, -\tau_i\}$)
    constitute the root system $\widetilde{A}_2$, which is impossible.
    For remaining pairs $\{u^{\nabla}_{\tau_k}, u^{\nabla}_{\tau_i}\}$, we have
    $-1 \leq u^{\nabla}_{\tau_k} - u^{\nabla}_{\tau_i} \leq 1$.
~\\

 \index{dual reflection $s^{*}_{\tau_i}$}
 \index{$s^{*}_{\tau_i}$, dual reflection}
\begin{proposition}
  \label{prop_link_diagr_conn}
  {\rm(i)} For dual reflections $s^{*}_{\tau_i}$, the following relations hold
  \begin{equation}
      \label{conn_dual_1}
        B_{\Gamma} s_{\tau_i} = s^*_{\tau_i} B_{\Gamma},
     \end{equation}
  \begin{equation}
     \label{conn_L_2}
    (s_{\tau_i}\gamma)^{\nabla} = s^{*}_{\tau_i}B_{\Gamma}\gamma_L = s^{*}_{\tau_i}\gamma^{\nabla}.
  \end{equation}
  {\rm(ii)} For $w^* \in W^{\vee}_S$ $(w^* := {}^t{w}^{-1})$, the dual element of $w \in W$, we have
    \begin{equation}
      \label{conn_L_3}
        (w\gamma)^{\nabla} = w^*{\gamma}^{\nabla} \quad(= {}^t{w}^{-1}{\gamma}^{\nabla}).
     \end{equation}
  {\rm(iii)} The following relations hold
    \begin{equation}
      \label{conn_L_5}
        \begin{split}
        & \mathscr{B}_{\Gamma} (s_{\tau_i} v)  =  \mathscr{B}_{\Gamma}(v)  \text{ for any } v \in L, \\
        & \mathscr{B}^{\vee}_{\Gamma}(s^*_{\tau_i} u^{\nabla}) =  \mathscr{B}^{\vee}_{\Gamma}(u^{\nabla}) \text{ for any } u^{\nabla} \in L^{\nabla}.
        \end{split}
     \end{equation}
     \end{proposition}

\PerfProof (i)  The equality \eqref{conn_dual_1} holds since the
following is true for
   any $u, v \in L$:
 \begin{equation*}
    \begin{split}
     & (s_{\tau_i}{u}, v)_{\Gamma} = ({u}, s_{\tau_i}{v})_{\Gamma}, \text{ i.e., }
     \langle B_{\Gamma} s_{\tau_i}{u}, v \rangle = \langle B_{\Gamma} {u}, s_{\tau_i}{v} \rangle =
     \langle s^*_{\tau_i}{B}_{\Gamma} {u}, v \rangle \text{, and } \\
     & \langle (B_{\Gamma} s_{\tau_i} - s^*_{\tau_i}{B}_{\Gamma}) {u}, v \rangle = 0.
    \end{split}
 \end{equation*}
~\\
Let us consider eq. \eqref{conn_L_2}. Since $(\tau_i, \mu) = 0$ for
any $\tau_i \in S$, and $s_{\tau_i}\mu = \mu$, we have
 \begin{equation*}
  \begin{split}
    (s_{\tau_i}\gamma)^{\nabla} =
      & \left (
      \begin{array}{c}
         (s_{\tau_i}\gamma, \tau_1) \\
         \dots \\
         (s_{\tau_i}\gamma, \tau_l) \\
      \end{array}
      \right ) =
      \left (
      \begin{array}{c}
         (s_{\tau_i}\gamma_L + \mu, \tau_1) \\
         \dots \\
         (s_{\tau_i}\gamma_L + \mu, \tau_l) \\
      \end{array}
      \right ) =
      \left (
      \begin{array}{c}
         (s_{\tau_i}\gamma_L, \tau_1) \\
         \dots \\
         (s_{\tau_i}\gamma_L, \tau_l) \\
      \end{array}
      \right ) =
      (s_{\tau_i}\gamma_L)^{\nabla}.
      \\
     \end{split}
  \end{equation*}
  Then, by \eqref{conn_dual_1} and \eqref{eq_labels_proj}
  \begin{equation*}
    (s_{\tau_i}\gamma)^{\nabla} =
    (s_{\tau_i}\gamma_L)^{\nabla} = B_{\Gamma}{s}_{\tau_i}\gamma_L =
    s^{*}_{\tau_i}B_{\Gamma}\gamma_L = s^{*}_{\tau_i}\gamma^{\nabla}.
  \end{equation*}
~\\

(ii)  Let $w = s_{\tau_1}s_{\tau_2}\dots{s}_{\tau_m}$ be a
decomposition of $w \in W$. Since  $s^*_{\tau} = {}^ts^{-1}_{\tau} =
{}^ts_{\tau}$,
 we deduce from \eqref{conn_L_2} the following:
\begin{equation*}
  \begin{split}
   & (w\gamma)^{\nabla} = (s_{\tau_1}s_{\tau_2}\dots{s}_{\tau_m}\gamma)^{\nabla} =
    s^*_{\tau_1}(s_{\tau_2}\dots{s}_{\tau_m}\gamma)^{\nabla} =
    s^*_{\tau_1}s^*_{\tau_2}({s}_{\tau_3}\dots{s}_{\tau_m}\gamma)^{\nabla} = \dots = \\
   & s^*_{\tau_1}{s}^*_{\tau_2}\dots{s}^*_{\tau_m}{\gamma}^{\nabla} =
     {}^t{(s_{\tau_m}\dots{s}_{\tau_2}{s}_{\tau_1})}{\gamma}^{\nabla} =
     {}^t{({s}_{\tau_1}s_{\tau_2}\dots{s}_{\tau_m})}^{-1}{\gamma}^{\nabla} =
     {w}^*{\gamma}^{\nabla}.
  \end{split}
\end{equation*}
~\\

 (iii)
 Further, \eqref{conn_L_5} holds since
 \begin{equation*}
    \begin{split}
    & \mathscr{B}_{\Gamma}(s_{\tau_i}{v}) = \langle B_{\Gamma} s_{\tau_i}{v}, s_{\tau_i}{v} \rangle =
     \langle s^*_{\tau_i}{B}_{\Gamma} {v}, s_{\tau_i}{v} \rangle =
     \langle {B}_{\Gamma} {v}, {v} \rangle = \mathscr{B}_{\Gamma}(v). \\
    \hspace{2.3cm} &
     \mathscr{B}_{\Gamma}^{\vee}(s^{*}_{\tau_i}{u^{\nabla}}) =
     \langle B_{\Gamma}^{\vee} s^{*}_{\tau_i}{u^{\nabla}}, s^{*}_{\tau_i}{u^{\nabla}} \rangle =
     \langle s_{\tau_i}{B}_{\Gamma}^{\vee} {u^{\nabla}}, s^{*}_{\tau_i}{u^{\nabla}} \rangle =
     \langle {B}_{\Gamma}^{\vee} {u^{\nabla}}, {u^{\nabla}} \rangle = \mathscr{B}_{\Gamma}^{\vee}(u^{\nabla}).
     \hspace{2.3cm} \qed \\
    \end{split}
 \end{equation*}
~\\
  The linkage system component $\mathscr{L}_{\Gamma'}(\Gamma)$,
  where $\Gamma <_D \Gamma'$ is a Dynkin extension, is the $W^{\vee}_S$-orbit on the set of linkage diagrams.
  We mentioned in \S\ref{sec_rat_p} that the rational number $p = \mathscr{B}^{\vee}_{\Gamma}(u^{\nabla})$
  is the \underline{invariant} of the linkage system component.
  We call this component {\it the extension of the Carter diagram $\Gamma$ by p} and
  denote this extension by $\{ \Gamma, p \}$.

\subsection{On $4$-cycles and $4$-cycles with a diagonal}

 In further considerations we need some facts on  $4$-cycles and $4$-cycles with a diagonal.

\subsubsection{How many endpoints may a linkage diagram have?}
 \index{bicolored decomposition}
 \index{$\alpha$-set, subset of roots corresponding to $w_1$}
 \index{$\beta$-set, subset of roots corresponding to $w_2$}
 \index{$\alpha$-endpoints}
 \index{$\beta$-endpoints}

 In what follows, we show that the number of endpoints in any linkage diagram is $\leq 6$
 (see Proposition \eqref{prop_numb_ep}(i)),  and in some cases this number is $\leq 4$
 (see Proposition \eqref{prop_numb_ep}(iii)).

\begin{proposition}
   \label{prop_numb_ep}
    Let $w = w_\alpha w_\beta$ be a bicolored decomposition
    of $w$ into the product of two involutions, and $\Gamma$ the Carter
    diagram corresponding to this decomposition.

    {\rm(i)} Let $\gamma^{\nabla}$ be any linkage diagram obtained from $\Gamma$.
    The number of $\alpha$-endpoints in $\gamma^{\nabla}$ $\leq 3$.
    The same applies to $\beta$-endpoints.
%%~\\

    {\rm(ii)} Let an $\alpha$-set in $\Gamma$ contain $3$ roots $\{\alpha_1, \alpha_2, \alpha_3\}$.
    Suppose there exist non-connected roots $\beta$ and $\gamma$ connected to every $\alpha_i$ from $\alpha$-set.
    Then roots $\{\alpha_1, \alpha_2, \alpha_3, \beta, \gamma\}$ are linearly dependent.
%% ~\\

    {\rm(iii)} Let $\{\alpha_1, \beta_1, \alpha_2, \beta_2\}$ constitute a square in $\Gamma$.
    Suppose $\gamma$ connected to all vertices of the square. Then roots
    $\{\alpha_1, \beta_1, \alpha_2, \beta_2, \gamma\}$ are linearly dependent.
    %% There does not exist
    %% a root $\gamma$ connected to all vertices of the square. %% \rm{(}\cite[Corollary $3.4{\rm(ii)}$]{St10}\rm{)}.
 \end{proposition}

\PerfProof
   (i) If the linkage diagram $\gamma^{\nabla}$ has $4$ $\alpha$-endpoints, then $\gamma^{\nabla}$
      contains the diagram
      $\widetilde{D}_4 = \{ \gamma, \alpha_1, \alpha_2, \alpha_3, \alpha_4 \}$:
      The vector\footnotemark[1] $\varphi = 2\gamma + \sum\limits_{i=1}^4\alpha_i$
      has zero length, since
\begin{equation}
  (\varphi, \varphi) = 4(\gamma, \gamma) + \sum\limits_{i=1}^4(\alpha_i, \alpha_i)  +
        4\sum\limits_{i=1}^4(\gamma, \alpha_i) = 4\cdot{2} + 4\cdot{2} - 16\cdot{1} = 0.
\end{equation}
   \footnotetext[1]{We denote vertices and the corresponding roots by the same letters.}
      Hence, $\varphi = 0$, contradicting the linear independence of roots
      $\{ \gamma, \alpha_1, \alpha_2, \alpha_3, \alpha_4 \}$.

\begin{figure}[h]
 \centering
 \includegraphics[scale=0.7]{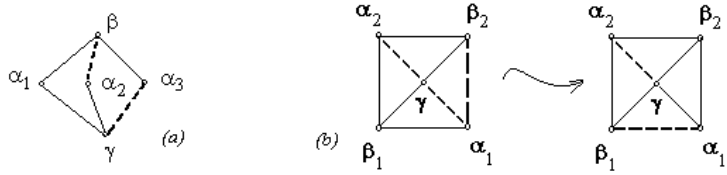}
 \caption{\hspace{3mm}Every cycle should contain an odd number of dotted edges,
          a case which cannot happen}
 %%%%%% The label must come after caption
 \label{3-4-cells}
 \end{figure}

   %% (ii)\footnotemark[2]
   (ii) Suppose there exist roots $\beta$ and $\gamma$ connected to every $\alpha_i$
  such that the vectors of the quintuple $\{\alpha_1, \alpha_2, \alpha_3, \beta, \gamma \}$
  are linearly independent, see Fig. \ref{3-4-cells}$(a)$.
  Then we have three cycles: $\{ \alpha_i, \gamma, \alpha_j, \beta \}$,  where  $1 \leq i < j  \leq 3$.
  Every cycle should contain an odd number of dotted edges.
  Let $n_1$, $n_2$, $n_3$ be the odd numbers of dotted edges in every cycle, therefore $n_1 + n_2 + n_3$ is odd.
  On the other hand, every dotted edge enters twice, so $n_1 + n_2 + n_3$ is even, which is a contradiction.
   %%\footnotetext[2]{For the reader's convenience, we give here a proof of (ii) and (iii)
   %%(originally it was presented in \cite[Corollary $3.4$]{St10}\rm{)}.}

  (iii)  Suppose $\gamma$ is connected to all vertices of the square %% $\{\alpha_1, \beta_1, \alpha_2, \beta_2\}$
  and $\{\alpha_1, \beta_1, \alpha_2, \beta_2, \gamma\}$
  are linearly independent.  Then we have $5$ cycles: Four triangles
  $\{ \alpha_i, \beta_j, \gamma \}$, where $i = 1,2$ and $j = 1,2$, and the square
  $\{\alpha_1, \beta_1, \alpha_2, \beta_2 \}$, see Fig. \ref{3-4-cells}$(b)$.
  Every cycle should contain an odd number of dotted edges.
  Let $n_1$, $n_2$, $n_3$, $n_4$, $n_5$ be the numbers of dotted edges in every cycle, therefore
  $n_1 + n_2 + n_3 + n_4 + n_5$ is odd. On the other hand, every dotted edge enters twice,
  so $n_1 + n_2 + n_3 + n_4 + n_5$ is even, which is a contradiction.
  For example, the left square in Fig. \ref{3-4-cells}$(b)$
  is transformed to the  right  one by the reflection $s_{\alpha_1}$,
  then the right square contains the cycle $\{\alpha_1, \beta_2, \gamma\}$
  with $3$ solid edges, i.e., the extended Dynkin diagram $\widetilde{A}_2$,
  contradicting Lemma \ref{lem_must_dotted}.
 \qed

\subsubsection{The diagonal in a square}
\begin{figure}[h]
\centering
\includegraphics[scale=0.8]{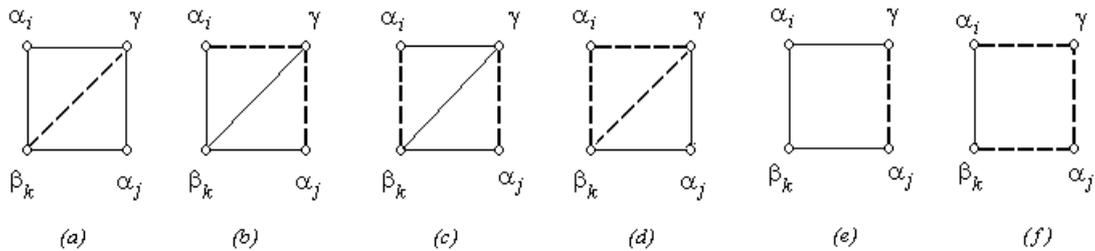}
 \caption{\hspace{3mm}Linkage diagrams containing a square.}
%%%%%% The label must come after caption
\label{square_n_diag}
\end{figure}

 The following proposition describes the structure of squares with and without diagonal.
 The squares arise during consideration of linkage diagrams near to the branch point of any Carter diagram,
 see Figs. \ref{E7a1_gamma_ij_3} - \ref{E7a3_gamma_ij_4}.
 The following statement is essential to the proof of Proposition \ref{prop_generating_diagr} determining the structure of loctets.
 \index{dotted ! - edge}
 \index{solid ! - edge}
 \index{dotted ! - diagonal}
 \index{solid ! - diagonal}
 \index{squares with a diagonal}
\begin{proposition}[On squares with a diagonal]
  \label{prop_diagonal}
   Let $\gamma$ form a linkage diagram containing
   the square $\{ \alpha_i, \beta_k, \alpha_j, \gamma \}$
   without the diagonal $\{ \alpha_i, \alpha_j \}$, i.e.,
 \begin{equation*}
   (\alpha_i, \beta_k) \neq 0, \quad  (\beta_k, \alpha_j) \neq 0, \quad
   (\alpha_j, \gamma) \neq 0, \quad  (\gamma, \alpha_i) \neq 0, \quad \text{ and } \quad  (\alpha_i, \alpha_j) = 0.
 \end{equation*}
   Let the roots $\{ \alpha_i, \beta_k, \alpha_j, \gamma \}$ be linearly independent.
   If the number of dotted edges in the square is even (resp. odd)  then the square has a diagonal
   (resp. has no diagonals).
  %% If there is an even number of dotted edges in the square, then the square has a diagonal.
  %% If there is an odd number of dotted edges in the square, then the square has no diagonals.
   Namely:

    {\rm (a)} If the square has no dotted edge, then $(\gamma, \beta_k) = 1$, i.e.,
    \underline{there exists the dotted diagonal $\{ \gamma, \beta_k \}$}, see Fig. \ref{square_n_diag},(a).

   {\rm (b)} If the square has two dotted edges $\{\gamma, \alpha_i\}$ and
   $\{\gamma, \alpha_j\}$ and remaining edges are solid, then $(\gamma, \beta_k) = -1$
   i.e., \underline{there exists the solid diagonal $\{ \gamma, \beta_k \}$}, see Fig. \ref{square_n_diag},(b).

   {\rm (c)} If the square has two dotted edges $\{\gamma, \alpha_j\}$ and
   $\{\beta_k, \alpha_i\}$  and remaining edges are solid, then $(\gamma, \beta_k) = -1$,
   i.e., \underline{there exists the solid diagonal  $\{ \gamma, \beta_k \}$},  see Fig. \ref{square_n_diag},(c).

   {\rm (d)} If the square has two dotted edges $\{\gamma, \alpha_i\}$ and $\{\beta_k, \alpha_i\}$
   and remaining edges are solid, then $(\gamma, \beta_k) = 1$,
   i.e., \underline{there exists the dotted diagonal  $\{ \gamma, \beta_k \}$},
   see Fig. \ref{square_n_diag},(d).

   {\rm (e)} If there is only one dotted edge $\{ \gamma, \alpha_j \}$ and remaining edges are
   solid, then $(\gamma, \beta_k) = 0$, i.e., \underline{there are no diagonals},
   see Fig. \ref{square_n_diag},(e).

   {\rm (f)} If there are three dotted edges
   $\{ \gamma, \alpha_i \}, \{ \gamma, \alpha_j \}, \{ \beta_k, \alpha_j \}$ and the remaining edge is
   solid, then $(\gamma, \beta_k) = 0$, i.e., \underline{there are no diagonals},
   see Fig. \ref{square_n_diag},(f).
\end{proposition}

 \index{linkage diagrams ! - containing a square}

  \PerfProof If there is no diagonal $\{ \gamma, \beta_k \}$ for one of cases (a),(b), (c) or (d),
   we get the extended Dynkin diagram $\widetilde{A}_3$ by following changes:
  \begin{equation*}
     (a) \quad \text{no changes,}
     \qquad (b) \quad \gamma \longrightarrow -\gamma,
     \qquad (c) \quad \gamma \longrightarrow -\gamma, \quad  \alpha_i \longrightarrow  -\alpha_i,
     \qquad (d) \quad \alpha_i \longrightarrow  -\alpha_i, \\
  \end{equation*}
  contradicting Lemma \ref{lem_must_dotted}.
  If the diagonal (dotted or solid) exists in cases (e) or (f)
  then one of the triangles obtained (if necessary, after the change $\gamma \longrightarrow -\gamma$) is
  the extended Dynkin diagram $\widetilde{A}_2$, contradicting Lemma \ref{lem_must_dotted}.
   \qed

\subsection{Loctets and unicolored linkage diagrams}
  \label{sec_loctets_diagr}
  \index{linkage diagrams ! - $\beta$-unicolored}
 In this section we give the complete description of linkage diagrams for every linkage system.
 It turns that each linkage diagram containing at least one non-zero $\alpha$-label $\alpha_1, \alpha_2$ or $\alpha_3$
 belongs to a certain $8$-cell linkage
 subsystem which we call the {\it loctet}, see Fig. \ref{fig_loctets}.
 Every linkage system is the union of several loctets and several $\beta$-unicolored
 linkage diagrams. For the exact description, see Tables \ref{tab_seed_linkages_6},
  Figs. \ref{D4a1_linkages}--\ref{D4_loctets}, \ref{Dk_al_linkages}, \ref{Dk_al_wind_rose},
  Theorem \ref{th_full_descr}.

\index{class of diagrams ! - $\mathsf{C4}$}
   \begin{proposition}
     Let $\Gamma$ be a Carter diagram in $\mathsf{C4}$ class, $\Gamma \neq D_4(a_1)$.
     Let $\gamma^{\nabla}$ be a linkage diagram from $\mathscr{L}(\Gamma)$,
       \begin{equation}
          \gamma^{\nabla} =
            \left (
            \begin{array}{c}
              a_1 \\
              a_2 \\
              a_3 \\
              \dots \\
              b_1 \\
              \dots \\
            \end{array}
            \right )
                =
            \left (
            \begin{array}{c}
              (\alpha_1, \gamma) \\
              (\alpha_2, \gamma) \\
              (\alpha_3, \gamma) \\
              \dots \\
              (\beta_1, \gamma) \\
              \dots \\
            \end{array}
            \right )
            \begin{array}{c}
              \alpha_1 \\
              \alpha_2 \\
              \alpha_3 \\
              \dots \\
              \beta_1 \\
              \dots \\
            \end{array}
       \end{equation}
       Among labels $a_i$ of the linkage diagram $\gamma^{\nabla}$, where $i=1,2,3$,
       at least one label $a_i$ is equal to $0$.
   \end{proposition}
 \PerfProof
   For $\Gamma \neq D_4(a_1)$, any
   Carter diagram out $\mathsf{C4}$ contains $D_5(a_1)$ as a subdiagram.
   In the diagram $D_5(a_1)$, the vertices  $\alpha_1, \alpha_2, \alpha_3$ are connected to $\beta_1$,
   see Fig. \ref{D5a1_D4_pattern}.
   Thus, no root $\gamma$ can be connected to all $\alpha_i$, where $i = 1,2,3$,
   otherwise we get the contradiction with Proposition \ref{prop_numb_ep}(ii).
 \qed
~\\

  The following proposition explains relations between linkage diagrams depicted
  in Fig. \ref{fig_loctets} and shows that every linkage diagram containing
  at least one non-zero $\alpha$-label belongs to one of the loctets $L_{12}$,  $L_{13}$, $L_{23}$.

\begin{figure}[h]
\centering
\includegraphics[scale=0.6]{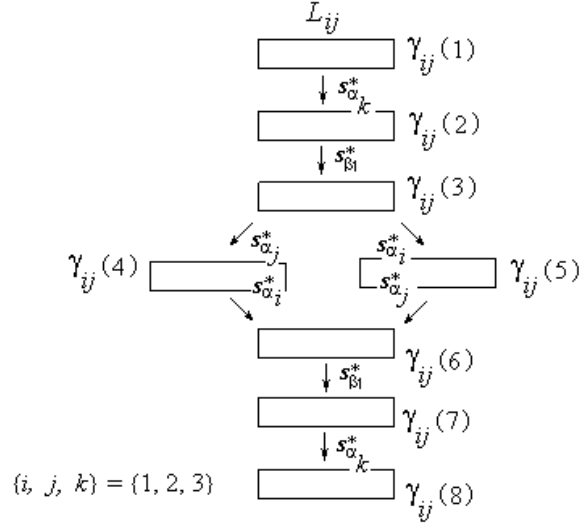}
 \caption{\hspace{3mm}{Linkage diagrams $\gamma^{\nabla}_{ij}(n)$, where $1 \le n \le 8$, $1 \le i < j \le 3$, in the loctet $L_{ij}$}}
%%%%%% The label must come after caption
\label{loctet_gamma_ij}
\end{figure}

  \index{label ! - linkage label vector}
  \index{linkage diagrams ! - in the loctet: $\gamma^{\nabla}_{ij}(n), 1 \leq n \leq 8$}
\begin{proposition}[Structure of the loctet]
  \label{prop_generating_diagr}
  {\rm(i)} The linkage label vector $\gamma^{\nabla}_{ij}(n)$ depicted in Fig. $\ref{loctet_gamma_ij}$
  (see also Fig. $\ref{fig_loctets}$) are connected by reflections $s^{*}_{\alpha_i}$, where $i = 1,2,3$,
  and reflection $s^{*}_{\beta_1}$ as follows:
\begin{equation}
  \label{refl_dual_linkages}
  \begin{split}
   & s^{*}_{\alpha_k}\gamma^{\nabla}_{ij}(8) = \gamma^{\nabla}_{ij}(7), \quad
     s^{*}_{\alpha_k}\gamma^{\nabla}_{ij}(1) = \gamma^{\nabla}_{ij}(2), \\
   & s^{*}_{\beta_1}\gamma^{\nabla}_{ij}(7) = \gamma^{\nabla}_{ij}(6), \quad
     s^{*}_{\beta_1}\gamma^{\nabla}_{ij}(2) = \gamma^{\nabla}_{ij}(3), \\
    s^{*}_{\alpha_i}\gamma^{\nabla}_{ij}(6) & =
     s^{*}_{\alpha_j}\gamma^{\nabla}_{ij}(3) = \gamma^{\nabla}_{ij}(4), \quad
     s^{*}_{\alpha_j}\gamma^{\nabla}_{ij}(6) =
     s^{*}_{\alpha_i}\gamma^{\nabla}_{ij}(3) = \gamma^{\nabla}_{ij}(5),
  \end{split}
\end{equation}
where $\{ i,j, k \} = \{1, 2, 3 \}$. Relations of the last line in
\eqref{refl_dual_linkages} hold up to permutation of indices $n = 4$
and $n = 5$ in  $\gamma^{\nabla}_{ij}(n)$.

  {\rm(ii)} If $\gamma^{\nabla}$ contains exactly two non-zero labels $a_i$, $a_j$ (corresponding to $\alpha_i, \alpha_j$),
  where $1 \leq i,j \leq 3$,
  then $\gamma^{\nabla}$ is one of the following linkage diagrams:
   \begin{equation*}
      \gamma^{\nabla}_{ij}(3), \quad \gamma^{\nabla}_{ij}(4), \quad
      \gamma^{\nabla}_{ij}(5), \quad \gamma^{\nabla}_{ij}(6).
   \end{equation*}

   {\rm(iii)} If $\gamma^{\nabla}$ contains exactly one non-zero label $a_i$ (corresponding to
   $\alpha_i$),  where $1 \leq i \leq 3$, then $\gamma^{\nabla}$ is one of the following linkage diagrams:
   \begin{equation*}
      \gamma^{\nabla}_{ij}(1), \quad \gamma^{\nabla}_{ij}(2), \quad
      \gamma^{\nabla}_{ij}(7), \quad \gamma^{\nabla}_{ij}(8).
   \end{equation*}
\end{proposition}

\PerfProof
  (i) Recall that vertices $\alpha_1$, $\alpha_2$, $\alpha_3$ are connected with $\beta_1$ by starlike shape,
  see Fig. \ref{D5a1_D4_pattern}.
  By \eqref{dual_refl} we have
  \begin{equation}
       \label{dual_refl_alpha_beta}
      \begin{split}
       (s^*_{\alpha_i}\gamma^{\nabla})_{\alpha_k} =
        \begin{cases}
            -\gamma^{\nabla}_{\alpha_i}  \text{ for } k = i, \\
            \hspace{3mm} \gamma^{\nabla}_{\alpha_k}  \text{ for } k \neq i, \\
        \end{cases} \quad
        (s^*_{\alpha_i}\gamma^{\nabla})_{\beta_k} =
        \begin{cases}
            \gamma^{\nabla}_{\beta_k} + \gamma^{\nabla}_{\alpha_i} \text{ for } (\beta_k, \alpha_i) = -1, \\
            \gamma^{\nabla}_{\beta_k} - \gamma^{\nabla}_{\alpha_i}  \text{ for } (\beta_k, \alpha_i) = 1, \\
            \gamma^{\nabla}_{\beta_k}  \text{ for } (\beta_k, \alpha_i) = 0.
        \end{cases} \\
          \\
       (s^*_{\beta_i}\gamma^{\nabla})_{\beta_k} =
        \begin{cases}
            -\gamma^{\nabla}_{\beta_i} \text{ for } k = i, \\
            \hspace{3mm} \gamma^{\nabla}_{\beta_k} \text{ for } k \neq i, \\
        \end{cases} \quad
        (s^*_{\beta_i}\gamma^{\nabla})_{\alpha_k} =
        \begin{cases}
            \gamma^{\nabla}_{\alpha_k} + \gamma^{\nabla}_{\beta_i} \text{ for } (\beta_i, \alpha_k) = -1, \\
            \gamma^{\nabla}_{\alpha_k} - \gamma^{\nabla}_{\beta_i} \text{ for } (\beta_i, \alpha_k) = 1, \\
            \gamma^{\nabla}_{\alpha_k}  \text{ for } (\beta_i, \alpha_k) = 0.
        \end{cases} \\
       \end{split}
    \end{equation}
  We show \eqref{refl_dual_linkages} only for $\gamma^{\nabla}_{ij}(n)$, where $n = 6, 7, 8$.
  Further, we use the \lq\lq{skew-symmetric relations\rq\rq~ for coordinates $\alpha_1$, $\alpha_2$, $\alpha_3$ and $\beta_1$
  for the following pairs of linkage diagrams $\{\gamma^{\nabla}_{ij}(1), \gamma^{\nabla}_{ij}(8)\}$,
  $\{\gamma^{\nabla}_{ij}(2), \gamma^{\nabla}_{ij}(7)\}$ and
  $\{\gamma^{\nabla}_{ij}(3), \gamma^{\nabla}_{ij}(6)\}$,  see Fig. \ref{fig_loctets}.
  One proves the remaining cases $n=1, 2, 3$ by changing the sign of coordinates $\alpha_1$, $\alpha_2$, $\alpha_3$ and $\beta_1$
  in $\gamma^{\nabla}_{ij}(n)$.
  Applying $s^*_{\alpha_k}$ to $\gamma^{\nabla}_{ij}(8)$, where $\{i, j, k\} = \{1, 2, 3\}$,
  we have the first line of \eqref{refl_dual_linkages} as follows:
\begin{equation*}
    %%\footnotesize
    %% \tiny
       \label{dual_refl_alpha_2}
    \begin{split}
     & \boxed{s^*_{\alpha_k}\gamma^{\nabla}_{ij}(8) = \gamma^{\nabla}_{ij}(7)} \\
       \\
    %%\tiny
      s^*_{\alpha_3}
      \left (
       \begin{array}{c}
          0 \\
          0 \\
          1 \\
          0 \\
       \dots \\
       \end{array} \right ) =
      \left (
       \begin{array}{c}
          0 \\
          0 \\
          -1 \\
          1 \\
       \dots \\
       \end{array} \right )
       \begin{array}{c}
          \alpha_1 \\
          \alpha_2 \\
          \alpha_3 \\
          \beta_1 \\
          \dots \\
       \end{array}, \quad
     s^*_{\alpha_2}
     & \left (
       \begin{array}{c}
          0 \\
          1 \\
          0 \\
          0 \\
       \dots \\
       \end{array} \right ) =
      \left (
       \begin{array}{c}
          0 \\
          -1 \\
          0 \\
          1 \\
       \dots \\
       \end{array} \right )
       \begin{array}{c}
          \alpha_1 \\
          \alpha_2 \\
          \alpha_3 \\
          \beta_1 \\
          \dots \\
       \end{array}, \quad
       s^*_{\alpha_1}
      \left (
       \begin{array}{c}
          1 \\
          0 \\
          0 \\
          0 \\
       \dots \\
       \end{array} \right ) =
      \left (
       \begin{array}{c}
          -1 \\
          0 \\
          0 \\
          1 \\
       \dots \\
       \end{array} \right )
       \begin{array}{c}
          \alpha_1 \\
          \alpha_2 \\
          \alpha_3 \\
          \beta_1 \\
          \dots \\
       \end{array}.
     \end{split}
  \end{equation*}
  Applying $s^*_{\beta_1}$ to $\gamma^{\nabla}_{ij}(7)$,
  we have the second line of \eqref{refl_dual_linkages}:
  \begin{equation*}
    %% \tiny
       \label{dual_refl_alpha_3}
    \begin{split}
     & \boxed{s^*_{\beta_1}\gamma^{\nabla}_{ij}(7) = \gamma^{\nabla}_{ij}(6)} \\
       \\
      s^*_{\beta_1}
      \left (
       \begin{array}{c}
          0 \\
          0 \\
          -1 \\
          1 \\
       \dots \\
       \end{array} \right ) =
      \left (
       \begin{array}{c}
          1 \\
          1 \\
          0 \\
          -1 \\
       \dots \\
       \end{array} \right )
       \begin{array}{c}
          \alpha_1 \\
          \alpha_2 \\
          \alpha_3 \\
          \beta_1 \\
          \dots \\
       \end{array}, \quad
      s^*_{\beta_1}
      & \left (
       \begin{array}{c}
          0 \\
          -1 \\
          0 \\
          1 \\
       \dots \\
       \end{array} \right ) =
      \left (
       \begin{array}{c}
          1 \\
          0 \\
          1 \\
          -1 \\
       \dots \\
       \end{array} \right )
       \begin{array}{c}
          \alpha_1 \\
          \alpha_2 \\
          \alpha_3 \\
          \beta_1 \\
          \dots \\
       \end{array}, \quad
       s^*_{\beta_1}
      \left (
       \begin{array}{c}
          -1 \\
          0 \\
          0 \\
          1 \\
       \dots \\
       \end{array} \right ) =
      \left (
       \begin{array}{c}
          0 \\
          1 \\
          1 \\
          -1 \\
       \dots \\
       \end{array} \right )
        \begin{array}{c}
          \alpha_1 \\
          \alpha_2 \\
          \alpha_3 \\
          \beta_1 \\
          \dots \\
       \end{array}.
     \end{split}
  \end{equation*}
  Applying $s^*_{\alpha_i}$, $s^*_{\alpha_j}$ to $\gamma^{\nabla}_{ij}(6)$
  we have the last line of \eqref{refl_dual_linkages}:
  \begin{equation*}
    %% \tiny
       \label{dual_refl_alpha_4}
    \begin{split}
     & \boxed{s^*_{\alpha_i}\gamma^{\nabla}_{ij}(6) = \gamma^{\nabla}_{ij}(4), \qquad
              s^*_{\alpha_j}\gamma^{\nabla}_{ij}(6) = \gamma^{\nabla}_{ij}(5)}   \\
       \\
      s^*_{\alpha_1}
      \left (
       \begin{array}{c}
          1 \\
          1 \\
          0 \\
          -1 \\
       \dots \\
       \end{array} \right ) & =
      \left (
       \begin{array}{c}
          -1 \\
          1 \\
          0 \\
          0 \\
       \dots \\
       \end{array} \right )
        \begin{array}{c}
          \alpha_1 \\
          \alpha_2 \\
          \alpha_3 \\
          \beta_1 \\
          \dots \\
       \end{array}, \quad
      s^*_{\alpha_2}
      \left (
       \begin{array}{c}
          1 \\
          1 \\
          0 \\
          -1 \\
       \dots \\
       \end{array} \right ) =
      \left (
       \begin{array}{c}
          1 \\
          -1 \\
          0 \\
          0 \\
       \dots \\
       \end{array} \right )
        \begin{array}{c}
          \alpha_1 \\
          \alpha_2 \\
          \alpha_3 \\
          \beta_1 \\
          \dots \\
       \end{array}, \quad
     \\ %%% end of 1-line
     \\
       s^*_{\alpha_1}
      \left (
       \begin{array}{c}
          1 \\
          0 \\
          1 \\
          -1 \\
       \dots \\
       \end{array} \right ) & =
      \left (
       \begin{array}{c}
          -1 \\
          0 \\
          1 \\
          0 \\
       \dots \\
       \end{array} \right )
        \begin{array}{c}
          \alpha_1 \\
          \alpha_2 \\
          \alpha_3 \\
          \beta_1 \\
          \dots \\
       \end{array}, \quad
      s^*_{\alpha_3}
      \left (
       \begin{array}{c}
          1 \\
          0 \\
          1 \\
          -1 \\
       \dots \\
       \end{array} \right ) =
      \left (
       \begin{array}{c}
          1 \\
          0 \\
          -1 \\
          0 \\
       \dots \\
       \end{array} \right )
        \begin{array}{c}
          \alpha_1 \\
          \alpha_2 \\
          \alpha_3 \\
          \beta_1 \\
          \dots \\
       \end{array},
       \\ %%% end of 2-line
       \\
     s^*_{\alpha_2}
      \left (
       \begin{array}{c}
          0 \\
          1 \\
          1 \\
          -1 \\
       \dots \\
       \end{array} \right ) & =
      \left (
       \begin{array}{c}
          0 \\
          -1 \\
          1 \\
          0 \\
       \dots \\
       \end{array} \right )
       \begin{array}{c}
          \alpha_1 \\
          \alpha_2 \\
          \alpha_3 \\
          \beta_1 \\
          \dots \\
       \end{array}, \quad
     s^*_{\alpha_3}
      \left (
       \begin{array}{c}
          0 \\
          1 \\
          1 \\
          -1 \\
       \dots \\
       \end{array} \right ) =
      \left (
       \begin{array}{c}
          0 \\
          1 \\
          -1 \\
          0 \\
       \dots \\
       \end{array} \right )
         \begin{array}{c}
          \alpha_1 \\
          \alpha_2 \\
          \alpha_3 \\
          \beta_1 \\
          \dots \\
       \end{array}.
    \end{split}
  \end{equation*}

 \index{linkage diagrams ! - examples, $\gamma^{\nabla}(3)$ for $E_7(a_1)$}
 \index{linkage diagrams ! - examples, $\gamma^{\nabla}(6)$ for $E_7(a_2)$}
  {\rm(ii)} Here, it suffices to prove that the label $b_1$ corresponding to the coordinate $\beta_1$
  is uniquely determined  by $\alpha_i, \alpha_j$. {\bf For the linkage diagram $\gamma^{\nabla}_{ij}(3)$,
  the statement follows from  Proposition \ref{prop_diagonal},(a), see Fig. \ref{square_n_diag},(a).}
  For example, for $E_7(a_1)$, the linkage diagrams $\gamma^{\nabla}_{ij}(3)$,
  where $\{ij\} \in \{ \{ 12 \}, \{ 13 \}, \{ 23 \} \}$,  depicted
\begin{figure}[H]
\centering
\includegraphics[scale=0.8]{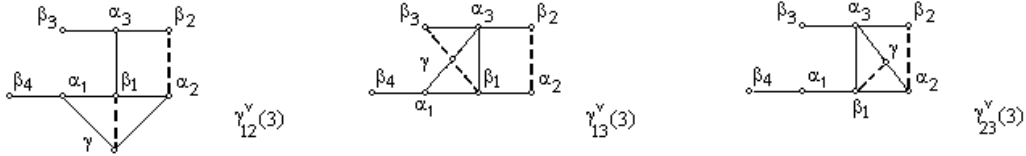}
 \caption{\hspace{3mm}The linkage diagrams $\gamma^{\nabla}_{ij}(3)$ for $E_7(a_1)$, loctets $L_{ij}^b$.}
%%%%%% The label must come after caption
\label{E7a1_gamma_ij_3}
\end{figure}
~\\
  in Fig. \ref{E7a1_gamma_ij_3}, see the linkage system $E_7(a_1)$, loctets $L_{ij}^b$  in Fig. \ref{E7a1_linkages}.
  {\bf For the linkage diagram $\gamma^{\nabla}_{ij}(6)$, the statement follows from
  Proposition \ref{prop_diagonal},(b), see Fig. \ref{square_n_diag},(b).}
  For example, for $E_7(a_2)$, the linkage diagrams $\gamma^{\nabla}_{ij}(6)$, where $\{ij\} \in \{\{ 12 \}, \{ 13 \}, \{ 23 \}\}$,
  depicted in Fig. \ref{E7a2_gamma_ij_6}, see the linkage system $E_7(a_2)$, loctets $L_{ij}^b$
  in Fig. \ref{E7a2_linkages}.
 \index{linkage diagrams ! - examples, $\gamma^{\nabla}(4)$ for $E_7(a_3)$}
\begin{figure}[H]
\centering
\includegraphics[scale=0.8]{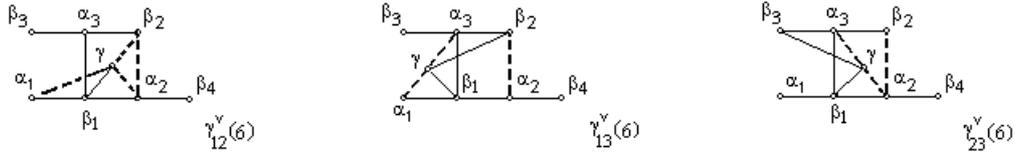}
 \caption{\hspace{3mm}The linkage diagrams $\gamma^{\nabla}_{ij}(6)$ for $E_7(a_2)$, loctets $L_{ij}^b$.}
%%%%%% The label must come after caption
\label{E7a2_gamma_ij_6}
\end{figure}
~\\
  {\bf For the linkage diagram $\gamma^{\nabla}_{ij}(4)$  and $\gamma^{\nabla}_{ij}(5)$,
  the statement follows from Proposition \ref{prop_diagonal},(e), see Fig. \ref{square_n_diag},(e).}
  For example, for $E_7(a_3)$,
  the linkage diagrams $\gamma^{\nabla}_{ij}(4)$, where $\{ij\} \in \{\{ 12 \}, \{ 13 \}, \{ 23 \}\}$,
  depicted in Fig. \ref{E7a3_gamma_ij_4}, see the linkage system $E_7(a_3)$, loctets $L_{ij}^b$
  in Fig. \ref{E7a3_linkages}.
\begin{figure}[H]
\centering
\includegraphics[scale=0.8]{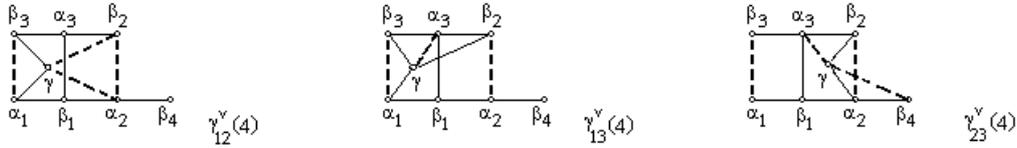}
 \caption{\hspace{3mm}The linkage diagrams $\gamma^{\nabla}_{ij}(4)$ for $E_7(a_3)$, loctets $L_{ij}^b$}
%%%%%% The label must come after caption
\label{E7a3_gamma_ij_4}
\end{figure}

 \index{dotted ! - edge}
 {\rm(iii)} The label $b_1$ corresponding to the coordinate $\beta_1$
  may take two values in $\{ -1, 0, 1 \}$ depending on the value of $a_i$,
  see Fig. \ref{fig_loctets}. Indeed, if $a_i = 1$, then $b_1 = (\gamma, \beta_1) \neq 1$, otherwise the triangle
  $\{\alpha_i, \beta_1, \gamma \}$ contains exactly two dotted edges, i.e., contains $\widetilde{A}_2$,
  contradicting Lemma \ref{lem_must_dotted}. Thus, $b_1 = -1$ or $b_1 = 0$.
  Respectively, we have linkage diagrams $\gamma^{\nabla}_{ij}(2)$ or $\gamma^{\nabla}_{ij}(8)$.
  If $a_i = -1$, then $b_1 = (\gamma, \beta_1) \neq -1$, otherwise the triangle
  $\{\alpha_i, \beta_1, \gamma \}$ does not contain any dotted edges. Thus, $b_1 = 1$ or $b_1 = 0$.
  Respectively, we have linkage diagrams $\gamma^{\nabla}_{ij}(7)$ or $\gamma^{\nabla}_{ij}(1)$:
\qed
~\\

 \index{class of diagrams ! - $\mathsf{DE4}$}
 \index{$\alpha$-label, coordinate out of $\alpha$-set}
 \index{label ! - $\alpha$-label}
 \begin{corollary}
   \label{corol_loctet}
    {\rm (i)} Any linkage diagram containing a non-zero $\alpha$-label $\alpha_1$, $\alpha_2$ or $\alpha_3$
    belongs to one of the loctets of the linkage system.

    {\rm(ii)} Any linkage diagram of the loctet uniquely determines the whole loctet.

    {\rm(iii)} If two loctets have one common linkage diagram, they coincide.

    {\rm(iv)} Every linkage diagram from the linkage system either belongs to one of the loctets
    or is $\beta$-unicolored.

 \end{corollary}
    \PerfProof  Statements (i) and (iv) follow from headings (ii), (iii)
     of Proposition \ref{prop_generating_diagr};
     statements (ii) and (iii) follow from heading (i) of Proposition \ref{prop_generating_diagr}.
    \qed
~\\

 The loctets of types
 $L_{12}$, $L_{13}$, $L_{23}$ are the {\bf main construction blocks} used for every
 linkage system, see all figures in
 Figs. \ref{D4a1_linkages}--\ref{D4_loctets}, \ref{Dk_al_linkages}, \ref{Dk_al_wind_rose}.
 By Proposition \ref{prop_generating_diagr}, any linkage diagram of a loctet
 gives rise to the whole loctet.
 By Theorem \ref{th_B_less_2}, to obtain all loctets associated with the given Carter diagram
 it suffices to find linkage diagrams $\gamma^{\nabla}_{ij}(n)$
 for a certain fixed $n \in \{1, 2, \dots, 8 \}$  satisfying the inequality:
\begin{equation*}
       \mathscr{B}^{\vee}_{\Gamma}(\gamma^{\nabla}_{ij}(n)) < 2.
\end{equation*}
The number of different loctets is defined by a number of different
linkage diagrams $\gamma^{\nabla}(n)$ for a given fixed $n$, where
$1 \leq n \leq 8$. In what follows, we enumerate loctets by linkage
diagrams $\gamma^{\nabla}(8)$. In \S\ref{sec_calc_gamma_8}, as an
example, we show how to calculate linkage diagrams
$\gamma^{\nabla}(8)$  of all loctets of $E_6(a_1)$. From Tables
\ref{sol_inequal_1}--\ref{sol_inequal_6} one can recover the
complete calculation of all linkage diagrams $\gamma^{\nabla}(8)$ of
all loctets for Carter diagrams $\Gamma \in \mathsf{C4} \coprod
\mathsf{DE4}$. The linkage diagrams $\gamma^{\nabla}_{ij}(6)$ for
every component and every loctet are listed in Table
\ref{tab_seed_linkages_6} for all linkage systems with $l < 8$.

 \index{$B_{\Gamma}$, partial Cartan matrix}
 \index{Cartan matrix ! - partial Cartan matrix $B_{\Gamma}$}
 The partial Cartan matrices $B_{\Gamma}$ and the inverse matrices $B^{-1}_{\Gamma}$
 for all Carter diagrams are presented in
 Tables \ref{tab_partial Cartan_1}--\ref{tab_partial Cartan_Al}.

 For Carter diagrams $A_l$, the technique of loctets does not work. In this case
 we construct the linkage system $\mathscr{L}(A_l)$ by induction, see \S\ref{sec_stratif_A}.

%% file: 6enumerat.tex
\section{\sc\bf Enumeration of linkage diagrams, loctets and linkage systems}
  \label{sec_enum_loctets}
  \index{class of diagrams ! - $\mathsf{C4}$}

 In this section we demonstrate some calculation examples of linkage diagrams, loctets and
 linkage systems. These calculations are based on findings of \S\ref{sec_loctets_diagr}.

\subsection{Calculation of linkage diagrams $\gamma^{\nabla}(8)$}
  \label{sec_calc_gamma_8}
  \index{loctet ! - $8$th linkage diagram}
  \index{linkage diagrams ! - calculation, $\gamma^{\nabla}(8)$}
   In this section, we consider only case $l < 8$. For case $D_l$, $l \geq 8$, see \S\ref{sec_D_lg8}, and for case
   $A_l$, $l \geq 8$, see \S\ref{sec_Al_lg8}.

 It seems a little easier to calculate the $8$th linkage diagram (we calculate it for every Carter diagram loctet) rather than to calculate any other linkage diagram of
 a loctet since $8$th linkage diagram contains $3$ zeroes among coordinates $\{ \alpha_1, \alpha_2, \alpha_3, \beta_1 \}$. We have
\begin{equation}
 \label{dual_linkage_1}
   \gamma^{\nabla}(8) =
 \left \{
  \begin{array}{ll}
     \hspace{-2mm}\{ a_1, a_2, a_3, 0,  b_2 \} & \text{ for } D_5(a_1),  \\
     \hspace{-2mm}\{ a_1, a_2, a_3, 0,  b_2, b_3\} & \text{ for } D_6(a_1), E_6(a_1), E_6(a_2),  \\
     \hspace{-2mm}\{ a_1, a_2, a_3, a_4, 0,  b_2 \} & \text{ for } D_6(a_2), \\
     \hspace{-2mm}\{ a_1, a_2, a_3, a_4, 0,  b_2, b_3\} & \text{ for } D_7(a_1), D_7(a_2), \\
     \hspace{-2mm}\{ a_1, a_2, a_3, 0,  b_2, b_3, b_4\} &  \text { for }  E_7(a_1)$,  $E_7(a_2)$, $E_7(a_3)$, $E_7(a_4),
  \end{array}
  \right . \\
\end{equation}
   where $a_i = 0,  a_j = 0, a_k = 1, b_1 = 0$ and
   $\{i, j, k \} = \{1, 2, 3 \}$  for type $L_{ij}$.
   The quadratic form $\mathscr{B}^{\vee}_{\Gamma}$ applied to the linkage label vector $\gamma^{\nabla}(8)$
 gives rise to a function $F$ of several variables $b_2, b_3, \dots$, where three dots mean $a_i$, $b_i$ for $i \geq 4$
 (if the $\alpha$-set (resp. $\beta$-set) contains $> 3$ labels), i.e.,
 \begin{equation*}
   F(b_2, b_3, \dots) := \mathscr{B}^{\vee}_{\Gamma}(\gamma^{\nabla}(8)).
 \end{equation*}
 Let $q(b_2, b_3, \dots)$ be the quadratic part of $F$ (containing only quadratic terms),
 $l(b_2, b_3, \dots)$ the linear part of $F$, and $f$ the free term:
 \begin{equation*}
        F(b_2, b_3, \dots) = q(b_2, b_3, \dots) + l(b_2, b_3, \dots) + f.
 \end{equation*}
 It is easily to see that $q(b_2, b_3, \dots)$ is the same for loctet types $L_{12}$, $L_{13}$, $L_{23}$.
 Consider the case, where the $\alpha$-set (resp. $\beta$-set) contains $\leq 3$ labels.
 The quadratic form $q(b_2, b_3)$ is determined by the principal submatrix of $B^{-1}_{\Gamma}$
 associated with coordinates $\beta_2$, $\beta_3$:
\begin{equation}
  \label{q_part}
  q(b_2, b_3) = d_{{\beta_2}{\beta_2}}b_2^2 + 2d_{{\beta_2}{\beta_3}}b_2{b_3} + d_{{\beta_3}{\beta_3}}b_3^2,
\end{equation}
 where $d_{ij}$ is the $\{i, j\}th$ element of the inverse matrix
 $B^{-1}_{\Gamma}$. The quadratic terms related to coordinates
 $\alpha_4$ or $\beta_4$ should be supplemented in the respective
 cases, see \eqref{dual_linkage_1}. The linear part
 $l(b_2, b_3)$ and the free term $f$ are as follows:
\begin{equation*}
  l(b_2, b_3) =
    \begin{cases}
      2(d_{{\alpha_1}{\beta_2}}b_2 + d_{{\alpha_1}{\beta_3}}b_3) \text{ for } L_{23}, \\
      2(d_{{\alpha_2}{\beta_2}}b_2 + d_{{\alpha_2}{\beta_3}}b_3) \text{ for } L_{13}, \\
      2(d_{{\alpha_3}{\beta_2}}b_2 + d_{{\alpha_3}{\beta_3}}b_3) \text{ for } L_{12}, \\
    \end{cases} \quad
  f =
    \begin{cases}
      d_{{\alpha_1}{\alpha_1}}  \text{ for } L_{23}, \\
      d_{{\alpha_2}{\alpha_2}}  \text{ for } L_{13}, \\
      d_{{\alpha_2}{\alpha_3}}  \text{ for } L_{12}. \\
    \end{cases}
\end{equation*}
  In section \S\ref{calc_examp_E6a1}, we give an example of the calculations for the case $E_6(a_1)$.
  By Tables \ref{sol_inequal_1}--\ref{sol_inequal_4}
  one can recover calculations for the remaining diagrams $\Gamma \in \mathsf{C4}$.

 \subsubsection{Calculation example for diagram $E_6(a_1)$}
   \label{calc_examp_E6a1}
   \index{linkage diagrams ! - examples for $E_6(a_1)$}
 The matrix $B_{\Gamma}^{-1}$ for $\Gamma = E_6(a_1)$ is given in Table \ref{tab_partial Cartan_1}.
 Here, by \eqref{q_part}:
 \begin{equation}
  \label{eq_q_E6a1}
   q(b_2, b_3) = \frac{4}{3}(b_2^2 + b_2{b_3} + b_3^2).
 \end{equation}
 In all cases below we use the fact that $b_2$, $b_3$ take values in $\{-1, 0, 1\}$, see Remark \ref{rem_values}.
~\\

  a) {\it Loctets $L_{12}$, $\gamma^{\nabla}_{12}(8) = \{0, 0, 1, 0, b_2, b_3 \}$.}
\begin{equation}
  \label{loctet_calc_E6a1}
  \begin{split}
    \mathscr{B}^{\vee}_{\Gamma}(\gamma^{\nabla}_{12}(8)) =
    &  \frac{1}{3}(10 + 2(4b_2 + 5b_3) + 4(b_2^2 + b_2{b_3} + b_3^2)) < 2, \quad \text {i.e.,} \\
    &   \frac{1}{3}(2(b_2 + b_3)^2 + 2(b_2 + 2)^2 + 2(b_3 + \frac{5}{2})^2 - \frac{21}{2}) < 2, \\
    &   \frac{2}{3}((b_2 + b_3)^2 + (b_2 + 2)^2 + (b_3 + \frac{5}{2})^2) < 2 + \frac{7}{2} = \frac{11}{2}, \\
    &  (b_2 + b_3)^2 + (b_2 + 2)^2  + (b_3 + \frac{5}{2})^2 < \frac{33}{4} = 8\frac{1}{4}. \\
   \end{split}
\end{equation}
 Recall that  $b_2, b_3 \in \{-1, 0, 1\}$. By \eqref{loctet_calc_E6a1} we have $b_2 \neq 1$, i.e.,
 $b_2 = -1$ or $b_2 = 0$. For $b_2 = 0$ we get
\begin{equation*}
     b_3^2 + (b_3 + \frac{5}{2})^2 < 4\frac{1}{4}.
\end{equation*}
  Only $b_3 = -1$ is suitable for this case. For $b_2 = -1$ we get
\begin{equation*}
     (b_3-1)^2 + (b_3 + \frac{5}{2})^2 < 7\frac{1}{4}.
\end{equation*}
 Again, $b_3 = -1$ is the single suitable solution.
 We get
 \begin{equation*}
   \gamma^{\nabla}_{12}(8) =  \{ 0, 0, 1, 0, 0, -1 \}  \text{ or } \gamma^{\nabla}_{12}(8) = \{ 0, 0, 1, 0, -1, -1 \}.
\end{equation*}
~\\

  b) {\it Loctets $L_{13}$, $\gamma^{\nabla}_{13}(8) = \{0, 1, 0, 0, b_2, b_3 \}$.}
\begin{equation*}
  \begin{split}
    \mathscr{B}^{\vee}_{\Gamma}(\gamma^{\nabla}_{13}(8)) =
    &  \frac{1}{3}(4 + 2(-b_2 + b_3) + 4(b_2^2 + b_2{b_3} + b_3^2)) < 2, \quad \text {i.e.,} \\
    &  \frac{1}{3}(2(b_2 + b_3) + 2(b_2 - \frac{1}{2})^2 + 2(b_3 + \frac{1}{2})^2 + 4 - \frac{1}{2} - \frac{1}{2})) < 2,  \\
    &  (b_2 + b_3)^2 + (b_2 - \frac{1}{2})^2  + (b_3 + \frac{1}{2})^2 < \frac{3}{2}. \\
   \end{split}
\end{equation*}
We get
\begin{equation*}
    \gamma^{\nabla}_{13}(8) =  \{ 0, 1, 0, 0, 0, 0 \} \text{ or } \gamma^{\nabla}_{13}(8) = \{ 0, 1, 0, 0,  1, -1 \}.
\end{equation*}
~\\

c) {\it Loctets $L_{23}$, $\gamma^{\nabla}_{23}(8) = \{1, 0, 0, 0, b_2, b_3 \}$.}
\begin{equation*}
  \begin{split}
    \mathscr{B}^{\vee}_{\Gamma}(\gamma^{\nabla}_{23}(8)) =
    &  \frac{1}{3}(4 + 2(b_2 + 2b_3) + 4(b_2^2 + b_2{b_3} + b_3^2)) < 2, \quad \text {i.e.,} \\
    &  (b_2 + b_3)^2 + (b_2 + \frac{1}{2})^2  + (b_3 + 1)^2 < \frac{9}{4}. \\
   \end{split}
\end{equation*}
 Here,
\begin{equation*}
   \gamma^{\nabla}_{23}(8) =  \{ 0, 1, 0, 0, 0, 0 \} \text{ or } \gamma^{\nabla}_{23}(8) = \{ 0, 1, 0, 0,  0, -1 \}.
\end{equation*}
~\\
 The corresponding $6$ loctets of the linkage system $E_6(a_1)$ are depicted in Fig. \ref{E6a1_linkages}.

\subsection{Calculation of the $\beta$-unicolored linkage diagrams}
  \label{sec_calc_homog}
  \index{linkage diagrams ! - calculation, $\beta$-unicolored}

 Now we consider $\beta$-unicolored linkage diagrams.
 Let  $a_i = (\alpha_i, \gamma)$ (resp. $b_i = (\beta_i, \gamma)$) be coordinates
 of linkage label vector $\gamma^{\nabla}$. For $\beta$-unicolored linkage diagram $\gamma^{\nabla}$,
 we have $a_i = 0$ for $i = 1,2,3$. In addition, we note that
 \begin{equation}
   \label{b1_is_0}
       b_1 = (\gamma, \beta_1) = 0,
 \end{equation}
 otherwise roots $\{\alpha_1, \alpha_2, \alpha_3, \beta_1, \gamma\}$ constitute the
 extended Dynkin diagram $\widetilde{D}_4$, which is impossible, see Lemma \ref{lem_non_ext_Dynkin}.
 Eq. \eqref{b1_is_0} holds for any Carter diagrams out $\mathsf{C4} \coprod \mathsf{DE4}$
 except for $D_4(a_1)$. Note that
 \begin{equation*}
   \mathscr{B}^{\vee}(\gamma^{\nabla}) = \mathscr{B}^{\vee}(-\gamma^{\nabla}) = q(b_2, b_3, \dots),
 \end{equation*}
 and solving the inequality $\mathscr{B}^{\vee}(\gamma^{\nabla}) < 2$
 we can assume that $b_2 > 0$, or $b_2 =0, b_3 > 0$, etc.
 We present here calculations only for the cases $E_6(a_1)$, $E_6(a_2)$ and $E_7(a_1)$.
 From Tables \ref{homog_inequal_1}--\ref{homog_inequal_3} one can recover calculations of the remaining cases.

\subsubsection{$\beta$-unicolored linkage diagrams in $\mathscr{L}(E_6(a_1))$ and $\mathscr{L}(E_6(a_2))$}
  \label{sec_homog_E6a1a2}
  \index{linkage diagrams ! - $\beta$-unicolored, in $\mathscr{L}(E_6(a_1))$ and $\mathscr{L}(E_6(a_2))$}

 For Carter diagrams $E_6(a_1)$ and $E_6(a_2)$, the $\beta$-unicolored linkage diagrams look as $\{0, 0, 0, 0, b_2, b_3\}$.
 For both diagrams $E_6(a_1)$ and $E_6(a_2)$, the $\beta$-unicolored linkage diagrams
 coincide since the principal $2\times2$ submatrices of $B^{-1}_{\Gamma}$ associated with coordinates $\beta_2$, $\beta_3$
 coincide and are as follows:
\begin{equation*}
  \label{homog_E6a1_E6a2}
   \frac{1}{3}
   \left [
  \begin{array}{cc}
     4 & 2 \\
     2 & 4 \\
  \end{array}
   \right ],
\end{equation*}
see Table \ref{tab_partial Cartan_1}. By \eqref{eq_q_E6a1} and Theorem \ref{th_B_less_2}, we have
\begin{equation}
  \label{homog_E6a1_E6a2_2}
    \begin{split}
    q(b_2, b_3) = & \frac{4}{3}(b_2^2 + b_3^2 + b_2{b_3}) =
                    \frac{4}{3}(\frac{1}{2}b_2^2 + \frac{1}{2}b_3^2 + \frac{1}{2}(b_2 + b_3)^2) < 2, \quad \text {i.e.,} \\
     & b_2^2 + b_3^2 + (b_2 + b_3)^2 < 3.
    \end{split}
\end{equation}
 There are exactly $6$ solutions of the inequality \eqref{homog_E6a1_E6a2_2}, the corresponding
linkage diagrams are:
\begin{equation}
  \label{homog_E6a1_E6a2_3}
   \begin{array}{lll}
      & \{ 0, 0, 0, 0, 0, 1 \},  & \{ 0, 0, 0, 0, 0, -1 \}, \\
      & \{ 0, 0, 0, 0, 1, 0 \},  & \{ 0, 0, 0, 0, -1, 0 \}, \\
      & \{ 0, 0, 0, 0, 1, -1 \},  & \{ 0, 0, 0, 0, -1, 1 \}. \\
   \end{array}
\end{equation}
 These $6$ linkage diagrams are located outside of the loctets in the linkage systems $\mathscr{L}(E_6(a_1))$, $\mathscr{L}(E_6(a_2))$,
 see Figs. \ref{E6a1_linkages}, \ref{E6a2_linkages}.

 \subsubsection{$\beta$-unicolored linkage diagrams in $\mathscr{L}(E_7(a_1))$}
 \index{linkage diagrams ! - $\beta$-unicolored, in $\mathscr{L}(E_7(a_1))$}

 Here, the $\beta$-unicolored linkage diagrams look as $\{0, 0, 0, 0, b_2, b_3, b_4\}$.
 The principal $3\times3$ submatrices  of $B^{-1}_{\Gamma}$ associated with coordinates $\beta_2$, $\beta_3$, $\beta_4$
 is as follows:
\begin{equation*}
  \label{homog_E7a1}
   \frac{1}{2}
   \left [
  \begin{array}{ccc}
     3 & 2 & 1 \\
     2 & 4 & 2 \\
     1 & 2 & 3 \\
  \end{array}
   \right ],
\end{equation*}
see Table \ref{tab_partial Cartan_2}. Then we have
 \begin{equation}
  \label{homog_E7a1_2}
   \begin{split}
  \mathscr{B}^{\vee}(\gamma^{\nabla}) = &
    \frac{1}{2}(3b_2^2 + 4b_3^2 + 3b_4^2 + 4b_2{b_3} + 2b_2{b_4} + 4b_3{b_4}) < 2, \quad
    \text {i.e.,}\\
  & 2(b_2 + b_3)^2 + 2(b_3 + b_4)^2 + (b_2 + b_4)^2 < 4.
  \end{split}
\end{equation}
 There are exactly $8$ solutions of the inequality \eqref{homog_E7a1_2}, the corresponding
  linkage diagrams are:
\begin{equation*}
  \label{homog_E7a1_3}
   \begin{array}{lll}
      & \{ 0, 0, 0, 0, 1, -1, 0 \},  & \{ 0, 0, 0, 0, -1, 1, 0 \}, \\
      & \{ 0, 0, 0, 0, 0, 1, -1 \},  & \{ 0, 0, 0, 0, 0, -1, 1 \}, \\
      & \{ 0, 0, 0, 0, 0, 0, 1 \},  & \{ 0, 0, 0, 0, 0, 0, -1 \}, \\
      & \{ 0, 0, 0, 0, 1, 0, 0 \},  & \{ 0, 0, 0, 0, -1, 0, 0 \}. \\
   \end{array}
\end{equation*}
 These $8$ linkage diagrams are located outside of the loctets in the linkage systems $\mathscr{L}(E_7(a_1))$,
 see Fig. \ref{E7a1_linkages}.

\subsection{Linkage systems for simply-laced Dynkin diagrams}

%% In this section we extend the previous results to simply-laced Dynkin diagrams.
%% First, we find simply-laced Dynkin diagrams such that each of them determines only one conjugacy class.

 \index{${\bf B}$, Cartan matrix associated with a Dynkin diagram}
 \index{Cartan matrix ! - ${\bf B}$}
 \index{$B_{\Gamma}$, partial Cartan matrix}
 \index{Cartan matrix ! - partial Cartan matrix $B_{\Gamma}$}
 \index{linkage diagrams ! - for Dynkin diagrams}
 In order to build the linkage systems for Dynkin diagrams $E_n$, $D_n$,
 we can use the technique of the partial Cartan matrix,
 linkage diagrams and loctets from \S\ref{sec_loctets_diagr}.
 Note that the partial Cartan matrix for Dynkin diagrams coincides with
 the usual Cartan matrix ${\bf B}$ associated with the given Dynkin diagram.
 Since $E_8$ does not have linkage diagrams, see Remark \ref{rem_l_less_8}, %% (i),
 we are interested only
 in $E_6$, $E_7$, $D_n$. In cases $E_6$, $E_7$, $D_5$, $D_6$, $D_7$, for the Cartan matrix ${\bf B}$
 and its inverse  ${\bf B}^{-1}$, see Table \ref{tab_partial Cartan_3}.
 One can recover the complete calculation  of linkage diagrams $\gamma^{\nabla}(8)$ of all loctets
 of $E_6$, $E_7$, $D_5$, $D_6$, $D_7$  by means of Tables \ref{sol_inequal_5}, \ref{sol_inequal_6}.
 The $\beta$-unicolored linkage diagrams look as follows:
\begin{equation*}
  \gamma^{\nabla} =
  \left \{
  \begin{array}{ll}
       \hspace{-2mm}\{0, 0, 0, 0, b_2, b_3 \} & \text{ for } E_6, \\
       \hspace{-2mm}\{0, 0, 0, a_4, 0, b_2, b_3 \} & \text{ for } E_7, \\
       \hspace{-2mm}\{0, 0, 0, 0, b_2 \} & \text{ for } D_5, \\
       \hspace{-2mm}\{0, 0, 0, a_4, 0, b_2 \} & \text{ for } D_6, \\
       \hspace{-2mm}\{0, 0, 0, a_4, 0, b_2, b_3 \} & \text{ for } D_7,
  \end{array}
  \right .
\end{equation*}
 see Tables \ref{homog_inequal_1}--\ref{homog_inequal_3}. The
 $\beta$-unicolored linkage diagrams are located outside of all loctets, see Figs. \ref{27_weight_diagr_E6__2comp}(top),
 \ref{E7pure_linkage_system}, \ref{D5pure_loctets}, \ref{D6pure_loctets}.

Note that for $E_6$, the principal matrix associated with coordinates $\beta_2$, $\beta_3$ coincide with
the principal matrix for the Carter diagrams $E_6(a_1)$, $E_6(a_2)$, see \S\ref{sec_homog_E6a1a2},
and, consequently, $\beta$-unicolored linkage diagrams coincide with these diagrams for
$E_6(a_1)$, $E_6(a_2)$, see $6$ solutions \eqref{homog_E6a1_E6a2_3} of the inequality \eqref{homog_E6a1_E6a2_2}.

  For $A_l$ the technique of loctets cannot be applied. The reason of this is lack of branch point.
  For $l = 5,6,7$, the linkage system $\mathscr{L}(A_l)$ consists of $A$-, $D$- and $E$-components,
  see \S\ref{sec_plis_Al}.
  For $l \leq 4$ and $l \geq 8$, the linkage system $\mathscr{L}(A_l)$ consists only of $A$- and $D$-components,
  see \S\ref{sec_Al_lg8}.  The linkage system components $\mathscr{L}_{A_{l+1}}(A_l)$ and
  $\mathscr{L}_{D_{l+1}}(A_l)$ are constructed by means of induction process, see \S\ref{sec_Al_lg8} and
  Fig. \ref{2_lines_A2nMin1} - Fig. \ref{A5_to_A6D6_pass}.

\subsubsection{Asymmetric relations between $A$-, $D$-, $E$-types of Dynkin diagrams}

 The relation between the Dynkin diagrams of $E$-type and $D$-type is asymmetric in the following sense:
 The root systems $E_n$ contains the root subsystem $D_{n-1}$ for $n = 6,7,8$;
 however, $E_6$ is not contained in $D_n$ for any n, see Remark \ref{rem_l_less_8}.

 Similarly, the relation between the Dynkin diagrams of $A$-type and Dynkin diagrams of $D$- and $E$-types is asymmetric,
 namely: The root systems of the $D$- and the $E$-type contain root subsystems $A_n$ for a certain $n$.
 The converse statement is not true: The root system of $A$-type does not contain root subsystems of
 $D$- and $E$-types, see Lemma \ref{lem_case_An}(ii).

 The asymmetric relations between $A$-, $D$-, $E$-types of Dynkin diagrams can be observed
 by means of the Table \ref{tab_val_Bmin1}. By abuse of language we can say that
 {\bf $E$-type contains $ADE$-types; $D$-type contains $AD$-types; $A$-type contains only $A$-type.}

%% file: 7weights.tex
\section{\sc\bf Theorem on coincidence of linkage and weight systems for Dynkin diagrams}
\subsection{Dominant weights and Dynkin labels}
  \label{sec_dominant_weight}

 \index{$\varPhi$, root system}
 \index{Dynkin labels}
 \index{label ! - Dynkin labels}
 \index{${\bf B}$, Cartan matrix associated with a Dynkin diagram}
 \index{Cartan matrix ! - ${\bf B}$}
 \index{$\Lambda$, weight lattice}
 \index{$\Lambda^{+}$, set of all dominant weights in $\Lambda$}
 \index{root system}
 \index{$\lambda^{\nabla}$, vector of Dynkin labels}
 \index{$\lambda$, weight in the weight lattice $\Lambda$}
 \index{$\overline{\omega}_i$, fundamental dominant weight}
 \index{weight ! - fundamental dominant weight}
 \index{weight ! - dominant weight}
 \index{weight ! - weight lattice}
 \index{dominant weight}

Let $\varPhi$ be the root system in the linear space $V$ with the Weyl group $W$
(associated with the given Dynkin diagram $\Gamma$). The set of vectors
$\Lambda$ = $\{\lambda \mid \lambda \in V\}$ such that
\begin{equation}
  \label{eq_weight_int_1}
   \langle \lambda, \alpha \rangle :=
   \frac{2(\lambda, \alpha)}{(\alpha, \alpha)} \in \mathbb{Z} \text { for any } \alpha \in \varPhi
\end{equation}
is said to be the {\it weight lattice} in $V$, and vectors $\lambda  \in \Lambda$ are called {\it weights}.
Clearly, the weight lattice contains the root system $\varPhi$, see \cite{H78, Bo02}.
Let $\Delta = \{\alpha_1,\dots \alpha_l\}$ be the set of positive simple roots of $\varPhi$,
and $\lambda$ a weight in $\Lambda$. If $(\lambda, \alpha_j) \geq 0$ for any $\alpha_j \in \Delta$,
the weight $\lambda$ is said to be {\it dominant}.  The set of all dominant weights is denoted by $\Lambda^+$,
this set lies in the fundamental Weyl chamber $\mathfrak{C}$.  The dominant weight $\overline{\omega}_i \in \Lambda^+$ is
said to be a {\it fundamental dominant weight} if
\begin{equation}
  \label{eq_weight_int_2}
   \frac{2(\overline{\omega}_i, \alpha_j)}{(\alpha_j, \alpha_j)} = \delta_{ij}.
\end{equation}
The coefficients $\{ \lambda_i \}$ of the decomposition of the
weight $\lambda \in \Lambda$ in the basis of fundamental weights $\overline{\omega}_i$ are said to be {\it Dynkin labels}.
In other words, for
\begin{equation}
  \label{eq_weight_int_3}
  \lambda = \sum\limits_{i=1}^l \lambda_i \overline{\omega}_i, \quad \text{we have} \quad
   \langle \lambda, \alpha_j \rangle = \lambda_j.
\end{equation}
For the simply-laced Dynkin diagrams\footnotemark[1], by \eqref{eq_weight_int_1} we have $(\lambda, \alpha_j) = \lambda_j$. In this case,
the {\it vector of Dynkin labels} looks as follows:
\footnotetext[1]{Recall that in this case, $(\alpha,\alpha) = 2$ for any simple root $\alpha$, see \eqref{eq_Kac}.}
\begin{equation}
  \label{eq_weight_int_4}
  \lambda^{\nabla} :=
  \left (
     \begin{array}{c}
        (\lambda, \alpha_1) \\
          \dots \\
        (\lambda, \alpha_l)
     \end{array}
  \right ) ~=~ {\bf B}\lambda,
\end{equation}
where ${\bf B}$ is the Cartan matrix associated with $\Gamma$.
Let $E$ be a module over the Weyl group $W$, $E^*$ be the dual space,
$W^{\vee}$ be the contragredient representation of $W$ in $E^*$.
The element $w^* \in W^{\vee}$ is defined as follows:
  \index{dual space $E^*$}
  \index{$E^*$, dual space to $E$, the representation space of $W^{\vee}$}
  \index{conjugation ! - of weights in the same $W$-orbit}
  \index{conjugate weight}
\begin{equation}
  \label{eq_weight_int_5}
     w^* := {}^t{w}^{-1}, \quad \text{ in particular, } \quad s^*_{\alpha} = {}^t{s_{\alpha}}.
\end{equation}
Then we have
\begin{equation}
  \label{eq_weight_int_6}
     w^*{\bf B} = {\bf B}w  ~\Longrightarrow~  {\bf B}w^{-1} = (w^*)^{-1}{\bf B} ~\Longrightarrow~
     w{\bf B}^{-1} = {\bf B}^{-1}w^*.   \\
\end{equation}
In particular,
\begin{equation}
  \label{eq_weight_int_6A}
    w^*\lambda^{\nabla} = w^*{\bf B}\lambda = {\bf B}w\lambda = (w\lambda)^{\nabla}, \\
\end{equation}
see Propositions \ref{prop_contagr_1}, \ref{prop_link_diagr_conn}.
By \eqref{eq_weight_int_1} $W$ acts on the weight lattice $\Lambda$.
By \eqref{eq_weight_int_6A} $W^{\vee}$ acts on the lattice spanned by vectors of Dynkin labels $\lambda^{\nabla}$,
where $\lambda \in \Lambda$.

 We say that a weight $\lambda$ is {\it conjugate} to a weight $\eta$, if they lie in the same orbit of $W$.
 The following proposition is one of central statements in the theory of root systems and weights.

 \begin{proposition}{\rm (}\cite[Ch.VI, $\S$1, $n^{\rm o}{10}$]{Bo02}{\rm )\footnotemark[2]}
   \label{prop_conjugate}
 For every weight $\lambda \in \Lambda$, there is one and only one fundamental weight $\overline{\omega}_i$
 conjugate to $\lambda$ .
  \qed
 \end{proposition}
 \index{Weyl group ! - $W$}
 \index{$W$, Weyl group}
 \index{weight system ! - (= weight diagram)}
 \index{weight system ! - of representations ${\bf 27}$ and $\overline{\bf 27}$ for $E_6$}
 \index{weight ! - conjugation}
 \index{fundamental weight}
 \index{fundamental representation}
 \index{$\mathfrak{g}$, simple Lie algebra}
 \index{highest weight of the representation}
 \index{weight ! - highest}
  \footnotetext[2]{This proposition follows from the fact the Weyl group $W$ acts simply transitively on the set of all Weyl chambers,
  see \cite[Ch.VI, $\S$1, $n^{\rm o}5$, Theorem 2]{Bo02}, \cite{Sn90}.}

 The fundamental weights $\overline{\omega}_i$, where $i=1,\dots,l$, generate $l$ non-intersecting $W$-orbits.
 These orbits are the components of the diagram we call {\it weight system}\footnotemark[3].
 The weight systems arise in the representation theory of semisimple Lie algebras, see \cite{Bo05, Sl81}.
 Let $\mathfrak{g}$ be the simple Lie algebra associated with the Dynkin diagram $\Gamma$.
 The irreducible finite-dimensional representation $V$ of $\mathfrak{g}$ is said to be {\it fundamental representation}
 if the highest weight of $V$ is a fundamental weight, see \cite{Bo05}.

  \footnotetext[3]{In the literature (see, for example, \cite{Va00}, \cite{Ch84}), the term
  \underline{weight diagram} is often used instead of the term \underline{weight system}. However, the term {\lq\lq}diagram{\rq\rq}
  is heavily overloaded in our context.} For weight systems of two dual $27$-dimensional fundamental representations
  of the semisimple Lie algebra $E_6$, see Fig. \ref{27_weight_diagr_E6__2comp}.
 In particle physics these representations are denoted by ${\bf 27}$ and $\overline{\bf 27}$, see \cite{Sl81}.

\subsection{Relationship between linkage system and weight system for Dynkin diagrams}

  Let $\overline{e}_i$ be the unit vector with a $1$ in the $i$th slot: $(\overline{e}_i)_j = \delta_{ij}$.

\begin{figure}[H]
\centering
\includegraphics[scale=0.4]{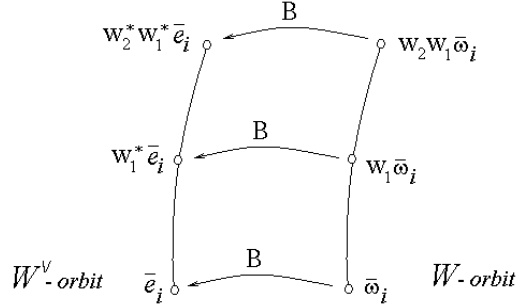}
\caption{
\hspace{3mm}The one-to-one correspondence between $W$-orbit and $W^{\vee}$-orbit}
%%%%%% The label must come after caption
\label{orbit_weights}
\end{figure}

 \begin{lemma}
   \label{lem_bijection_orbits}
 {\rm(i)} There is a one-to-one correspondence between the $W$-orbit of
  the fundamental weight $\overline{\omega}_i$ and the $W^{\vee}$-orbit of the unit vector $\overline{e}_i$.

 {\rm(ii)}
  The correspondence between $W$- and $W^{\vee}$-orbits is carried out by Cartan matrix ${\bf B}$ as follows:
  \begin{equation}
    \label{eq_corres_orbits}
    \begin{split}
      & {\bf B}(\overline{\omega}_i) = \overline{e}_i, \\
      & {\bf B}(w_1\overline{\omega}_i) = w^*_1\overline{e}_i, \\
      & {\bf B}(w_2{w}_1\overline{\omega}_i) = w^*_2w^*_1\overline{e}_i,
    \end{split}
  \end{equation}
 where $w_1$ and $w_2$ are any elements in $W$, see Fig. $\ref{orbit_weights}$.
 \end{lemma}

  \PerfProof
  By \eqref{eq_weight_int_2}
 \begin{equation}
  \label{eq_weight_int_7}
    \overline{\omega}_i = {\bf B}^{-1} \overline{e}_i.
\end{equation}
  By \eqref{eq_weight_int_6}, we have
 \begin{equation}
  \label{eq_relat_lw_2}
    w\overline{\omega}_i = w{\bf B}^{-1} \overline{e}_i = {\bf B}^{-1}(w^{*}\overline{e}_i).
\end{equation}
 Relation \eqref{eq_relat_lw_2} yields the desired correspondence.
 The first and second relations from \eqref{eq_corres_orbits} follow from \eqref{eq_weight_int_7} and \eqref{eq_relat_lw_2}.
 Further, by \eqref{eq_weight_int_6}
 \begin{equation*}
   {\bf B}(w_2{w}_1\overline{\omega}_i) = w^*_2{\bf B}w_1\overline{\omega}_i = w^*_2{w}^{*}_1{\bf B}\overline{\omega}_i =
   w^*_2{w}^{*}_1\overline{e}_i.
 \end{equation*}
 \qed

  \begin{theorem}[On linkage systems for Dynkin diagrams]
   \label{th_linkages_n_weights}
  Let $\Gamma$ be a simply-laced Dynkin diagram,  $\mathfrak{g}$ the simple Lie algebra associated with $\Gamma$.
  Every $A$-, $D$- or $E$-component of the linkage system $\mathscr{L}(\Gamma)$ coincides with
  a weight system of one of fundamental representations of $\mathfrak{g}$.
  \end{theorem}

  \PerfProof
  The fact $\gamma^{\nabla} \in \mathscr{L}(\Gamma)$ means that
   \begin{equation*}
      \gamma^{\nabla} =
        \left (
          \begin{array}{c}
            (\gamma, \tau_1) \\
            \dots \\
            (\gamma, \tau_l) \\
          \end{array}
        \right )
  \end{equation*}
  for a certain root $\gamma$ not lying in the space $L$ spanned by the $\Gamma$-associated root subset
  $S = \{\tau_1,\dots,\tau_l\}$, see \eqref{dual_gamma}.
  In addition,
  \begin{equation*}
      \gamma = \gamma_L + \mu, \qquad \gamma^{\nabla} = {\bf B}\gamma_L.
  \end{equation*}
  where $\gamma_L$ is the projection of $\gamma$ on the linear subspace $L = [\tau_1,\dots,\tau_l]$
  and $\mu$ depends on the type of the extension, see \S\ref{sec_inverse_qvadr} and Proposition \ref{prop_unique_val}.
  It suffices to prove that any $\gamma^{\nabla}$
  lies in the $W^{\vee}$-orbit of a certain unit vector $\overline{e}_i$.
  For every simply root $\tau_i \in \varPhi$,
  the label $\gamma^{\nabla}_i = (\gamma, \tau_i) = (\gamma_L, \tau_i)$ takes value in $\{-1, 0, 1\}$,
  see \S\ref{sec_linkage}. Then $(\gamma_L, \alpha) \in \mathbb{Z}$ for any root $\alpha \in \varPhi$.
  By \eqref{eq_weight_int_1}, since $(\alpha, \alpha) = 2$ we have
  $\langle \gamma_L, \alpha \rangle  \in \mathbb{Z}$ for any root $\alpha \in \varPhi$.
  Thus, the \underline{projection $\gamma_L$ of any linkage root $\gamma$ is a weight}.
  Then, by Proposition \ref{prop_conjugate}
  $\gamma_L$ is conjugate to some fundamental weight $\overline{\omega}_i$:
  \begin{equation}
    \label{eq_orbit_omega}
      \gamma_L = w\overline{\omega}_i.
  \end{equation}
  Finally, by \eqref{eq_relat_lw_2}
  \begin{equation*}
     \gamma^{\nabla} = {\bf B}\gamma_L = {\bf B}w\overline{\omega}_i = {\bf B}{\bf B}^{-1}w^*\overline{e}_i = w^*{\overline{e}_i},
  \end{equation*}
  i.e., $\gamma^{\nabla}$ lies in the $W^{\vee}$-orbit of $\overline{e}_i$:
  \begin{equation}
    \label{eq_gamma_nab_in_orbit}
     \gamma^{\nabla} =  w^*{\overline{e}_i},
  \end{equation}
  and the weight $w\overline{\omega}_i$ (on the $W$-orbit of the
  fundamental weight $\overline{\omega}_i$) uniquely corresponds to $\gamma^{\nabla}$:
  \begin{equation*}
     w\overline{\omega}_i = {\bf B}^{-1}\gamma^{\nabla}.
  \end{equation*}
\qed

  \begin{corollary}
    \label{col_one_unicolored}
    For any simply-laced Dynkin diagram $\Gamma$, every $A$-, $D$- or $E$-component of
    the linkage system $\mathscr{L}(\Gamma)$ contains a unique $\beta$-unicolored or $\alpha$-unicolored
    linkage diagram $\gamma^{\nabla}$ such that $\gamma_L = {\bf B}^{-1}\gamma^{\nabla}$ coincides with one of
    fundamental weights $\overline{\omega}_i$.
  \end{corollary}

     \PerfProof  Let us take $\gamma^{\nabla} =  \overline{e}_i$, see \eqref{eq_gamma_nab_in_orbit}.
      The vector $\overline{e}_i$ has only one non-zero coordinate, i.e., $\overline{e}_i$ is $\beta$-unicolored or $\alpha$-unicolored.
      By Theorem \ref{th_linkages_n_weights} and Proposition \ref{prop_conjugate} $\gamma^{\nabla}$ is unique.
      \qed

  \begin{remark}{\rm
     If $\Gamma$ is a Carter diagram but not a Dynkin diagram, Corollary \ref{col_one_unicolored} does not
     hold. For example, for the linkage system $\mathscr{L}(D_6(a_1))$, the first $64$-element component
     does not contain such linkage diagrams at all, see Fig. \ref{D6a1_linkages}(left),  and the second
     $64$-element component contains $2$ such linkage diagrams:
     \begin{equation*}
        \gamma^{\nabla}_{13}(8) = \{0, 1, 0, 0, 0, 0\} \in L_{13}^c, \quad
        \gamma^{\nabla}_{12}(8) = \{0, 0, 1, 0, 0, 0\} \in L_{12}^c,
     \end{equation*}
     see Fig. \ref{D6a1_linkages}(right).
     }
  \end{remark}

%% file: 8D_partialDlak.tex
%%~\\
\section{\sc\bf $D$-type linkage systems}
  \label{sec_stratif_D}

 \subsection{The linkage systems  $\mathscr{L}(D_l)$ and $\mathscr{L}(D_l(a_k))$ for $l \geq 8$}
    \label{sec_D_lg8}
    \index{linkage system  ! - $\mathscr{L}(D_l)$ and $\mathscr{L}(D_l(a_k))$ for $l \geq 8$}
    \index{Dynkin extension ! - $D_l <_D D_{l+1}$ and $D_l(a_k) <_D D_{l+1}$ for $l \geq 8$}
    \index{$\Gamma <_D \Gamma'$, Dynkin extension}
    \index{root stratum}

\begin{figure}[H]
\centering
\includegraphics[scale=0.85]{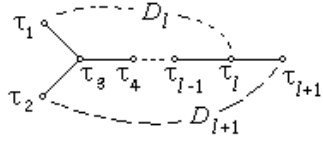}
\caption{
\hspace{3mm}The Dynkin extension $D_l <_D D_{l+1}$}
%%%%%% The label must come after caption
\label{Dl_ext_Dl1}
\end{figure}

 In this section, we assume that $\Gamma$ is one of the Carter diagrams $D_l$ or $D_l(a_k)$,
 see Fig. \ref{Dl_ext_Dl1}$(a)$, $(b)$. Let $S$ be a $\Gamma$-associated root subset.
\begin{lemma}
  \label{lem_coinc_linkage}
    Let $D_{l+1}$ be the Dynkin extension of the Carter diagram $D_l$, and $\tau_{l+1}$ a
    simple positive root from the root stratum $\varPhi(D_{l+1})\backslash\varPhi(S)$:
 \begin{equation}
      \varPhi(S) = \{ \tau_1, \dots, \tau_l \}, \quad
      \varPhi(D_{l+1}) = \{ \tau_1, \dots, \tau_l, \tau_{l+1} \}.
 \end{equation}
    Let $\varphi$ be a positive root in the root stratum $\varPhi(D_{l+1})\backslash\varPhi(S)$,
    and $\mu_{max}$ be maximal root in $\varPhi(D_{l+1})$.

     {\rm (i)} The vector
   \begin{equation}
     \label{eq_Dl_ext1}
      \delta = \mu_{max} - \varphi + \tau_{l+1}
   \end{equation}
   is also a root in $\varPhi(D_{l+1})\backslash\varPhi(S)$.

  \index{$L$, linear space ! $= [\tau_1,\dots,\tau_l]$}
  \index{quadratic form ! - $\mathscr{B}$, quadratic Tits form associated with the Cartan matrix ${\bf B}$}
  \index{$\mathscr{B}$, quadratic Tits form associated with the Cartan matrix ${\bf B}$}
    %% Let $L = [\tau_1, \dots, \tau_l]$.
    {\rm (ii)}
    The linkage label vectors $\varphi^{\nabla}$ and $-\delta^{\nabla}$ coincide:
   \begin{equation}
     \label{eq_Dl_ext2}
      \delta^{\nabla} = -\varphi^{\nabla}.
   \end{equation}
\end{lemma}

  \PerfProof  Let $\mathscr{B}$ be the quadratic Tits form associated with $D_{l+1}$.

  (i)  By \eqref{eq_Kac}
  it suffices to prove that
  \begin{equation}
    \label{eq_Kac_Dlak}
    \mathscr{B}(\delta) = 2.
  \end{equation}
   If $\varphi = \mu_{max}$ (resp. $\tau_{l+1}$), then $\delta = \tau_{l+1}$ (resp. $\mu_{max}$).
   In both cases, $\delta$ is the root in $\varPhi(D_{l+1})\backslash\varPhi(S)$. Suppose $\varphi \neq \mu_{max}, \tau_{l+1}$.
   We need to prove that
   \begin{equation}
     \label{eq_Dl_ext3}
     \begin{split}
       & \mathscr{B}(\mu_{max} - \varphi + \tau_{l+1}) = 2, \text{ i.e.}, \\
       & \mathscr{B}(\mu_{max}) + \mathscr{B}(\varphi) + \mathscr{B}(\tau_{l+1}) +
           2(\mu_{max}, \tau_{l+1}) - 2(\varphi, \mu_{max} + \tau_{l+1}) = 2.
     \end{split}
   \end{equation}
   Since $\varphi$, $\tau_{l+1}$, $\mu_{max}$ are roots,
   we have  $\mathscr{B}(\varphi) = \mathscr{B}(\tau_{l+1}) = \mathscr{B}(\mu_{max}) = 2$.
   Then the last equality in \eqref{eq_Dl_ext3} is equivalent to
  \begin{equation*}
     %%\label{eq_Dl_ext4}
     \begin{split}
       & 4 +  2(\mu_{max}, \tau_{l+1}) - 2(\varphi, \mu_{max} + \tau_{l+1}) = 0, \text{ i.e., } \\
       & (\varphi, \mu_{max} + \tau_{l+1}) - (\mu_{max}, \tau_{l+1}) = 2.
     \end{split}
   \end{equation*}
   Further, $(\mu_{max}, \tau_i) = 0$ for any $i \neq l$. In particular, $(\mu_{max}, \tau_{l+1}) = 0$ and
   it suffices to prove that
   \begin{equation}
     \label{eq_Dl_ext5}
       (\varphi, \mu_{max} + \tau_{l+1}) = 2 \text{ for any } \varphi \in \varPhi(D_{l+1})\backslash\varPhi(S).
   \end{equation}
   We have
   \begin{equation}
      \label{eq_Dl_ext6}
     (\gamma, \mu_{max} + \tau_{l+1}) =
       \begin{cases}
           2 \text{ for } \gamma = \tau_{l+1}, \text{ since } \mu_{max} \perp \tau_{l+1}, \\
           0 \text{ for } \gamma = \tau_l, \text{ since } (\tau_l, \mu_{max}) = 1, (\tau_l, \tau_{l+1}) = -1, \\
           0 \text{ for } \gamma = \tau_i, \text{ where } i < l.
       \end{cases}
   \end{equation}

   The fact that $\varphi \in \varPhi(D_{l+1})\backslash\varPhi(S)$
   and $\varphi$ is positive means
   that $\tau_{l+1}$ enters with coefficient $1$ into the decomposition of $\varphi$ with relation to $\{\tau_1, \dots, \tau_{l+1}\}$.
   Then, by \eqref{eq_Dl_ext6} the relation \eqref{eq_Dl_ext5} holds for
   any root $\varphi \in \varPhi(D_{l+1})\backslash\varPhi(S)$.
   Therefore, \eqref{eq_Kac_Dlak} holds.  In other words, $\delta$ is also a root.
   By \eqref{eq_Dl_ext1} $\tau_{l+1}$ also enters with coefficient $1$ into the decomposition of $\delta$
   with relation to $\{\tau_1, \dots, \tau_{l+1}\}$.
   Thus, $\delta \in \varPhi(D_{l+1})\backslash\varPhi(S)$.
~\\

   (ii) By \eqref{eq_Dl_ext1} we have
   \begin{equation*}
       (\delta, \tau_i) = (\mu_{max} + \tau_{l+1}, \tau_i) - (\varphi, \tau_i), \text{ where } 1 \leq i \leq l.
   \end{equation*}
   By \eqref{eq_Dl_ext6} for $i \leq l$, we have
    \begin{equation}
     \label{ex_def_beta}
        \delta^{\nabla}_i = (\delta, \tau_i) =  -(\varphi, \tau_i) = -\varphi^{\nabla}_i, \text{ where } 1 \leq i \leq l.
   \end{equation}
   %% Since $\delta^{\nabla}$ has only coordinates $\delta^{\nabla}_i$ for $1 \leq i \leq l$,
   %% the coordinate $\tau_{l+1}$ does not matter.
   %% Thus,
   Since $L = [\tau_1, \dots, \tau_l]$, we have $\delta^{\nabla} = -\varphi^{\nabla}.$
  \qed
~\\

  Now, we can calculate the size of the linkage systems $\mathscr{L}(D_l)$ and $\mathscr{L}(D_l(a_k))$ for $l \geq 8$.
  For sizes of the linkage systems $\mathscr{L}(D_5)$,  $\mathscr{L}(D_6)$, $\mathscr{L}(D_7)$,
  and $\mathscr{L}(D_5(a_k))$,  $\mathscr{L}(D_6(a_k))$, $\mathscr{L}(D_7(a_k))$,
  see \S\ref{sec_plis_Dl_Dlak}.

\begin{figure}[H]
\centering
\includegraphics[scale=0.5, angle=270]{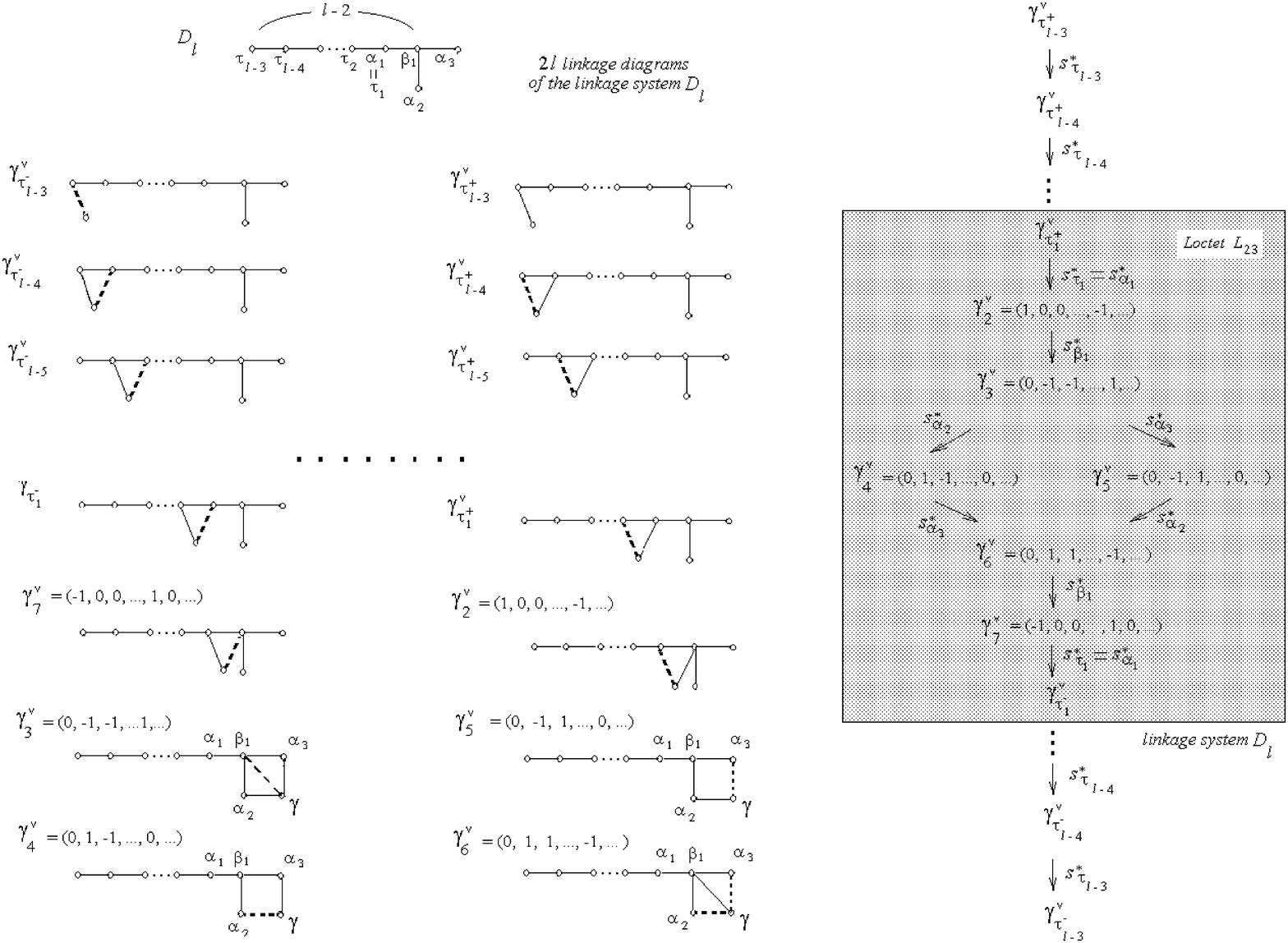}
\caption{\hspace{3mm}The linkage system $D_l$ for $l > 7$, $2l$ linkages.
 The single loctet $L_{23}$ is depicted in the shaded area}
%%%%%% The label must come after caption
\label{Dlpu_linkages}
\end{figure}

  \begin{corollary}
    \label{cor_size_pls}
    {\rm(i)}
      For any $l$, the size of the linkage system component $\mathscr{L}_{D_{l+1}}(D_l)$
      (resp. $\mathscr{L}_{D_{l+1}}(D_l(a_k))$, where $1 \leq k \leq \big [\frac{l-2}{2} \big ]$),
      is equal to $2l$.

    {\rm(ii)}
      For $l=4$ and $l \geq 8$, the size of the linkage system $\mathscr{L}(D_l)$
      (resp. $\mathscr{L}(D_l(a_k))$, where $1 \leq k \leq \big [\frac{l-2}{2} \big ]$), is equal to $2l$.
  \end{corollary}
    \PerfProof (i) The number of roots of $\varPhi(D_l)$ is $2l(l-1)$, see \cite[Table IV]{Bo02}.
    Then by \eqref{eq_regext4} and Corollary \ref{corol_Estim_2} we have
   \begin{equation*}
        \mid \mathscr{L}_{D_{l+1}}(D_l(a_k)) \mid ~=~ \mid \mathscr{L}_{D_{l+1}}(D_l) \mid ~\leq~ 2(l+1)l - 2l(l-1) = 4l. \\
   \end{equation*}
   By Lemma \ref{lem_coinc_linkage}(ii), this size is twice less:
    \begin{equation*}
        \mid \mathscr{L}_{D_{l+1}}(D_l(a_k)) \mid ~=~ \mid \mathscr{L}_{D_{l+1}}(D_l) \mid ~\leq~ 2l.
   \end{equation*}
%%~\\
      In Fig. \ref{Dlpu_linkages}, we present $2l$ linkage label vectors for the Carter diagram $D_l$.
      So, we conclude that there are \underline{exactly $2l$ linkage label vectors} in
      the $D$-component $\mathscr{L}_{D_{l+1}}(D_l)$ and, therefore, in the $D$-component $\mathscr{L}_{D_{l+1}}(D_l(a_k))$.
~\\

  (ii)
      For $l=4$ and $l \geq 8$, we have
         $\mathscr{L}(D_l) = \mathscr{L}_{D_{l+1}}(D_l)$ and
         $\mathscr{L}(D_l(a_k)) = \mathscr{L}_{D_{l+1}}(D_l(a_k))$.
      So, $\mid \mathscr{L}(D_l) \mid = 2l$ and $\mid \mathscr{L}(D_l(a_k)) \mid = 2l$.
      For $l = 5,6,7$, this is not so, since
       \begin{equation*}
         \mathscr{L}(D_l) = \mathscr{L}_{E_{l+1}}(D_l) \cup  \mathscr{L}_{D_{l+1}}(D_l)
           \text{ and }
         \mathscr{L}(D_l(a_k)) = \mathscr{L}_{E_{l+1}}(D_l(a_k)) \cup  \mathscr{L}_{D_{l+1}}(D_l(a_k)). \qed
      \end{equation*}

 \subsection{The linkage systems  $\mathscr{L}(D_l)$ and $\mathscr{L}(D_l(a_k))$ for $l=5,6,7$}
    \label{sec_plis_Dl_Dlak}
    \index{linkage system  ! - $\mathscr{L}(D_l)$ and $\mathscr{L}(D_l(a_k))$ for $l < 8$}
    \index{Dynkin extension ! - $D_l <_D D_{l+1}$ and $D_l(a_k) <_D D_{l+1}$ for $l < 8$}
    \index{$\Gamma <_D \Gamma'$, Dynkin extension}

 For $l=5, 6, 7$, in addition to Dynkin extensions $D_l <_D D_{l+1}$, there are Dynkin extensions of $E$-type,
 namely
 \begin{equation*}
    D_5 <_D  E_6,  \quad D_6 <_D E_7, \quad D_7 <_D  E_8.
 \end{equation*}
 The corresponding $E$-components are $\mathscr{L}_{E_6}(D_{5})$,  $\mathscr{L}_{E_7}(D_{6})$, and $\mathscr{L}_{E_8}(D_{7})$.

 \subsubsection{The $D$- and $E$-components $\mathscr{L}_{D_6}(D_5)$, $\mathscr{L}_{E_6}(D_5)$,
   $\mathscr{L}_{D_6}(D_5(a_1))$ and $\mathscr{L}_{E_6}(D_5(a_1))$}
   \index{linkage system  ! - $\mathscr{L}(D_5)$ and $\mathscr{L}(D_5(a_1))$}

 By \eqref{eq_regext4} and Corollary \ref{corol_Estim_2} we have
 \begin{equation*}
   \mid \mathscr{L}_{E_6}(D_5(a_1)) \mid = \mid \mathscr{L}_{E_6}(D_5) \mid ~\leq~ |\varPhi(E_6)| - |\varPhi(D_5)| = 72 - 40 = 32.
 \end{equation*}
~\\
 The linkage system $\mathscr{L}(D_{5})$ (resp. $\mathscr{L}(D_5(a_1)$) consists of two parts:
 \begin{equation}
   \label{eq_plises_D5}
   \begin{split}
     & \boxed{ \mathscr{L}(D_5) = \mathscr{L}_{D_6}(D_5) ~\cup~  \mathscr{L}_{E_6}(D_5)}, \\
     & \boxed{ \mathscr{L}(D_5(a_1)) = \mathscr{L}_{D_6}(D_5(a_1)) ~\cup~  \mathscr{L}_{E_6}(D_5(a_1))}.
   \end{split}
 \end{equation}
~\\
 By Corollary \ref{cor_size_pls}(i), we have
   \index{linkage diagrams ! - in $\mathscr{L}(D_5)$, $\mathscr{L}(D_5(a_1))$}
 \begin{equation}
  \begin{split}
   & \mid \mathscr{L}_{D_6}(D_5) \mid = \mid \mathscr{L}_{D_6}(D_5(a_1)) \mid = 2\times5 =  10, \\
   & \mid \mathscr{L}(D_5) \mid  ~\leq~ 10 + 32 = 42, \\
   & \mid \mathscr{L}(D_5(a_1)) \mid ~\leq~ 10 + 32 = 42.
  \end{split}
 \end{equation}
There are $42$ linkage diagrams of $\mathscr{L}(D_{5})$ presented in Fig. \ref{D5pure_loctets};
so we conclude that there are \underline{exactly $42$ diagrams in $\mathscr{L}(D_5)$},
and, therefore, in $\mathscr{L}(D_5(a_1))$, see Fig. \ref{D5a1_linkages}.

\subsubsection{The $D$- and $E$-components $\mathscr{L}_{D_7}(D_6)$, $\mathscr{L}_{D_7}(D_6(a_k))$,
   $\mathscr{L}_{E_7}(D_6)$ and $\mathscr{L}_{E_7}(D_6(a_k))$}
   \index{linkage system  ! - $\mathscr{L}(D_6)$ and $\mathscr{L}(D_6(a_k))$}

 Let $S_1$ (resp. $S_2$) be $D_6(a_1)$-associated (resp. $D_6(a_2)$-associated) root subset.
 As above, by \eqref{eq_regext4} and Corollary \ref{corol_Estim_2} we have
 \begin{equation*}
     \mid \mathscr{L}_{E_7}(D_6(a_2)) \mid =
     \mid \mathscr{L}_{E_7}(D_6(a_1)) \mid =
     \mid \mathscr{L}_{E_7}(D_6) \mid ~\leq~ \mid \varPhi(E_7) \mid -  \mid \varPhi(D_6) \mid = 126 - 60 = 66. \\
 \end{equation*}

 \begin{table}[H]
  %%\footnotesize
  \centering
  \renewcommand{\arraystretch}{1.3}
  \begin{tabular} {|c|c|c|c|}
  \hline
       & $\varphi \in \varPhi(D_8)\backslash\varPhi(D_7)$ &  $\delta \in \varPhi(D_8)\backslash\varPhi(D_7)$
       & $\delta^{\nabla} = -\varphi^{\nabla} \in \mathscr{L}_{D_8}(D_7)$ \\
    \hline
     1 & $\begin{array}{ccccccc}
            0 & 0 & 0 & 0 & 0 & 0 & 1 \\
              & 0
         \end{array}$ &
         $\begin{array}{ccccccc}
            1 & 2 & 2 & 2 & 2 & 2 & 1 \\
              & 1
         \end{array}$ &
          $\begin{array}{ccccccc}
             0 & 0 & 0 & 0 & 0 & -1 \\
              & 0
         \end{array}$ \\
    \hline
     2 & $\begin{array}{ccccccc}
            0 & 0 & 0 & 0 & 0 & 1 & 1 \\
              & 0
         \end{array}$ &
         $\begin{array}{ccccccc}
            1 & 2 & 2 & 2 & 2 & 1 & 1 \\
              & 1
         \end{array}$ &
          $\begin{array}{ccccccc}
             0 & 0 & 0 & 0 & -1 & 1 \\
              & 0
         \end{array}$ \\
    \hline
     3 & $\begin{array}{ccccccc}
            0 & 0 & 0 & 0 & 1 & 1 & 1 \\
              & 0
         \end{array}$ &
         $\begin{array}{ccccccc}
            1 & 2 & 2 & 2 & 1 & 1 & 1 \\
              & 1
         \end{array}$ &
          $\begin{array}{ccccccc}
             0 & 0 & 0 & -1 & 1 & 0 \\
              & 0
         \end{array}$ \\
    \hline
     4 & $\begin{array}{ccccccc}
            0 & 0 & 0 & 1 & 1 & 1 & 1 \\
              & 0
         \end{array}$ &
         $\begin{array}{ccccccc}
            1 & 2 & 2 & 1 & 1 & 1 & 1 \\
              & 1
         \end{array}$ &
          $\begin{array}{ccccccc}
             0 & 0 & -1 & 1 & 0 & 0 \\
              & 0
         \end{array}$ \\
    \hline
     5 & $\begin{array}{ccccccc}
            0 & 0 & 1 & 1 & 1 & 1 & 1 \\
              & 0
         \end{array}$ &
         $\begin{array}{ccccccc}
            1 & 2 & 1 & 1 & 1 & 1 & 1 \\
              & 1
         \end{array}$ &
          $\begin{array}{ccccccc}
             0 & -1 & 1 & 0 & 0 & 0 \\
              & 0
         \end{array}$ \\
    \hline
     6 & $\begin{array}{ccccccc}
            0 & 1 & 1 & 1 & 1 & 1 & 1 \\
              & 0
         \end{array}$ &
         $\begin{array}{ccccccc}
            1 & 1 & 1 & 1 & 1 & 1 & 1 \\
              & 1
         \end{array}$ &
          $\begin{array}{ccccccc}
             -1 & 1 & 0 & 0 & 0 & 0 \\
              & -1
         \end{array}$ \\
    \hline
     7 & $\begin{array}{ccccccc}
            1 & 1 & 1 & 1 & 1 & 1 & 1 \\
              & 0
         \end{array}$ &
         $\begin{array}{ccccccc}
            0 & 1 & 1 & 1 & 1 & 1 & 1 \\
              & 1
         \end{array}$ &
          $\begin{array}{ccccccc}
             1 & 0 & 0 & 0 & 0 & 0 \\
              & -1
         \end{array}$ \\
    \hline
       \end{tabular}
  \vspace{2mm}
  \caption{\hspace{3mm} $7$ pairs of positive roots $\varphi, \delta \in D_8$ such that $\delta^{\nabla} = -\varphi^{\nabla}$.
     There are also $7$ pairs of negative roots $\varphi, \delta \in D_8$ such that $\delta^{\nabla} = -\varphi^{\nabla}$}
  \label{tab_pairs_roots_D8}
  \end{table}
~\\
 The linkage system $\mathscr{L}(D_6)$, (resp. $\mathscr{L}(D_6(a_1))$ and $\mathscr{L}(D_6(a_2))$) consists of two parts:
 \begin{equation}
   \label{eq_plises_D6}
  \begin{split}
    & \boxed{ \mathscr{L}(D_6) = \mathscr{L}_{D_7}(D_6) ~\cup~  \mathscr{L}_{E_7}(D_6)}, \\
    & \boxed{ \mathscr{L}(D_6(a_1)) = \mathscr{L}_{D_7}(D_6(a_1)) ~\cup~  \mathscr{L}_{E_7}(D_6(a_1))}, \\
    & \boxed{ \mathscr{L}(D_6(a_2)) = \mathscr{L}_{D_7}(D_6(a_2)) ~\cup~  \mathscr{L}_{E_7}(D_6(a_2))}.
   \end{split}
 \end{equation}
~\\
 By Corollary \ref{cor_size_pls}(i), we have
  \begin{equation*}
    \mid \mathscr{L}_{D_7}(D_6) \mid ~=~ \mid \mathscr{L}_{D_7}(D_6(a_1)) \mid ~=~
                 \mid \mathscr{L}_{D_7}(D_6(a_2)) \mid ~=~ 2\times6 = 12. \\
   \end{equation*}
 By \eqref{eq_plises_D6},
  \begin{equation*}
   \mid \mathscr{L}(D_6(a_1)) \mid ~=~  \mid \mathscr{L}(D_6(a_2)) \mid ~=~ \mid \mathscr{L}(D_6) \mid ~\leq~ 12 + 66 = 78.
   \end{equation*}
  Let the coordinates of roots of $E_7$ be as follows
  \begin{equation*}
   \begin{array}{cccccc}
     \tau_6 & \tau_1 & \tau_2 & \tau_3 & \tau_4 & \tau_5 \\
            &   & \beta_2 \\
    \end{array}
   \end{equation*}
~\\
Consider the maximal (resp. minimal) root in $E_7$:
  \begin{equation*}
   \mu_{max} = \left (
   \begin{array}{cccccc}
     2 & 3 & 4 & 3 & 2 & 1 \\
       &   & 2 \\
   \end{array}
     \right ), \quad
   \mu_{min} = -\left (
   \begin{array}{cccccc}
     2 & 3 & 4 & 3 & 2 & 1 \\
       &   & 2 \\
   \end{array}
     \right ).
  \end{equation*}
  Roots $\pm\mu_{max}$ are orthogonal to any simple root except for  $\tau_6$.
  Further, $\pm\mu_{max} \in \varPhi(E_7) \backslash \varPhi(D_6)$ and $\pm\mu_{max}\not\in \varPhi(D_6)$,
  where $\varPhi(D_6)$ is spanned by $\{ \tau_1, \tau_2, \tau_3, \tau_4,  \tau_5, \beta_2 \}$.
  Vectors $\pm\mu_{max}^{\nabla}$ from $\mathscr{L}_{E_7}(D_6)$ are zero, i.e.,
 \begin{equation*}
   \mid \mathscr{L}(D_6(a_1)) \mid ~=~  \mid \mathscr{L}(D_6(a_2)) \mid ~=~ \mid \mathscr{L}(D_6) \mid ~\leq~  78 - 2 = 76.
 \end{equation*}

  There are $76$ linkage diagrams of $\mathscr{L}(D_6)$
  (resp.  $\mathscr{L}(D_6(a_1))$, resp.  $\mathscr{L}(D_6(a_2))$) presented in Fig. \ref{D6pure_loctets}
  (resp. Fig. \ref{D6a1_linkages},  resp. Fig. \ref{D6a2_linkages}).
  We conclude that there are \underline{exactly $76$ linkage diagrams} in each of linkage systems
  $\mathscr{L}(D_6)$, $\mathscr{L}(D_6(a_1))$, and $\mathscr{L}(D_6(a_2))$.

  \subsubsection{The $D$- and $E$-components $\mathscr{L}_{D_8}(D_7)$, $\mathscr{L}_{D_8}(D_7(a_k))$, $\mathscr{L}_{E_8}(D_7)$,
   and $\mathscr{L}_{E_8}(D_7(a_k))$}
   \index{linkage system  ! - $\mathscr{L}(D_7)$ and $\mathscr{L}(D_7(a_k))$}
 By \eqref{eq_regext4} and Corollary \ref{corol_Estim_2}
 \begin{equation}
  \begin{split}
   \label{eq_estim_D7_1}
      \mid \mathscr{L}_{E_8}(D_7(a_2)) \mid =  \mid \mathscr{L}_{E_8}(D_7(a_1)) \mid = &
      \mid \mathscr{L}_{E_8}(D_7) \mid ~\leq~ \\
      & \mid \varPhi(E_8) \mid - \mid \varPhi(D_7) \mid = 240 - 84 = 156. \\
  \end{split}
 \end{equation}

 \begin{table}[H]
  %%\footnotesize
  \centering
  \renewcommand{\arraystretch}{1.3}
  \begin{tabular} {|c|c|c|c|}
  \hline
       & $\eta \in \varPhi(E_8)\backslash\varPhi(D_7)$ &  $\lambda \in \varPhi(E_8)\backslash\varPhi(D_7)$
       & $\eta^{\nabla} = -\lambda^{\nabla}\in \mathscr{L}_{E_8}(D_7)$ \\
    \hline
     1 & $\begin{array}{ccccccc}
            2 & 3 & 4 & 3 & 2 & 1 & 0 \\
              &  &  2
         \end{array}$ &
         $\begin{array}{ccccccc}
            2 & 4 & 6 & 5 & 4 & 3 & 2 \\
              &  & 3
         \end{array}$ &
          $\begin{array}{ccccccc}
             0 & 0 & 0 & 0 & 0 & -1 \\
              & 0
         \end{array}$ \\
    \hline
     2 & $\begin{array}{ccccccc}
            2 & 3 & 4 & 3 & 2 & 1 & 1 \\
              & & 2
         \end{array}$ &
         $\begin{array}{ccccccc}
            2 & 4 & 6 & 5 & 4 & 3 & 1 \\
              &  & 3
         \end{array}$ &
          $\begin{array}{ccccccc}
             0 & 0 & 0 & 0 & -1 & 1 \\
              & 0
         \end{array}$ \\
    \hline
     3 & $\begin{array}{ccccccc}
            2 & 3 & 4 & 3 & 2 & 2 & 1 \\
              & & 2
         \end{array}$ &
         $\begin{array}{ccccccc}
            2 & 4 & 6 & 5 & 4 & 2 & 1 \\
              & & 3
         \end{array}$ &
          $\begin{array}{ccccccc}
             0 & 0 & 0 & -1 & 1 & 0 \\
              & 0
         \end{array}$ \\
    \hline
     4 & $\begin{array}{ccccccc}
            2 & 3 & 4 & 3 & 3 & 2 & 1 \\
              & & 2
         \end{array}$ &
         $\begin{array}{ccccccc}
            2 & 4 & 6 & 5 & 3 & 2 & 1 \\
              & & 3
         \end{array}$ &
          $\begin{array}{ccccccc}
             0 & 0 & -1 & 1 & 0 & 0 \\
              & 0
         \end{array}$ \\
    \hline
     5 & $\begin{array}{ccccccc}
            2 & 3 & 4 & 4 & 3 & 2 & 1 \\
              & & 2
         \end{array}$ &
         $\begin{array}{ccccccc}
            2 & 4 & 6 & 4 & 3 & 2 & 1 \\
              & & 3
         \end{array}$ &
          $\begin{array}{ccccccc}
             0 & -1 & 1 & 0 & 0 & 0 \\
              & 0
         \end{array}$ \\
    \hline
     6 & $\begin{array}{ccccccc}
            2 & 3 & 5 & 4 & 3 & 2 & 1 \\
              & & 2
         \end{array}$ &
         $\begin{array}{ccccccc}
            2 & 4 & 5 & 4 & 3 & 2 & 1 \\
              & & 3
         \end{array}$ &
          $\begin{array}{ccccccc}
             -1 & 1 & 0 & 0 & 0 & 0 \\
              & -1
         \end{array}$ \\
    \hline
     7 & $\begin{array}{ccccccc}
            2 & 3 & 5 & 4 & 3 & 2 & 1 \\
              & & 3
         \end{array}$ &
         $\begin{array}{ccccccc}
            2 & 4 & 5 & 4 & 3 & 2 & 1 \\
              & & 2
         \end{array}$ &
          $\begin{array}{ccccccc}
             -1 & 0 & 0 & 0 & 0 & 0 \\
                & 1
         \end{array}$ \\
    \hline
       \end{tabular}
  \vspace{2mm}
  \caption{\hspace{3mm} $7$ pairs of positive roots $\eta, \lambda \in \varPhi(E_8)$ such that $\eta^{\nabla} = -\lambda^{\nabla}$.
     There are also $7$ pairs of negative roots $\eta, \lambda \in \varPhi(E_8)$ such that $\eta^{\nabla} = -\lambda^{\nabla}$}
  \label{tab_pairs_roots_E8}
  \end{table}

\begin{figure}[H]
\centering
\includegraphics[scale=0.4]{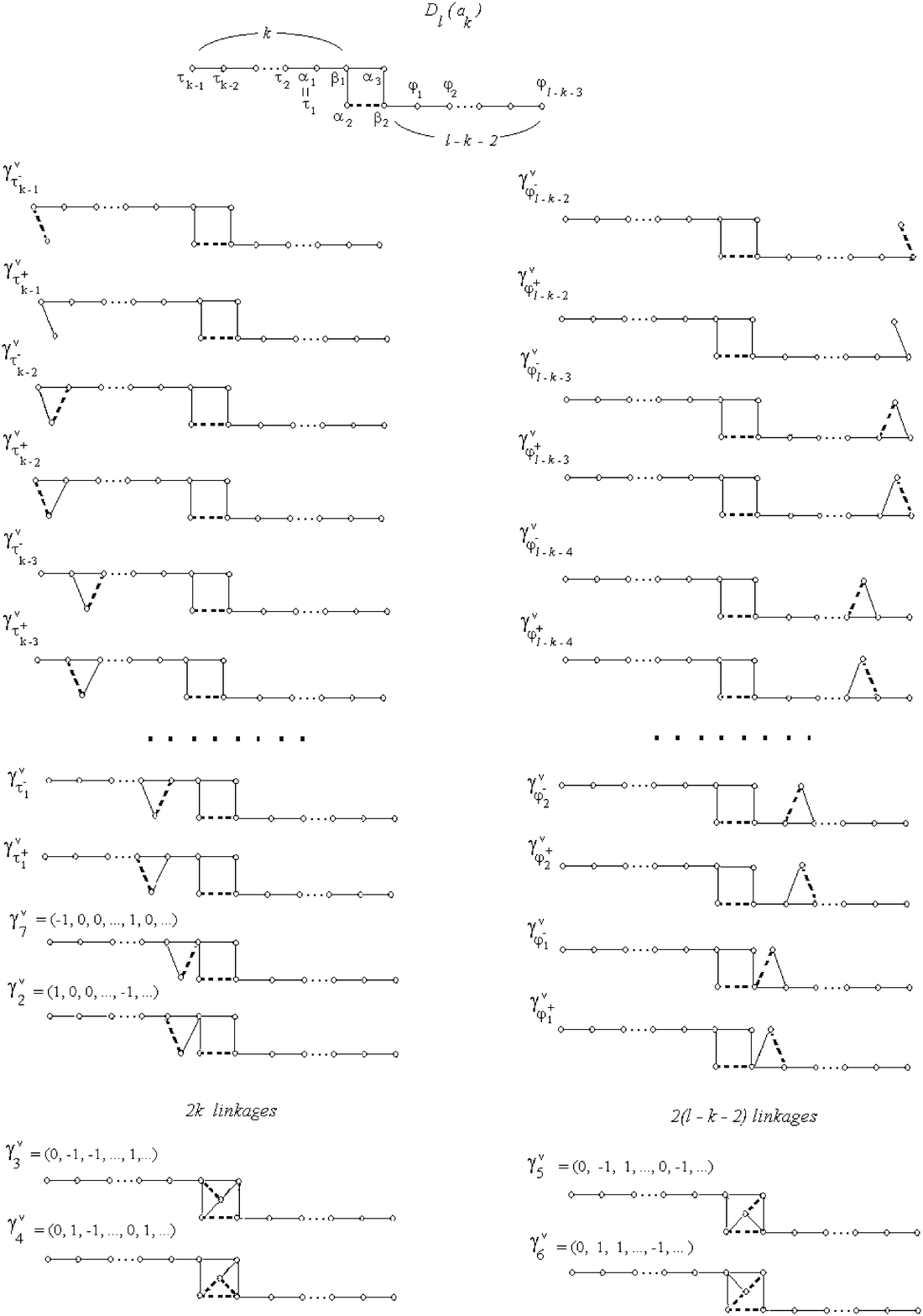}
\caption{\hspace{3mm}$D_l(a_k)$ for $l > 7$, $1$ loctet, $2l$ linkage diagrams}
%%%%%% The label must come after caption
\label{Dk_al_linkages}
\end{figure}

\begin{figure}[H]
\centering
\includegraphics[scale=0.5]{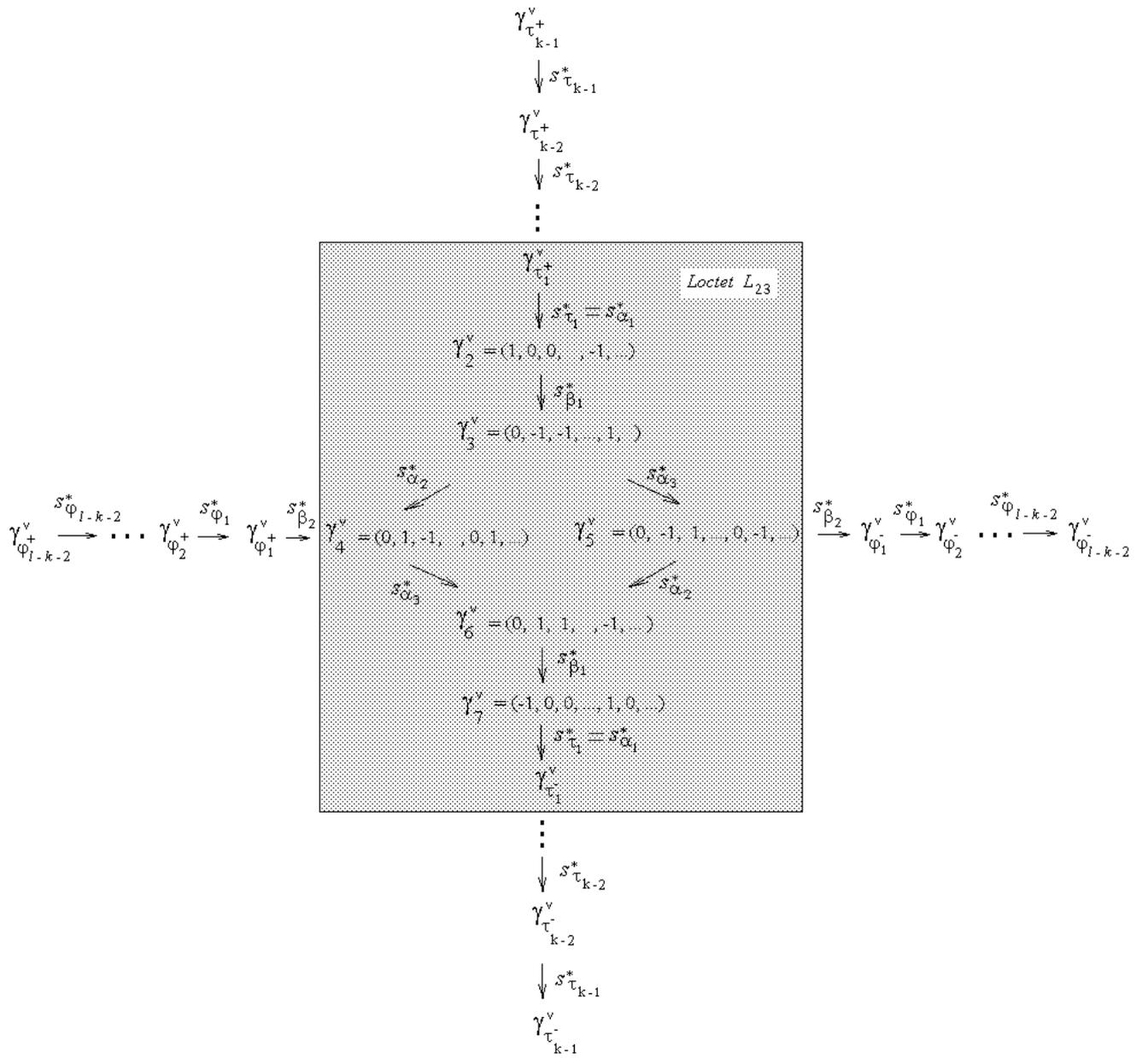}
\caption{\hspace{3mm}The linkage system $D_l(a_k)$ for $l > 7$ (wind rose of linkages).
 The single loctet $L_{23}$ is depicted in the shaded area}
%%%%%% The label must come after caption
\label{Dk_al_wind_rose}
\end{figure}

~\\
~\\

 \index{$L$, linear space ! - spanned by the $\Gamma$-associated root subset}
 \index{linkage label vectors ! - for $D_7$}
 \index{$L^{\nabla}$, space of linkage labels}
  Let  $L$ be the linear space spanned by the simple roots
  $\{\alpha_1, \alpha_2, \alpha_3, \alpha_4, \beta_1, \beta_2, \beta_3\}$ of  the root system $\varPhi(D_7)$,
  and $L^{\nabla}$ is the space of linkage labels, see \S\ref{sec_proj_linkage}.
  %% and $L'$ be the linear space spanned by $L$ and the simple root $\beta_4 \in \varPhi(E_8)$, linearly independent of
  %% the space $L$.
  We have
\begin{equation*}
  \begin{split}
   & \mathscr{L}_{D_8}(D_7)  =  \{ \gamma^{\nabla} | \gamma \in \varPhi(D_8)\backslash\varPhi(D_7) \} \subset L^{\nabla}, \\
   & \mathscr{L}_{E_8}(D_7)  = \{ \gamma^{\nabla} | \gamma \in \varPhi(E_8)\backslash\varPhi(D_7) \} \subset L^{\nabla}.
   \end{split}
 \end{equation*}
  For each of the seven pairs of $\{\eta, \lambda\}$ of positive roots in $\varPhi(E_8)$
  given in Table \ref{tab_pairs_roots_E8},  we have $\eta^{\nabla} = -\lambda^{\nabla}$.
  There are also seven pairs of negative roots such that every pair gives the same linkage label vector in  $\varPhi(E_8)\backslash\varPhi(D_7)$.
  Thus, we have to subtract $14$ in the estimation \eqref{eq_estim_D7_1} for $\mathscr{L}_{E_8}(D_7)$:
 \begin{equation}
    \label{eq_estim_D7_2}
   \mid \mathscr{L}_{E_8}(D_7(a_2)) \mid =  \mid \mathscr{L}_{E_8}(D_7(a_1)) \mid =
   \mid \mathscr{L}_{E_8}(D_7) \mid ~\leq~ 156 - 14 = 142.
 \end{equation}
  In addition, we have here a new phenomenon: components $\mathscr{L}_{E_8}(D_7)$ and
  $\mathscr{L}_{D_8}(D_7)$ overlap. Namely, roots of Tables \ref{tab_pairs_roots_D8} and
  \ref{tab_pairs_roots_E8} yield the same linkage label vectors in $L^{\nabla}$, i.e.,
  \begin{equation}
   \label{eq_incl_1}
     \mathscr{L}_{D_8}(D_7) \subset  \mathscr{L}_{E_8}(D_7).
  \end{equation}
 Hence,
 \begin{equation}
   \label{eq_estim_D7_3}
      \mathscr{L}(D_7) = \mathscr{L}_{D_8}(D_7) \cup \mathscr{L}_{E_8}(D_7) = \mathscr{L}_{E_8}(D_7), \text{ i.e., }
         \mid \mathscr{L}(D_7) \mid =  \mid \mathscr{L}_{E_8}(D_7) \mid.
   \end{equation}
 ~\\
 Thus, by \eqref{eq_estim_D7_2}, \eqref{eq_estim_D7_3}
  \begin{equation}
     \label{eq_estim_D7_5}
     \mid \mathscr{L}(D_7(a_1)) \mid ~=~ \mid \mathscr{L}(D_7(a_2)) ~=~ \mid \mathscr{L}(D_7) \mid ~\leq~ 142.
  \end{equation}
 In Figs. \ref{D7pu_loctets_comp1}, \ref{D7pu_loctets_comp2}, \ref{D7a1_D7a2_D7pu_loctets_comp3}$(a)$
 we have $2\times64 + 14 = 142$ linkage label vectors  in the linkage system $\mathscr{L}(D_7)$.
 Then we conclude that there are \underline{exactly $142$ elements} in each of linkage systems $\mathscr{L}(D_7)$},
 $\mathscr{L}(D_7(a_1))$ and $\mathscr{L}(D_7(a_2))$.

%% file: 8E_partialElak.tex
~\\
\section{\sc\bf $E$-type linkage systems}
  \label{sec_stratif_E}

\subsection{The linkage systems  $\mathscr{L}(E_l)$ and $\mathscr{L}(E_l(a_k))$ for $l = 6,7$}
    \label{sec_plis_El_Elak}
 \index{linkage system  ! - $\mathscr{L}(E_l)$ and $\mathscr{L}(E_l(a_k))$}

 In this section, we suppose that $\Gamma$ is one of the Carter diagrams
 \begin{equation}
  \label{eq_only_E6_E7_types}
   \begin{split}
     & E_6, E_6(a_k), \text{ where } k = 1,2, \\
     & E_7, E_7(a_k), \text{ where } k = 1,2,3,4,
   \end{split}
 \end{equation}
 see Table \ref{tab_part_root_syst}.
 Let $S$ be a $\Gamma$-associated root subset.

 \subsubsection{The linkage systems $\mathscr{L}(E_6)$, $\mathscr{L}(E_6(a_k))$ for $k=1,2$}
  \label{sec_E6_E6ak}
 Since the root system $\varPhi(E_6)$ is not contained in the root system $\varPhi(D_n)$, see Lemma \ref{lem_E6_not_in_Dn},
 then by means of Theorem \ref{prop_bijection_root_syst} and
 bijective maps of Table \ref{tab_part_root_syst} we deduce
 that $\varPhi(E_6(a_k))\footnotemark[1]$, where $k=1,2$, are also not contained in $\varPhi(D_n)$.
 \footnotetext[1]{By abuse of notation, we write $\varPhi(E_6(a_k))$ meaning the partial root system $\varPhi(S)$,
 where $S$ is the $E_6(a_k)$-associated subset;
 we write $\varPhi(E_7(a_k))$ meaning the partial root system $\varPhi(S)$, where $S$ is the $E_7(a_k)$-associated subset.}

 Since $E_6 \subset E_7, E_7(a_1), E_7(a_2)$ and $E_6(a_1) \subset E_7(a_3), E_7(a_4)$, we derive
 that $\varPhi(E_7)$ and $\varPhi(E_7(a_k))$ are not contained in $\varPhi(D_n)$ for $k=1,2,3,4$.
 Thus, for Carter diagrams \eqref{eq_only_E6_E7_types}, a Dynkin extension
 can be obtained \underline{only by means of $E_{l+1}$}:
 \begin{equation*}
   \begin{split}
     & E_6 <_D E_7, \quad E_6(a_k) <_D E_7, \text{ where } k = 1,2 \\
     %% & \\
     & E_7 <_D E_8, \quad E_7(a_k) <_D E_8,  \text{ where } k = 1,2,3,4.
   \end{split}
 \end{equation*}
~\\
 Then the linkage systems for $E_l$ (resp. $E_l(a_k)$) are
 \begin{equation}
   \label{eq_linksyst_E6E7}
   \begin{split}
     & \mathscr{L}(E_6) = \mathscr{L}_{E_7}(E_6), \quad \mathscr{L}(E_6(a_k)) = \mathscr{L}_{E_7}(E_6(a_k)),
        \quad k = 1,2, \\
     & \mathscr{L}(E_7) = \mathscr{L}_{E_8}(E_7), \quad \mathscr{L}(E_7(a_k)) = \mathscr{L}_{E_8}(E_7(a_k)),
        \quad k = 1,2,3,4.
   \end{split}
 \end{equation}
 By \eqref{eq_regext4} and Corollary \ref{corol_Estim_2}
 \begin{equation}
   \label{eq_linksyst_E6E7_1}
     \mid \mathscr{L}(E_6(a_1)) \mid ~=~ \mid \mathscr{L}(E_6(a_2)) \mid ~=~
     \mid \mathscr{L}(E_6) \mid ~\leq~ \mid \varPhi(E_7) \mid - \mid \varPhi(E_6) \mid = 126 - 72 = 54.
 \end{equation}

 In Figs. \ref{E6a1_linkages}, \ref{E6a2_linkages},
 \ref{27_weight_diagr_E6__2comp}(top)
 we have $54$ linkage label vectors in the linkage system $\mathscr{L}(E_6)$, (resp. $\mathscr{L}(E_6(a_1))$, resp. $\mathscr{L}(E_6(a_2))$).
 We conclude that there are \underline{exactly $54$ elements} in each of linkage systems
 $\mathscr{L}(E_6)$}, $\mathscr{L}(E_6(a_1))$, $\mathscr{L}(E_6(a_2))$.

 \subsubsection{The linkage systems $\mathscr{L}(E_7)$, $\mathscr{L}(E_7(a_k))$ for $k=1,2,3,4$}
  Throughout this subsection we assume that $k=1,2,3,4$.
 \begin{lemma}
   \label{lem_compl_roots_E8}
   \index{Dynkin extension ! - $E_7 <_D E_8$,  $E_7(a_k) <_D E_8$}
   \index{$\Gamma <_D \Gamma'$, Dynkin extension}
    \index{root stratum}

   Consider the Dynkin extension $E_7 <_D E_8$. Let $S = \{\tau_1,\dots,\tau_7\}$ be $E_7$-associated root subset.
   Let $\varphi$ be a positive root in the root stratum $\varPhi(E_8)\backslash\varPhi(S)$,
   and $\mu_{max}$ be the maximal root in $\varPhi(E_8)$:
  \begin{equation}
   \label{eq_max_in_E8}
    \varphi =
    \begin{array}{ccccccc}
      \tau_1 & \tau_3 & \tau_4 & \tau_5 & \tau_6 & \tau_7 & \tau_8 \\
             & & \tau_2 \\
    \end{array} \qquad
    \mu_{max} = \begin{array}{ccccccc}
      2 & 4 & 6 & 5 & 4 & 3 & 2 \\
             & & 3 \\
    \end{array}
  \end{equation}

     {\rm (i)} The vector
   \begin{equation}
     \label{eq_E8_ext1}
      \delta = \mu_{max} - \varphi,
   \end{equation}
   where $\varphi \neq \mu_{max}$, is also the root in $\varPhi(E_8)\backslash\varPhi(S)$.

 \index{$L$, linear space ! - spanned by the $\Gamma$-associated root subset}
 \index{$L^{\nabla}$, space of linkage labels}
    %% Let $L = [\tau_1,\dots,\tau_7]$.
    {\rm (ii)}
    The linkage label vectors $\varphi^{\nabla}$ and $-\delta^{\nabla}$ coincide:
   \begin{equation}
     \label{eq_Dl_ext2}
      \delta^{\nabla} = -\varphi^{\nabla}.
   \end{equation}
 \end{lemma}

 \PerfProof \rm (i)  As above, in Lemma \ref{lem_coinc_linkage}, we need to prove that
  \begin{equation}
     \label{eq_E8_ext3}
      %%\begin{split}
       \mathscr{B}(\delta) = 2, \text{ i.e., }
       \mathscr{B}(\mu_{max}) + \mathscr{B}(\varphi) - 2(\varphi, \mu_{max}) = 2. \\
      %%\end{split}
   \end{equation}
 Eq. \eqref{eq_E8_ext3} is equivalent to
  \begin{equation}
     \label{eq_E8_ext4}
       (\varphi, \mu_{max}) = 1.
   \end{equation}
~\\
 We have
  \begin{equation}
      \label{eq_E8_ext5}
     (\gamma, \mu_{max}) =
       \begin{cases}
           1 \text{ for } \gamma = \tau_8, \\
           0 \text{ for } \gamma = \tau_i, \text{ where } i < 8.
       \end{cases}
   \end{equation}

  The fact that $\varphi \in \varPhi(E_8)\backslash\varPhi(S)$, where $\varphi \neq \mu_{max}$
  and $\varphi$ is positive, means
  that $\tau_8$ enters with coefficient $1$ into the decomposition of $\varphi$ with relation to $\{\tau_1, \dots, \tau_8\}$,
  see \cite[Table VII]{Bo02}.
  Then, by \eqref{eq_E8_ext5} the relation \eqref{eq_E8_ext4} holds for
  any root $\varphi \in \varPhi(E_8)\backslash\varPhi(S)$.
  Therefore, \eqref{eq_E8_ext3} holds, i.e., $\delta$ is also a root.
  By \cite[Table VII]{Bo02} the coordinate $\tau_8$ of $\mu_{max}$ is $2$.
  By \eqref{eq_E8_ext1} $\tau_8$ also enters with coefficient $1$ into the decomposition of $\delta$
  with relation to $\{\tau_1, \dots, \tau_8\}$.
  Thus $\delta$ also belongs to $\varPhi(E_8)\backslash\varPhi(S)$.
~\\

  \rm (ii)
  For $\Gamma = E_7$,  by \eqref{eq_E8_ext5} for $\tau_i \neq  \tau_8$, we have
  \index{ $\varPhi(S)$, partial root system}
  \index{root system ! - partial}
  \index{partial root system}
  \index{linkage label vector}
     \begin{equation}
     \label{eq_E8_ext6}
        \delta^{\nabla}_i = (\delta, \tau_i) = (\mu_{max} - \varphi, \tau_i) =  -(\varphi, \tau_i) = -\varphi^{\nabla}_i,
         \text{ where } \delta \in \varPhi(E_8)\backslash\varPhi(E_7).
   \end{equation}
 Since $L = [\tau_1,\dots,\tau_7]$. we have $\delta^{\nabla} = -\varphi^{\nabla}$.
  \qed

  \begin{corollary}
    \label{cor_size_E7}
      The size of the linkage system $\mathscr{L}(E_7)$
      (resp. $\mathscr{L}(E_7(a_k)$) is equal to $56$.
  \end{corollary}

  \PerfProof
 By \eqref{eq_regext4} and Corollary \ref{corol_Estim_2} we have
 \begin{equation}
   \label{eq_linksyst_E7_1}
       \mid \mathscr{L}(E_7(a_k)) \mid ~=~ \mid \mathscr{L}(E_7) \mid ~\leq~ \mid \varPhi(E_8) \mid - \mid \varPhi(E_7) \mid = 240 - 126 = 114.
 \end{equation}
 Further, consider the maximal root $\mu_{max}$  and minimal root $\mu_{min} = -\mu_{max}$  in $\mathscr{L}(E_8)$:
  \begin{equation*}
 \pm{\mu}_{max} =
 \pm\begin{array}{ccccccc}
      2 & 4 & 6 & 5 & 4 & 3 & 2 \\
        &   & 3
 \end{array}.
  \end{equation*}
 We have $\pm{\mu}_{max} \in \varPhi(E_8) \backslash \varPhi(E_7)$ and $\pm{\mu}_{max} \not\in \varPhi(E_7)$.
 Roots $\pm{\mu}_{max}$ are orthogonal to all $\tau_i$ except for $\tau_8$, see \eqref{eq_max_in_E8}.
 Then vectors $\pm\mu_{max}^{\nabla}$ are zero. Thus, by \eqref{eq_linksyst_E7_1} we have
 \begin{equation}
   \label{eq_linksyst_E7_3}
      \mid \mathscr{L}(E_7(a_k)) \mid ~=~ \mid \mathscr{L}(E_7) \mid ~\leq~ 112.
 \end{equation}
 Finally, by Lemma \ref{lem_compl_roots_E8} roots $\varphi$ and $\varphi - \mu_{max}$
 give the same linkage label vectors. Therefore,
 \index{linkage label vectors ! - in $\mathscr{L}(E_7)$}
 \begin{equation}
   \label{eq_linksyst_E7_4}
      \mid \mathscr{L}(E_7(a_k)) \mid ~=~ \mid \mathscr{L}(E_7) \mid ~\leq~ 56.
 \end{equation}

 In Figs. \ref{E7a1_linkages}, \ref{E7a2_linkages}, \ref{E7a3_linkages}, \ref{E7a4_linkages},
 \ref{E7pure_linkage_system}
 we have $56$ linkage label vectors in the linkage system $\mathscr{L}(E_7)$,
 (resp. $\mathscr{L}(E_7(a_k))$ for $k = 1,2,3,4$).
 We conclude that there are \underline{exactly $56$ elements} in each of linkage systems $\mathscr{L}(E_7)$ and $\mathscr{L}(E_7(a_k))$
 for $k = 1,2,3,4$.
 \qed

%% file: 8A_partialAl.tex
~\\
\section{\sc\bf $A$-type linkage systems}
 \label{sec_stratif_A}

\subsection{The linkage systems  $\mathscr{L}(A_l)$ for $l = 5,6,7$}
  \label{sec_plis_Al}
  \index{linkage system  ! - $\mathscr{L}(A_l)$ for $l = 5,6,7$}

 In this section, we suppose that $\Gamma = A_l$ with $l = 5,6,7$.
 Let $S$ be a $\Gamma$-associated root subset. First, we have
 \begin{equation}
   \label{eq_first_Al}
   \begin{split}
    & \mathscr{L}(A_l) = \mathscr{L}_{A_{l+1}}(A_l) \cup \mathscr{L}_{D_{l+1}}(A_l)  \text{ for } l \neq 5,6,7, \\
    & \mathscr{L}(A_l) = \mathscr{L}_{A_{l+1}}(A_l) \cup \mathscr{L}_{D_{l+1}}(A_l) \cup \mathscr{L}_{E_{l+1}}(A_l)
           \text{ for } l = 5,6,7. \\
   \end{split}
 \end{equation}

 By \eqref{eq_regext4} the sizes of $A$- $D$- and $E$-components are estimated as follows:
 \begin{equation}
   \label{eq_first_Al_2}
   \Small
   \begin{split}
    & \mid \mathscr{L}_{A_{l+1}}(A_l) \mid ~\leq~ \mid \varPhi(A_{l+1}) \mid - \mid \varPhi(A_l) \mid =
              (l+1)(l+2) - l(l+1) = 2(l+1), \\
    & \mid \mathscr{L}_{D_{l+1}}(A_l) \mid ~\leq~ \mid \varPhi(D_{l+1}) \mid - \mid \varPhi(A_l) \mid =
              2l(l+1) - l(l+1) = l(l+1), \\
    & \mid \mathscr{L}_{E_6}(A_5) \mid ~\leq~ \mid \varPhi(E_6) \mid - \mid \varPhi(A_5) \mid = 72 - 30 = 42, \\
    & \mid \mathscr{L}_{E_7}(A_6) \mid ~\leq~ \mid \varPhi(E_7) \mid - \mid \varPhi(A_6) \mid = 126 - 42 = 84, \\
    & \mid \mathscr{L}_{E_8}(A_7) \mid ~\leq~ \mid \varPhi(E_8) \mid - \mid \varPhi(A_7) \mid = 240 - 56 = 184. \\
   \end{split}
 \end{equation}

\subsubsection{The linkage systems  $\mathscr{L}(A_3)$}
 By \eqref{eq_first_Al} and \eqref{eq_first_Al_2} we have

 \begin{equation}
   \label{eq_first_Al_3}
      \mathscr{L}(A_3) = \mathscr{L}_{A_4}(A_3) \cup \mathscr{L}_{D_4}(A_3), \text{ where }
      \mid \mathscr{L}_{A_4}(A_3) \mid \leq 8, \quad  \mid \mathscr{L}_{D_4}(A_3) \mid \leq 12.
 \end{equation}

 \begin{figure}[H]
 \centering
 \includegraphics[scale=0.56]{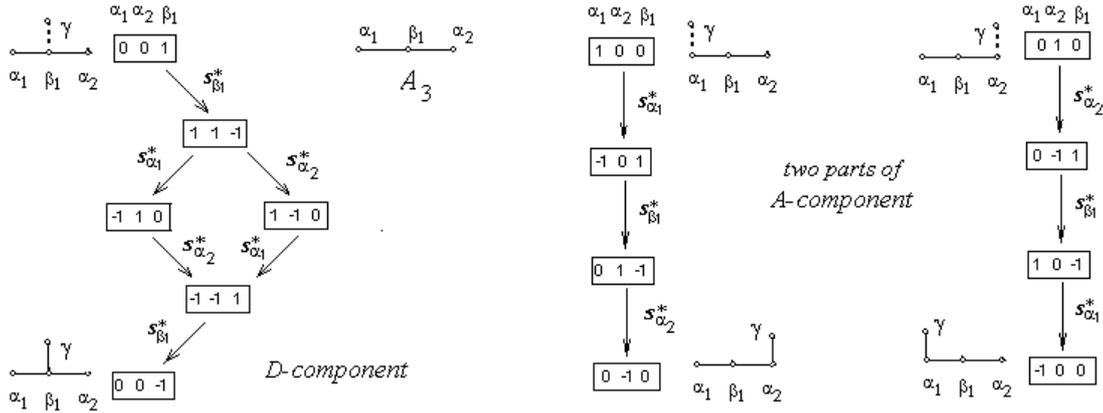}
 \caption{
 \hspace{3mm}The single part of the $D$-component and two parts of the $A$-component
  of the linkage system $\mathscr{L}(A_3)$}
 %%%%%% The label must come after caption
 \label{A3_to_A4D4}
 \end{figure}

 In Fig. \ref{A3_to_A4D4}, we have $8$ linkage label vectors in the $A$-component of the
 linkage system $\mathscr{L}(A_3)$. By \eqref{eq_first_Al_3} we conclude that there are
 \underline{exactly $8$ linkages in the $A$-component} $\mathscr{L}_{A_4}(A_3)$.
 For the component $\mathscr{L}_{D_4}(A_3)$, there is a different picture:
 By Table \ref{tab_pairs_D4overA4} the size of $\mathscr{L}_{A_4}(A_3)$ is twice less than the estimate in \eqref{eq_first_Al_3},
 i.e., $\mid \mathscr{L}_{D_4}(A_3) \mid \leq 6$. Thus, $\mid \mathscr{L}(A_3) \mid = 6 + 8 = 14$,
 see Fig. \ref{A3_to_A4D4}.

 \begin{table}[H]
  \footnotesize
  \centering
  \renewcommand{\arraystretch}{1.3}
  \begin{tabular} {|c|c|c|c|}
  \hline
       & $\varphi \in \varPhi(D_4)\backslash\varPhi(A_3)$ &  $\delta \in \varPhi(D_4)\backslash\varPhi(A_3)$
       & $\delta^{\nabla} = -\varphi^{\nabla} \in \mathscr{L}_{D_4}(A_3)$ \\
    \hline
     1 & $\begin{array}{ccccccc}
            0 & 1 & 1 \\
              & 1
         \end{array}$ &
         $\begin{array}{ccccccc}
            1 & 1 & 0  \\
              & 1
         \end{array}$ &
          $\begin{array}{ccccccc}
             \alpha_1 & \beta_1 & \alpha_2 \\
              -1 & 0 & 1  \\
               %% & & \\
         \end{array}$ \\
    \hline
     2 & $\begin{array}{ccccccc}
            1 & 1 & 1 \\
              & 1
         \end{array}$ &
         $\begin{array}{ccccccc}
            0 & 1 & 0  \\
              & 1
         \end{array}$ &
          $\begin{array}{ccccccc}
           \alpha_1 & \beta_1 & \alpha_2 \\
             1 & -1 & 1  \\
              %% & & \\
         \end{array}$ \\
    \hline
     3 & $\begin{array}{ccccccc}
            1 & 2 & 1 \\
              & 1
         \end{array}$ &
         $\begin{array}{ccccccc}
            0 & 0 & 0  \\
              & 1
         \end{array}$ &
          $\begin{array}{ccccccc}
            \alpha_1 & \beta_1 & \alpha_2 \\
             0 & 1 & 0  \\
              %%  & & \\
         \end{array}$  \\
    \hline
     4 & $\begin{array}{ccccccc}
            0 & -1 & -1 \\
              & -1
         \end{array}$ &
         $\begin{array}{ccccccc}
            -1 & -1 & 0  \\
              & -1
         \end{array}$ &
          $\begin{array}{ccccccc}
             \alpha_1 & \beta_1 & \alpha_2 \\
               1 & 0 & -1  \\
               %% & & \\
         \end{array}$ \\
    \hline
     5 & $\begin{array}{ccccccc}
            -1 & -1 & -1 \\
              & -1
         \end{array}$ &
         $\begin{array}{ccccccc}
            0 & -1 & 0  \\
              & -1
         \end{array}$ &
          $\begin{array}{ccccccc}
           \alpha_1 & \beta_1 & \alpha_2 \\
             -1 & 1 & -1  \\
              %% & & \\
         \end{array}$ \\
    \hline
     6 & $\begin{array}{ccccccc}
            -1 & -2 & -1 \\
               & -1
         \end{array}$ &
         $\begin{array}{ccccccc}
            0 & 0 & 0  \\
              & -1
         \end{array}$ &
          $\begin{array}{ccccccc}
            \alpha_1 & \beta_1 & \alpha_2 \\
             0 & -1 & 0  \\
              %%  & & \\
         \end{array}$  \\
    \hline
  \end{tabular}
  \vspace{2mm}
  \caption{\hspace{3mm} $6$ pairs of roots $\varphi, \delta \in \varPhi(D_4)\backslash\varPhi(A_3)$ such that
    $\delta^{\nabla} = -\varphi^{\nabla}$}
    \label{tab_pairs_D4overA4}
  \end{table}

\subsubsection{The linkage systems  $\mathscr{L}(A_4)$}
 By \eqref{eq_first_Al} and \eqref{eq_first_Al_2} we have

 \begin{equation}
   \label{eq_first_Al_4}
      \mathscr{L}(A_4) = \mathscr{L}_{A_5}(A_4) \cup \mathscr{L}_{D_5}(A_4), \text{ where }
      \mid \mathscr{L}_{A_5}(A_4) \mid \leq 10, \quad  \mid \mathscr{L}_{D_5}(A_4) \mid \leq 20.
 \end{equation}

  The linkage system $\mathscr{L}(A_4)$ contains $4$ parts. One of the two parts of the $D$-component
  (resp. $A$-component) of the linkage system $\mathscr{L}(A_4)$ corresponds to the positive roots
   of $\varPhi(D_5)\backslash\varPhi(A_4)$ (resp. $\varPhi(A_5)\backslash\varPhi(A_4)$),
  another part of the $D$-component (resp. $A$-component) of the linkage system $\mathscr{L}(A_4)$
  corresponds to the negative roots of
  $\varPhi(D_5)\backslash\varPhi(A_4)$ (resp. $\varPhi(A_5)\backslash\varPhi(A_4)$).
  In all, \underline{there are $10 + 20 = 30$ elements in $\mathscr{L}(A_4)$}.

 \begin{figure}[H]
 \centering
 \includegraphics[scale=0.55]{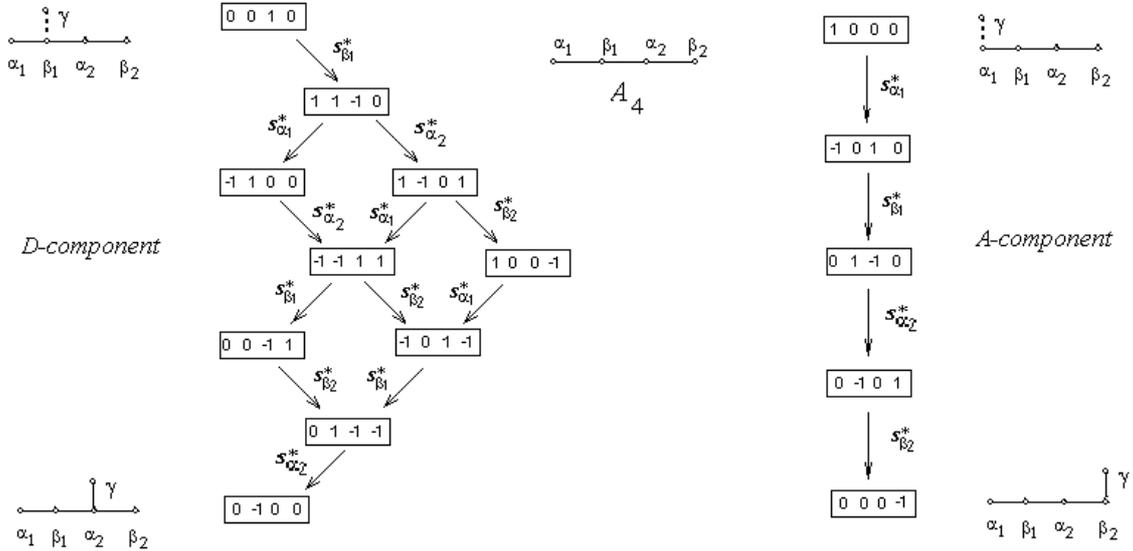}
 \caption{
 \hspace{3mm}One of the two parts of the $D$-component and one of two parts of the $A$-component
  of the linkage system $\mathscr{L}(A_4)$}
 %%%%%% The label must come after caption
 \label{A4_to_A5D5}
 \end{figure}

\subsubsection{The linkage systems  $\mathscr{L}(A_5)$}
 By \eqref{eq_first_Al} and \eqref{eq_first_Al_2} we have

 \begin{equation}
   \label{eq_first_A5_1}
    \begin{split}
      & \mathscr{L}(A_5) = \mathscr{L}_{A_6}(A_5) \cup \mathscr{L}_{D_6}(A_5) \cup \mathscr{L}_{E_6}(A_5), \text{ where } \\
      & \mid \mathscr{L}_{A_6}(A_5) \mid ~\leq~ 12, \quad  \mid \mathscr{L}_{D_6}(A_5) \mid ~\leq~ 30, \quad
      \mid \mathscr{L}_{E_6}(A_5) \mid ~\leq~ 72 - 30 = 42.
    \end{split}
 \end{equation}
 Let the coordinates of roots of $E_6$ be as follows:
  \begin{equation*}
   \begin{array}{cccccc}
     \tau_1 & \tau_3 & \tau_4 & \tau_5 & \tau_6 \\
            &   & \tau_2 \\
    \end{array}
   \end{equation*}
 Consider the maximal and minimal roots in $E_6$:
  \begin{equation*}
   \mu_{max} = \left (
   \begin{array}{cccccc}
     1 & 2 & 3 & 2 & 1 \\
       &   & 2 \\
   \end{array}
     \right ), \quad
   \mu_{min} = -\left (
   \begin{array}{cccccc}
     1 & 2 & 3 & 2 & 1 \\
       &   & 2 \\
   \end{array}
     \right ).
  \end{equation*}

  The roots $\mu_{max}$ (resp. $\mu_{min}$) is orthogonal to any simple root except for  $\tau_2$.
  The vectors $\pm\mu_{max}^{\nabla}$ are zero.
  \index{linkage label vectors ! - $\mu^{\nabla}_{max}$ and $\mu^{\nabla}_{min}$ for $E_6$}
  Then by \eqref{eq_first_A5_1}
 \begin{equation}
   \label{eq_first_A5_2}
   \mid \mathscr{L}_{E_6}(A_5) \mid ~\leq~  42 - 2 = 40.
  \end{equation}

  \begin{lemma}
   \label{lem_compl_roots_A5}
   Let $E_6$ be the Dynkin extension of the Carter diagram $\Gamma = A_5$,
   $\varphi$ be a certain positive root in $\varPhi(E_6)\backslash\varPhi(A_5)$.
   \index{Dynkin extension ! - $A_5 <_D E_6$}
   \index{$\Gamma <_D \Gamma'$, Dynkin extension}

     {\rm (i)} The vector
   \begin{equation}
     \label{eq_A5_ext1}
      \delta = \mu_{max} - \varphi,
   \end{equation}
   where $\varphi \neq \mu_{max}$, is also a root in $\varPhi(E_6)\backslash\varPhi(A_5)$.
 ~\\

 \index{$L$, linear space ! - spanned by the $\Gamma$-associated root subset}
  %% Let $L = [\tau_1,\dots,\tau_5]$.
    {\rm (ii)}
    The linkage label vectors $\varphi^{\nabla}$ and $-\delta^{\nabla}$ coincide:
   \begin{equation}
     \label{eq_A5_ext2}
      \delta^{\nabla} = -\varphi^{\nabla}.
   \end{equation}
 \end{lemma}

 \PerfProof  \rm (i) As above, in Lemmas \ref{lem_coinc_linkage} and \ref{lem_compl_roots_E8}, we need to prove that
  \begin{equation}
     \label{eq_A5_ext3}
        \mathscr{B}(\mu_{max}) + \mathscr{B}(\varphi) - 2(\varphi, \mu_{max}) = 2. \\
   \end{equation}
In other words, we need to prove that
  \begin{equation}
     \label{eq_A5_ext4}
       (\varphi, \mu_{max}) = 1.
   \end{equation}
~\\
 We have
  \begin{equation}
      \label{eq_A5_ext5}
     (\gamma, \mu_{max}) =
       \begin{cases}
           1 \text{ for } \gamma = \tau_2, \\
           0 \text{ for } \gamma = \tau_i, \text{ where } i \neq 2.
       \end{cases}
   \end{equation}
  The fact that $\varphi \in \varPhi(E_6)\backslash\varPhi(A_5)$, $\varphi \neq \mu_{max}$ and $\varphi$ is positive means
  that $\tau_2$ enters with coefficient $1$ into the decomposition of $\varphi$ with relation to $\{\tau_1, \dots, \tau_6\}$.
  By \eqref{eq_A5_ext5} the relation \eqref{eq_A5_ext4} holds for
  any root $\varphi \in \varPhi(E_6)\backslash\varPhi(A_5)$.
  Therefore, \eqref{eq_A5_ext3} holds, i.e., $\delta$ is also a root.
  By \cite[Table V]{Bo02} the coordinate $\tau_2$ of $\mu_{max}$ is $2$.
  By \eqref{eq_A5_ext1} $\tau_2$ also enters with coefficient $1$ into the decomposition of $\delta$
  with relation to $\{\tau_1, \dots, \tau_6\}$.
  Thus, $\delta$ also belongs to $\varPhi(E_6)\backslash\varPhi(A_5)$.

 \rm (ii) By \eqref{eq_A5_ext5} for $\varphi \neq \tau_2$, we have
    \begin{equation}
     \label{eq_A5_ext6}
        \delta^{\nabla}_i = (\delta, \tau_i) =  -(\varphi, \tau_i) = -\varphi^{\nabla}_i.
   \end{equation}
 Since $L = [\tau_1,\tau_3,\tau_4,\tau_5,\tau_6]$, we have $\delta^{\nabla} = -\varphi^{\nabla}.$
 \qed
~\\

\begin{figure}[H]
 \centering
 \includegraphics[scale=0.55]{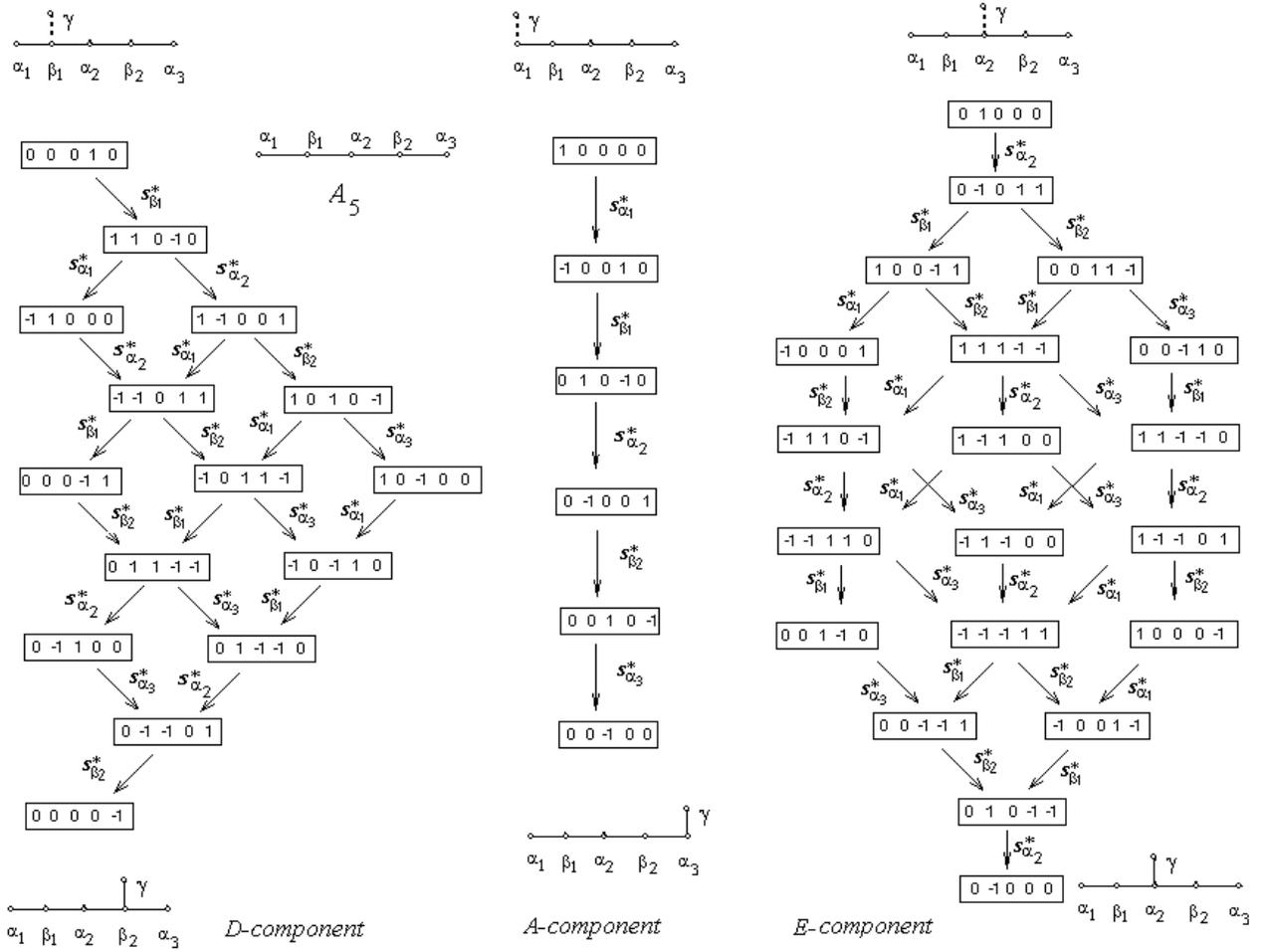}
 \caption{
 \hspace{3mm}One of the two parts of the $D$-component, one of the two parts of the $A$-component, and
 the $E$-component of the linkage system $\mathscr{L}(A_5)$}
 %%%%%% The label must come after caption
 \label{A5_to_A6D6E6}
 \end{figure}

 By \eqref{eq_first_A5_1}, \eqref{eq_first_A5_2}
 the two parts of the $A$-component contain $2 \times 6 = 12$ elements; the two parts of the $D$-component
 contain $2 \times 15 = 30$ elements, see Fig. \ref{A5_to_A6D6E6}. By Lemma \ref{lem_compl_roots_A5} and
 eq. \eqref{eq_first_A5_2} we see that
 the $E$-component contains $\leq 40/2 = 20$ elements. In Fig. \ref{A5_to_A6D6E6}, we see that
 the $E$-component contains  exactly $20$ elements. In all, the linkage system
  \underline{$\mathscr{L}(A_5)$ contains $12 + 30 + 20 = 62$ elements.}

\subsubsection{The linkage systems  $\mathscr{L}(A_6)$}
 By \eqref{eq_first_Al} and \eqref{eq_first_Al_2} we have
  \index{linkage diagrams ! - in $\mathscr{L}(A_6)$}

 \begin{equation}
   \label{eq_first_A6_1}
    %%\Small
    \begin{split}
      & \mathscr{L}(A_6) = \mathscr{L}_{A_7}(A_6) \cup \mathscr{L}_{D_7}(A_6) \cup \mathscr{L}_{E_7}(A_6), \text{ where } \\
      & \mid \mathscr{L}_{A_7}(A_6) \mid ~\leq~ 14, \quad  \mid \mathscr{L}_{D_7}(A_6) \mid ~\leq~ 42, \quad
      \mid \mathscr{L}_{E_7}(A_6) \mid ~\leq~ 126 - 42 = 84.
    \end{split}
 \end{equation}
~\\
 However, there are $7$ pairs of positive (resp. $7$ pairs of negative) roots  $\{\eta, \lambda\}$,
 where $\eta \in \varPhi(E_7)\backslash\varPhi(A_6)$ and $\lambda \in \varPhi(A_7)\backslash\varPhi(A_6)$
 such that the corresponding linkage label vectors coincide up-to-sign:
 $$
   \eta^{\nabla} = -\lambda^{\nabla},
 $$
 see Table \ref{tab_pairs_roots_A6ext}.
 In other words, we have to subtract $14$ roots from the estimate of $\mid \mathscr{L}_{E_7}(A_6) \mid$:
 \begin{equation}
   \label{eq_first_A6_2}
   %% \Small
      \mid \mathscr{L}_{A_7}(A_6) \mid ~\leq~ 14, \quad  \mid \mathscr{L}_{D_7}(A_6) \mid ~\leq~ 42, \quad
      \mid \mathscr{L}_{E_7}(A_6) \mid ~\leq~ 70.
 \end{equation}
There are $2 \times 7 + 2 \times 21 + 2 \times 35 = 14 + 42 + 70 = 126$ linkage diagrams
in $\mathscr{L}(A_6)$ presented in Fig. \ref{A6_to_A7D7E7}. We conclude that there are
\underline{exactly $126$ linkage diagrams in $\mathscr{L}(A_6)$}.

\subsubsection{The linkage systems  $\mathscr{L}(A_7)$}
  \index{root system}
 By \eqref{eq_first_Al} and \eqref{eq_first_Al_2} we have
 \begin{equation}
   \label{eq_first_A7_1}
    \begin{split}
      & \mathscr{L}(A_7) = \mathscr{L}_{A_8}(A_7) \cup \mathscr{L}_{D_8}(A_7) \cup \mathscr{L}_{E_8}(A_7), \text{ where } \\
      & \mid \mathscr{L}_{A_8}(A_7) \mid ~\leq~ 16, \quad  \mid \mathscr{L}_{D_8}(A_7) \mid ~\leq~ 56, \quad
      \mid \mathscr{L}_{E_8}(A_7) \mid ~\leq~ 240 - 56 = 184.
    \end{split}
 \end{equation}

  \begin{table}[H]
  %%\footnotesize
  \centering
  \renewcommand{\arraystretch}{1.3}
  \begin{tabular} {|c|c|c|c|}
  \hline
       & $\eta \in \varPhi(E_7)\backslash\varPhi(A_6)$ &  $\lambda \in \varPhi(A_7)\backslash\varPhi(A_6)$
       & $\eta^{\nabla} = -\lambda^{\nabla}  \in \mathscr{L}(A_6)$ \\
    \hline
     1 & $\begin{array}{ccccccc}
          \alpha_1 & \beta_1 & \alpha_2 & \beta_2 & \alpha_3 & \beta_3 \\
            1 & 2 & 3 & 2 & 1 & 0 \\
              &  &  2
         \end{array}$ &
         $\begin{array}{ccccccc}
          & \alpha_1 & \beta_1 & \alpha_2 & \beta_2 & \alpha_3 & \beta_3 \\
            1 & 1 & 1 & 1 & 1 & 1 & 1 \\
         \end{array}$ &
          $\begin{array}{ccccccc}
          \alpha_1 & \beta_1 & \alpha_2 & \beta_2 & \alpha_3 & \beta_3 \\
             0 & 0 & 0 & 0 & 0 & -1 \\
         \end{array}$ \\
    \hline
     2 & $\begin{array}{ccccccc}
            1 & \hspace{2mm}2 & \hspace{2mm}3 & \hspace{2mm}2 & \hspace{2mm}1 & \hspace{2mm}1 \\
              &  &  \hspace{2mm}2
         \end{array}$ &
         $\begin{array}{ccccccc}
            1 & \hspace{2mm}1 & \hspace{2mm}1 & \hspace{2mm}1 & \hspace{2mm}1 & \hspace{2mm}1 & \hspace{2mm}0 \\
         \end{array}$ &
          $\begin{array}{ccccccc}
             \hspace{2mm}0 & \hspace{2mm}0 & \hspace{2mm}0 & \hspace{2mm}0 & \hspace{2mm}-1 & \hspace{2mm}1 \\
         \end{array}$ \\
    \hline
     3 & $\begin{array}{ccccccc}
            1 & \hspace{2mm}2 & \hspace{2mm}3 & \hspace{2mm}2 & \hspace{2mm}2 & \hspace{2mm}1 \\
              &  &  \hspace{2mm}2
         \end{array}$ &
         $\begin{array}{ccccccc}
            1 & \hspace{2mm}1 & \hspace{2mm}1 & \hspace{2mm}1 & \hspace{2mm}1 & \hspace{2mm}0 & \hspace{2mm}0 \\
            %%1 & 1 & 1 & 1 & 1 & 0 & 0 \\
         \end{array}$ &
          $\begin{array}{ccccccc}
             \hspace{2mm}0 & \hspace{2mm}0 & \hspace{2mm}0 & \hspace{2mm}-1 & \hspace{2mm}1 & \hspace{2mm}0 \\
         \end{array}$ \\
    \hline
     4 & $\begin{array}{ccccccc}
            1 & \hspace{2mm}2 & \hspace{2mm}3 & \hspace{2mm}3 & \hspace{2mm}2 & \hspace{2mm}1 \\
              &  &  \hspace{2mm}2
         \end{array}$ &
         $\begin{array}{ccccccc}
            1 & \hspace{2mm}1 & \hspace{2mm}1 & \hspace{2mm}1 & \hspace{2mm}0 & \hspace{2mm}0 & \hspace{2mm}0 \\
            %%1 & 1 & 1 & 1 & 0 & 0 & 0 \\
         \end{array}$ &
          $\begin{array}{ccccccc}
             \hspace{2mm}0 & \hspace{2mm}0 & \hspace{2mm}-1 & \hspace{2mm}1 & \hspace{2mm}0 & \hspace{2mm}0 \\
         \end{array}$ \\
    \hline
     5 & $\begin{array}{ccccccc}
            1 & \hspace{2mm}2 & \hspace{2mm}4 & \hspace{2mm}3 & \hspace{2mm}2 & \hspace{2mm}1 \\
              &  &  \hspace{2mm}2
         \end{array}$ &
         $\begin{array}{ccccccc}
            1 & \hspace{2mm}1 & \hspace{2mm}1 & \hspace{2mm}0 & \hspace{2mm}0 & \hspace{2mm}0 & \hspace{2mm}0 \\
           %% 1 & 1 & 1 & 0 & 0 & 0 & 0 \\
         \end{array}$ &
          $\begin{array}{ccccccc}
             \hspace{2mm}0 & \hspace{2mm}-1 & \hspace{2mm}1 & \hspace{2mm}0 & \hspace{2mm}0 & \hspace{2mm}0 \\
         \end{array}$ \\
    \hline
     6 & $\begin{array}{ccccccc}
            1 & \hspace{2mm}3 & \hspace{2mm}4 & \hspace{2mm}3 & \hspace{2mm}2 & \hspace{2mm}1 \\
              &  &  \hspace{2mm}2
         \end{array}$ &
         $\begin{array}{ccccccc}
            1 & \hspace{2mm}1 & \hspace{2mm}0 & \hspace{2mm}0 & \hspace{2mm}0 & \hspace{2mm}0 & \hspace{2mm}0 \\
           %% 1 & 1 & 0 & 0 & 0 & 0 & 0 \\
         \end{array}$ &
          $\begin{array}{ccccccc}
             \hspace{2mm}-1 & \hspace{2mm}1 & \hspace{2mm}0 & \hspace{2mm}0 & \hspace{2mm}0 & \hspace{2mm}0 \\
         \end{array}$ \\
    \hline
     7 & $\begin{array}{ccccccc}
            2 & \hspace{2mm}3 & \hspace{2mm}4 & \hspace{2mm}3 & \hspace{2mm}2 & \hspace{2mm}1 \\
              &  &  \hspace{2mm}2
         \end{array}$ &
         $\begin{array}{ccccccc}
            1 & \hspace{2mm}0 & \hspace{2mm}0 & \hspace{2mm}0 & \hspace{2mm}0 & \hspace{2mm}0 & \hspace{2mm}0 \\
            %%1 & 0 & 0 & 0 & 0 & 0 & 0 \\
         \end{array}$ &
          $\begin{array}{ccccccc}
             \hspace{2mm}-1 & \hspace{2mm}0 & \hspace{2mm}0 & \hspace{2mm}0 & \hspace{2mm}0 & \hspace{2mm}0 \\
         \end{array}$ \\
    \hline
       \end{tabular}
  \vspace{2mm}
  \caption{\hspace{3mm} $7$ pairs of positive roots $\{\eta, \lambda\}$ such that $\eta^{\nabla} = -\lambda^{\nabla}$.
     There are also $7$ pairs of negative roots $\{\eta, \lambda\}$ such that $\eta^{\nabla} = -\lambda^{\nabla}$}
  \label{tab_pairs_roots_A6ext}
  \end{table}
~\\
 However, there are $8$ pairs of positive (resp. $8$ pairs of negative) roots  $\{\eta, \lambda\}$,
 where $\eta \in \varPhi(E_8)\backslash\varPhi(A_7)$ and $\lambda \in \varPhi(A_8)\backslash\varPhi(A_7)$
 such that the corresponding  linkage label vectors coincide up-to-sign:
 $$
   \eta^{\nabla} = -\lambda^{\nabla}.
 $$
 The pairs $\{\eta_i, \lambda_i\}$, where $1 \leq i \leq 8$ are depicted in Table \ref{tab_pairs_roots_A7extE8}.
 In other words, we have to subtract $16$ linkage label vectors from the estimate of $\mid \mathscr{L}_{E_8}(A_7) \mid$:
 \begin{equation}
   \label{eq_first_A7_2}
      \mid \mathscr{L}_{E_8}(A_7) \mid ~\leq~ 184 - 16 = 168.
 \end{equation}

\begin{remark}
{\rm
For any root $\gamma$ from the root system $\varPhi(A_l)$, we consider the sequence of coordinates
of the linkage label vector $\gamma^{\nabla}$
ordered as they lie in the diagram $A_l$, see $A$-components in Figs. \ref{A6_to_A7D7E7}, \ref{A7_to_A8D8_1}, i.e.,
 \index{linear order for the linkage diagram}
 \index{bicolored order}
\begin{equation}
  \label{lin_order}
    \gamma^{\nabla} = \{ a_1, b_1, a_2, b_2, \dots, a_i, b_i, \dots \}, \text{ where } a_i = (\alpha_i, \gamma), b_i = (\beta_i, \gamma).
\end{equation}
Then we say that the linkage diagram $\gamma^{\nabla}$ is given in a {\it linear order}. Note
that usually we consider coordinates of the linkage label vector $\gamma^{\nabla}$ in a {\it bicolored order} given as follows:
\begin{equation*}
    \gamma^{\nabla} = \{ a_1, a_2, \dots, b_1, b_2, \dots \}.
\end{equation*}
In Table \ref{tab_pairs_roots_A7extE8}, the linkage label vectors $\lambda^{\nabla}_i$ are given in
a linear order. The same vectors in the $A$-component in Fig. \ref{A7_to_A8D8_1} are given in the bicolored order.
\qed
}
\end{remark}

\begin{figure}[H]
\centering
 \vspace{1cm}
\includegraphics[scale=0.6]{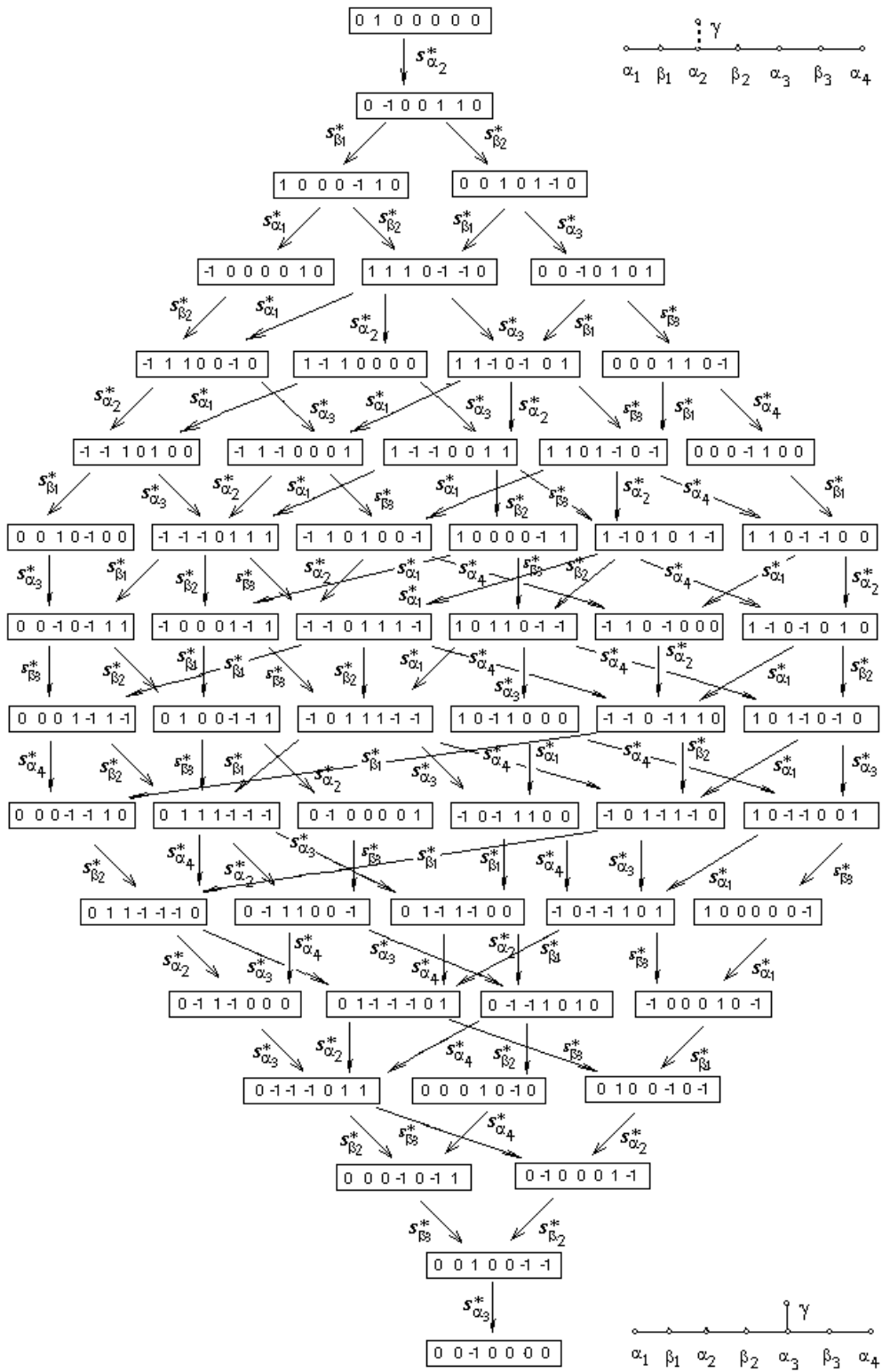}
\caption{\hspace{3mm} The linkage system $A_7$ (cont). One of the two $56$-element parts
 of the $E$-component}
%%%%%% The label must come after caption
\label{A7_to_E8}
\end{figure}

 \index{Weyl group ! - $W$}
 \index{$W$, Weyl group}
 \index{space of linkage labels $L^{\nabla}$ for $\varPhi(A_7)$}
 In addition, there are $28$ pairs of positive (resp. $28$ pairs of negative) roots  $\{\eta, \lambda\}$,
 where $\eta \in \varPhi(E_8)\backslash\varPhi(A_7)$ and $\lambda \in \varPhi(D_8)\backslash\varPhi(A_7)$
 such that the corresponding linkage label vectors coincide: $\eta^{\nabla} = \lambda^{\nabla}$.
 We do not present all of the $28$ pairs $\{\eta_i, \lambda_i\}$:
 It suffices to produce one pair, the remaining $27$ pairs are obtained by action of the Weyl group. Namely,
 let us take roots $\eta_1 \in \varPhi(E_8)\backslash\varPhi(A_7)$ and $\lambda_1 \in \varPhi(D_8)\backslash\varPhi(A_7)$
 such that $\eta_1^{\nabla} = \lambda_1^{\nabla}$;
 let $w \in W$, where $W$ is the Weyl group associated with $\varPhi(A_7)$.
 Though $w\eta_1$ and $w\lambda_1$ lie in different root systems,
 $w^*\eta^{\nabla}_1$ and $w^*\lambda^{\nabla}_1$ coincide:
 \begin{equation*}
    (w\eta_1)^{\nabla} = w^*\eta^{\nabla}_1 = w^*\lambda^{\nabla}_1 = (w\lambda_1)^{\nabla},
 \end{equation*}
 see Proposition \ref{prop_link_diagr_conn}(ii).
 Consider the following pair $\{\eta_1, \lambda_1\}$:
 \begin{equation*}
  \Small
  %%\begin{split}
    \eta_1 =
   \begin{array}{ccccccc}
      2 & 4 & 5 & 4 & 3 & 2 & 1 \\
      & & 2
   \end{array}  \in \varPhi(E_8)\backslash\varPhi(A_7), \qquad
   %% \\ \\
    \lambda_1 =
   \begin{array}{ccccccc}
      0 & 0 & 0 & 0 & 0 & 0 & 0 \\
      &  -1
   \end{array} \in \varPhi(D_8)\backslash\varPhi(A_7).
 %% \end{split}
 \end{equation*}
 Then
  \begin{equation*}
  \Small
    \eta^{\nabla}_1 = \lambda^{\nabla}_1 =
  \begin{cases}
   \begin{array}{ccccccc}
      0 & 1 & 0 & 0 & 0 & 0 & 0 \\
      \alpha_1 & \beta_1 & \alpha_2 & \beta_2 & \alpha_3 & \beta_3 & \alpha_4
   \end{array} \text{ in linear order,}
   \\ \\
   \begin{array}{ccccccc}
      0 & 0 & 0 & 0 & 1 & 0 & 0 \\
      \alpha_1 & \alpha_2 & \alpha_3 & \alpha_4 & \beta_1 & \beta_2 & \beta_3
   \end{array} \text{ in bicolored order.}
  \end{cases}
 \end{equation*}
The linkage label vector $\lambda^{\nabla}_1$ is the top element
in one of two $28$-element parts of the $D$-component of $\mathscr{L}(A_7)$, see Fig. \ref{A7_to_A8D8_1}.
Therefore, we have to subtract $2 \times 28 = 56$ linkage label vectors
from the estimate of $\mid \mathscr{L}_{E_8}(A_7) \mid$ in \eqref{eq_first_A7_2}:
 \begin{equation}
   \label{eq_first_A7_3}
      \mid \mathscr{L}_{E_8}(A_7) \mid ~\leq~ 168 - 56 = 112.
 \end{equation}
By \eqref{eq_first_A7_1}, \eqref{eq_first_A7_3}, we have
\begin{equation}
   \label{eq_first_A7_4}
      \mid \mathscr{L}_{A_8}(A_7) \mid ~\leq~ 16, \quad  \mid \mathscr{L}_{D_8}(A_7) \mid ~\leq~ 56, \quad
      \mid \mathscr{L}_{E_8}(A_7) \mid ~\leq~ 112.
 \end{equation}
In accordance to Figs. \ref{A7_to_A8D8_1} and \ref{A7_to_E8} we have
\underline{exactly $16 + 56 + 112 = 184$ linkages in $\mathscr{L}(A_7)$}.

%%\clearpage
%%~\\
%%~\\
  \begin{table}[H]
  \footnotesize
  \centering
  \renewcommand{\arraystretch}{1.5}
  \begin{tabular}{|c|c|c|c|}
  \hline
       & $\eta \in \varPhi(E_8)\backslash\varPhi(A_7)$ &
         $\lambda \in \varPhi(A_8)\backslash\varPhi(A_7)$ &
         $\lambda^{\nabla} = \eta^{\nabla}$ \\
   \hline
     & $\begin{array}{ccccccc}
          \alpha_1  \hspace{2mm}\beta_1  \hspace{2mm}\alpha_2  \hspace{2mm}\beta_2  \hspace{2mm}\alpha_3
                    \hspace{2mm}\beta_3   \hspace{2mm}\alpha_4 \\
      \end{array}$ &
     $\begin{array}{cccccccc}
          & \hspace{2mm}\alpha_1  \hspace{2mm}\beta_1  \hspace{2mm}\alpha_2  \hspace{2mm}\beta_2  \hspace{2mm}\alpha_3
                    \hspace{2mm}\beta_3   \hspace{2mm}\alpha_4 \\
      \end{array}$ &
     $\begin{array}{ccccccc}
           \alpha_1  \hspace{2mm}\beta_1  \hspace{2mm}\alpha_2  \hspace{2mm}\beta_2  \hspace{2mm}\alpha_3
                    \hspace{2mm}\beta_3   \hspace{2mm}\alpha_4 \\
        \end{array}$  \\
    \hline
     1 & $\begin{array}{ccccccc}
            1 & 3 & 5 & 4 & 3 & 2 & 1 \\
              &  &  3
         \end{array}$ &
         $\begin{array}{cccccccc}
            1 & 0 & 0 & 0 & 0 & 0 & 0 & 0\\
         \end{array}$ &
         $\begin{array}{ccccccc}
            -1 & 0 & 0 & 0 & 0 & 0 & 0 \\
         \end{array}$  \\
    \hline
     2 & $\begin{array}{ccccccc}
            2 & 3 & 5 & 4 & 3 & 2 & 1 \\
              &  &  3
         \end{array}$ &
         $\begin{array}{cccccccc}
            1 & 1 & 0 & 0 & 0 & 0 & 0 & 0\\
         \end{array}$ &
         $\begin{array}{ccccccc}
             1 & -1 & 0 & 0 & 0 & 0 & 0 \\
         \end{array}$ \\
    \hline
     3 & $\begin{array}{ccccccc}
            2 & 4 & 5 & 4 & 3 & 2 & 1 \\
              &  &  3
         \end{array}$ &
         $\begin{array}{cccccccc}
            1 & 1 & 1 & 0 & 0 & 0 & 0 & 0\\
         \end{array}$ &
         $\begin{array}{ccccccc}
            0 & 1 & -1 & 0 & 0 & 0 & 0 \\
         \end{array}$ \\
    \hline
     4 & $\begin{array}{ccccccc}
            2 & 4 & 6 & 4 & 3 & 2 & 1 \\
              &  &  3
         \end{array}$ &
         $\begin{array}{cccccccc}
            1 & 1 & 1 & 1 & 0 & 0 & 0 & 0\\
         \end{array}$ &
         $\begin{array}{ccccccc}
            0 & 0 & 1 & -1 & 0 & 0 & 0 \\
         \end{array}$ \\
    \hline
     5 & $\begin{array}{ccccccc}
            2 & 4 & 6 & 5 & 3 & 2 & 1 \\
              &  &  3
         \end{array}$ &
         $\begin{array}{cccccccc}
            1 & 1 & 1 & 1 & 1 & 0 & 0 & 0\\
         \end{array}$ &
         $\begin{array}{ccccccc}
            0 & 0 & 0 & 1 & -1 & 0 & 0 \\
         \end{array}$ \\
    \hline
     6 & $\begin{array}{ccccccc}
            2 & 4 & 6 & 5 & 4 & 2 & 1 \\
              &  &  3
         \end{array}$ &
         $\begin{array}{cccccccc}
            1 & 1 & 1 & 1 & 1 & 1 & 0 & 0\\
         \end{array}$ &
         $\begin{array}{ccccccc}
            0 & 0 & 0 & 0 & 1 & -1 & 0 \\
         \end{array}$ \\
    \hline
     7 & $\begin{array}{ccccccc}
            2 & 4 & 6 & 5 & 4 & 3 & 1 \\
              &  &  3
         \end{array}$ &
         $\begin{array}{cccccccc}
            1 & 1 & 1 & 1 & 1 & 1 & 1 & 0\\
         \end{array}$ &
         $\begin{array}{ccccccc}
            0 & 0 & 0 & 0 & 0 & 1 & -1 \\
         \end{array}$ \\
     \hline
     8 & $\begin{array}{ccccccc}
            2 & 4 & 6 & 5 & 4 & 3 & 2 \\
              &  &  3
         \end{array}$ &
         $\begin{array}{cccccccc}
            1 & 1 & 1 & 1 & 1 & 1 & 1 & 1\\
         \end{array}$ &
         $\begin{array}{ccccccc}
            0 & 0 & 0 & 0 & 0 & 0 & 1 \\
         \end{array}$ \\
    \hline
       \end{tabular}
  \vspace{2mm}
  \caption{\hspace{3mm} $8$ pairs of positive roots $\{\eta, \lambda\}$ such that $\eta^{\nabla} = \lambda^{\nabla}$.
     There are also $8$ pairs of negative roots $\{\eta, \lambda\}$ such that $\eta^{\nabla} = \lambda^{\nabla}$}
  \label{tab_pairs_roots_A7extE8}
  \end{table}

\subsection{The linkage systems  $\mathscr{L}(A_l)$ for $l \geq 8$}
    \label{sec_Al_lg8}
    \index{Dynkin extension ! - $A_l <_D D_{l+1}$}
    \index{$\Gamma <_D \Gamma'$, Dynkin extension}
    \index{linkage system  ! - $\mathscr{L}(A_l)$ for $l \geq 8$}

\begin{figure}[H]
\centering
 \vspace{1cm}
\includegraphics[scale=0.47]{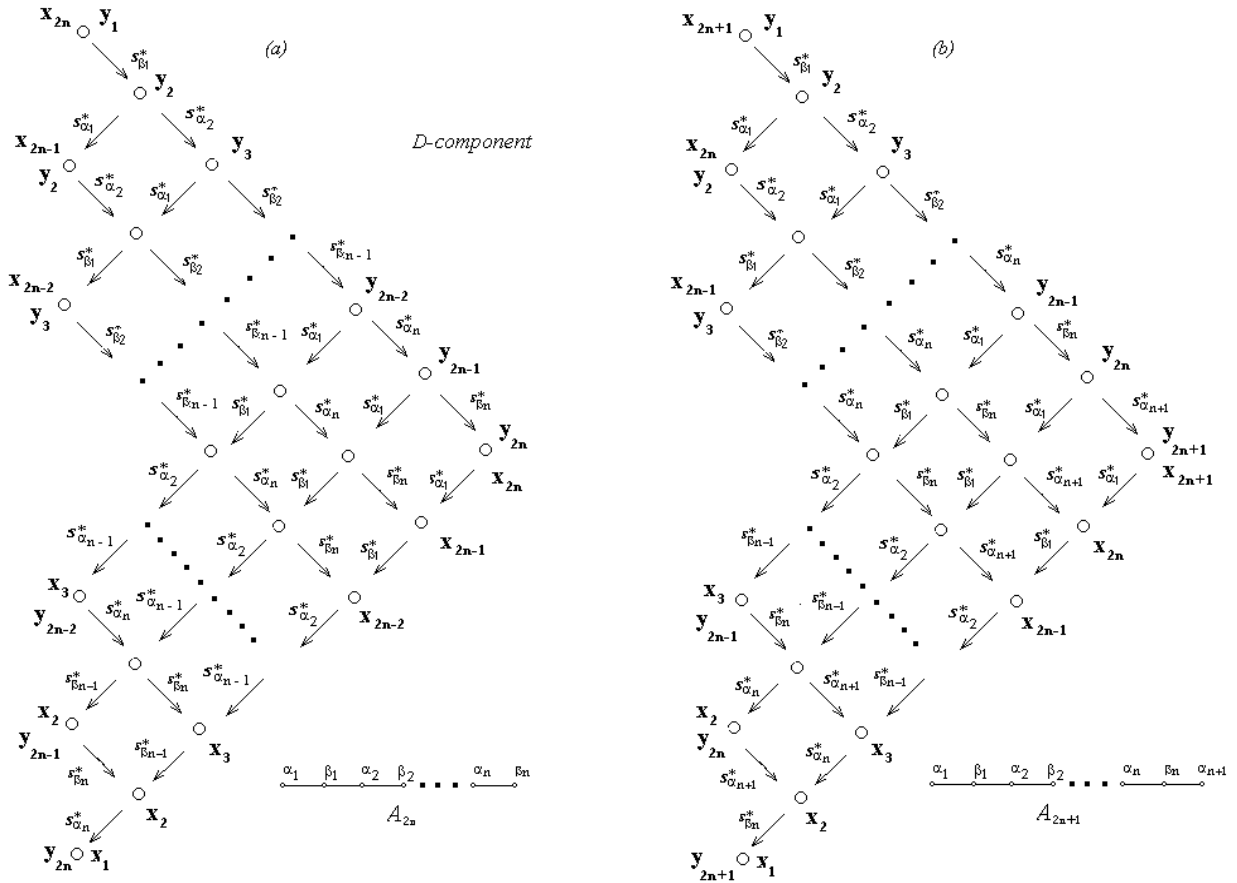}
\caption{\hspace{3mm} Coordinates $\{x_i, y_j\}$ in the $D$-component
of the linkage system $A_{2n+1}$}
%%%%%% The label must come after caption
\label{Al_coord}
\end{figure}

 For convenience, we introduce coordinates $\{x_i, y_j\}$, where $1 \leq i,j \leq l$, for every
 linkage diagram in the $D$-component of the linkage system
 $\mathscr{L}_{D_{l+1}}(A_l)$.
 In Fig. \ref{Al_coord}, the path $y_i$ (resp. $x_i$), where $1 \leq i \leq l$, consists of reflections acting
 in the direction from Nord-West to South-East (resp. from Nord-East to South-West).
 The linkage diagram lying in the intersection node of paths $x_i$ and $y_j$  has coordinates $\{x_i, y_j\}$.

\begin{figure}[H]
\centering
\includegraphics[scale=0.8]{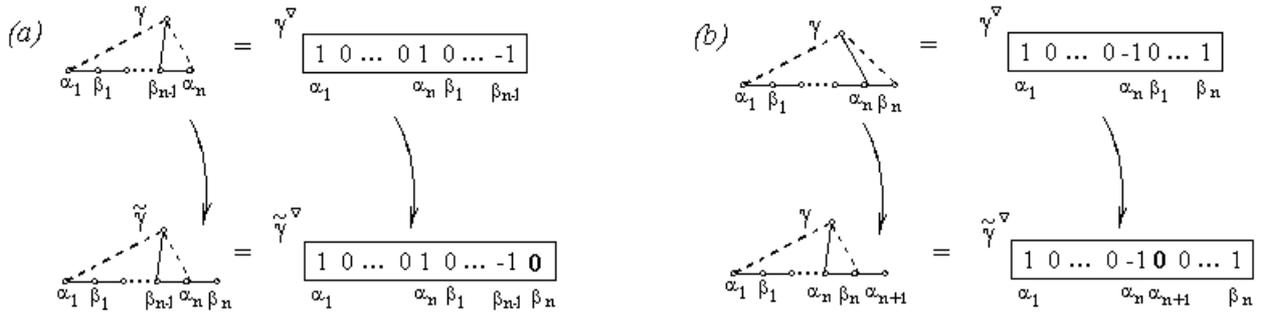}
\caption{\hspace{3mm} Extension $\gamma^{\nabla} \rightarrow  \widetilde\gamma^{\nabla}$
of linkage diagrams $(a)$ and $(b)$ and equivalent linkage label vectors. Additional coordinates are in bold}
%%%%%% The label must come after caption
\label{one_diagr_An}
\end{figure}

\begin{figure}[H]
\centering
 \vspace{1cm}
\includegraphics[scale=0.5]{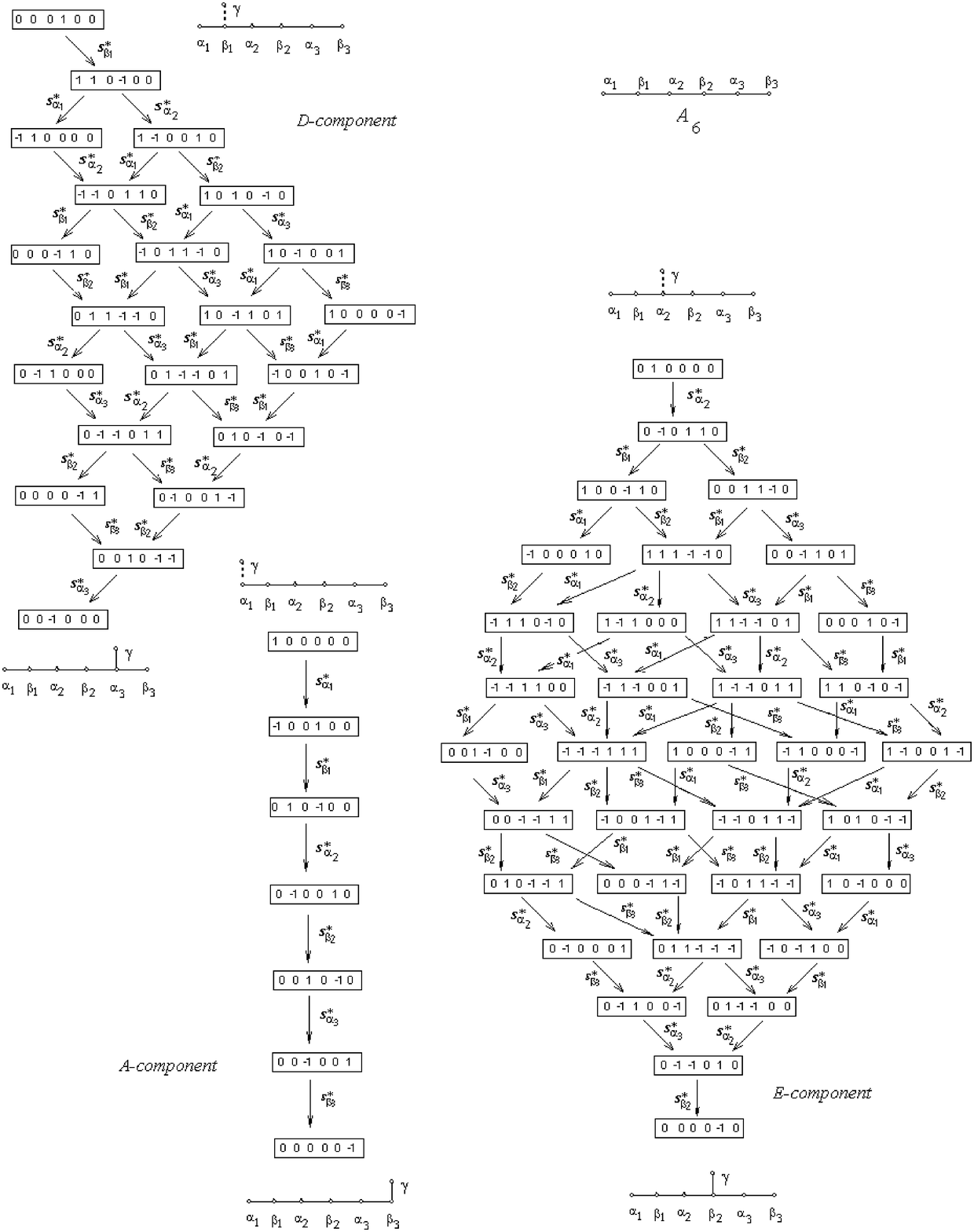}
\caption{\hspace{3mm}The linkage system $A_6$ has $6$ components: The two parts
 of the $A$-component, the two $21$-element parts of the $D$-component, and the two $35$-element parts of the $E$-component}
%%%%%% The label must come after caption
\label{A6_to_A7D7E7}
\end{figure}

\begin{figure}[H]
\centering
 %%\vspace{1cm}
\includegraphics[scale=0.5]{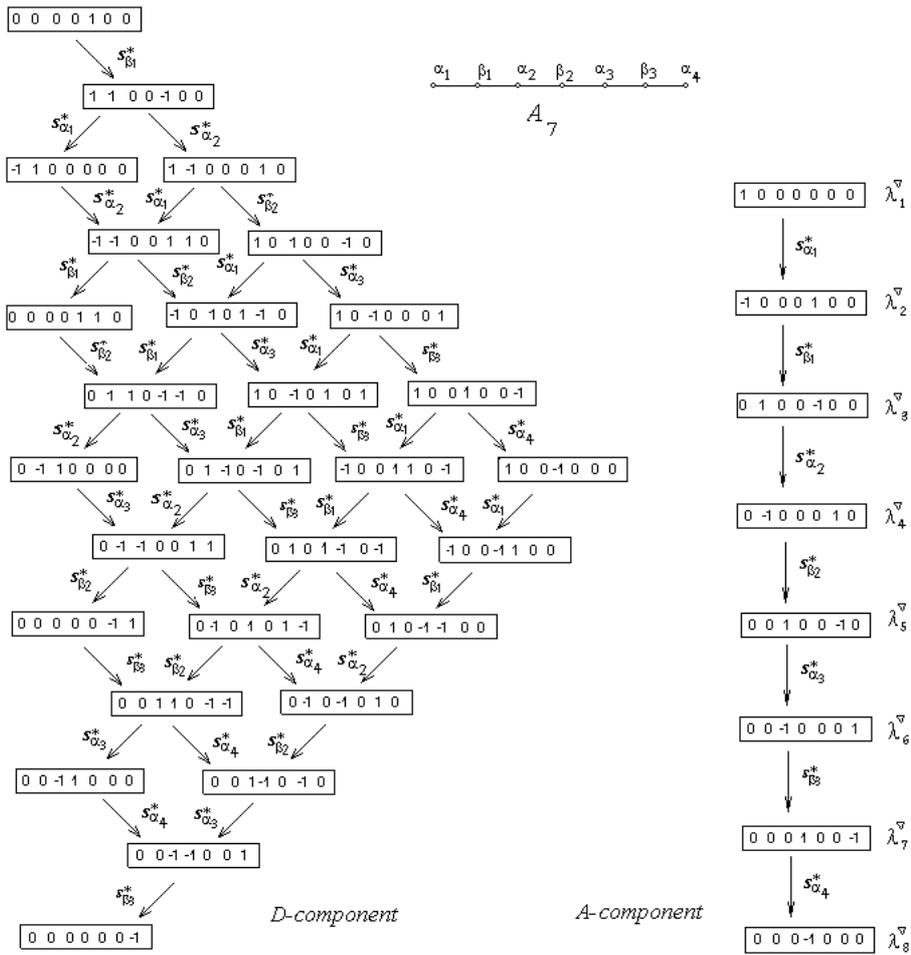} %% 0.5
\vspace{1cm}
\caption{\hspace{3mm} The linkage system $A_7$. One of the two $8$-element parts
 of the $A$-component (see Table \ref{tab_pairs_roots_A7extE8}) and one of the two $28$-element parts of the $D$-component}
%%%%%% The label must come after caption
\label{A7_to_A8D8_1}
\end{figure}

\begin{figure}[H]
\centering
 \vspace{1cm}
\includegraphics[scale=0.7]{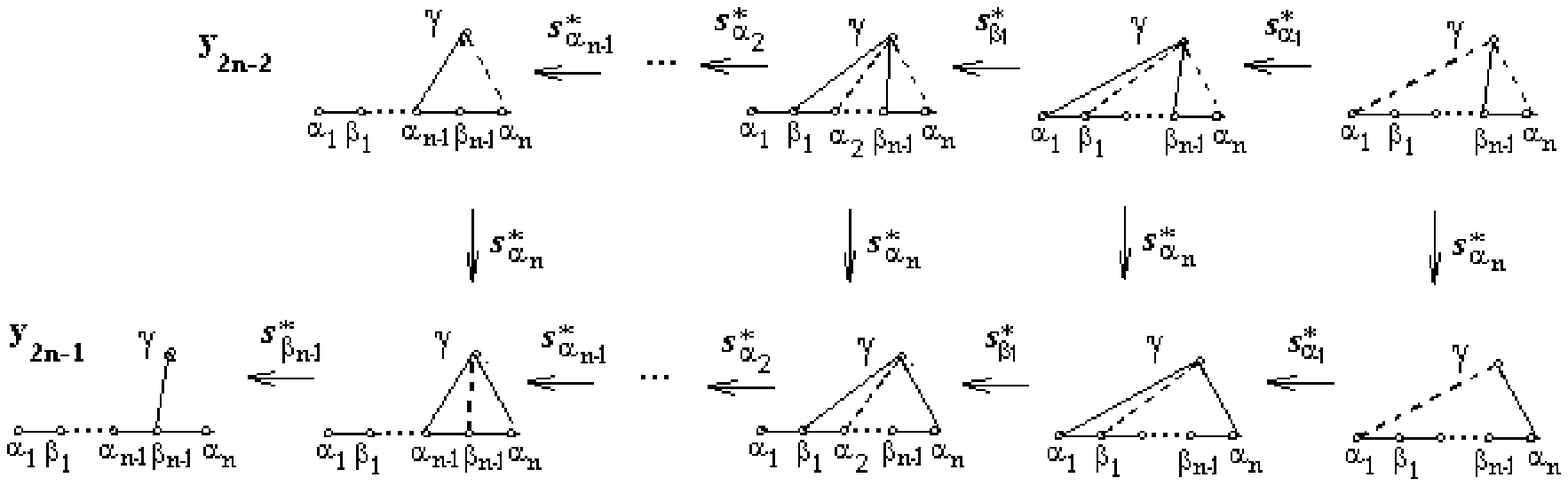}
\caption{\hspace{3mm} The last $2$ lines of the $D$-component of $A_{2n-1}$}
%%%%%% The label must come after caption
\label{2_lines_A2nMin1}
 \vspace{1cm}
\includegraphics[scale=0.7]{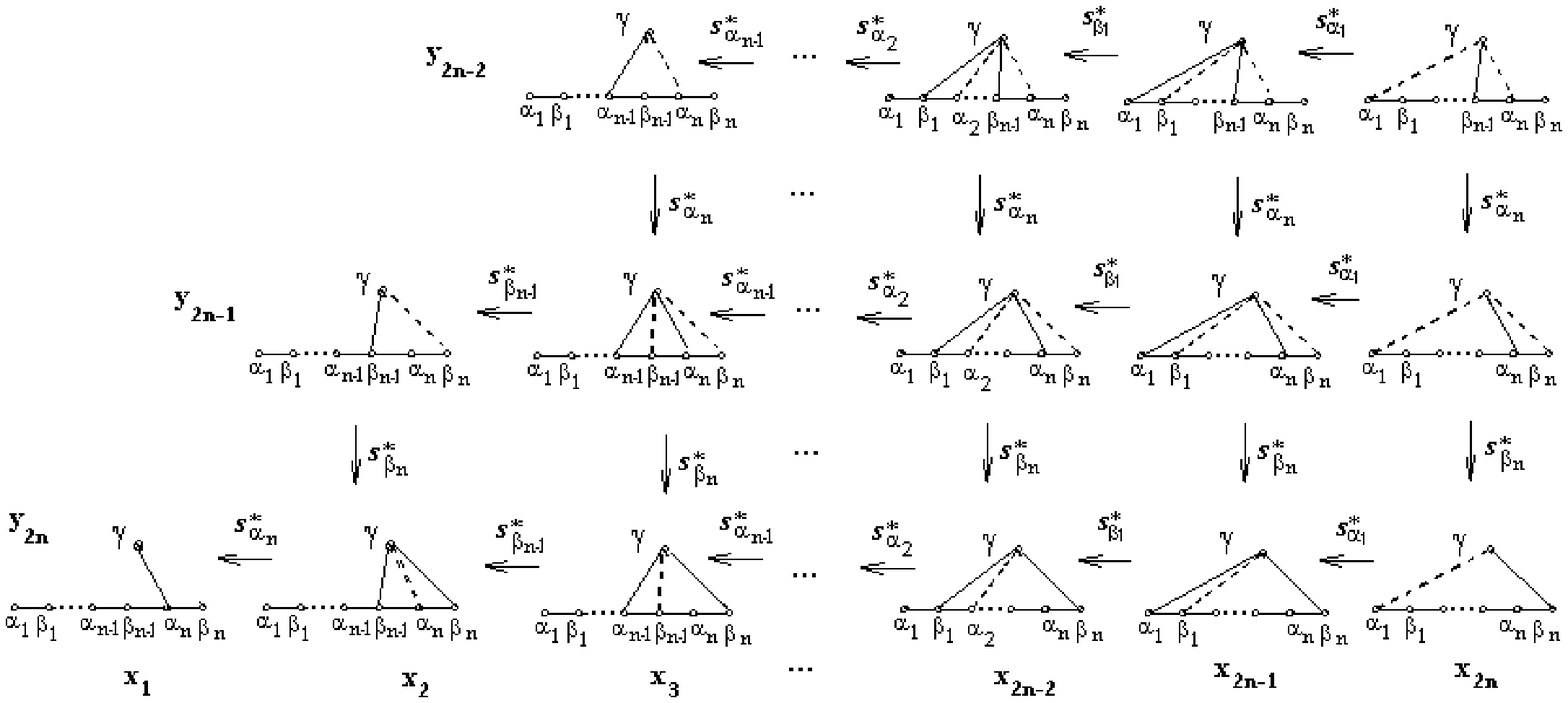}
\caption{\hspace{3mm} The last $3$ lines of the $D$-component of $A_{2n}$}
%%%%%% The label must come after caption
\label{3_lines_A2n}
 \vspace{1cm}
\includegraphics[scale=0.7]{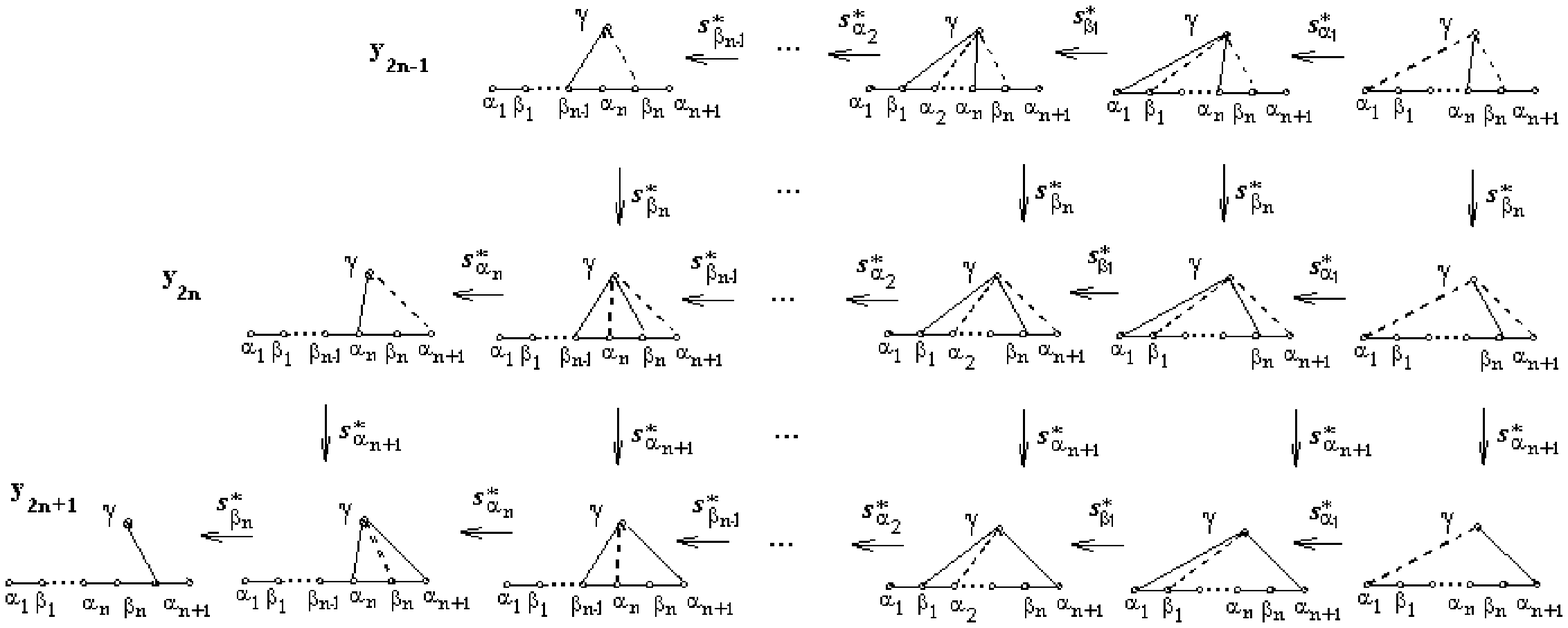}
\caption{\hspace{3mm} The last $3$ lines of the $D$-component of $A_{2n+1}$}
%%%%%% The label must come after caption
\label{3_lines_A2nPl1}
\end{figure}

 \begin{proposition}[On the $D$-component of the linkage systems  $\mathscr{L}(A_l)$]
  \label{ptop_Dcomp}
   Consider the Dynkin extension $A_l <_D D_{l+1}$.

  ${\rm (i)}$
   For $l > 3$, the $D$-component of the linkage systems  $\mathscr{L}(A_l)$ consists of two non-connected parts $\mathscr{L}^+$ and
   $\mathscr{L}^-$ corresponding to the positive (resp. negative) linkage roots in the root stratum $\varPhi(D_{l+1})\backslash\varPhi(A_l)$.
   For $l = 3$, components $\mathscr{L}^+$ and $\mathscr{L}^-$ coincide, see Fig. \ref{A3_to_A4D4}.

   For $l=2n$ (resp. $l=2n+1$), where $n \geq 2$ (resp. $n \geq 1$), reflections acting on linkage diagrams of
   $\mathscr{L}^+$ are given in Fig. \ref{Al_coord}$(a)$ (resp. Fig. \ref{Al_coord}$(b)$);
   the linkage diagram $\gamma^{\nabla}$ in the top of $\mathscr{L}^+$ (vertex $\{x_{2n}, y_1\}$ (resp. $\{x_{2n+1}, y_1\}$))
   depicted in Fig. \ref{A2n_top_root}$(a)$ (resp. $(b)$).

 \begin{figure}[H]
 \centering
 \includegraphics[scale=0.6]{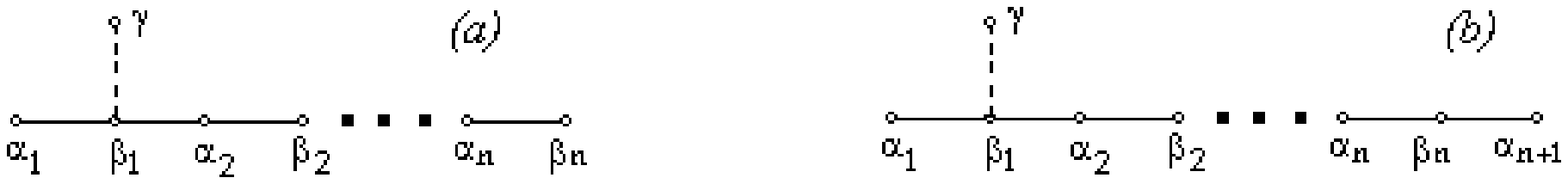}  \\
 \begin{equation*}
 \begin{array}{lllllllll}
  \gamma^{\nabla} = &\{0,~~ \dots, &0,~~~&1,~~0,\dots, &0 \} \qquad\qquad
  \gamma^{\nabla} = &\{0,~~ \dots, &0,~~~&1,~~0, \dots, &0 \}\\
    &~\alpha_1 &\alpha_n &\beta_1 &\beta_n &~\alpha_1 &\alpha_{n+1} &\beta_1 &\beta_n
 \end{array}
 \end{equation*}
 \caption{}
 \label{A2n_top_root}
 \end{figure}
   Remaining linkage diagrams of $\mathscr{L}^+$ are constructed from $\gamma$ by means of reflections
\begin{equation}
  \begin{split}
   & s^*_{\alpha_1}, s^*_{\beta_1}, s^*_{\alpha_2}, s^*_{\beta_2}, \dots, s^*_{\alpha_n}, s^*_{\beta_n} \text{ for } l = 2n, \text{ or }, \\
   & s^*_{\alpha_1}, s^*_{\beta_1}, s^*_{\alpha_2}, s^*_{\beta_2}, \dots, s^*_{\alpha_n}, s^*_{\beta_n}, s^*_{\alpha_{n+1}}
   \text{ for } l = 2n + 1.
  \end{split}
\end{equation}
   see Fig. \ref{Al_coord}.
   The three last lines of $\mathscr{L}^+$  corresponding to coordinates
   $y_{2n-2}$, $y_{2n-1}$ and $y_{2n}$ (resp. $y_{2n-1}$, $y_{2n}$ and $y_{2n+1}$ )
   are given in Fig. \ref{3_lines_A2n} (resp. Fig. \ref{3_lines_A2nPl1}).

   (For $3 \leq l \leq 7$, see Figs. \ref{A3_to_A4D4} - \ref{A7_to_A8D8_1}).
~\\

  ${\rm (ii)}$
   For $l > 3$, we have
 \begin{equation}
   \mid \mathscr{L}_{D_{l+1}(A_l)} \mid = l(l+1).
 \end{equation}
 \end{proposition}

\begin{figure}[H]
\centering
\includegraphics[scale=0.42]{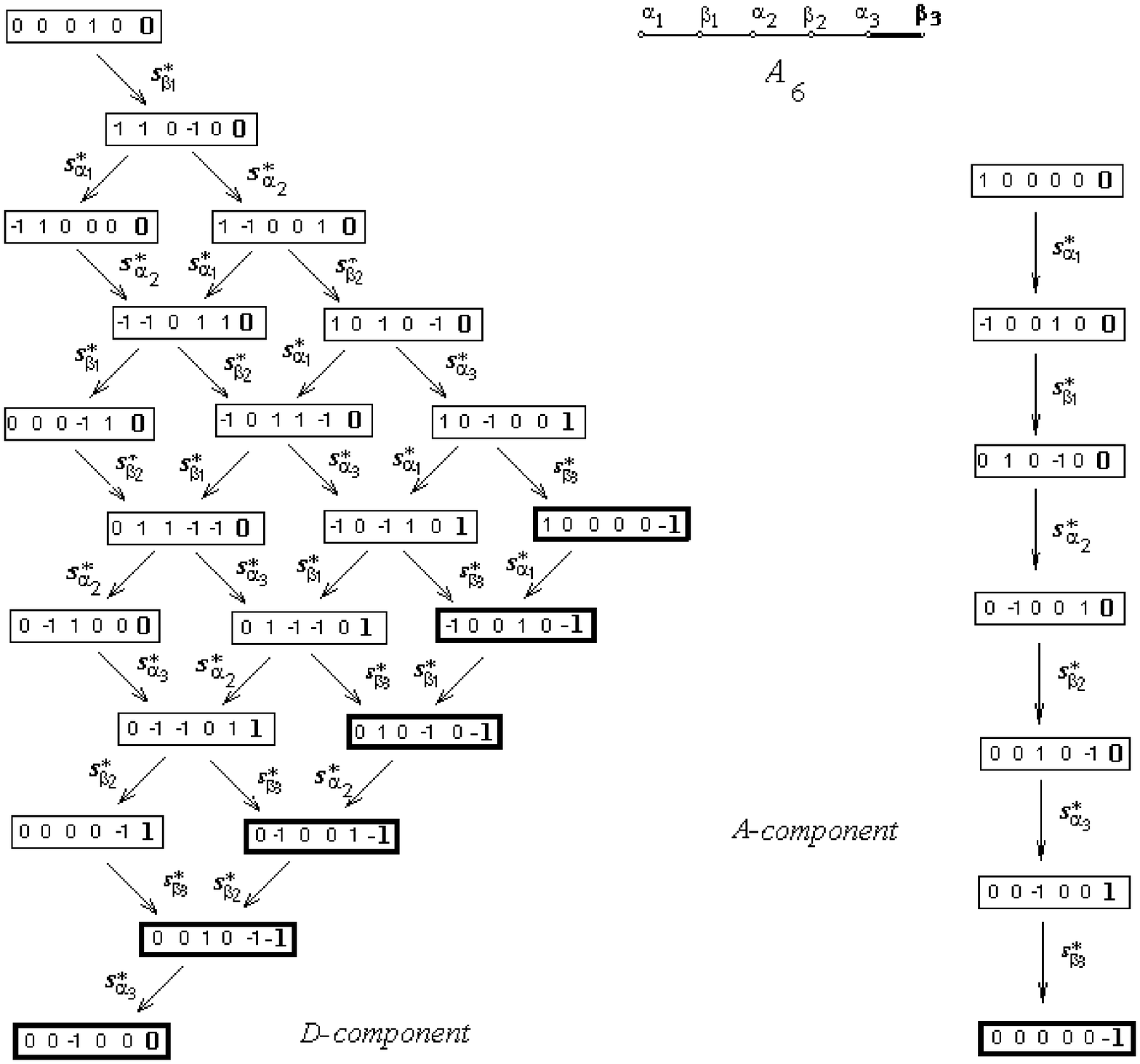}
\vspace{1cm}
\includegraphics[scale=0.48]{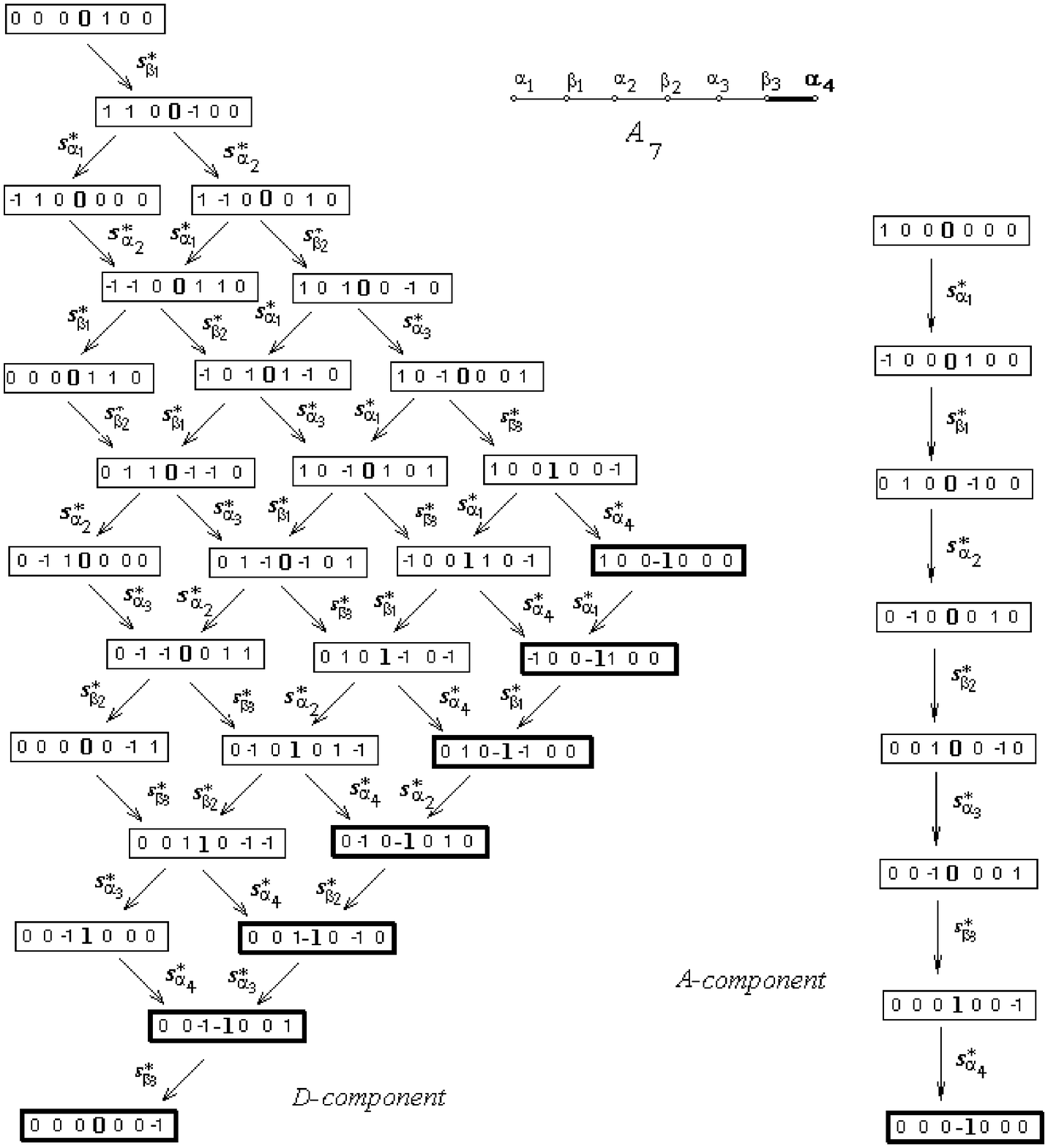}
\caption{\hspace{3mm}The transition
  $\mathscr{L}(A_5)  \Longrightarrow \mathscr{L}(A_6)$ and $\mathscr{L}(A_6)  \Longrightarrow \mathscr{L}(A_7)$.
  The new coordinate $\beta_3$ for the case $\mathscr{L}(A_6)$ (resp. $\alpha_4$ for the case $\mathscr{L}(A_7)$)
  and new linkage diagrams are in bold}
%%%%%% The label must come after caption
\label{A5_to_A6D6_pass}
\end{figure}

\begin{figure}[H]
\centering
\includegraphics[scale=0.6]{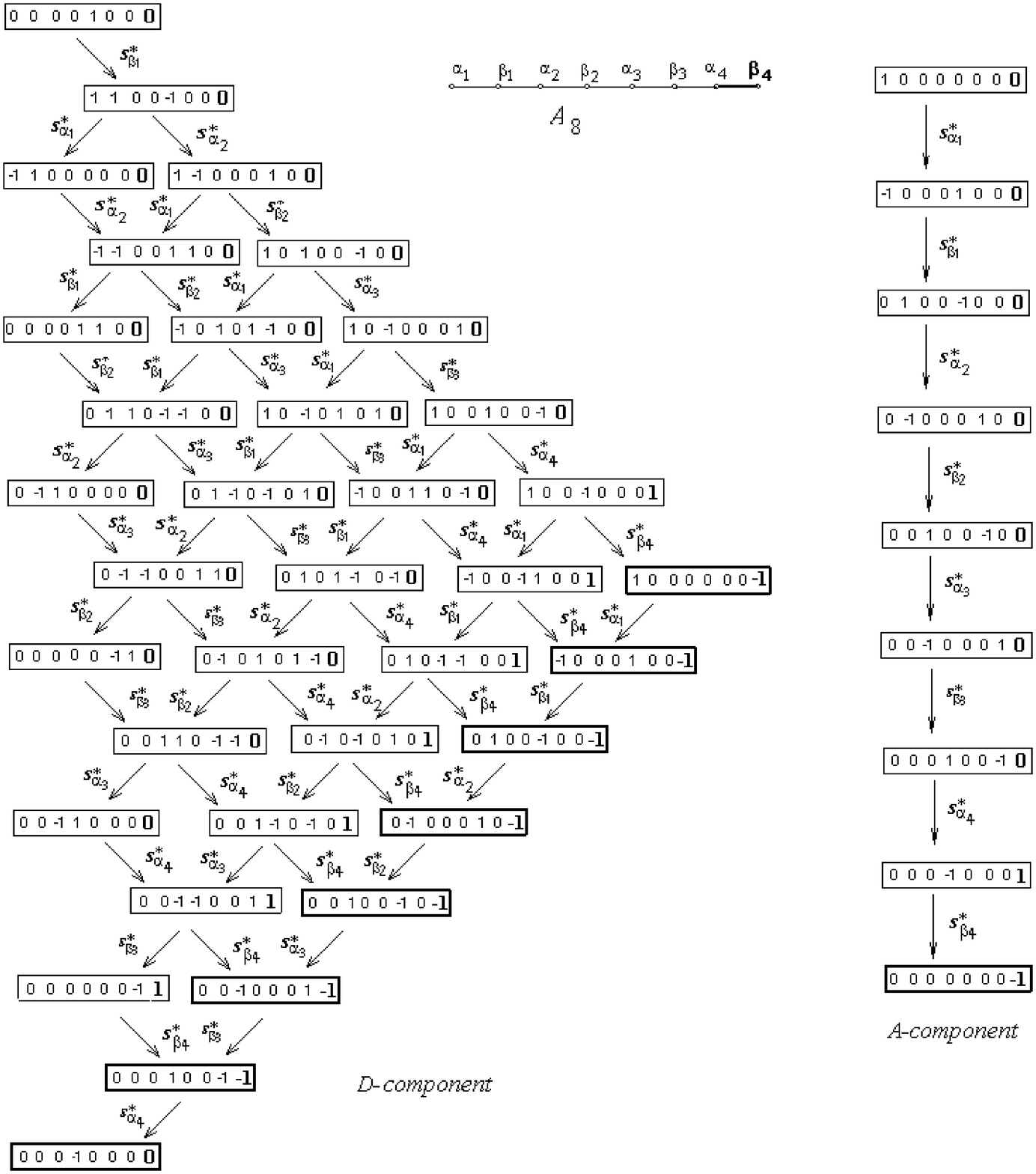}
\caption{\hspace{3mm}The transition from $A$-components and
 $D$-components in the linkage system $A_7$ to the linkage system $A_8$.
 The new coordinate $\beta_4$ and new linkage diagrams are marked in bold}
%%%%%% The label must come after caption
%%\label{A7_to_A8D8}
\label{A8_to_A9D9_mono}
\end{figure}

   \PerfProof (i) For transitions
  $\mathscr{L}(A_5)  \Rightarrow \mathscr{L}(A_6)$ and $\mathscr{L}(A_6)  \Rightarrow \mathscr{L}(A_7)$,
  see Fig. \ref{A5_to_A6D6_pass}. For the transition $\mathscr{L}(A_7)  \Rightarrow \mathscr{L}(A_8)$,
  see Fig. \ref{A8_to_A9D9_mono}. We prove the proposition by induction.
~\\

  \underline{Transition $2n-1 \Rightarrow 2n$.} Consider the triangle
  $\Delta_{2n-2} = \{\{y_1, x_{2n}\}, \{y_{2n-2},x_{2n}\}, \{y_{2n-2},x_3\}\}$ is the linkage subsystem
  containing linkage diagrams with the coordinate $y_i$ for $i \leq 2n-2$, see Fig. \ref{Al_coord}$(a)$.
  To transit from $2n-1$ to $2n$ coordinates,
  we extend any linkage diagram  $\gamma^{\nabla}$ in the triangle $\Delta_{2n-2}$ and in the line $y_{2n-2}$
  in Fig. \ref{2_lines_A2nMin1} by the vertex $\beta_n$ non-connected with $\gamma$. In other words, we have
  $\gamma^{\nabla}_{\beta_n} = (\gamma, \beta_n) = 0,$ see Figs. \ref{A8_to_A9D9_mono}, \ref{Al_coord}$(a)$.
   Remember that any linkage diagram is equivalent to a certain
   linkage label vector, see \S\ref{sec_linkage_diagr_0}. For the extension of the linkage diagram and
   the equivalent linkage label vectors, see Fig. \ref{one_diagr_An}$(a)$.
~\\
  By induction hypothesis,
  any reflection $s^*_{\varphi}$, which acts within the triangle $\Delta_{2n-2}$,
  differs from $s^*_{\alpha_n}$ and $s^*_{\beta_n}$, see Fig. \ref{Al_coord}.
  So $s^*_{\varphi}$ \underline{do not touch the coordinate $\gamma^{\nabla}_{\beta_n}$}.
  We put $\gamma^{\nabla}_{\beta_n} = 0$ to every linkage $\gamma^{\nabla}$ from $\Delta_{2n-2}$ and
  the subsystem given by the triangle $\Delta_{2n-2}$ is saved in the transition $2n-1 \Rightarrow 2n$.
  We see this process of the extension by the vertex $\beta_n$ in lines $y_{2n-2}$ of diagrams in Fig. \ref{2_lines_A2nMin1}
  and Fig. \ref{3_lines_A2n}.

~\\
  By induction hypothesis, we have $\gamma^{\nabla}_{\alpha_n} = 1$
  for any linkage diagram $\gamma^{\nabla}$ in the line $y_{2n-2}$ in $\Delta_{2n-2}$, see Fig. \ref{2_lines_A2nMin1}.
  The same holds for
  extensions of these diagrams in the line $y_{2n-2}$ in Fig. \ref{3_lines_A2n}.  Consider now
  the last lines of the linkage system $\mathscr{L}(A_{2n})$ in  Fig. \ref{3_lines_A2n}.
  By acting $s^*_{\alpha_n}$ on every linkage diagram $\gamma^{\nabla}$ in the line $y_{2n-2}$, we get linkage diagrams
  with $\gamma^{\nabla}_{\beta_n} = 1$  in the line $y_{2n-1}$.
  Really, by \eqref{dual_refl} coordinates
  $\gamma^{\nabla}_{\alpha_n}$ and $\gamma^{\nabla}_{\beta_n}$
  for any linkage diagram $\gamma^{\nabla}$ in the line $y_{2n-2}$ are transformed as follows:
  \begin{equation*}
     (s^*_{\alpha_n}\gamma^{\nabla})_{\beta_n} = \gamma^{\nabla}_{\beta_n} + \gamma^{\nabla}_{\alpha_n}, \qquad
     (s^*_{\alpha_n}\gamma^{\nabla})_{\alpha_n} = -\gamma^{\nabla}_{\alpha_n},
  \end{equation*}
  see \eqref{dual_refl_alpha_beta}.
  The linkage diagram $\{x_2, y_{2n-1}\}$ in the line $y_{2n-1}$ in Fig. \ref{3_lines_A2n}
  is obtained by acting $s^*_{\beta_{n-1}}$ on the linkage diagram $\{x_3, y_{2n-1}\}$.
  Further, by acting $s^*_{\beta_n}$ on every linkage diagram
  (except for linkage diagrams $\{x_2, y_{2n}\}$ and $\{x_1, y_{2n}\}$)
  in the line $y_{2n-1}$ we get the last line with $\gamma^{\nabla}_{\alpha_n} = 0$ and $\gamma^{\nabla}_{\beta_n} = -1$.
  The linkage diagram $\{x_2, y_{2n}\}$ (resp. $\{x_1, y_{2n}\}$) is obtained
  by acting $s^*_{\beta_{n-1}}$ (resp. $s^*_{\alpha_n}$) on the linkage diagram $\{x_3, y_{2n}\}$ (resp. $\{x_2, y_{2n}\}$).
~\\

  \underline{Transition $2n \Rightarrow 2n+1$} is considered in a similar way.
~\\

  (ii) It is easy to see from Fig. \ref{Al_coord} that
  \begin{equation*}
    \mid \mathscr{L}^+ \mid =  \mid \mathscr{L}^+ \mid = \frac{l(l+1)}{2}.
  \end{equation*}
  Since $\mathscr{L}_{D_{l+1}(A_l)} = \mathscr{L}^+ \coprod \mathscr{L}^-$ the statement is proved.

\qed

 \begin{proposition}[On the $A$-component of the linkage systems  $\mathscr{L}(A_l)$]
   The $A$-component of the linkage systems  $\mathscr{L}(A_l)$ consists of two non-connected parts $\mathscr{L}^+$ and
   $\mathscr{L}^-$ corresponding to the positive (resp. negative) linkage roots in the Dynkin extension $A_l <_D A_{l+1}$.
   The number of linkage diagrams in the $A$-component of the linkage systems  $\mathscr{L}(A_l)$ is equal to $2(l+1)$.
 \end{proposition}

 \PerfProof This statement is also proved by induction. Consider, the \underline{transition ${2n-1} \Rightarrow 2n$}.
 For $n = 3$ this process is shown in Fig. \ref{A5_to_A6D6_pass}, see $A_5 \Rightarrow A_6$.
 The first $2n-1$ linkage diagrams are extended by the vertex $\gamma^{\nabla}_{\beta_{n}} = 0$.
 Under the action of $s^*_{\alpha_n}$ the ($2n$)th linkage diagram is extended by the coordinate $\gamma^{\nabla}_{\beta_n} = 1$.
 Under the action of $s^*_{\beta_n}$ the ($2n+1$)th linkage diagram is added,
 the coordinate $\gamma^{\nabla}_{\beta_{n}}$ of this linkage diagram is equal $-1$.
 \qed

\clearpage

\clearpage

%% file: AB_treesMatrix.tex
\clearpage
\section{\sc\bf Some properties of Carter and connection diagrams}
  \label{sec_Dynkin_diagr}

\subsection{Similarity of Carter diagrams} %%%  operating on connection diagrams}
  \label{seq_sim_diagram}
 Talking about a certain diagram $\Gamma$ we actually have in
 mind a set of roots with orthogonality relations as it is prescribed by the diagram $\Gamma$. %% diagram = set of roots
 We try to find some common properties of sets of roots (from the root systems associated with
 the simple Lie algebras) and diagrams associated with these sets. %% we study sets of toots
 These diagrams are not necessarily Dynkin diagrams since sets of roots we study
 are not necessarily sets of simple roots and are not
 root subsystems. %% diagrams are not necessarily Dynkin diagram
 We use the term \lq{\lq}Dynkin diagram\rq\rq to describe connected sets of
 linearly independent \underline{simple roots} in the root system. %% why we use term Dynkin diagram
 Similarly, \lq{\lq}Carter diagrams{\rq\rq} describe connected sets of linearly independent roots,
 not necessarily simple, and such that any cycle is even. %% even cycles

\begin{figure}[h]
\centering
\includegraphics[scale=0.7]{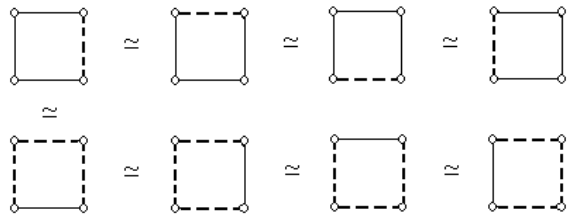}
\caption{\hspace{3mm}Eight similar $4$-cycles equivalent to $D_4(a_1)$}
%%%%%% The label must come after caption
\label{8-equiv_4cycles}
\end{figure}

 \index{similarity of Carter diagrams}
 \index{connection diagram}
 Two connection diagrams obtained from each other by a sequence
 of reflections \eqref{eq_1_equiv}, are said to be {\it similar} connection diagrams, see Fig. \ref{8-equiv_4cycles}.
\begin{equation}
  \label{eq_1_equiv}
   \alpha \longmapsto -\alpha.
\end{equation}
 A transformation of connection diagrams obtained by a sequence of reflections \eqref{eq_1_equiv}
 is said to be a {\it similarity transformation} or {\it similarity}.
\index{similar Carter diagrams}

 By applying similarity \eqref{eq_1_equiv} any solid edge with an endpoint vertex being $\alpha$
 can be changed to a dotted one and vice versa; this does not change, however, the corresponding reflection:

\begin{equation*}
   s_{\alpha} = s_{-\alpha}.
\end{equation*}

  \begin{remark}[On trees]
   \label{rem_tree}
  {\rm For the set
  $\{ \alpha_1, \dots, \alpha_i, \alpha_{i+1}, \dots \alpha_n \}$
  forming a tree, we may assume that, up to the similarity, all non-zero inner products
  $(\alpha_i, \alpha_j)$ are negative. Indeed, if $(\alpha_i, \alpha_j) >  0$,
  we apply similarity transformation $\alpha_j \longmapsto -\alpha_j$, consider all
  inner products $(\alpha_k, \alpha_j) > 0$ and repeat similarity transformations
  $\alpha_k  \longmapsto -\alpha_k$ if necessary. This process converges
  since the diagram is a tree.
  }
 \end{remark}

  \subsection{The ratio of lengths of roots}
   \label{sec_trees}
   \index{obtuse angle between roots}
  Let $\Gamma$ be a Dynkin diagram, and
  $\sqrt{t}$ be the ratio of the length of any long root to the length of any short root.
  The inner product between two long roots is
$$
   (\alpha, \beta) = \sqrt{t}\cdot\sqrt{t}\cdot\cos(\widehat{\alpha, \beta}) =
   \sqrt{t}\cdot\sqrt{t}\cdot \big ( \pm\frac{1}{2} \big ) =
   \pm\displaystyle\frac{t}{2}.
$$
 By Remark \ref{rem_tree}, we may put $(\alpha, \beta) = -\displaystyle\frac{t}{2}$.
The inner product between two short roots is
$$
   (\alpha, \beta) = \cos(\widehat{\alpha, \beta}) = \pm\displaystyle\frac{1}{2}.
$$
 Again, by Remark \ref{rem_tree}, we may put $(\alpha, \beta) =  -\displaystyle\frac{1}{2}$.
 The inner product $(\alpha, \beta)$ between roots of different lengths is
$$
   (\alpha, \beta) = 1\cdot\sqrt{t}\cdot\cos(\widehat{\alpha, \beta}) =
   1\cdot\sqrt{t}\cdot \big ( \pm\frac{\sqrt{t}}{2} \big ) =
   \pm\displaystyle\frac{t}{2}.
$$
As above, we choose the obtuse angle and put
  $(\alpha, \beta) = -\displaystyle\frac{t}{2}$.
~\\
We can summarize:
\begin{equation}
 \label{eq_all_inner_prod}
  (\alpha, \beta) =
  \begin{cases}
    -\frac{1}{2} \quad  \text{for} & \norm{\alpha} = \norm{\beta} = 1, \\
    \\
    -1 \quad  \text{for} & \norm{\alpha}  = \norm{\beta} =  2, \quad \text{or} \quad \norm{\alpha}  = 1, ~\norm{\beta} = 2, \\
    \\
    -\frac{3}{2} \quad  \text{for} & \norm{\alpha} = \norm{\beta} = 3, \quad \text{or} \quad \norm{\alpha}  = 1, ~\norm{\beta} = 3, \\
  \end{cases}
\end{equation}
where all angles $\widehat{\alpha, \beta}$ are obtuse.

  \subsection{Basic lemmas}
    \label{sec_basic_lem}
  \begin{lemma}
    \label{lem_must_dotted}
  There is no root subset (in the root system associated with a Dynkin diagram)
  forming a simply-laced cycle containing only solid edges.
  Every cycle in the Carter diagram or in the connection diagram contains at least one solid edge
  and at least one dotted edge.
  \end{lemma}
  \PerfProof Suppose a subset $S = \{ \alpha_1, \dots,  \alpha_n\} \subset \varPhi $
 forms a cycle containing only solid edges. Consider the vector
 \begin{equation*}
     v = \sum\limits_{i=1}^n \alpha_i.
 \end{equation*}
 The value of the quadratic Tits form $\mathscr{B}$ (see \cite{St08}) on $v$ is equal to
 \begin{equation*}
    \mathscr{B}(v) =  \sum\limits_{i \in \Gamma_0} 1 - \sum\limits_{i \in \Gamma_1} 1 =  n - n = 0,
 \end{equation*}
 where $\Gamma_0$ (resp. $\Gamma_1$) is the set of all vertices (resp. edges)
 of the diagram $\Gamma$ associated with $S$.
 Therefore, $v = 0$ and elements of the root subset $S$ are linearly dependent. \qed
 \\

   The following proposition is true only for trees.
\begin{lemma}[Lemma 8, \cite{Ca72}]
  \label{lem_non_ext_Dynkin}
  Let $S = \{ \alpha_1, \dots,  \alpha_n\} $ be a subset of linearly independent (not necessarily simple)
  roots of a certain root system $\varPhi$, and $\Gamma$ be the Dynkin diagram corresponding to $\varPhi$.
  Let $\Gamma_S$  the Carter diagram or the connection diagram  associated with $S$.
  If $\Gamma_S$ is a tree, then $\Gamma_S$ is a Dynkin diagram.

\end{lemma}
\PerfProof
 If $\Gamma_S$ is not a Dynkin diagram, then $\Gamma_S$ contains an
 extended Dynkin diagram $\widetilde{\Gamma}$ as a subdiagram. Since $\Gamma_S$
 is a tree, we can turn all dotted edges to solid ones\footnotemark[1],
 see Remark \ref{rem_tree}.
\footnotetext[1]{This fact is not true for cycles, since by Lemma
\ref{lem_must_dotted} we cannot
 eliminate all dotted edges.}
 Further, we consider the vector
 \begin{equation}
    \label{eq_lin_dep}
        v = \sum\limits_{i \in \widetilde{\Gamma}_0}{t_i\alpha_i},
 \end{equation}
 where $\widetilde{\Gamma}_0$ is the set of all vertices of $\widetilde{\Gamma}$, and $t_i$ (where $i \in \widetilde{\Gamma}_0$)
 are the coefficients of the nil-root, see \cite{Kac80}. Let the remaining coefficients corresponding
 to  $\Gamma_S \backslash  \widetilde{\Gamma}$ be equal to $0$.
 Let $\mathscr{B}$ be the positive definite quadratic Tits form (see \cite{St08}) associated with the diagram $\Gamma$,
 and $(\cdot\hspace{0.7mm},\cdot)$ the symmetric bilinear form associated with $\mathscr{B}$.
 Let  $\{\delta_i \mid i \in \widetilde{\Gamma}_0\}$
 be the set of simple roots associated with vertices $\widetilde{\Gamma}_0$.
 For all $i,j \in \widetilde{\Gamma}_0$, we have $(\alpha_i, \alpha_j) = (\delta_i, \delta_j)$,
 since this value is described by edges of $\widetilde{\Gamma}$. Therefore,
 \begin{equation*}
    \mathscr{B}(v) =
    \sum\limits_{i,j \in \widetilde{\Gamma}_0}t_i{t}_j(\alpha_i, \alpha_j) =
    \sum\limits_{i,j \in \widetilde{\Gamma}_0}t_i{t}_j(\delta_i, \delta_j) =
    \mathscr{B}(\sum\limits_{i \in \widetilde{\Gamma}_0}{t_i\delta_i}) = 0.
 \end{equation*}
 Since $\mathscr{B}$ is a positive definite form, we have $v = 0$, i.e.,
 vectors $\alpha_i$ are linearly dependent. This contradicts the definition of the set $S$.
 \qed

 \index{similarity transformation for Carter diagrams}
\begin{lemma}
  \label{lem_case_An} The root system $A_n$ does not contain $\Gamma$-associated root subsets for

  {\rm ~(i)} $\Gamma$ is a cycle of length $\geq 4$,

  {\rm (ii)} $\Gamma$ is $D_4$.
\end{lemma}
 \PerfProof
   (i) Recall that any root in $A_n$ is of the form $\pm(e_i - e_j)$, where $1 \leq i < j \leq n+1$.
 Then, up to the similarity $\alpha \longmapsto -\alpha$,
 a cycle of roots is of one of the followings forms:
 \begin{equation*}
  %%\label{eq_An_1}
  \begin{split}
   & \{ e_{i_1} - e_{i_2},  e_{i_2} - e_{i_3}, \dots,
    e_{i_{k-1}} - e_{i_k}, e_{i_{k}} - e_{i_1} \}, \\
   & \{ e_{i_1} - e_{i_2},  e_{i_2} - e_{i_3}, \dots,
    e_{i_{k-1}} - e_{i_k}, -(e_{i_{k}} - e_{i_1}) \}.
   \end{split}
 \end{equation*}
 In the first case, the sum of all these roots is equal to $0$, and roots are linearly dependent.
 In the second case, the sum of the $k-1$ first roots is equal to the last one, and
 roots are also linearly dependent. Thus, for $A_n$, there are no cycles of linearly independent roots.

 Note that $A_n$ do has cycles of length $3$ containing linearly independent roots, for example:
\begin{equation*}
    \{ e_{i_1} - e_{i_2},  e_{i_2} - e_{i_3}, e_{i_2} - e_{i_4}\}. \\
 \end{equation*}

  (ii) If a certain root subset contains $D_4$-associated root subset then it contains one of the following $4$-element subsets:
 \begin{equation*}
  \begin{split}
   & S_1 = \{ e_{i_1} - e_{i_2},  e_{i_2} - e_{i_3}, e_{i_3} - e_{i_4}, e_{i_2} - e_{i_5} \} \\
   & S_2 = \{ e_{i_1} - e_{i_2},  e_{i_2} - e_{i_3}, e_{i_3} - e_{i_4}, e_{i_3} - e_{i_5} \} \\
   \end{split}
 \end{equation*}
  In both cases $e_{i_2} - e_{i_3}$ is the branch point. In the first (resp. second) case
  $e_{i_2} - e_{i_5}$ (resp. $e_{i_3} - e_{i_5}$) is connected to
  $e_{i_1} - e_{i_2}$ and $e_{i_2} - e_{i_3}$ (resp. $e_{i_2} - e_{i_3}$ and $e_{i_3} - e_{i_4}$).
  Then $\{ e_{i_1} - e_{i_2},  e_{i_2} - e_{i_3},  e_{i_2} - e_{i_5}\}$ (resp. $\{ e_{i_2} - e_{i_3},  e_{i_3} - e_{i_4},  e_{i_3} - e_{i_5}\}$)
  forms the $3$-cycle and $S_1$ (resp. $S_2$) is not $D_4$.
  \qed

 \begin{lemma}
   \label{lem_E6_not_in_Dn}
  The root systems $E_6$, $E_7$, $E_8$ are not contained in the root system $D_n$.
 \end{lemma}

 \PerfProof
 Roots of the root systems $E_6$ and $D_n$ are as follows:
 \begin{equation}
 \begin{split}
   & D_n \quad
    \begin{cases}
      \pm{e_i}  \pm{e_j} \quad (1 \le  i  < j  \le l),
    \end{cases} \\
   \\
   & E_6 \quad
    \begin{cases}
    \begin{split}
     & \pm{e_i}  \pm{e_j} \quad (1 \le  i  < j  \le 5), \\
     \\
     & \pm{\frac{1}{2}} \left ( e_8 - e_7 - e_6 + \sum\limits_{i=1}^5(-1)^{\nu(i)}e_i \right ),
     \text{ where } \sum\limits_{i=1}^5\nu(i) \text{ is even, }
    \end{split}
   \end{cases} \\
  \end{split}
 \end{equation}
 see \cite[Tables IV and V]{Bo02}.
 It is clear that some roots of $E_6$ cannot be obtained as roots of $D_n$.
 Thus, $\varPhi(E_6) \not\subset \varPhi(D_n)$
 and, therefore, also $\varPhi(E_7) \not\subset \varPhi(D_n)$ and $\varPhi(E_8) \not\subset \varPhi(D_n)$.
 \qed

\begin{figure}[h]
\centering
\includegraphics[scale=0.62]{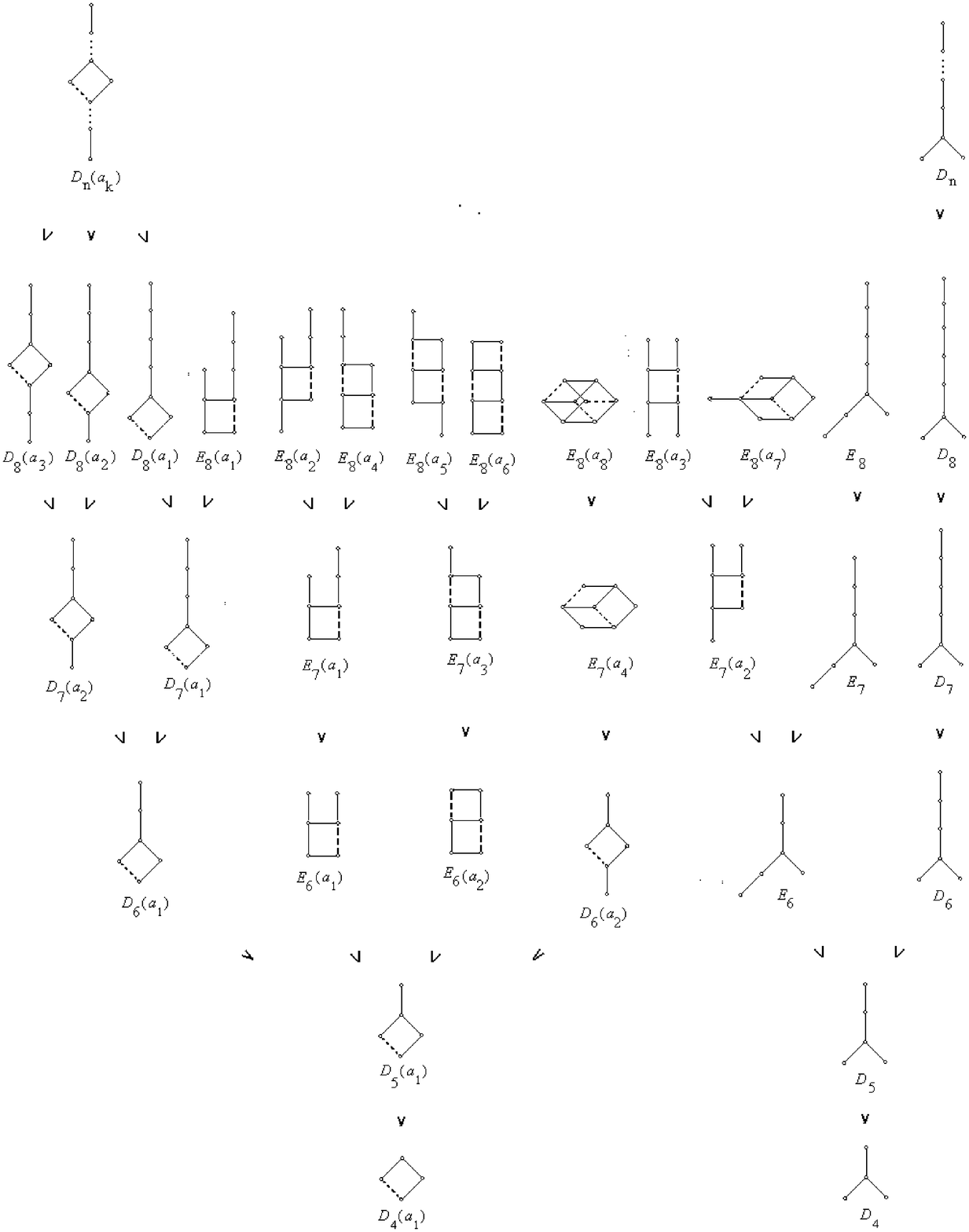}
\caption{\hspace{3mm}The ordered tree of Carter diagrams from
$\mathsf{C4} \coprod \mathsf{DE4}$}
%%%%%% The label must come after caption
\label{diagram_tree}
\end{figure}

 \begin{lemma}
   \label{lem_E6ak_not_in_Dn}
  The partial root systems
   \begin{equation}
    \label{eq_E678ak}
    \left \{
     \begin{array}{ll}
         \hspace{-1mm}E_6(a_k) & \text{for } k = 1,2, \\
         \hspace{-1mm}E_7(a_k) & \text{for } k = 1,2,3,4, \\
         \hspace{-1mm}E_8(a_k) & \text{for } 1 \leq k \leq 8
     \end{array}
    \right .
   \end{equation}
  are not contained in the root system $D_n$.
 \end{lemma}
  \PerfProof %% As in Lemma \ref{lem_E6_not_in_Dn},
  Since any $E_7(a_k)$ and $E_8(a_k)$ in \eqref{eq_E678ak} contains $E_6(a_1)$ or $E_6(a_2)$ as a subdiagram,
  see Fig. \ref{diagram_tree}, it suffices to prove that $E_6(a_1) \not\subset D_n$ and $E_6(a_2) \not\subset D_n$.
  %% the statement for $E_6(a_1)$ and $E_6(a_2)$.
  This follows from Corollary \ref{corol_biject}(i) and Lemma \ref{lem_E6_not_in_Dn}.
  %%% prop_bijection_root_syst
  \qed

\clearpage
\subsection{The ordered tree of Carter diagrams}
~\\
~\\

\subsection{$\Gamma$-associated root subsets and conjugacy classes}

\subsubsection{Two $\Gamma$-associated conjugacy classes}
   \label{sec_two_diff_calsses}
\begin{figure}[h]
\includegraphics[scale=1.8]{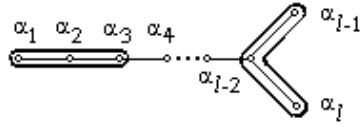}
 \caption{\hspace{3mm}Elements $s_{\alpha_1}s_{\alpha_3}s_{\alpha_2}$ and
     $s_{\alpha_{l-1}}s_{\alpha_l}s_{\alpha_{l-2}}$ are not conjugate}
%%%%%% The label must come after caption
\label{Dl_w1_w2}
\end{figure}

 There exist $\Gamma$-associated elements
 $w_1$ and $w_2$ such that $w_1 \not\simeq w_2$. For example, the Carter diagram $A_3$
 determines two different conjugacy classes in $D_l$, see Fig. \ref{Dl_w1_w2}; for details, see \cite[\S{B.2.2}]{St10}.

\index{conjugate sets} \index{non-conjugate $\Gamma$-associated
sets}

 \subsubsection{Two non-conjugate $\Gamma$-associated sets}
   Let $S_1 = \{\varphi_1,\dots,\varphi_n\}$ and $S_2 = \{\delta_1,\dots,\delta_n\}$
   be two $\Gamma$-associated sets of roots. The sets $S_1$ and $S_2$ are
   said to be {\it conjugate} if there exists an element $T \in W$
   such that $T: \varphi_i \longmapsto \delta_i$ for $i = 1,\dots,n$. In this case, we write
 \begin{equation*}
    S_1 \simeq S_2 \text{ and } T{S}_1 = S_2.
 \end{equation*}
   Let $w_1$ (resp. $w_2$) be any $S_1$-associated (resp. $S_2$-associated) element.
   If $S_1 \simeq S_2$, then $w_1 \simeq w_2$.
~\\

\begin{figure}[h]
 \centering
\includegraphics[scale=0.9]{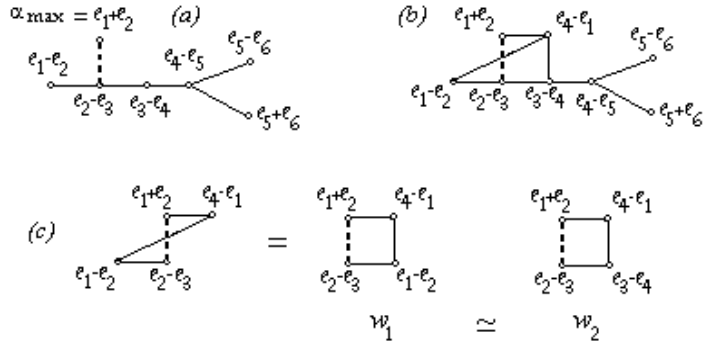}
\caption{\hspace{3mm}Equivalence of the $\mathcal{C}_1$-associated
element $w_1$ and
         the $\mathcal{C}_2$-associated element $w_2$}
%%%%%% The label must come after caption
\label{example_D6_2sq}
\end{figure}

   There exist, however, conjugate elements $w_1$ and $w_2$ such that $S_1 \not\simeq S_2$.
   Consider two $4$-cycles in $D_6$:
 \begin{equation*}
   \begin{split}
     & \mathcal{C}_1 = \{ e_1 + e_2,  e_4 - e_1, e_1 - e_2, e_2 - e_3  \}, \\
     & \mathcal{C}_2 = \{ e_1 + e_2,  e_4 - e_1, e_3 - e_4, e_2 - e_3  \}.
   \end{split}
 \end{equation*}
 These sets are non-conjugate: $\mathcal{C}_1 \not\simeq \mathcal{C}_2$, see Fig. \ref{example_D6_2sq}
 and \S\ref{sec_example_4cycles},
 but the $\mathcal{C}_1$-associated element $ w_1 = s_{e_1 + e_2}s_{e_1 - e_2}s_{e_4 - e_1}s_{e_2 - e_3}$
 and the $\mathcal{C}_2$-associated element $w_2 = s_{e_1 + e_2}s_{e_3 - e_4}s_{e_4 - e_1}s_{e_2 - e_3}$ are
 conjugate.
 %% : $w_1 \simeq w_2$.

\subsubsection{Example of equivalent $4$-cycles}
  \label{sec_example_4cycles}
 Let us take the diagram $D_6$ with simple roots:
 \begin{equation}
   \begin{split}
     & \alpha_1 = e_1 - e_2, \quad \alpha_2 = e_2 - e_3, \quad \alpha_3 = e_3 - e_4, \\
     & \alpha_4 = e_4 - e_5, \quad \alpha_5 = e_5 - e_6, \quad \alpha_2 = e_5 + e_6,
   \end{split}
 \end{equation}
 and the maximal root $\alpha = e_1 + e_2$, see Fig. \ref{example_D6_2sq}$(a)$.
 Consider the following two $4$-cycles:
 \begin{equation}
   \begin{split}
     & \mathcal{C}_1 = \{ e_1 + e_2,  e_4 - e_1, e_1 - e_2, e_2 - e_3  \}, \\
     & \mathcal{C}_2 = \{ e_1 + e_2,  e_4 - e_1, e_3 - e_4, e_2 - e_3  \}.
   \end{split}
 \end{equation}

 Since the $4$-index dipole can be mapped only onto
 a $4$-index dipole, see \cite[Lemma B.1]{St10}, and $\mathcal{C}_2$ consists of two $4$-index dipoles,
 it follows that $\mathcal{C}_1$ and $\mathcal{C}_2$ can not be conjugate.
 However, $\mathcal{C}_1$ and $\mathcal{C}_2$
 have the \underline{common second dipole}  $\{ e_4 - e_1, e_2 - e_3\}$.
 Moreover, let $w_1$ (resp. $w_2$) be $\mathcal{C}_1$-associated (resp. $\mathcal{C}_2$-associated):
 \begin{equation}
   \begin{split}
     & w_1 = s_{e_1 + e_2}s_{e_1 - e_2}s_{e_4 - e_1}s_{e_2 - e_3}, \\
     & w_2 = s_{e_1 + e_2}s_{e_3 - e_4}s_{e_4 - e_1}s_{e_2 - e_3},
   \end{split}
 \end{equation}
then elements $w_1$ and $w_2$ are conjugate:
 \begin{equation}
   \begin{split}
     w_1 = & s_{e_1 + e_2}s_{e_1 - e_2}s_{e_4 - e_1}s_{e_2 - e_3} =
             s_{e_1 + e_2}(s_{e_1 - e_2}s_{e_2 - e_3})s_{e_4 - e_1} = \\
           & s_{e_1 + e_2}s_{e_2 - e_3}s_{(e_2 - e_3) + (e_1 - e_2)}s_{e_4 - e_1} =
             s_{e_1 + e_2}s_{e_2 - e_3}(s_{e_1 - e_3}s_{e_4 - e_1}) = \\
           & s_{e_1 + e_2}s_{e_2 - e_3}s_{e_4 - e_1}s_{(e_4 - e_1) + (e_1 - e_3)} =
             s_{e_1 + e_2}s_{e_2 - e_3}s_{e_4 - e_1}s_{e_4 - e_3} \stackrel{s_{e_4 - e_3}}{\simeq}   \\
           & s_{e_1 + e_2}s_{e_3 - e_4}s_{e_4 - e_1}s_{e_2 - e_3} = w_2.  \qed
   \end{split}
 \end{equation}
\clearpage

%%\section{\sc\bf The partial Cartan matrix $B_{\Gamma}$ and the inverse matrix $B^{-1}_{\Gamma}$}
%% \label{sec_inv_matr}

\subsection{The diagonal elements of  $B^{-1}_{\Gamma}$ for $A_l$, $D_l$, $D_l(a_k)$}
  \index{$b^{\vee}_{\eta, \eta}$ (diagonal element of $B^{-1}_{\Gamma}$)}

\begin{figure}[h]
\centering
\includegraphics[scale=0.5]{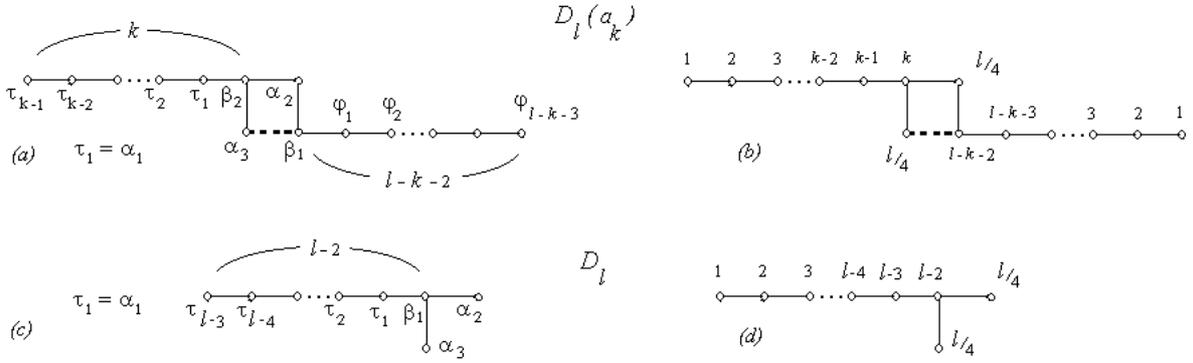}
\caption{\hspace{3mm}The numerical labels in the right hand side
 are the diagonal elements of $B^{-1}_{\Gamma}$}
%%%%%% The label must come after caption
\label{Dl_ak_Carter_diagr_in_prop}
\end{figure}

 Let $\eta$ be one of vertices of $\Gamma = D_l(a_k)$ (resp. $\Gamma = D_l$) such that $\eta \neq \alpha_2, \alpha_3$,
 see Fig. \ref{Dl_ak_Carter_diagr_in_prop},
 \begin{equation}
        \label{diag_elem_1}
    \Small %% \footnotesize
        \eta \in
        \begin{cases}
           \{ \tau_1, \dots, \tau_{k-1}, \varphi_1, \dots, \varphi_{l-k-3},
          \alpha_2, \alpha_3, \beta_1, \beta_2\}  & \text{ for } \Gamma = D_l(a_k),  \\
           \{ \tau_1, \dots, \tau_{l-3},
          \alpha_2, \alpha_3, \beta_1, \beta_2\}  & \text{ for } \Gamma = D_l,
        \end{cases}
 \end{equation}
 Removing the vertex $\eta$ with its bonds from $\Gamma$ we get the diagram $\Gamma'$ which is decomposed,
 except for $\eta = \tau_{k-1}$ and $\eta = \varphi_{l-k-3}$,
 into the  union of two connected subdiagrams:
  \begin{equation}
    \label{eq_one_of_diagr}
     \Small
     \Gamma' =
        \begin{cases}
           A_{d-1} \oplus D_{l-p-1},  & \text{ where } d = l-k-2 \text{~(resp. } d = k), \\
                                  & \text{ for } \eta = \beta_1 ~(\text{resp. } \eta = \beta_2), \\
           A_{d-1} \oplus D_{l-p-1(a_i)},  & \text{ where } d = k-i, \text{~(resp. } d = l-k-2-i), \\
                & \text{ for } \eta = \tau_i ~(1 \leq i \leq k-1) ~(\text{resp. }
                               \eta = \varphi_i ~(1 \leq i \leq l-k-3)).  \\
        \end{cases}
  \end{equation}
  For $\eta = \tau_{k-1}$ (resp. $\eta = \varphi_{l-k-3}$),
  we have $d = 1$ and $\Gamma' = D_{l-1}(a_{k-1})$ (resp. $\Gamma' = D_{l-1}(a_{l-k-3})$).
  In Fig. \ref{Dl_ak_Carter_diagr_in_prop}$(b)$ and $(d)$, we associate
  the numerical label $d$ with the corresponding vertex $\eta$.

 \begin{proposition}
   \label{prop_determinants}
    {\rm (i)} The determinant of the partial Cartan matrix $B_{\Gamma}$ is as follows:
  \begin{equation}
    \Small
    \label{eq_det_AD}
      \det{B_{\Gamma}} =
        \begin{cases}
           l + 1 & \text{for } A_l, \text{ where } l \geq 2; \text{ here } B_{\Gamma} = {\bf B}, \\
           4     & \text{for } D_l \text{ and }  D_l(a_k), \text{ where } l \geq 4. \\
        \end{cases}
  \end{equation}

    {\rm  (ii)} Let $b^{\vee}_{\eta,\eta}$ be diagonal elements of $B_{\Gamma}^{-1}$,
     where $\eta$ is given by \eqref{diag_elem_1}.
     For $\Gamma = D_l(a_k)$ or $\Gamma = D_l$,
      \begin{equation}
        \label{bi_eq_1}
  \Small
     b^{\vee}_{\eta,\eta} =
     \begin{cases}
         \displaystyle\frac{l}{4} \text{ for } \eta = \alpha_2, \text{ or } \eta = \alpha_3, \\
         d \text{ for } \eta \text{ given in }  \eqref{diag_elem_1},
     \end{cases}
     \end{equation}
 where $d$ is given in the vertex $\eta$
 in Fig. \ref{Dl_ak_Carter_diagr_in_prop}$(b)$ and $(d)$, and by the relation \eqref{eq_one_of_diagr}.

   {\rm  (iii)} For $\Gamma = A_l$, let $b^{\vee}_{\eta,\eta}$ be diagonal elements of $B_{\Gamma}^{-1}$, where
 \begin{equation}
  \Small
   \label{seq_vetices}
   \eta \in
  \begin{cases}
     \alpha_1, \beta_1, \alpha_2, \beta_2, \dots, \alpha_n, \beta_n & \text{ for } l = 2n,  \\
     \alpha_1, \beta_1, \alpha_2, \beta_2, \dots, \alpha_n, \beta_n, \alpha_{n+1} & \text{ for } l = 2n + 1,
  \end{cases}
 \end{equation}
 see Fig. $\ref{fig_A_l}$. Then
   \begin{equation}
     \label{eq_diag_elem_0}
      b^{\vee}_{\eta,\eta} =  d - \frac{d^2}{l+1},
   \end{equation}
 where $d$ is the sequential number of the vertex in eq. \eqref{seq_vetices} or in Fig. $\ref{fig_A_l}$.

\begin{figure}[h]
\centering
\includegraphics[scale=0.8]{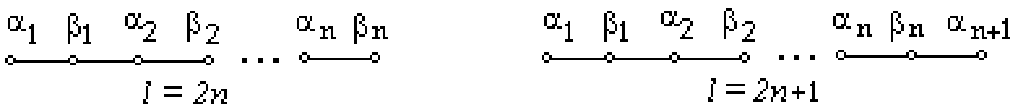}
\caption{\hspace{3mm}The Dynkin diagrams $A_l$}
%%%%%% The label must come after caption
\label{fig_A_l}
\end{figure}

\end{proposition}

\PerfProof  (i) This statement is easily verified for $A_1$, $A_2$,
$D_4$, $D_5$ and $D_4(a_1)$.
   By induction
 \begin{equation*}
  \Small
   \begin{split}
    & \det{{\bf B}(A_{l+1})} = 2\det{{\bf B}(A_{l})} -  \det{{\bf B}(A_{l-1})} =
        2(l+1) - l = l + 2, \\
    & \det{{\bf B}(D_{l+1})} = 2\det{{\bf B}(D_{l})} -  \det{{\bf B}(D_{l-1})} = 4,
   \end{split}
 \end{equation*}
 where  ${\bf B}(A_{l})$ (resp. ${\bf B}(D_{l})$)  is the Cartan matrix for $A_{l}$ (resp. $D_{l}$).

 To get the determinant of the partial Cartan matrix $B_{\Gamma}$ for the diagram $\Gamma = D_{l+1}(a_k)$,
 we expand $\det B_{\Gamma}(D_{l+1}(a_k))$ with respect
 to the minors corresponding to the $i$-th line (associated with the vertex $i$ of $\Gamma$).
 By induction, we have
 \begin{equation*}
  \Small
   \begin{array}{lll}
    & \det B_{\Gamma}(D_{l+1}(a_k)) = 2\det B_{\Gamma}(D_{l}(a_{k-1})) -  \det B_{\Gamma}(D_{l-1}(a_{k-2})) = 4
       & \text{ for } k > 2, \\
    & \det B_{\Gamma}(D_{l+1}(a_{l - k - 1})) = 2\det B_{\Gamma}(D_{l}(a_{l - k - 2})) -
                                                   \det B_{\Gamma}(D_{l-1}(a_{l - k - 3})) = 4
       & \text{ for } l - k  > 3, \\
    & \det B_{\Gamma}(D_{l+1}(a_2)) = 2\det B_{\Gamma}(D_{l}(a_1)) -  \det {\bf B}(D_{l-1}) = 4
       & \text{ for } k = 2 \text{ or } l - k = 3.
    \end{array}
 \end{equation*}
 In the case $k < 2$ and $l - k  < 3$, we have $l < 5$, i.e., $\Gamma = D_4(a_1)$.
~\\

  (ii) For $\Gamma = D_l(a_k)$, we have
   \begin{equation}
  \Small
    b^{\vee}_{\eta,\eta} =
    \begin{cases}
    ~ \displaystyle\frac{\det {B_{\Gamma}(D_{l}(a_k))}}{\det {B_{\Gamma}(D_{l+1}(a_k))}} = \frac{4}{4} = 1
        & \text{ for } $d = 1$, \\
    & \\
    ~ \displaystyle\frac{\det{\bf B}(A_{d-1}) \det B_{\Gamma}(D_{l-d}(a_k))}
                {\det B_{\Gamma}(D_{l+1}(a_k))} = \frac{d \cdot 4}{4} = d
        & \text{ for } d > 1, \\
    & \\
    ~  \displaystyle\frac{\det {{\bf B}(A_{l-1})}}{\det {B_{\Gamma}(D_{l+1}(a_k))}} = \frac{l}{4}
       & \text{ for } \eta = \alpha_2, \alpha_3.
     \end{cases}
   \end{equation}
   For $\Gamma = D_l$, we have
   \begin{equation}
  \Small
    b^{\vee}_{\eta,\eta} =
    \begin{cases}
    ~ \displaystyle\frac{\det {{\bf B}(D_{l})}}{\det {{\bf B}(D_{l+1})}} = \frac{4}{4} = 1
        & \text{ for } $d = 1$, \\
    & \\
    ~ \displaystyle\frac{\det{\bf B}(A_{d-1}) \det {\bf B}(D_{l-d})}
                {\det {\bf B}(D_{l+1})} = \frac{d \cdot 4}{4} = d
        & \text{ for } d > 1, \\
    & \\
    ~  \displaystyle\frac{\det {{\bf B}(A_{l-1})}}{\det {{\bf B}(D_{l+1})}} = \frac{l}{4}
        & \text{ for } \eta = \alpha_2, \alpha_3.
     \end{cases}
   \end{equation}
~\\

  (iii) For $\Gamma = A_l$, we have:
  \begin{equation}
   \Small
        b^{\vee}_{\eta,\eta} = \frac{\det \overline{B}_{\eta,\eta}}{\det {\bf B}(A_l)},
  \end{equation}
  where  $\overline{B}_{\eta,\eta}$ is the matrix obtained from $B_{\Gamma}$ by deleting the $\eta$th column and $\eta$th row.
  The matrix $\overline{B}_{\eta,\eta}$ splits into the direct sum:
  \begin{equation*}
        \overline{B}_{\eta,\eta} = {\bf B}(A_{d-1}) \oplus {\bf B}(A_{l - d}).
  \end{equation*}
   Hence,
  \begin{equation}
    \label{eq_diag_elem_1}
        b^{\vee}_{\eta,\eta} = \frac{\det {\bf B}(A_{d-1}) \times \det {\bf B}(A_{l-d})}{\det {\bf B}(A_l)}.
  \end{equation}
  By \eqref{eq_det_AD} $\det {\bf B}(A_{l-1}) = l$ and by \eqref{eq_diag_elem_1} we have
  \begin{equation*}
        b^{\vee}_{\eta,\eta} = \frac{d \times (l + 1 -d)}{l + 1} =  d - \frac{d^2}{l+1}. \qed
  \end{equation*}
  ~\\

   For expressions ${\bf B^{-1}}$ for the Dynkin diagram $A_l$, see \cite[p. 295]{OV90}.
   We give the partial  Cartan matrices and its inverse for
   Carter diagrams $E_l(a_k)$, (where $l < 7$), for $A_l$ (where $l \leq 8$),
   for $D_l(a_k)$ (where $l \leq 8$) in
   Tables \ref{tab_partial Cartan_1}, \ref{tab_partial Cartan_2}, \ref{tab_partial Cartan_2a}, \ref{tab_partial Cartan_3},
   \ref{tab_partial Cartan_Al}.

  \index{$B_{\Gamma}$, partial Cartan matrix}
  \index{Cartan matrix ! - partial Cartan matrix $B_{\Gamma}$}

\subsection{Simply extendable Carter diagrams}
  \label{sect_extendable}
  \index{simply linked vertex}
  \index{simply linkable vertex}
  \index{simply extendable Carter diagram}
  \index{simply extended diagram}

\begin{figure}[H]
\centering
\includegraphics[scale=0.7]{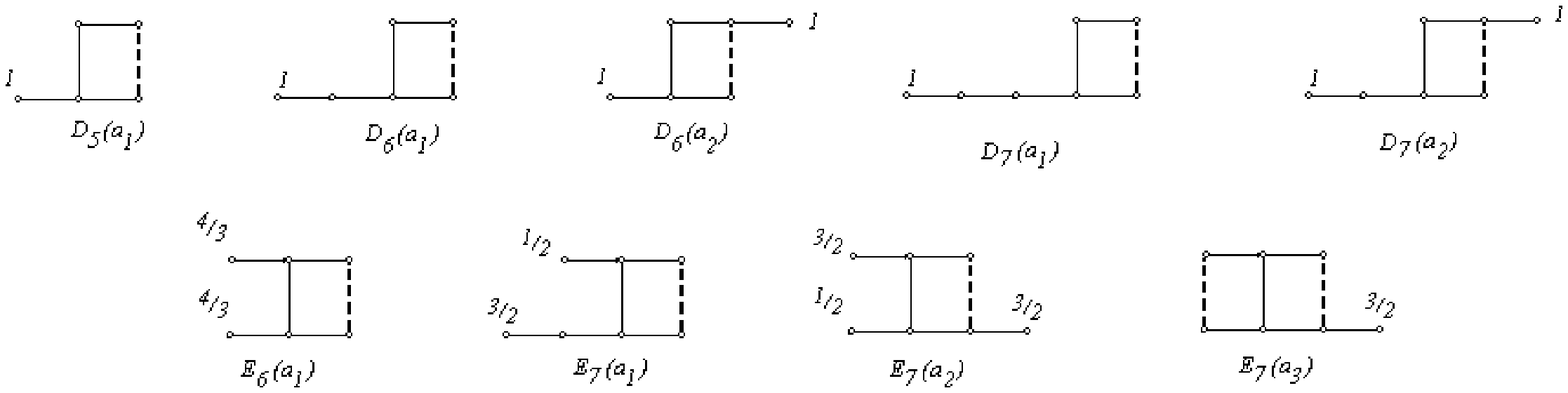}
\caption{
\hspace{3mm}Simply linkable vertices are marked by numerical values, which are the diagonal elements
of $B^{-1}_{\Gamma}$, see Tables \ref{tab_partial Cartan_1} - \ref{tab_partial Cartan_3}.}
%%%%%% The label must come after caption
\label{examples_ext}
\end{figure}

  We say that the Carter diagram $\Gamma$ is {\it simply extendable in the vertex $\tau_p$} if
  the new diagram $\widetilde\Gamma$ is obtained by adding an extra vertex $\tau_{l+1}$ and only one connection edge
  $\{\tau_p, \tau_{l+1}\}$, and $\widetilde\Gamma$ is also the Carter diagram.
  In this case, the extra vertex $\tau_{l+1}$ is called {\it simply linked},
  the vertex $\tau_p$ is called {\it simply linkable},
  and $\widetilde\Gamma$ is called {\it simply extended}, see Fig. \ref{examples_ext}.
  We will show that extensibility in the vertex is closely
  associated with the value of the diagonal element $b^{\vee}_{\tau_p,\tau_p}$  of the matrix $B_{\Gamma}^{-1}$,
  the inverse of the partial Cartan matrix $B_{\Gamma}$.

 \index{$B_{\Gamma}$, partial Cartan matrix}
 \index{Cartan matrix ! - partial Cartan matrix $B_{\Gamma}$}
 \index{$B^{-1}_{\Gamma}$, inverse of the partial Cartan matrix}
 \index{Cartan matrix ! - inverse $B^{-1}_{\Gamma}$ of the partial Cartan matrix}
 \begin{proposition}
  \label{prop_crit_sim_ext}
   The Carter diagram $\Gamma$ is simply extendable in the vertex $\tau_p$ if and only if
  \begin{equation}
    b^{\vee}_{\tau_p,\tau_p} < 2,
  \end{equation}
   where $b^{\vee}_{\tau_p,\tau_p}$ is the diagonal element (corresponding to the vertex $\tau_p$)
   of the matrix $B_{\Gamma}^{-1}$.
 \end{proposition}

 \PerfProof This is a direct consequence of Theorem \ref{th_B_less_2}. Indeed, the linkage label vector
 $\gamma^{\nabla}$ corresponding to the simply linked vertex $\tau_p$ is the vector with $1$ or $-1$ in the place $\tau_p$
 and zeros in remaining places. Then $\mathscr{B}^{\vee}_{\Gamma}(\gamma^{\nabla}) = b^{\vee}_{\tau_p,\tau_p}$. \qed

\begin{remark}
 \label{rem_val_le_2}
 {\rm
  For $D_l(a_k)$ (resp. $D_l$), where $l \geq 8$, we have $b^{\vee}_{\eta,\eta} < 2$ only for the endpoints $\tau_{k-1}$ and
  $\varphi_{l - k -3}$ (resp. only for the endpoint $\tau_{l-3}$), see Fig. \ref{Dl_ak_Carter_diagr_in_prop}$(a)$,$(b)$
  (resp. Fig. \ref{Dl_ak_Carter_diagr_in_prop}$(c)$,$(d)$).
  For $D_4(a_k)$, $D_5(a_k)$, $D_6(a_k)$, $D_7(a_k)$ and $D_4$, $D_5$, $D_6$, $D_7$,
  we also have $b^{\vee}_{\eta,\eta} < 2$ for $\eta = \alpha_2, \alpha_3$.
  }
\end{remark}

\subsubsection{Simple extensions for the Carter diagram $A_l$}

\begin{proposition}
  \label{prop_simple_ex}

   {\rm(i)} For $l \leq 4$  or $l \geq 8$, the simply linkable vertices $\tau$ are only vertices with
   $d \in \{1, 2, l-1, l\}$, see Table $\ref{tab_extend_diagr}$.

   {\rm(ii)} For $A_5$, $A_6$ and $A_7$, except for $d  \in \{1, 2, l-1, l\}$,
   there are some additional simply linkable vertices $\tau$, see Table $\ref{tab_extendable_567}$
   and Table $\ref{tab_partial Cartan_Al}$.
\end{proposition}

 \begin{table}[H]
  \centering
  \renewcommand{\arraystretch}{1.7}
  %% {\rm
  \begin{tabular} {|c|c|c|c|}
  \hline
  $d$   & Simply linkable vertex $\tau$ &  Simply extended diagram    \cr
     \hline
     1      & $\alpha_1$    & $\begin{array}{c} \includegraphics[scale=0.6]{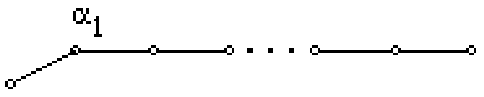}  \quad A_{l+1} \end{array}$  \\
     \hline
     2      & $\beta_1$     & $\begin{array}{c} \includegraphics[scale=0.6]{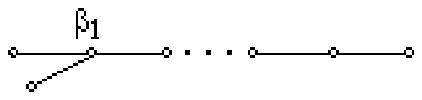} \qquad  D_{l+1} \end{array}$    \\
     \hline
     $l-1$  &
             $\tau_{l-1} = \begin{cases} \alpha_n \text{ for } l = 2n \\
                                          \beta_n \text{ for } l = 2n+1 \end{cases}$     &
                                             $\begin{array}{c} \includegraphics[scale=0.6]{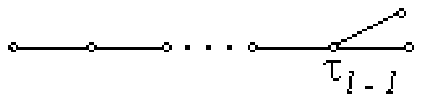} \qquad
                                                  A_{l+1}\end{array}$   \\
     \hline
     $l$    & $\tau_l =
                 \begin{cases}
                     \beta_n  \quad \text{ for } l = 2n \\
                     \alpha_{n+1} \text{ for } l = 2n+1
                 \end{cases} $  &
                 $\begin{array}{c} \includegraphics[scale=0.6]{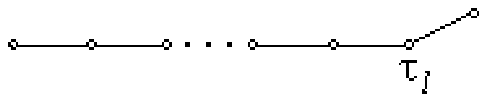} \quad  D_{l+1}\end{array}$     \\
                 %% don't confuse Al_ext_1 (odin) Al_ext_l (el)
     \hline
   \end{tabular}
    %% }
  \vspace{2mm}
  \caption{\small\hspace{3mm} The simply extendable diagrams for $A_l$, where $l \leq 4$  or $l \geq 8$}
  \label{tab_extend_diagr}
\end{table}

\begin{table}[H]
  \centering
  \renewcommand{\arraystretch}{1.7}
  %% {\rm
  \begin{tabular} {|c|c|c|c|}
   \hline
     $\Gamma$  & Simply linkable vertex $\tau$  & $b^{\vee}_{\eta,\eta}$  &  Simply extended diagram \\
   \hline
    $A_5$     &  $\alpha_2$          & $b^{\vee}_{\alpha_2,\alpha_2}  = \frac{3}{2}$
      & $\begin{array}{c} \includegraphics[scale=0.6]{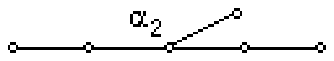}  \qquad \qquad E_6 \end{array}$  \\
   \hline
    $A_6$     &  $\alpha_2$, $\beta_2$   & $b^{\vee}_{\alpha_2,\alpha_2}  = b^{\vee}_{\beta_2,\beta_2} = \frac{12}{7}$
      & $\begin{array}{c}
          \begin{array}{c}
            \includegraphics[scale=0.6]{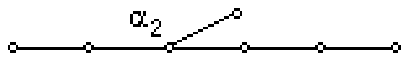} \\
            \includegraphics[scale=0.6]{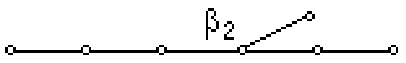}
           \end{array}
              \qquad E_7 \end{array}$  \\
   \hline
    $A_7$     &  $\alpha_2$, $\alpha_3$  & $b^{\vee}_{\alpha_2,\alpha_2}  = b^{\vee}_{\alpha_3,\alpha_3} = \frac{15}{8}$
       & $\begin{array}{c}
          \begin{array}{c}
            \includegraphics[scale=0.6]{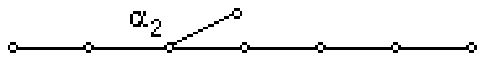} \\
            \includegraphics[scale=0.6]{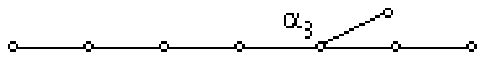}
           \end{array}
            \quad E_8 \end{array}$  \\
   \hline
  \end{tabular}
   %% }
  \vspace{2mm}
  \caption{\small\hspace{3mm} Additional simply extendable diagrams for $A_5$, $A_6$, $A_7$}
  \label{tab_extendable_567}
\end{table}

  \PerfProof
  (i) The statement holds for $l \leq 4$ since in this case any $d$ is from $\{1, 2, l-1, l\}$,
    and the simple extensions in these cases are as follows:
      $A_3$ for $l = 2$; $A_4$ or $D_4$ for $l = 3$; $A_5$ or $D_5$ for $l = 4$.

     Let $l \geq 8$. By Proposition \ref{prop_crit_sim_ext}, the Carter diagram $\Gamma$ is
     simply extendable in the vertex $\tau_d$ if and only if $b^{\vee}_{\eta,\eta} < 2$.
     By eq. \eqref{eq_diag_elem_0} from Proposition \ref{prop_determinants},  we have to solve the inequality
   \begin{equation*}
     \label{eq_diag_elem_2}
     \Small
          d - \frac{d^2}{l+1} < 2, \text{ i.e., } \quad   d^2 - d(l+1) +2(l+1)^2 > 0,  \\
   \end{equation*}
   so,
   \begin{equation*}
     \Small
          d  < d_- \text{ or } d > d_+, \text{ where }
           d_{\pm} = \frac{(l + 1) \pm \sqrt{(l+1)^2 - 8(l+1)}}{2}.
   \end{equation*}
Since $l \geq 8$, we have $l - 5  ~\leq~  \sqrt{(l+1)^2 - 8(l+1)}$. Then
   \begin{equation}
     \label{eq_diag_elem_3}
     \Small
     \begin{split}
          &  l-2 ~\leq~ \frac{l+1 + \sqrt{(l+1)^2 - 8(l+1)}}{2} = d_{+} < d, \\
          & 3 ~\geq~ \frac{l+1 - \sqrt{(l+1)^2 - 8(l+1)}}{2} = d_{-} > d,
     \end{split}
   \end{equation}
 By, \eqref{eq_diag_elem_3}
   \begin{equation*}
       \Small
          d < d_- \leq 3 \text{ and } d > d_+ \geq l - 2.
   \end{equation*}
   In other words,
   \begin{equation*}
      \Small
          d \in \{ 1, 2, l - 1, l \},
   \end{equation*}
   see Table \ref{tab_extend_diagr}.
  ~\\

  (ii) This statement is immediately verified by Table \ref{tab_partial Cartan_Al}. \qed

\clearpage
\subsection{The partial Cartan matrix $B_{\Gamma}$ and the inverse matrix $B^{-1}_{\Gamma}$}
  \label{sec_inv_matr}
~\\
 \begin{table}[H]
\tiny
  \centering
  \renewcommand{\arraystretch}{1.5}
  \begin{tabular} {|c|c|c|}
  \hline
      The Carter   & The partial Cartan         &   The inverse          \cr
        diagram    &  matrix $B_{\Gamma}$              &   matrix $B^{-1}_{\Gamma}$    \\
  \hline  %1
     & & \\
     $\begin{array}{c} \includegraphics[scale=0.6]{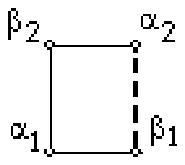} \\
        {\bf D_4(a_1)} \end{array}$  &
     $\begin{array}{c}
       \left [
   \begin{array}{cccc}
     2 & 0 &   -1 & -1  \\
     0 & 2 &   1 & -1 \\
     -1 & 1 &  2 & 0  \\
     -1 & -1 &   0 & 2  \\
  \end{array}
  \right ]
     \end{array}$   &
    $\begin{array}{c}
    \frac{1}{2} \left [ \begin{array}{cccc}
     2 & 0 &   1 & 1  \\
     0 & 2 &   -1 & 1 \\
     1 & -1 &  2 & 0  \\
     1 &  1 &   0 & 2  \\
  \end{array} \right ]
     \end{array}$    \\
    & & \\
    \hline
    & & \\
     $\begin{array}{c} \includegraphics[scale=0.6]{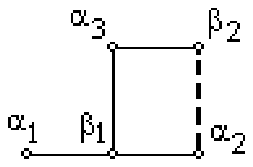} \\
        {\bf D_5(a_1) = D_5(a_2)} \end{array}$ &
     $\begin{array}{c}
       \left [
   \begin{array}{ccccc}
     2 & 0 &   0 & -1 & 0  \\
     0 & 2 &   0 & -1 & 1\\
     0 & 0 &   2 & -1 & -1  \\
     -1 & -1 &  -1 & 2 & 0  \\
      0 &  1 &   -1  & 0 & 2 \\
  \end{array}
  \right ]
     \end{array}$   &
      $\begin{array}{c}
    \frac{1}{4} \left [ \begin{array}{ccccc}
     4 & 2 &   2 & 4 & 0  \\
     2 & 5 &   1 & 4 & -2\\
     2 & 1 &   5  & 4 & 2  \\
     4 & 4 &   4 & 8 & 0  \\
     0 & -2 &  2  & 0 & 4 \\
  \end{array} \right ]
     \end{array}$    \\
    & & \\
    \hline
    & & \\
     $\begin{array}{c} \includegraphics[scale=0.6]{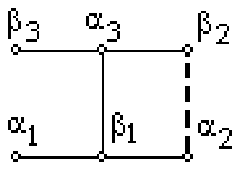} \\
        {\bf E_6(a_1)} \end{array}$ &
     $\begin{array}{c}
      \left [
   \begin{array}{cccccc}
     2 & 0 &   0 & -1 & 0  & 0 \\
     0 & 2 &   0 & -1 & 1  & 0\\
     0 & 0 &   2 & -1 & -1 & -1 \\
     -1 & -1 &  -1  & 2 & 0 & 0 \\
      0 &  1 &  -1  & 0 & 2 & 0 \\
      0 &  0 &  -1   & 0 & 0 & 2 \\
  \end{array}
  \right ]
     \end{array}$   &
     $\begin{array}{c}
      \frac{1}{3}\left [
   \begin{array}{cccccc}
     4 & 2 & 4 &  5 &  1 & 2 \\
     2 & 4 & 2 &  4 & -1 & 1 \\
     4 & 2 & 10 &  8 & 4 & 5 \\
     5 & 4 & 8 &  10 & 2 & 4 \\
     1 & -1 & 4 & 2 &  4 & 2 \\
     2 & 1 & 5 &  4 &  2 & 4 \\
  \end{array}
  \right ]
     \end{array}$   \\
    & & \\
    \hline
    & & \\
     $\begin{array}{c} \includegraphics[scale=0.6]{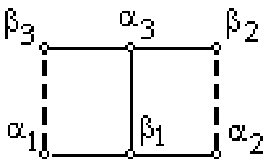} \\
        {\bf E_6(a_2)} \end{array}$ &
      $\begin{array}{c}
      \left [
   \begin{array}{cccccc}
     2 & 0 &   0 & -1 & 0  & 1 \\
     0 & 2 &   0 & -1 & 1  & 0\\
     0 & 0 &   2 & -1 & -1 & -1 \\
     -1 & -1 &  -1  & 2 & 0 & 0 \\
      0 &  1 &  -1  & 0 & 2 & 0 \\
      1 &  0 &  -1   & 0 & 0 & 2 \\
  \end{array}
  \right ]
     \end{array}$   &
           $\begin{array}{c}
      \frac{1}{3}\left [
   \begin{array}{cccccc}
     4 & 2 & 0 &  3 & -1 & -2 \\
     2 & 4 & 0 &  3 & -2 & -1 \\
     0 & 0 & 6 &  3 & 3 &  3 \\
     3 &  3  &  3 &  6 & 0 & 0 \\
     -1 & -2 &  3 &  0 & 4 & 2 \\
     -2 & -1 &  3 &  0 & 2 & 4 \\
  \end{array}
  \right ]
     \end{array}$  \\
    & & \\
    \hline  %1
     & & \\
     $\begin{array}{c} \includegraphics[scale=0.6]{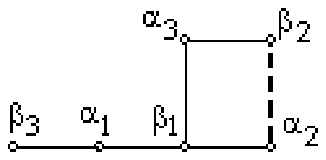} \\
        {\bf D_6(a_1) = D_6(a_3)} \end{array}$  &
     $\begin{array}{c}
       \left [
   \begin{array}{cccccc}
     2 & 0 &   0 & -1 & 0  & -1 \\
     0 & 2 &   0 & -1 & 1  & 0\\
     0 & 0 &   2 & -1 & -1 & 0 \\
     -1 & -1 &  -1  & 2 & 0 & 0 \\
      0 &  1 &  -1  & 0 & 2 & 0 \\
     -1 &  0 &  0   & 0 & 0 & 2 \\
  \end{array}
  \right ]
     \end{array}$   &
    $\begin{array}{c}
    \frac{1}{2} \left [
       \begin{array}{cccccc}
     4 & 2 &  2 &   4  & 0 & 2 \\
     2 & 3 &  1 &   3 & -1 & 1  \\
     2 & 1 &  3 &   3 &  1 & 1 \\
     4 & 3 &  3 &   6 &  0 & 2 \\
     0 &  -1 & 1 &  0 &  2 & 0 \\
     2 &  1 &  1 &  2 &  0 & 2 \\
  \end{array}
 \right ]
     \end{array}$    \\
    & & \\
     \hline
     & & \\
     $\begin{array}{c} \includegraphics[scale=0.6]{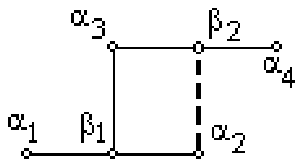} \\
        {\bf D_6(a_2)} \end{array}$  &
     $\begin{array}{c}
       \left [
   \begin{array}{cccccc}
      2 & 0  &   0  &  0  & -1  & 0 \\
      0 & 2  &   0  &  0  & -1  & 1 \\
      0 & 0  &   2  &  0  & -1  & -1 \\
      0 & 0  &   0  &  2  &  0  & -1 \\
     -1 & -1 &  -1  &  0  &  2  & 0 \\
      0 &  1 &  -1  & -1  &  0  & 2 \\
  \end{array}
  \right ]
     \end{array}$   &
    $\begin{array}{c}
    \frac{1}{2} \left [
       \begin{array}{cccccc}
     2 & 1 &  1  &  0  & 2 & 0 \\
     1 & 3 &  0  &  -1  & 2 & -2  \\
     1 & 0 &  3 &   1 &  2  &  2 \\
     0 &  -1 &  1 &  2 & 0  & 2 \\
     2 &  2 &  2 &  0 & 4 & 0 \\
     0 & -2 &  2 &  2 & 0 & 4 \\
  \end{array}
 \right ]
     \end{array}$    \\
    & & \\
     \hline
\end{tabular}
  \vspace{2mm}
  \caption{\small\hspace{3mm} The partial Cartan matrix $B_{\Gamma}$ and the inverse matrix $B^{-1}_{\Gamma}$ for Carter
   diagrams with the number of vertices $l < 7$}
  \label{tab_partial Cartan_1}
\end{table}

 \begin{table}[H]
\tiny
  \centering
  \renewcommand{\arraystretch}{1.4}
  \begin{tabular} {|c|c|c|}
  \hline
      The Carter   & The partial Cartan         &   The inverse          \cr
        diagram    &  matrix $B_{\Gamma}$              &   matrix $B^{-1}_{\Gamma}$    \\
     \hline
     & & \\
          $\begin{array}{c} \includegraphics[scale=0.6]{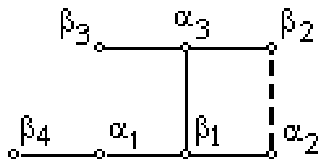} \\
        {\bf E_7(a_1)} \end{array}$  &
 $\begin{array}{c}
      \left [
   \begin{array}{ccccccc}
      2 & 0 &   0 & -1 & 0  & 0 & -1 \\
     0 & 2 &   0 & -1 & 1   & 0 & 0 \\
     0 & 0 &   2 & -1 & -1  & -1 & 0 \\
     -1 & -1 &  -1  & 2 & 0 & 0 & 0 \\
      0 &  1 &  -1  & 0 & 2 & 0 & 0 \\
      0 &  0 &  -1  & 0 & 0 & 2 & 0 \\
     -1 &  0 &  0   & 0 & 0 & 0 & 2 \\
  \end{array}
  \right ]
     \end{array}$   &
   $\begin{array}{c}
      \frac{1}{2}\left [
   \begin{array}{ccccccc}
      8 & 4 &  8 &  10 & 2 & 4 & 4 \\
      4 & 4 &  4 &  6 &  0 & 2 & 2 \\
      8 & 4 & 12 &  12 & 4 & 6 &  4 \\
      10 & 6 & 12 & 15 & 3 & 6  & 5 \\
      2 &  0 & 4 &  3 & 3 & 2  & 1 \\
      4 & 2 &  6 &  6 & 2 & 4  & 2 \\
      4 & 2 &  4 &  5 & 1 & 2  & 3 \\
  \end{array}
  \right ]
     \end{array}$   \\
     & & \\
     \hline
    & & \\
     $\begin{array}{c} \includegraphics[scale=0.6]{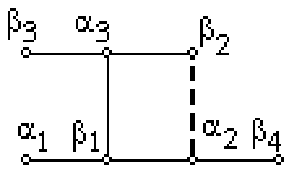} \\
        {\bf E_7(a_2)}
     \end{array}$ &
     $\begin{array}{c}
      \left [
   \begin{array}{ccccccc}
     2 & 0 &   0 & -1 &  0  & 0 & 0 \\
     0 & 2 &   0 & -1 &  1  & 0 & -1 \\
     0 & 0 &   2 & -1 & -1  & -1 & 0 \\
     -1 & -1 &  -1  & 2 & 0 & 0 & 0 \\
      0 &  1 &  -1  & 0 & 2 & 0 & 0 \\
      0 &  0 &  -1  & 0 & 0 & 2 & 0 \\
      0 & -1 &   0  & 0 & 0 & 0 & 2 \\
  \end{array}
  \right ]
     \end{array}$   &
     $\begin{array}{c}
      \frac{1}{2}\left [
   \begin{array}{ccccccc}
      4 & 4 &  4 &  6 & 0 & 2 & 2 \\
      4 & 8 &  4 &  8 & -2 & 2 & 4 \\
      4 & 4 &  8 &  8 & 2 & 4 &  2 \\
      6 & 8 & 8 &  12 & 0 & 4  & 4 \\
      0 & -2 & 2 &  0 & 3 & 1  & -1 \\
      2 &  2 & 4 &  4 & 1 & 3  & 1 \\
      2 &  4 & 2 &  4 & -1 & 1  & 3 \\
  \end{array}
  \right ]
     \end{array}$   \\
     & & \\
     \hline
    & & \\
     $\begin{array}{c} \includegraphics[scale=0.6]{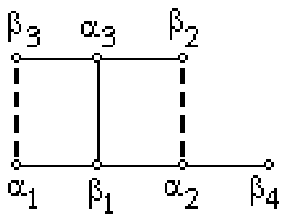} \\
       {\bf E_7(a_3)} \end{array}$ &
 $\begin{array}{c}
      \left [
   \begin{array}{ccccccc}
      2 & 0 & 0 &   -1 & 0 & 1 & 0 \\
     0 & 2 & 0 &    -1 & 1 & 0 & -1 \\
     0 & 0 & 2 &    -1 & -1 & -1 & 0 \\
     -1 & -1 & -1 &   2 & 0 & 0 & 0 \\
     0 &  1 & -1 &    0 & 2 & 0 & 0 \\
     1 &  0 &  -1 &    0 & 0 & 2 & 0 \\
     0 & -1 & 0 &     0 & 0 & 0 & 2 \\
  \end{array}
  \right ]
   \end{array}$  &
   $\begin{array}{c}
      \frac{1}{2}\left [
   \begin{array}{ccccccc}
      4 & 4 &  0 &  4 & -2 & -2 & 2 \\
      4 & 8 &  0 &  6 & -4 & -2 & 4 \\
      0 & 0 &  4 &  2 & 2 &  2 &  0 \\
      4 & 6 & 2 &  7 & -2 & -1  & 3 \\
      -2 & -4 & 2 &  -2 & 4 & 2  & -2 \\
      -2 & -2 & 2 &  -1 & 2 & 3  & -1 \\
      2 & 4 &  0 &  3 & -2 & -1  & 3 \\
  \end{array}
  \right ]
     \end{array}$   \\
     & & \\
     \hline
    & & \\
     $\begin{array}{c} \includegraphics[scale=0.6]{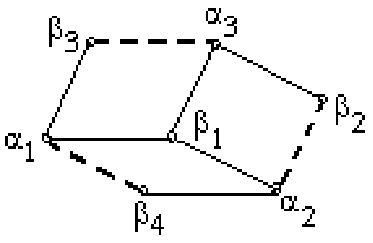} \\
        {\bf E_7(a_4)}  \end{array}$ &
      $\begin{array}{c}
      \left [
   \begin{array}{ccccccc}
      2 & 0 & 0 &   -1 & 0 & -1 & 1 \\
     0 & 2 & 0 &     -1 & 1 & 0 & -1 \\
     0 & 0 & 2 &     -1 & -1 & 1 & 0 \\
     -1 & -1 & -1 &   2 & 0 & 0 & 0 \\
     0 &  1 & -1 &    0 & 2 & 0 & 0 \\
     -1 & 0 &  1 &    0 & 0 & 2 & 0 \\
     1 & -1 & 0 &    0 & 0 & 0 & 2 \\
  \end{array}
  \right ]
     \end{array}$   &
      $\begin{array}{c}
      \frac{1}{2}\left [
   \begin{array}{ccccccc}
      4 & 0 &  0 &  2 & 0 & 2 & -2 \\
      0 & 4 &  0 &  2 & -2 & 0 & 2 \\
      0 & 0 &  4 &  2 & 2 &  -2 &  0 \\
      2 & 2 & 2 &  4 & 0 & 0  & 0 \\
      0 & -2 & 2 &  0 & 3 & -1  & -1 \\
      2 & 0 &  -2 & 0 & -1 & 3  & -1 \\
      -2 & 2 & 0 &  0 & -1 & -1  & 3 \\
  \end{array}
  \right ]
     \end{array}$   \\
     & & \\
     \hline
    & & \\
     $\begin{array}{c} \includegraphics[scale=0.6]{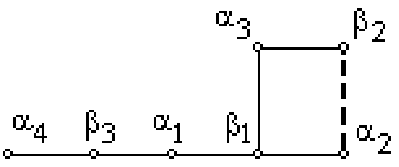} \\
        {\bf D_7(a_1) = D_7(a_4)} \end{array}$  &
     $\begin{array}{c}
       \left [
   \begin{array}{ccccccc}
      2  & 0 &  0 &  0 & -1 &  0  & -1 \\
      0  & 2 &  0 &  0 & -1 &  1  & 0 \\
      0  & 0 &  2 &  0 & -1 & -1  & 0 \\
      0  & 0 &  0 &  2 &  0 &  0  & -1\\
     -1 & -1 & -1 &  0 &  2 &  0  & 0 \\
      0 &  1 & -1 &  0 &  0 &  2  &  0 \\
     -1 &  0 &  0 & -1 &  0 &  0  &  2 \\
  \end{array}
  \right ]
     \end{array}$   &
    $\begin{array}{c}
    \frac{1}{4} \left [
       \begin{array}{ccccccc}
      12  & 6 &  6 &  4  & 12 & 0 & 8 \\
      6   & 7 &  3 &   2 & 8 &  -2 & 4 \\
      6   & 3 &  7 &  2 &  8 &  2 & 4 \\
      4   & 2 &  2 &  4 &  4 &  0  & 4\\
      12  &  8 &  8 &  4 &  16 &  0  & 8 \\
      0  &  -2 &  2 &  0 &  0 &  4 & 0 \\
      8   &  4 &  4 &  4 &  8 &  0 &  8 \\
  \end{array}
 \right ]
     \end{array}$    \\
     & & \\
     \hline
    & & \\
     $\begin{array}{c} \includegraphics[scale=0.6]{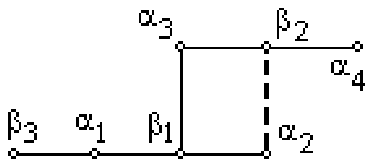} \\
        {\bf D_7(a_2) = D_7(a_3)} \end{array}$  &
     $\begin{array}{c}
       \left [
  \begin{array}{ccccccc}
      2  & 0 &  0 &  0 & -1 &  0 & -1 \\
      0  & 2 &  0 &  0 & -1 &  1 & 0 \\
      0  & 0 &  2 &  0 & -1 & -1 & 0 \\
      0  & 0 &  0 &  2 &  0 & -1 & 0\\
     -1 & -1 & -1 &  0 &  2 &  0  & 0 \\
      0 &  1 & -1 &  -1 &  0 &  2  &  0 \\
     -1 &  0 & 0 &  0 &  0 & 0 &  2 \\
  \end{array}
  \right ]
     \end{array}$   &
    $\begin{array}{c}
    \frac{1}{4} \left [
       \begin{array}{ccccccc}
     8 & 4 & 4 &  0 & 8 & 0 & 4 \\
     4  & 7 & 1 & -2 & 6 & -4 & 2 \\
     4  & 1 & 7 &  2 & 6 & 4 & 2 \\
     0  & -2 & 2 & 4 & 0 & 4 & 0 \\
     8  & 6 & 6 &  0 & 12 & 0 & 4 \\
     0  & -4 & 4 & 4 & 0 & 8 & 0 \\
     4  & 2 & 2 & 0 & 4 & 0 &  4 \\
  \end{array}
 \right ]
     \end{array}$    \\
    & & \\
     %%\hline
     \hline
\end{tabular}
  \vspace{2mm}
  \caption{\small\hspace{3mm} (cont.) The partial Cartan matrix $B_{\Gamma}$
   and the inverse matrix $B^{-1}_{\Gamma}$ for Carter
   diagrams with the number of vertices $l = 7$}
  \label{tab_partial Cartan_2}
\end{table}

 \index{${\bf B}$, Cartan matrix associated with a Dynkin diagram}
 \index{Cartan matrix ! - ${\bf B}$}

 \begin{table}[H]
\tiny
  \centering
  \renewcommand{\arraystretch}{1.5}
  \begin{tabular} {|c|c|c|}
  \hline
      The Carter   & The partial Cartan         &   The inverse          \cr
        diagram    &  matrix $B_{\Gamma}$              &   matrix $B^{-1}_{\Gamma}$    \\
  \hline
      & & \\
     $\begin{array}{c} \includegraphics[scale=0.6]{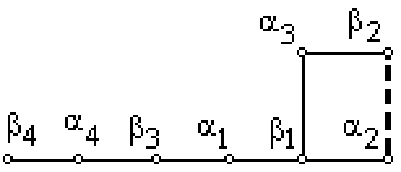} \\
        {\bf D_8(a_1) = D_8(a_5)} \end{array}$  &
     $\begin{array}{c}
       \left [
  \begin{array}{cccccccc}
      2  & 0 &  0 &  0 & -1 &  0 & -1 & 0 \\
      0  & 2 &  0 &  0 & -1 &  1 & 0 & 0 \\
      0  & 0 &  2 &  0 & -1 & -1 & 0 & 0 \\
      0  & 0 &  0 &  2 &  0 & 0  &-1 & -1 \\
     -1 & -1 & -1 &  0 &  2 & 0  & 0 & 0 \\
      0 &  1 & -1 &  0 &  0 & 2  & 0 & 0 \\
     -1 &  0 & 0 &  -1 &  0 & 0  & 2 & 0 \\
      0 &  0 &  0 & -1 &  0 & 0  & 0 & 2 \\
  \end{array}
  \right ]
     \end{array}$   &
    $\begin{array}{c}
    \frac{1}{2} \left [
       \begin{array}{cccccccc}
     8  & 4 & 4 & 4 & 8 & 0 & 6 & 2 \\
     4  & 4 & 2 & 2 & 5 & -1 & 3 & 1 \\
     4  & 2 & 4 & 2 & 5 & 1 & 3 & 1 \\
     4  & 2 & 2 & 4 & 4 & 0 & 4 & 2 \\
     8  & 5 & 5 & 4 & 10 & 0 & 6 & 2 \\
     0  & -1 & 1 & 0 & 0 & 2 & 0 & 0 \\
     6  & 3 & 3 & 4 & 6 & 0 &  6 & 2 \\
     2  & 1 & 1 & 2 & 2 & 0 & 2  & 2 \\
  \end{array}
 \right ]
     \end{array}$    \\
    & & \\
     \hline
      & & \\
     $\begin{array}{c} \includegraphics[scale=0.6]{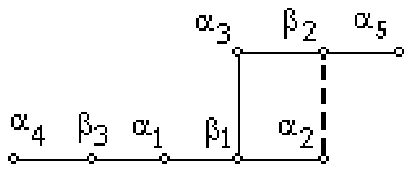} \\
        {\bf D_8(a_2) = D_8(a_4)} \end{array}$  &
     $\begin{array}{c}
       \left [
  \begin{array}{cccccccc}
      2  & 0 &  0 &  0 &  0 & -1 &  0 & -1 \\
      0  & 2 &  0 &  0 &  0 & -1 &  1 & 0  \\
      0  & 0 &  2 &  0 &  0 & -1 & -1 & 0  \\
      0  & 0 &  0 &  2 &  0 & 0  & 0 & -1  \\
      0  & 0 &  0 &  0 &  2 & 0 & -1 & 0 \\
     -1 & -1 & -1 &  0 &  0 & 2 & 0  & 0  \\
      0 &  1 & -1 &  0 & -1 & 0 & 2  & 0  \\
     -1 &  0 & 0 &  -1 &  0 & 0 & 0  & 2  \\
  \end{array}
  \right ]
     \end{array}$   &
    $\begin{array}{c}
    \frac{1}{2} \left [
       \begin{array}{cccccccc}
     6  & 3 & 3 & 2 & 0 & 6 & 0 & 4 \\
     3  & 4 & 1 & 1 & -1 & 4 & -2 & 2 \\
     3  & 1 & 4 & 1 & 1 & 4 & 2 & 2 \\
     2  & 1 & 1 & 2 & 0 & 2 & 0 & 2 \\
     0  & -1 & 1 & 0 & 2 & 0 & 2 & 0 \\
     6  & 4 & 4 & 2 & 0 & 8 & 0 & 4 \\
     0  & -2 & 2 & 0 & 2 & 0 &  4 & 0 \\
     4  & 2 & 2 & 2 & 0 & 4 & 0  & 4 \\
  \end{array}
 \right ]
     \end{array}$    \\
    & & \\
     \hline
      & & \\
     $\begin{array}{c} \includegraphics[scale=0.6]{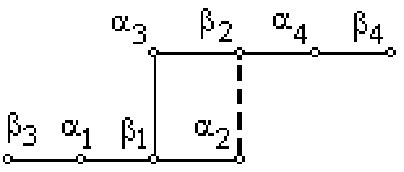} \\
        {\bf D_8(a_3)} \end{array}$  &
     $\begin{array}{c}
       \left [
  \begin{array}{cccccccc}
      2  & 0 &  0 &  0 & -1 &  0 & -1 & 0 \\
      0  & 2 &  0 &  0 & -1 &  1 & 0 & 0 \\
      0  & 0 &  2 &  0 & -1 & -1 & 0 & 0 \\
      0  & 0 &  0 &  2 &  0 & -1 & 0 & -1 \\
     -1 & -1 & -1 &  0 &  2 & 0  & 0 & 0 \\
      0 &  1 & -1 & -1 &  0 & 2  & 0 & 0 \\
     -1 &  0 & 0  &  0 &  0 & 0  & 2 & 0 \\
      0 &  0 &  0 & -1 &  0 & 0  & 0 & 2 \\
  \end{array}
  \right ]
     \end{array}$   &
    $\begin{array}{c}
    \frac{1}{2} \left [
       \begin{array}{cccccccc}
     4  & 2 & 2 & 0 & 4 & 0 & 2 & 0 \\
     2  & 4 & 0 & -2 & 3 & -3 & 1 & -1 \\
     2  & 0 & 4 & 2 & 3 & 3 & 1 & 1 \\
     0  & -2 & 2 & 4 & 0 & 4 & 0 & 2 \\
     4  & 3 & 3 & 0 & 6 & 0 & 2 & 0 \\
     0  & -3 & 3 & 4 & 0 & 6 & 0 & 2 \\
     2  & 1 & 1 & 0 & 2 & 0 &  2 & 0 \\
     0  & -1 & 1 & 2 & 2 & 0 & 2 & 0 \\
  \end{array}
 \right ]
     \end{array}$    \\
    & & \\
     \hline
      & & \\
     $\begin{array}{c} \includegraphics[scale=0.6]{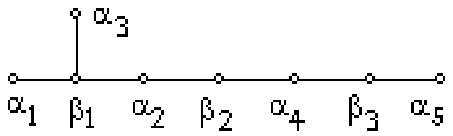} \\
        {\bf D_8} \end{array}$  &
     $\begin{array}{c}
       \left [
  \begin{array}{cccccccc}
      2  & 0 &  0 &  0 &  0 & -1 &  0 & 0 \\
      0  & 2 &  0 &  0 &  0 & -1 & -1 & 0 \\
      0  & 0 &  2 &  0 &  0 & -1 &  0 & 0 \\
      0  & 0 &  0 &  2 &  0 & 0  & -1 & -1 \\
      0  & 0 &  0 &  0 &  2 & 0  & 0 & -1 \\
     -1 & -1 & -1 &  0 &  0 & 2  & 0 & 0 \\
      0 & -1 &  0 &  -1 &  0 & 0  & 2  & 0 \\
      0 &  0 &  0 &  -1 & -1 & 0  & 0 & 2 \\
  \end{array}
  \right ]
     \end{array}$   &
    $\begin{array}{c}
    \frac{1}{2} \left [
       \begin{array}{cccccccc}
     8  & 4 & 4 & 4 & 8 & 0 & 6 & 2 \\
     4  & 4 & 2 & 2 & 5 & -1 & 3 & 1 \\
     4  & 2 & 4 & 2 & 5 & 1 & 3 & 1 \\
     4  & 2 & 2 & 4 & 4 & 0 & 4 & 2 \\
     8  & 5 & 5 & 4 & 10 & 0 & 6 & 2 \\
     0  & -1 & 1 & 0 & 0 & 2 & 0 & 0 \\
     6  & 3 & 3 & 4 & 6 & 0 &  6 & 2 \\
     2  & 1 & 1 & 2 & 2 & 0 & 2  & 2 \\
  \end{array}
 \right ]
     \end{array}$    \\
    & & \\
     \hline
\end{tabular}
  \vspace{2mm}
  \caption{\small\hspace{3mm} (cont.) The partial Cartan matrix $B_{\Gamma}$
   and the inverse matrix $B^{-1}_{\Gamma}$ for Carter
   diagrams $D_8(a_1)$, $D_8(a_2)$, $D_8(a_3)$, $D_8$}
  \label{tab_partial Cartan_2a}
\end{table}

 \begin{table}[H]
\tiny
  \centering
  \renewcommand{\arraystretch}{1.5}
  \begin{tabular} {|c|c|c|}
  \hline
      The Carter   & The Cartan         &   The inverse          \cr
        diagram    &  matrix ${\bf B}$              &   matrix ${\bf B}^{-1}$    \\
  \hline
      & & \\
     $\begin{array}{c} \includegraphics[scale=0.6]{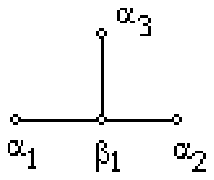} \\
        {\bf D_4} \end{array}$ &
     $\begin{array}{c}
       \left [
   \begin{array}{cccc}
     2 & 0 & 0  &   -1  \\
     0 & 2 & 0  &   -1  \\
     0 & 0 & 2  &   -1  \\
    -1 & -1 & -1 &   2  \\
  \end{array}
  \right ]
     \end{array}$   &
      $\begin{array}{c}
    \frac{1}{2} \left [ \begin{array}{cccc}
     2 & 1 & 1 & 2  \\
     1 & 2 & 1 & 2  \\
     1 & 1 & 2 & 2  \\
     2 & 2 & 2 & 4  \\
  \end{array} \right ]
     \end{array}$    \\
    & & \\
  \hline  %1
    & & \\
     $\begin{array}{c} \includegraphics[scale=0.6]{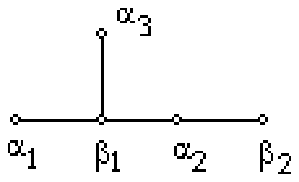} \\
        {\bf D_5} \end{array}$ &
     $\begin{array}{c}
       \left [
   \begin{array}{ccccc}
     2 & 0 & 0  &   -1 &  0 \\
     0 & 2 & 0  &   -1 & -1 \\
     0 & 0 & 2  &   -1 &  0 \\
    -1 & -1 & -1 &   2 &  0 \\
     0 & -1 & 0  &   0 & 2  \\
  \end{array}
  \right ]
     \end{array}$   &
      $\begin{array}{c}
    \frac{1}{4} \left [ \begin{array}{ccccc}
     5 & 4 & 3 & 6 & 2  \\
     4 & 8 & 4 & 8 & 4  \\
     3 & 4 & 5 & 6 & 2  \\
     6 & 8 & 6 & 12 & 4 \\
     2 & 4 & 2 & 4 & 4 \\
  \end{array} \right ]
     \end{array}$    \\
    & & \\
  \hline  %2
    & & \\
     $\begin{array}{c} \includegraphics[scale=0.6]{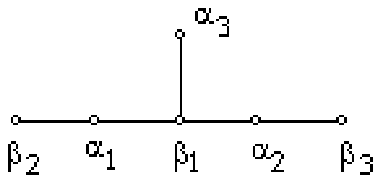} \\
        {\bf E_6} \end{array}$ &
     $\begin{array}{c}
       \left [
   \begin{array}{cccccc}
     2 & 0 & 0  &   -1 & -1 & 0  \\
     0 & 2 & 0  &   -1 &  0 & -1 \\
     0 & 0 & 2  &   -1 &  0 & 0  \\
    -1 & -1 & -1 &   2 & 0 & 0  \\
    -1 &  0 & 0  &   0 & 2 & 0 \\
      0 & -1 & 0 &   0 & 0 & 2 \\
  \end{array}
  \right ]
     \end{array}$   &
      $\begin{array}{c}
    \frac{1}{3} \left [ \begin{array}{cccccc}
     10 & 8 & 6 & 12 & 5 & 4  \\
     8 & 10 & 6 & 12 & 4 & 5 \\
     6 & 6 &  6  & 9 & 3 & 3  \\
    12 & 12 &  9 & 18 & 6 & 6 \\
     5 & 4 &  3  & 6 & 4 & 2 \\
     4 & 5 & 3 & 6 & 2 & 4 \\
  \end{array} \right ]
     \end{array}$    \\
    & & \\
  \hline  %3
    & & \\
     $\begin{array}{c} \includegraphics[scale=0.6]{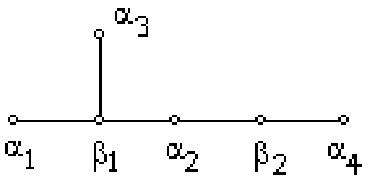} \\
        {\bf D_6} \end{array}$ &
     $\begin{array}{c}
       \left [
   \begin{array}{cccccc}
     2 & 0 & 0  &  0 & -1 &  0 \\
     0 & 2 & 0  &  0 & -1 & -1 \\
     0 & 0 & 2  &  0 & -1 &  0 \\
     0 & 0 & 0  &  2 &  0 & -1 \\
    -1 & -1 & -1 & 0 &  2 &  0 \\
     0 & -1 & 0  & -1 & 0 & 2  \\
  \end{array}
  \right ]
     \end{array}$   &
      $\begin{array}{c}
    \frac{1}{2} \left [ \begin{array}{cccccc}
     3 & 3 & 2 & 1 & 4 & 2 \\
     3 & 6 & 3 & 2 & 6 & 4 \\
     2 & 3 & 3 & 1 & 4 & 2 \\
     1 & 2 & 1 & 2 & 2 & 2 \\
     4 & 6 & 4 & 2 & 8 & 4 \\
     2 & 4 & 2 & 2 & 4 & 4 \\
  \end{array} \right ]
     \end{array}$    \\
    & & \\
     \hline %4
    & & \\
     $\begin{array}{c} \includegraphics[scale=0.6]{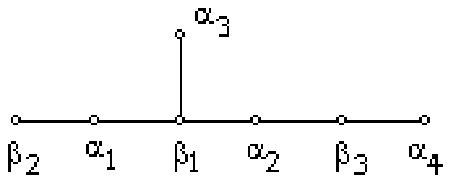} \\
        {\bf E_7} \end{array}$ &
     $\begin{array}{c}
       \left [
   \begin{array}{ccccccc}
     2 & 0 & 0  &  0 & -1 & -1 & 0  \\
     0 & 2 & 0  &  0 & -1 &  0 & -1 \\
     0 & 0 & 2  &  0 & -1 &  0 & 0  \\
     0 & 0 & 0  &  2 &  0 &  0 & -1 \\
    -1 & -1 & -1 & 0 &  2 &  0 & 0  \\
    -1 &  0 &  0 & 0 &  0 &  2 & 0 \\
     0 & -1 &  0 & -1 &  0 &  0 & 2 \\
  \end{array}
  \right ]
     \end{array}$   &
      $\begin{array}{c}
    \frac{1}{2} \left [ \begin{array}{ccccccc}
     12 & 12 & 8 & 4  & 16 & 6 & 8  \\
     12 & 15 & 9 & 5  & 18 & 6 & 10 \\
     8  & 9  & 7 & 3  & 12 & 4 & 6 \\
     4  & 5  & 3 & 3  & 6  & 2 & 4 \\
     16 & 18 & 12 & 6 & 24 & 8 & 12 \\
     6  & 6  & 4 & 2  & 8  & 4 & 4 \\
     8  & 10 & 6 & 4 & 12  & 4 & 8 \\
  \end{array} \right ]
     \end{array}$    \\
    & & \\
     \hline  %5
    & & \\
     $\begin{array}{c} \includegraphics[scale=0.6]{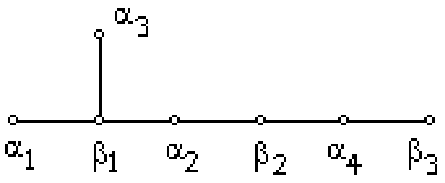} \\
        {\bf D_7} \end{array}$ &
     $\begin{array}{c}
       \left [
   \begin{array}{ccccccc}
     2 & 0 & 0  &  0 & -1 &  0 & 0 \\
     0 & 2 & 0  &  0 & -1 & -1 & 0 \\
     0 & 0 & 2  &  0 & -1 &  0 & 0 \\
     0 & 0 & 0  &  2 &  0 & -1 & -1 \\
    -1 & -1 & -1 & 0 &  2 &  0 & 0 \\
     0 & -1 & 0  & -1 & 0 &  2 & 0 \\
     0 & 0  & 0  & -1 & 0 &  0 & 2 \\
  \end{array}
  \right ]
     \end{array}$   &
      $\begin{array}{c}
    \frac{1}{4} \left [ \begin{array}{ccccccc}
     7 & 8  & 5 & 4 & 10 & 6 & 2 \\
     8 & 16 & 8 & 8 & 16 & 12 & 4 \\
     5 & 8  & 7 & 4 & 10 & 6 & 2 \\
     4 & 8  & 4 & 8 & 8  & 8  & 4 \\
    10 & 16 & 10 & 8 & 20 & 12 & 4 \\
     6 & 12 & 6 & 8 & 12 & 12 & 4 \\
     2 & 4  &  2 & 4 & 4 & 4  & 4 \\
  \end{array} \right ]
     \end{array}$    \\
    & & \\
     \hline
\end{tabular}
  \vspace{2mm}
  \caption{\small\hspace{3mm} (cont.) The Cartan matrix ${\bf B}$ and the inverse matrix ${\bf B}^{-1}$
  for the conjugacy classes $E_6$, $E_7$, $D_4$, $D_5$, $D_6$, $D_7$}
  \label{tab_partial Cartan_3}
\end{table}

 \begin{table}[H]
 \tiny
  \centering
  \renewcommand{\arraystretch}{1.5}
  \begin{tabular} {|c|c|c|c|}
  \hline
             & The Cartan  &  The inverse        &  S.e. \cr
     $A_l$   & matrix ${\bf B}$        &  matrix ${\bf B}^{-1}$  &  in $\tau$   \\
             %%&                     &                     &  \\
         %%% S.e. - simple extension
   \hline
    &  &  &  \\
     $\begin{array}{c} \includegraphics[scale=0.6]{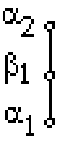}
        \\ A_3 \end{array}$  &
      $\left [
   \begin{array}{ccc}
    2  &  0  &  -1 \\
    0  &  2  &  -1 \\
   -1  & -1  &  2  \\
  \end{array}
  \right ]$
   $\begin{array}{c}
    \alpha_1 \\
    \alpha_2 \\
    \beta_1  \\
    \end{array}$
     &
     $\frac{1}{4}\left [
     \begin{array}{ccc}
        3  &  1  & 2 \\
        1  &  3  & 2 \\
        2  &  2  & 4 \\
      \end{array}
      \right ]$
   $\begin{array}{c}
    \alpha_1 \\
    \alpha_2 \\
    \beta_1  \\
    \end{array}$
    &
      $\begin{array}{c}
        A_4 \\
        A_4 \\
        D_4 \\
      \end{array}$ \\
     &  &  &  \\
     \hline
    &  & &  \\
     $\begin{array}{c} \includegraphics[scale=0.6]{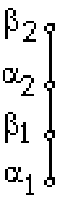}
       \\ A_4 \end{array}$  &
      $\left [
   \begin{array}{cccc}
    2  &  0  &  -1  &  0 \\
    0  &  2  &  -1  &  -1 \\
   -1  & -1  &  2  &  0 \\
    0  & -1  &  0  & 2 \\
  \end{array}
  \right ]$
  $\begin{array}{c}
    \alpha_1 \\
    \alpha_2 \\
    \beta_1 \\
    \beta_2\\
    \end{array}$ &
     $ \frac{1}{5}\left [
     \begin{array}{cccc}
        4  &  2  & 3  &  1 \\
        2  &  6  & 4  &  3 \\
        3  &  4  & 6  &  2 \\
        1  &  3  & 2  &  4 \\
  \end{array} \right ]$
       $\begin{array}{c}
    \alpha_1 \\
    \alpha_2 \\
    \beta_1 \\
    \beta_2\\
    \end{array}$ &
      $\begin{array}{c}
        A_5 \\
        D_5 \\
        D_5 \\
        A_5 \\
      \end{array}$
           \\
      &  &  &  \\
     \hline
    &  & & \\
     $\begin{array}{c} \includegraphics[scale=0.6]{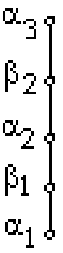}
       \\ A_5 \end{array}$  &
      $\left [
   \begin{array}{ccccc}
    2  &  0  &  0 & -1  &  0 \\
    0  &  2  &  0 & -1  &  -1 \\
    0  &  0  &  2 &  0  &  -1 \\
   -1  & -1  &  0 &  2  &  0 \\
    0  & -1  & -1  &  0  & 2 \\
  \end{array}
  \right ]$
          $\begin{array}{c}
    \alpha_1 \\
    \alpha_2 \\
    \alpha_3 \\
    \beta_1 \\
    \beta_2\\
    \end{array}$ &
     $ \frac{1}{6}\left [
     \begin{array}{ccccc}
              5   &  3  &  1   &  4 & 2 \\
              3   &  9  &  3   &  6 & 6 \\
              1   &  3  &  5   &  2 & 4 \\
              4   &  6  &  2   &  8 & 4 \\
              2   &  6  &  4   &  4 & 8 \\
  \end{array}
  \right ]$
     $\begin{array}{c}
    \alpha_1 \\
    \alpha_2 \\
    \alpha_3 \\
    \beta_1 \\
    \beta_2\\
    \end{array}$ &
      $\begin{array}{c}
        A_6 \\
        E_6 \\
        A_6 \\
        D_6 \\
        D_6 \\
      \end{array}$ \\
      &  &  &   \\
     \hline
      &  &  &   \\
     $\begin{array}{c} \includegraphics[scale=0.6]{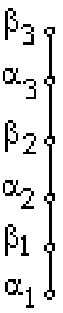}
        \\ A_6 \end{array}$  &
      $\left [
   \begin{array}{cccccc}
    2  &  0  &  0 & -1  &  0  & 0 \\
    0  &  2  &  0 & -1  &  -1 & 0 \\
    0  &  0  &  2 &  0  &  -1 & -1 \\
   -1  & -1  &  0 &  2  &  0  & 0 \\
    0  & -1  & -1  &  0 &  2  & 0 \\
    0  & 0   & -1 &   0 &  0  & 2 \\
  \end{array}
  \right ]$
     $\begin{array}{c}
    \alpha_1 \\
    \alpha_2 \\
    \alpha_3 \\
    \beta_1 \\
    \beta_2\\
    \beta_3 \\
    \end{array}$ &
     $ \frac{1}{7}\left [
     \begin{array}{cccccc}
              6   &  4  &  1   &  5 & 3 & 1 \\
              4   &  12  & 6   &  8 & 9 & 3 \\
              2   &  6  &  10  &  4 & 8 & 5 \\
              5   &  8  &  4   &  10 & 6 & 2 \\
              3   &  9  &  8   &  6 & 12 & 4 \\
              1   &  3  &  5   &  2 & 4 & 6 \\
  \end{array}
  \right ]$
   $\begin{array}{c}
    \alpha_1 \\
    \alpha_2 \\
    \alpha_3 \\
    \beta_1 \\
    \beta_2\\
    \beta_3 \\
    \end{array}$  &
      $\begin{array}{c}
        A_7 \\
        E_7 \\
        D_7 \\
        D_7 \\
        E_7 \\
        A_7 \\
      \end{array}$ \\
      &  &  &   \\
     \hline
      &  &  &   \\
     $\begin{array}{c} \includegraphics[scale=0.6]{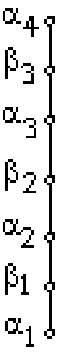}
        \\ A_7 \end{array}$  &
      $\left [
   \begin{array}{ccccccc}
    2  &  0  &  0 & 0 & -1  &  0  & 0 \\
    0  &  2  &  0 & 0 & -1  &  -1 & 0 \\
    0  &  0  &  2 & 0 &  0  &  -1 & -1 \\
    0  &  0  &  0 & 2 &  0  &  0  & -1 \\
   -1  & -1  &  0 & 0 &  2  &  0  & 0 \\
    0  & -1  & -1 & 0 &  0  &  2  & 0 \\
    0  & 0   & -1 & -1 &  0 &  0  & 2 \\
  \end{array}
  \right ]$
     $\begin{array}{c}
    \alpha_1 \\
    \alpha_2 \\
    \alpha_3 \\
    \alpha_4 \\
    \beta_1 \\
    \beta_2\\
    \beta_3 \\
    \end{array}$ &
     $ \frac{1}{8}\left [
     \begin{array}{ccccccc}
              7   &  5  &  3  & 1  &  6  & 4 & 2 \\
              5   &  15 &  9  & 3  &  10 & 12 & 6 \\
              3   &  9  &  15 & 5  &  6  & 12 & 10 \\
              1   &  3  &  5  & 7  &  2  & 4  & 6 \\
              6   &  10 &  6  & 2  &  12 & 8  & 4 \\
              4   &  12  & 12 & 4  &  8 & \fbox{16} & 8 \\
              2   &  6  &  10 & 6  &  4 & 8 & 12 \\
  \end{array}
  \right ]$
     $\begin{array}{c}
    \alpha_1 \\
    \alpha_2 \\
    \alpha_3 \\
    \alpha_4 \\
    \beta_1 \\
    \beta_2\\
    \beta_3 \\
    \end{array}$ &
      $\begin{array}{c}
        A_8 \\
        E_8 \\
        E_8 \\
        A_8 \\
        D_8 \\
         ~  \\
        D_8 \\
      \end{array}$ \\
      &  &  &   \\
     \hline
     $\begin{array}{c} \includegraphics[scale=0.6]{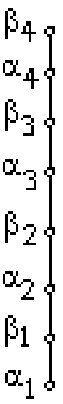}
        \\ A_8 \end{array}$  & $\left [
   \begin{array}{cccccccc}
    2  &  0  &  0 & 0 & -1  &  0  &  0 &  0 \\
    0  &  2  &  0 & 0 & -1  &  -1 &  0 &  0 \\
    0  &  0  &  2 & 0 &  0  &  -1 & -1 &  0 \\
    0  &  0  &  0 & 2 &  0  &  0  & -1 & -1 \\
   -1  & -1  &  0 & 0 &  2  &  0  &  0 &  0 \\
    0  & -1  & -1 & 0 &  0  &  2  &  0 &  0 \\
    0  & 0   & -1 & -1 &  0 &  0  &  2 &  0 \\
    0  & 0   &  0 & -1 &  0 &  0  &  0 &  2 \\
  \end{array} \right ]$
     &
     $ \frac{1}{9}\left [
     \begin{array}{cccccccc}
              8   &         6  &  4   &   2  &  7  &  5 & 3 & 1 \\
              6   &  \fbox{18} &  12  &  6  &  12 & 15  & 9 & 3 \\
              4   &  12  &  \fbox{20} &  10  &  8  & 16 & 15 & 5 \\
              2   &  6   &     10     & 14   &  4  & 8  & 12 & 7 \\
              7   &  12  &  8  & 4    &  14  & 10  & 6  & 2  \\
              5   &  15  &  16 & 8    &  10  & \fbox{20} & 12 & 4 \\
              3   &  9   &  15 & 12   &  6   & 12 & \fbox{18} & 6 \\
              1   &  3   &  5  & 7    &  2   & 4  & 6 & 8 \\
  \end{array}
   \right ]$
     &
      $\begin{array}{c}
        A_8 \\
        ~ \\
        ~ \\
        D_8 \\
        D_8 \\
        ~ \\
        ~ \\
        A_8 \\
      \end{array}$ \\
      &  &  &   \\
     \hline
  \end{tabular}
  \vspace{2mm}
  \caption{\small\hspace{3mm} (cont.) The Cartan matrix ${\bf B}$ and its inverse matrix ${\bf B}^{-1}$
    for the Carter diagram $A_l$. The diagonal elements $b^{\vee}_{\tau, \tau}$ in ${\bf B}^{-1}$ correspond to
    simple extension in the point $\tau$ if $b^{\vee}_{\tau, \tau} < 2$. In frames: $b^{\vee}_{\tau, \tau} \geq 2$.}
  \label{tab_partial Cartan_Al}
\end{table}

%% file: C_linkageDiagr.tex
\clearpage
\section{\sc\bf Linkage diagrams $\gamma^{\nabla}(8)$ and inequality $\mathscr{B}^{\vee}_{\Gamma}(\gamma^{\nabla}) < 2$}
   \label{sec_solutions}

\subsection{The linkage diagrams $\gamma^{\nabla}_{ij}(8)$
 and solutions of inequality $\mathscr{B}^{\vee}_{\Gamma}(\gamma^{\nabla}_{ij}(8)) < 2$}
~\\

\begin{table}[H]
  \centering
  \renewcommand{\arraystretch}{1.4}
  \small
  \begin{tabular} {||c|c|c|l|l||}
  \hline \hline
     The Carter &  $L_{ij}$  & $\mathscr{B}^{\vee}_{\Gamma}(\gamma^{\nabla}_{ij}(8)) = p < 2$
                                & The linkage diagram  & $p$ \cr
      diagram &    &            &  \qquad \qquad $\gamma^{\nabla}_{ij}(8)$  &  \\
    \hline  %%% E_6(a_1)
                  &  &   $q = 4(b_2^2 + b_2{b_3} + b_3^2)$ & & \cr
                  \cline{3-3}
                  &  $L_{12}$
                  &  $\frac{1}{3}(10 + 2(4b_2 + 5b_3) + q) < 2$, or
                  &  $\{ 0, 0, 1, 0, 0, -1 \}$  &  \cr
                  &
                  &  $(b_2 + b_3)^2 + (b_2 + 2)^2  + (b_3 + \frac{5}{2})^2 < \frac{33}{4}$
                  &  $\{ 0, 0, 1, 0, -1, -1 \}$ & \cr
                  \cline{2-4}
    $E_6(a_1)$    &  $L_{13}$
                  &  $\frac{1}{3}(4 + 2(-b_2 + b_3) + q) < 2$, or
                  &  $\{ 0, 1, 0, 0, 0, 0 \}$  & $\frac{4}{3}$ \cr
                  &
                  &  $(b_2 + b_3)^2 + (b_2 - \frac{1}{2})^2  + (b_3 + \frac{1}{2})^2 < \frac{3}{2}$
                  &  $\{ 0, 1, 0, 0,  1, -1 \}$ &   \cr
                  \cline{2-4}
                  &  $L_{23}$
                  &  $\frac{1}{3}(4 + 2(b_2 + 2b_3) + q) < 2$,
                  &  $\{ 1, 0, 0, 0, 0, 0 \}$  &  \cr
                  &
                  &  $(b_2 + b_3)^2 + (b_2 + \frac{1}{2})^2  + (b_3 + 1)^2 < \frac{9}{4}$
                  & $\{ 1, 0, 0, 0,  0, -1 \}$ & \\
    \hline\hline   %%% E_6(a_2)
                  &  &   $q = 4(b_2^2 + b_2{b_3} + b_3^2)$ & & \cr
                  \cline{3-3}
                  &   $L_{12}$
                  &  $\frac{1}{3}(6 + 2(3b_2 + 3b_3) + q) < 2$, or
                  &  $\{ 0, 0, 1, 0, 0, -1 \}$  &  \cr
                  &
                  & $(b_2 + b_3)^2 + (b_2 + \frac{3}{2})^2  + (b_3 + \frac{3}{2})^2 < \frac{9}{2}$
                  &  $\{ 0, 0, 1, 0, -1, 0 \}$ & \cr
                  \cline{2-4}
       $E_6(a_2)$ &   $L_{13}$
                  &   $\frac{1}{3}(4 + 2(-2b_2 - b_3) + q < 2$, or
                  &  $\{ 0, 1, 0, 0, 0, 0 \}$  & $\frac{4}{3}$ \cr
                  &
                  & $(b_2 + b_3)^2 + (b_2 - 1)^2  + (b_3 - \frac{1}{2})^2 < \frac{9}{4}$
                  &  $\{ 0, 1, 0, 0,  1, 0 \}$ &  \cr
                  \cline{2-4}
                  &   $L_{23}$
                  & $\frac{1}{3}(4 + 2(-b_2 - 2b_3) + q) < 2$, or
                  & $\{ 1, 0, 0, 0, 0, 0 \}$ & \cr
                  &
                  & $(b_2 + b_3)^2 + (b_2 - \frac{1}{2})^2  + (b_3 - 1)^2 < \frac{9}{4}$
                  & $\{ 1, 0, 0, 0,  0, 1 \}$ &  \\
      \hline\hline   %%% E_7(a_1)
                  &  &   $q = 3b_2^2 + 4b_3^2 + 3b_4^2 + 4b_2b_3 + 2b_2b_4 + 4b_3b_4$ & & \cr
                  \cline{3-3}
                  &   $L_{12}$
                  &   $\frac{1}{2}(12 + 2(4b_2 + 6b_3 + 4b_4) +  q) < 2$, or
                  & $\{ 0, 0, 1, 0, 0, -1, -1 \}$ & \cr
                  &
                  & $(b_2 + b_4 + 1)^2 +
                    2(b_3 + b_4 + \frac{3}{2})^2  + 2(b_2 + b_3 + \frac{3}{2})^2 < 2$
                  & $\{ 0, 0, 1, 0, -1, -1, 0 \}$ & \cr
                  \cline{2-4}
     $E_7(a_1)$   & $L_{13}$
                  & $\frac{1}{2}(4 + 2(2b_3 + 2b_4) + q) < 2$, or
                  & $\{ 0, 1, 0, 0, 1, -1, 0 \}$ &  $\frac{3}{2}$\cr
                  &
                  & $(b_2 + b_4)^2 + 2(b_3 + b_4 + 1)^2 + 2(b_2 + b_3)^2 < 2$
                  & $\{ 0, 1, 0, 0, 0, 0, -1 \}$ & \cr
                  \cline{2-4}
                  & $L_{23}$
                  & $\frac{1}{2}(8 + 2(2b_2 + 4b_3 + 4b_4) + q) <2$, or
                  & $\{ 1, 0, 0, 0, 0, 0, -1 \}$ & \cr
                  &
                  & $(b_2 + b_4 + 1)^2 +
                     2(b_3 + b_4 + \frac{3}{2})^2  + 2(b_2 + b_3 + \frac{1}{2})^2 < 2$
                  & $\{ 1, 0, 0, 0, 0, -1, -1 \}$ & \\
      \hline\hline   %%% E_7(a_2)
                  &  &   $q = 3b_2^2 + 3b_3^2 + 3b_4^2 + 2b_2b_3 - 2b_2b_4 + 2b_3b_4$ & & \cr
                  \cline{3-3}
                  & $L_{12}$
                  & $\frac{1}{2}(8 + 2(2b_2 + 4b_3 + 2b_4) +  q) < 2$, or
                  & $\{ 0, 0, 1, 0, 0, -1, 0 \}$ & \cr
                  &
                  & $(b_2 + b_3)^2 + (b_3 + b_4)^2 + (b_2 - b_4)^2 + (b_2 + 2)^2 \quad\quad\quad$
                  & $\{ 0, 0, 1, 0, -1, -1, -1 \}$ & \cr
                  &
                  & $\quad\quad\quad\quad\quad\quad\quad + (b_3 + 4)^2 + (b_4 + 2)^2 < 20$ & & \cr
                  \cline{2-4}
       $E_7(a_2)$ & $L_{13}$
                  & $\frac{1}{2}(8 + 2(-2b_2 + 2b_3 + 4b_4) + q) < 2$, or
                  & $\{ 0, 1, 0, 0, 0, 0, -1 \}$ &  $\frac{3}{2}$\cr
                  &
                  & $(b_2 + b_3)^2 + (b_3 + b_4)^2 + (b_2 - b_4)^2 + (b_2 - 2)^2  \quad\quad\quad$
                  & $\{ 0, 1, 0, 0, 1, -1, -1 \}$ & \cr
                  &
                  & $\quad\quad\quad\quad\quad\quad\quad + (b_3 + 2)^2 + (b_4 + 4)^2 < 20$ & & \cr
                  \cline{2-4}
                  & $L_{23}$
                  & $\frac{1}{2}(4 + 2(2b_3 + 2b_4) + q) < 2$, or
                  & $\{ 1, 0, 0, 0, 0, 0, -1 \}$ & \cr
                  &
                  & $(b_2 + b_3)^2 + (b_3 + b_4)^2 + (b_2 - b_4)^2 + (b_3 + 2)^2  \quad\quad\quad$
                  & $\{ 1, 0, 0, 0, 0, -1, 0 \}$ & \cr
                  &
                  & $\quad\quad\quad\quad\quad\quad\quad + (b_4 + 2)^2 + b_2^2 < 8$ & & \cr
       \hline\hline   %%%
\end{tabular}
  \vspace{2mm}
  \caption{\hspace{3mm}The linkage diagrams $\gamma^{\nabla}_{ij}(8)$ obtained as solutions of inequality
    $\mathscr{B}^{\vee}_{\Gamma}(\gamma^{\nabla}_{ij}(8)) < 2$}
  \label{sol_inequal_1}
  \end{table}

\begin{table}[H]
  \centering
  \renewcommand{\arraystretch}{1.6}
  \small
  \begin{tabular} {||c|c|c|l|l||}
  \hline \hline
      The Carter &  $L_{ij}$  & $\mathscr{B}^{\vee}_{\Gamma}(\gamma^{\nabla}_{ij}(8)) = p < 2$
                                & The linkage diagram  & $p$ \cr
      diagram &    &            &  \qquad \qquad $\gamma^{\nabla}_{ij}(8)$  &  \\
    \hline  %%% E_7(a_3)
                  &  &   $q = 4b_2^2 + 3b_3^2 + 3b_4^2 + 4b_2{b_3} - 4b_2{b_4} - 2b_3{b_4}$ & & \cr
                  \cline{3-3}
                  &  $L_{12}$
                  &  $\frac{1}{2}(4 + 2(2b_2 + 2b_3) + q) < 2$, or
                  &  $\{ 0, 0, 1, 0, 0, -1, 0 \}$  &  \cr
                  &
                  &  $2(b_2 + b_3 + 1)^2 + (b_3 - b_4)^2 + 2(b_4 - b_2)^2 < 2$
                  &  $\{ 0, 0, 1, 0, -1, 0, -1 \}$ & \cr
                  \cline{2-4}
    $E_7(a_3)$    &  $L_{13}$
                  &  $\frac{1}{2}(8 + 2(-4b_2 -2b_3 + 4b_4) + q) < 2$, or
                  &  $\{ 0, 1, 0, 0, 0, 0, -1 \}$  & $\frac{3}{2}$ \cr
                  &
                  &  $2(b_2 + b_3 - \frac{1}{2})^2 +
                       (b_3 - b_4 - 1)^2  + 2(b_4 -b_2 + \frac{3}{2})^2 < 2$
                  &  $\{ 0, 1, 0, 0, 1,  0, -1 \}$ &   \cr
                  \cline{2-4}
                  &  $L_{23}$
                  &  $\frac{1}{2}(4 + 2(-2b_2 - 2b_3 + 2b_4) + q) < 2$, or
                  &  $\{ 1, 0, 0, 0, 0, 1, 0 \}$  &  \cr
                  &
                  &  $2(b_2 + b_3 - \frac{1}{2})^2 +
                       (b_3 - b_4 - 1)^2 + 2(b_4 - b_2 + \frac{1}{2})^2 < 2$
                  &  $\{ 1, 0, 0, 0, 0, 0, -1 \}$ & \\
     \hline  %%% E_7(a_4)
                  &  &   $q = 3b_2^2 + 3b_3^2 + 3b_4^2 - 2b_2{b_3} - 2b_2{b_4} - 2b_3{b_4}$ & & \cr
                  \cline{3-3}
                  &  $L_{12}$
                  &  $\frac{1}{2}(4 + 2(2b_2 - 2b_3) + q) < 2$, or
                  &  $\{ 0, 0, 1, 0, 0, 1, 0 \}$  &  \cr
                  &
                  &  $(b_2 - b_3)^2 + (b_3 - b_4)^2 + (b_2 - b_4)^2\quad\quad\quad\quad$
                  &  $\{ 0, 0, 1, 0, -1, 0, 0 \}$ & \cr
                  &
                  &  $\quad\quad\quad\quad+ (b_2 + 2)^2  + (b_3 - 2)^2 + b_4^2 < 2$
                  &  & \cr
                  \cline{2-4}
    $E_7(a_4)$    &  $L_{13}$
                  &  $\frac{1}{2}(4 + 2(-2b_2 + 2b_4) + q) < 2$, or
                  &  $\{ 0, 1, 0, 0, 1, 0, 0 \}$  & $\frac{3}{2}$ \cr
                  &
                  &  $(b_2 - b_3)^2 + (b_3 - b_4)^2 + (b_2 - b_4)^2\quad\quad\quad\quad$
                  &  $\{ 0, 1, 0, 0, 0, 0, -1 \}$ & \cr
                  &
                  &  $\quad\quad\quad\quad+ (b_2 - 2)^2  + b_3^2 + (b_4 + 2)^2 < 2$
                  &  & \cr
                  \cline{2-4}
                  &  $L_{23}$
                  &  $\frac{1}{2}(4 + 2(2b_3 - 2b_4) + q) < 2$, or
                  &  $\{ 1, 0, 0, 0, 0, 0, 1 \}$  &  \cr
                  &
                  &  $(b_2 - b_3)^2 + (b_3 - b_4)^2 + (b_2 - b_4)^2\quad\quad\quad\quad$
                  &  $\{ 1, 0, 0, 0, 0, -1, 0 \}$ & \cr
                  &
                  &  $\quad\quad\quad\quad+ b_2^2  + (b_3 + 2)^2 + (b_4 - 2)^2 < 2$
                  &  & \\
    \hline\hline  %%% D_5(a_1)
                  & $L_{12}$
                  & $\frac{1}{4}(5 + 2(2b_2) + 4b^2) < 2$, or
                  & $\{ 0, 0, 1, 0, 0 \}$ & \cr
                  &
                  & $(b_2 + \frac{1}{2})^2 < 1$
                  & $\{ 0, 0, 1, 0, -1 \}$ & $\frac{5}{4}$ \cr
                  \cline{2-4}
       $D_5(a_1)$ & $L_{13}$
                  & $\frac{1}{4}(5 - 2(2b_2) +  4b^2) < 2$, or
                  & $\{ 0, 1, 0, 0, 0 \}$ & \cr
                  &
                  & $(b_2 - \frac{1}{2})^2 < 1$
                  & $\{ 0, 1, 0, 0, 1 \}$ & \cr
                  \cline{2-5}
                  & $L_{23}$
                  & $\frac{1}{4}(4 +  4b^2) < 2$, or $b_2^2 < 1$
                  & $\{ 1, 0, 0, 0, 0 \}$ & $1$ \\
       \hline\hline  %%% D_6(a_1)
                  &  &   $q = 2b_2^2 + 2b_3^2$ & & \cr
                  \cline{3-3}
                  & $L_{12}$
                  & $\frac{1}{2}(3 + 2(b_2 + b_3) +  q) < 2$, or
                  & $\{ 0, 0, 1, 0, 0, 0 \}$ & \cr
                  &
                  & $(b_2 + \frac{1}{2})^2 + (b_3 + \frac{1}{2})^2 < 1$
                  & $\{ 0, 0, 1, 0, -1, 0 \}$ &  \cr
                  & & & $\{ 0, 0, 1, 0, 0, -1 \}$ & \cr
                  & & & $\{ 0, 0, 1, 0, -1, -1 \}$ & \cr
                  \cline{2-4}
       $D_6(a_1)$ & $L_{13}$
                  & $\frac{1}{2}(3 + 2(-b_2 + b_3) +  q) < 2$, or
                  & $\{ 0, 1, 0, 0, 0, 0 \}$ & $\frac{3}{2}$ \cr
                  &
                  & $(b_2 - \frac{1}{2})^2 + (b_3 + \frac{1}{2})^2 < 1$
                  & $\{ 0, 1, 0, 0, 1, 0\}$ & \cr
                  & & & $\{ 0, 1, 0, 0, -1 \}$ & \cr
                  & & & $\{ 0, 1, 0, 1, -1 \}$ & \cr
                  \cline{2-5}
                  & $L_{23}$
                  & $\frac{1}{2}(4 + 2b_3 + q) < 2$, or $b_2^2 + (b_3 + 1)^2 < 1$
                  & $\{ 1, 0, 0, 0, -1 \}$ & $1$ \\
    \hline\hline
    \end{tabular}
  \vspace{2mm}
  \caption{\hspace{3mm}(cont.) The linkage diagrams $\gamma^{\nabla}_{ij}(8)$
  obtained as solutions of inequality  $\mathscr{B}^{\vee}_{\Gamma}(\gamma^{\nabla}_{ij}(8)) < 2$}
  \label{sol_inequal_2}
  \end{table}

  \begin{table}[H]
  \centering
  \renewcommand{\arraystretch}{1.7}
  \small
  \begin{tabular} {||c|c|c|l|l||}
  \hline \hline
      The Carter &  $L_{ij}$  & $\mathscr{B}^{\vee}_{\Gamma}(\gamma^{\nabla}_{ij}(8)) = p < 2$
                                & The linkage diagram  & $p$ \cr
      diagram &    &            &  \qquad \qquad $\gamma^{\nabla}_{ij}(8)$  &  \\
     \hline\hline  %%% D_6(a_2)
                  &  &   $q = 2a_4^2 + 4b_2^2 + 4a_4{b_2}$ & & \cr
                  \cline{3-3}
                  & $L_{12}$
                  & $\frac{1}{2}(3 + 2(a_4 + 2b_2) +  q) < 2$, or
                  & $\{ 0, 0, 1, 0, 0, 0 \}$ & \cr
                  &
                  & $(a_4 + 2b_2 + 1)^2 + a_4^2 < 2$
                  & $\{ 0, 0, 1, 0, 0, -1 \}$ &  \cr
                  & & & $\{ 0, 0, 1, 1, 0, -1 \}$ & \cr
                  & & & $\{ 0, 0, 1, -1, 0, 0 \}$ & \cr
                  \cline{2-4}
       $D_6(a_2)$ & $L_{13}$
                  & $\frac{1}{2}(3 + 2(-a_4 - 2b_2) +  q) < 2$, or
                  & $\{ 0, 1, 0, 0, 0, 0 \}$ & $\frac{3}{2}$ \cr
                  &
                  & $(a_4 + 2b_2 - 1)^2 + a_4^2 < 2$
                  & $\{ 0, 1, 0, 0, 0, 1 \}$ &  \cr
                  & & & $\{ 0, 1, 0, -1, 0, 1 \}$ & \cr
                  & & & $\{ 0, 1, 0,  1, 0, 0 \}$ & \cr
                  \cline{2-5}
                  & $L_{23}$
                  & $\frac{1}{2}(2 +  q) < 2$, or $(a_4 + 2b_2)^2 + a_4^2 < 1$
                  & $\{ 1, 0, 0, 0, 0, 0 \}$ & $1$ \\
     \hline\hline  %%% D_7(a_1)
                  &  &   $q = 4a_4^2 + 4b_2^2 + 8b_3^2 + 8a_4{b_3}$ & & \cr
                  \cline{3-3}
                  & $L_{12}$
                  & $\frac{1}{4}(7 + 2(2a_4 + 2b_2 + 4b_3) +  q) < 2$, or
                  & $\{ 0, 0, 1, 1, 0, 0, -1 \}$ & \cr
                  &
                  & $(a_4 + b_3 + \frac{1}{2})^2 + (b_2 + \frac{1}{2})^2 + (b_3 + \frac{1}{2})^2 < 1$
                  & $\{ 0, 0, 1, 1, 0, -1, -1 \}$ &  \cr
                  & & & $\{ 0, 0, 1, 0, 0, 0, -1 \}$ & \cr
                  & & & $\{ 0, 0, 1, 0, 0, -1, 0 \}$ & \cr
                  & & & $\{ 0, 0, 1, 0, 0, 0,  0 \}$ & \cr
                  & & & $\{ 0, 0, 1, 0, 0, -1, -1 \}$ & \cr
                  & & & $\{ 0, 0, 1, -1, 0, 0, 0 \}$ & \cr
                  & & & $\{ 0, 0, 1, -1, 0, -1, 0 \}$ & \cr
                  \cline{2-4}
         $D_7(a_1)$ & $L_{13}$
                  & $\frac{1}{4}(7 + 2(2a_4 - 2b_2 + 4b_3) +  q) < 2$, or
                  & $\{ 0, 1, 0, 1, 0, 0, -1 \}$ & $\frac{7}{4}$ \cr
                  &
                  & $(a_4 + b_3 + \frac{1}{2})^2 + (b_2 - \frac{1}{2})^2 + (b_3 + \frac{1}{2})^2 < 1$
                  & $\{ 0, 1, 0, 1, 0,  1, -1 \}$ &  \cr
                  & & & $\{ 0, 1, 0, 0, 0, 0, -1 \}$ & \cr
                  & & & $\{ 0, 1, 0, 0, 0, 1, 0 \}$ & \cr
                  & & & $\{ 0, 1, 0, 0, 0, 0,  0 \}$ & \cr
                  & & & $\{ 0, 1, 0, 0, 0, 1, -1 \}$ & \cr
                  & & & $\{ 0, 1, 0, -1, 0, 0, 0 \}$ & \cr
                  & & & $\{ 0, 1, 0, -1, 0, 1, 0 \}$ & \cr
                  \cline{2-5}
                  & $L_{23}$
                  & $\frac{1}{4}(12 + 2(4a_4 + 8b_3) +  q) < 2$, or
                  & $\{ 1, 0, 0, 0, 0, 0, -1 \}$ & $1$ \cr
                  & & $(a_4 + b_3 + 1)^2 + (b_3 + 1)^2 + b_2^2 < 1$
                  & &  \\
    \hline\hline
    \end{tabular}
  \vspace{2mm}
  \caption{\hspace{3mm}(cont.) The linkage diagrams $\gamma^{\nabla}_{ij}(8)$
   obtained as solutions of inequality $\mathscr{B}^{\vee}_{\Gamma}(\gamma^{\nabla}_{ij}(8)) < 2$}
  \label{sol_inequal_3}
  \end{table}

  \begin{table}[H]
  \centering
  \renewcommand{\arraystretch}{1.7}
  \small
  \begin{tabular} {||c|c|c|l|l||}
  \hline \hline
      The Carter &  $L_{ij}$  & $\mathscr{B}^{\vee}_{\Gamma}(\gamma^{\nabla}_{ij}(8)) = p < 2$
                                & The linkage diagram  & $p$ \cr
      diagram &    &            &  \qquad \qquad $\gamma^{\nabla}_{ij}(8)$  &  \\
     \hline\hline  %%% D_7(a_2)
                  &  &   $q = 4a_4^2 + 8b_2^2 + 4b_3^2 + 8a_4{b_2}$ & & \cr
                  \cline{3-3}
                  & $L_{12}$
                  & $\frac{1}{4}(7 + 2(2a_4 + 4b_2 + 2b_3) +  q) < 2$, or
                  & $\{ 0, 0, 1, 1, 0, -1, 0 \}$ & \cr
                  &
                  & $(a_4 + b_3 + \frac{1}{2})^2 + (b_2 + \frac{1}{2})^2 + (b_3 + \frac{1}{2})^2 < 1$
                  &     $\{ 0, 0, 1, 1, 0, -1, -1 \}$ &  \cr
                  & & & $\{ 0, 0, 1, 0, 0, -1, 0 \}$ & \cr
                  & & & $\{ 0, 0, 1, 0, 0, 0, -1 \}$ & \cr
                  & & & $\{ 0, 0, 1, 0, 0,  0, 0 \}$ & \cr
                  & & & $\{ 0, 0, 1, 0, 0, -1, -1 \}$ & \cr
                  & & & $\{ 0, 0, 1, -1, 0, 0,  0 \}$ & \cr
                  & & & $\{ 0, 0, 1, -1, 0, 0, -1 \}$ & \cr
                  \cline{2-4}
         $D_7(a_2)$ & $L_{13}$
                 & $\frac{1}{4}(7 + 2(-2a_4 - 4b_2 + 2b_3) +  q) < 2$, or
                  & $\{ 0, 1, 0, 0, 0, 0, 0 \}$ & $\frac{7}{4}$ \cr
                  &
                  & $(a_4 + b_2 - \frac{1}{2})^2 + (b_2 - \frac{1}{2})^2 + (b_3 + \frac{1}{2})^2 < 1$
                  & $\{ 0, 1, 0, 1, 0,  0, 0 \}$ &  \cr
                  & & & $\{ 0, 1, 0, 0, 0, 1, 0 \}$ & \cr
                  & & & $\{ 0, 1, 0, -1, 0, 1,  0 \}$ & \cr
                  & & & $\{ 0, 1, 0, 0, 0, 0, -1 \}$ & \cr
                  & & & $\{ 0, 1, 0, 1, 0, 0, -1 \}$ & \cr
                  & & & $\{ 0, 1, 0, 0, 0, 1, -1 \}$ & \cr
                  & & & $\{ 0, 1, 0, -1, 0, 1, -1 \}$ & \cr
                  \cline{2-5}
                  & $L_{23}$
                  & $\frac{1}{4}(8 + 8b_3) +  q) < 2$, or
                  & $\{ 1, 0, 0, 0, 0, 0, -1 \}$ & $1$ \cr
                  & & $(a_4 + b_2)^2 + (b_3 + 1)^2 + b_2^2 < 1$
                  & &  \\
    \hline\hline
    \end{tabular}
  \vspace{2mm}
  \caption{\hspace{3mm}(cont.) The linkage diagrams $\gamma^{\nabla}_{ij}(8)$
  obtained as solutions of inequality $\mathscr{B}^{\vee}_{\Gamma}(\gamma^{\nabla}_{ij}(8)) < 2$}
  \label{sol_inequal_4}
  \end{table}

  %% ============================ E6, E7, D5, D6 ==================

  \begin{table}[H]
  \centering
  \renewcommand{\arraystretch}{1.4}
  \small
  \begin{tabular} {||c|c|c|l|l||}
  \hline \hline
     The Carter &  $L_{ij}$  & $\mathscr{B}^{\vee}_{\Gamma}(\gamma^{\nabla}_{ij}(8)) = p < 2$
                                & The linkage diagram  & $p$ \cr
      diagram &    &            &  \qquad \qquad $\gamma^{\nabla}_{ij}(8)$  &  \\
    \hline  %%% E_6
                  &  &   $q = 4b_2^2 + 4b_2{b_3} + 4b_3^2$ & & \cr
                  \cline{3-3}
                  &  $L_{12}$
                  &  $\frac{1}{3}(6 + 2(3b_2 + 3b_3) + q) < 2$, or
                  &  $\{ 0, 0, 1, 0, 0, -1 \}$  &  \cr
                  &
                  &  $(b_2 + b_3)^2 + (b_2 + \frac{3}{2})^2  + (b_3 + \frac{3}{2})^2 < \frac{9}{2}$
                  &  $\{ 0, 0, 1, 0, -1, 0 \}$ & \cr
                  \cline{2-4}
    $E_6$         &  $L_{13}$
                  &  $\frac{1}{3}(10 + 2(4b_2 + 5b_3) + q) < 2$, or
                  &  $\{ 0, 1, 0, 0, 0, -1 \}$  & $\frac{4}{3}$ \cr
                  &
                  &  $(b_2 + b_3)^2 + (b_2 + 2)^2  + (b_3 + \frac{5}{2})^2 < \frac{33}{4}$
                  &  $\{ 0, 1, 0, 0,  -1, -1 \}$ &   \cr
                  \cline{2-4}
                  &  $L_{23}$
                  &  $\frac{1}{3}(10 + 2(5b_2 + 4b_3) + q) < 2$,
                  &  $\{ 1, 0, 0, 0, -1, 0 \}$  &  \cr
                  &
                  &  $(b_2 + b_3)^2 + (b_2 + \frac{5}{2})^2  + (b_3 + 2)^2 < \frac{33}{4}$
                  & $\{ 1, 0, 0, 0,  -1, -1 \}$ & \\
    \hline\hline   %%% E_7
                  &  &   $q = 3a_4^2 + 4b_2^2 + 8b_3^2 + 4a_4{b}_2 + 8b_3{a}_4 + 8b_2{b_3}$ & & \cr
                  \cline{3-3}
                  &   $L_{12}$
                  &  $\frac{1}{2}(7 + 2(3a_4 + 4b_2 + 6b_3) + q) < 2$, or
                  &  $\{ 0, 0, 1, 0, 0, 0, -1 \}$  &  \cr
                  &
                  & $(a_4 + 2b_2 + 2b_3 + 2)^2 + (a_4 + 2b_3 + 1)^2  + a_4^2 < 2$
                  &  $\{ 0, 0, 1, 0, 0, -1, 0 \}$ & \cr
                  \cline{2-4}
       $E_7$ &   $L_{13}$
                  &   $\frac{1}{2}(15 + 2(5a_4 + 6b_2 + 10b_3) + q < 2$, or
                  &  $\{ 0, 1, 0, 0, 0, -1, -1 \}$  & $\frac{3}{2}$ \cr
                  &
                  & $(a_4 + 2b_2 + 2b_3 + 3)^2 + (a_4 + 2b_3 + 2)^2  + a_4^2 < 2$
                  &  $\{ 0, 1, 0, 0, 0, 0, -1 \}$ &  \cr
                  \cline{2-4}
                  &   $L_{23}$
                  & $\frac{1}{3}(12 + 2(4a_4 + 6b_2 + 8b_3) + q) < 2$, or
                  & $\{ 1, 0, 0, 1, 0, -1, -1 \}$ & \cr
                  &
                  &  $(a_4 + 2b_2 + 2b_3 + 3)^2 + (a_4 + 2b_3 + 1)^2  + a_4^2 < 2$
                  & $\{ 1, 0, 0, -1,  0, -1, 0 \}$ &  \\
      \hline\hline   %%% D_5
                  &  &   $q = 4b_2^2$ & & \cr
                  \cline{3-3}
                  &   $L_{12}$
                  &   $\frac{1}{4}(5 + 2b_2 + (2 + 4b_2)b_2) < 2$, or
                  & $\{ 0, 0, 1, 0, -1 \}$ & $\frac{5}{4}$ \cr
                  &
                  & $(2b_2 + 1)^2 < 4$
                  & $\{ 0, 0, 1, 0, 0\}$ & \cr
                  \cline{2-5}
     $D_5$        & $L_{13}$
                  & $\frac{1}{4}(8 + 4b_2 + (4 + 4b_2)b_2) < 2$, or
                  & $\{ 0, 1, 0, 0, -1 \}$ &  $1$\cr
                  &
                  & $(b_2 + 1)^2 < 1$
                  &  & \cr
                  \cline{2-5}
                  & $L_{23}$
                  & $\frac{1}{4}(5 + 2b_2 + (2 + 4b_2)b_2) <2$, or
                  & $\{ 1, 0, 0, 0, -1 \}$ & $\frac{5}{4}$ \cr
                  &
                  & $(2b_2 + 1)^2 < 4$
                  & $\{ 1, 0, 0, 0, 0 \}$ & \\
      \hline\hline   %%% D6
                  &  &   $q = 4b_2^2 + 4b_2{a}_4 + 2b_3b_4$ & & \cr
                  \cline{3-3}
                  & $L_{12}$
                  & $\frac{1}{2}(3 + 2(a_4 + 2b_2) +  q) < 2$, or
                  & $\{ 0, 0, 1, 0, 0, 0 \}$ & $\frac{3}{2}$ \cr
                  & & & $\{ 0, 0, 1, 0, 0, -1 \}$ & \cr
                  &
                  & $(2b_2 + a_4 + 1)^2 + a_4^2 < 2$
                  & $\{ 0, 0, 1, -1, 0, 0 \}$ & \cr
                  & & & $\{ 0, 0, 1, 1, 0, -1 \}$ & \cr
                  \cline{2-5}
       $D_6$ & $L_{13}$
                  & $\frac{1}{2}(6 + 2(2a_4 + 4b_2) + q) < 2$, or
                  & $\{ 0, 1, 0, 0, -1 \}$ &  $1$\cr
                  &
                  & $(2b_2 + a_4 + 2)^2 + a_4^2 < 2$
                  & $\{ 0, 1, 0, 0, -1 \}$ & \cr
                  \cline{2-5}
                  & $L_{23}$
                  & $\frac{1}{2}(3 + 2(a_4 + 2b_2) +  q) < 2$, or
                  & $\{ 1, 0, 0, 0, 0, 0 \}$ & $\frac{3}{2}$ \cr
                  & & & $\{ 1, 0, 0, -1, 0, 0 \}$ & \cr
                  &
                  & $(2b_2 + a_4 + 1)^2 + a_4^2 < 2$
                  & $\{ 1, 0, 0, 0, 0, -1 \}$ & \cr
                  & & & $\{ 1, 0, 0, 1, 0, -1 \}$ & \cr
       \hline\hline   %%%
\end{tabular}
  \vspace{2mm}
  \caption{\hspace{3mm}(cont.) The linkage diagrams $\gamma^{\nabla}_{ij}(8)$ obtained
  as solutions of inequality $\mathscr{B}^{\vee}_{\Gamma}(\gamma^{\nabla}_{ij}(8)) < 2$}
  \label{sol_inequal_5}
  \end{table}

  \begin{table}[H]
  \centering
  \renewcommand{\arraystretch}{1.7}
  \small
  \begin{tabular} {||c|c|c|l|l||}
  \hline \hline
      The Carter &  $L_{ij}$  & $\mathscr{B}^{\vee}_{\Gamma}(\gamma^{\nabla}_{ij}(8)) = p < 2$
                                & The linkage diagram  & $p$ \cr
      diagram &    &            &  \qquad \qquad $\gamma^{\nabla}_{ij}(8)$  &  \\
     \hline\hline  %%% D_7
                  &  &   $q = 8a_4^2 + 12b_2^2 + 4b_3^2 + 16b_2{a}_4 + 8b_2{b}_3 + 8b_3{a}_4$ & & \cr
                  \cline{3-3}
                  & $L_{12}$
                  & $\frac{1}{4}(7 + 4(2a_4 + 3b_2 + b_3) +  q) < 2$, or
                  & $\{ 0, 0, 1, 0, 0, 0, 0 \}$ & \cr
                  &
                  & $(a_4 + b_2 + b_3 + \frac{1}{2})^2 + (a_4 + b_2 + \frac{1}{2})^2 + (b_2 + \frac{1}{2})^2 < 1$
                  &     $\{ 0, 0, 1, 0, 0, -1, 1 \}$ &  \cr
                  & & & $\{ 0, 0, 1, 1, 0, -1, -1 \}$ & \cr
                  & & & $\{ 0, 0, 1, -1, 0, 0, 0 \}$ & \cr
                  %%----------------------------------------
                  & & & $\{ 0, 0, 1, 1, 0,  -1, 0 \}$ & $\frac{7}{4}$ \cr
                  & & & $\{ 0, 0, 1, -1, 0, 0, 1 \}$ & \cr
                  & & & $\{ 0, 0, 1, 0, 0, 0,  -1 \}$ & \cr
                  & & & $\{ 0, 0, 1, 0, 0, -1, 0 \}$ & \cr
                  \cline{2-5}
         $D_7$   & $L_{13}$
                 & $\frac{1}{4}(16 + 2(8a_4 +12b_2 + 4b_3) +  q) < 2$, or
                  & $\{ 0, 1, 0, 0, 0, -1, 0 \}$ & $1$ \cr
                  &
                  & $(a_4 + b_2 + b_3 + 1)^2 + (a_4 + b_2 + 1)^2 + (b_2 + 1)^2 < 1$ & & \cr
                  \cline{2-5}
                  & $L_{23}$
                  & $\frac{1}{4}(7 + 4(2a_4 + 3b_2 + b_3) +  q) < 2$, or
                      & $\{ 1, 0, 0, 0, 0, 0, 0 \}$ & $1$ \cr
                  &
                  & $(a_4 + b_2 + b_3 + \frac{1}{2})^2 + (a_4 + b_2 + \frac{1}{2})^2 + (b_2 + \frac{1}{2})^2 < 1$
                  &     $\{ 1, 0, 0, 0, 0, -1, 1 \}$ &  \cr
                  & & & $\{ 1, 0, 0, 1, 0, -1, -1 \}$ & \cr
                  & & & $\{ 1, 0, 0, -1, 0, 0, 0 \}$ & \cr
                  %%----------------------------------------
                  & & & $\{ 1, 0, 0, 1, 0,  -1, 0 \}$ & $\frac{7}{4}$ \cr
                  & & & $\{ 1, 0, 0, -1, 0, 0, 1 \}$ & \cr
                  & & & $\{ 1, 0, 0, 0, 0, 0,  -1 \}$ & \cr
                  & & & $\{ 1, 0, 0, 0, 0, -1, 0 \}$ & \cr
    \hline\hline
    \end{tabular}
  \vspace{2mm}
  \caption{\hspace{3mm}(cont.) The linkage diagrams $\gamma^{\nabla}_{ij}(8)$
  obtained as solutions of inequality $\mathscr{B}^{\vee}_{\Gamma}(\gamma^{\nabla}_{ij}(8)) < 2$}
  \label{sol_inequal_6}
  \end{table}

\newpage
\subsection{$\beta$-unicolored linkage diagrams.
 Solutions of inequality $\mathscr{B}^{\vee}_{\Gamma}(\gamma^{\nabla}) < 2$}
   \index{linkage diagrams ! - $\beta$-unicolored}
~\\

\begin{table}[H]
  \centering
  \renewcommand{\arraystretch}{1.6}
  \small
  \begin{tabular} {||c|l|l|l||}
  \hline \hline
      Diagram & \multicolumn{2}{c|}{$\mathscr{B}^{\vee}_{\Gamma}(\gamma^{\nabla}) = p < 2$,}  & $p$ \cr
              & \multicolumn{2}{c|}{$\gamma^{\nabla}$ -  $\beta$-unicolored linkage diagram} &  \\
    \hline\hline  %%% E_6(a_1)
                  &  \multicolumn{2}{c|}{$\frac{1}{3}(4b_2^2 + 4b_3^2 + 4b_2{b_3}) < 2$, or} &  \cr
    $E_6(a_1)$, $E_6(a_2)$
                  &  \multicolumn{2}{c|}{$(b_2 + b_3)^2 + b_2^2  + b_3^2 < 3$} & \cr
    \cline{2-3}
                  &  $\{ 0, 0, 0, 0, 0, 1 \}$ \quad\quad\quad\quad\quad\quad & $\{ 0, 0, 0, 0, 0, -1 \}$ & $\frac{4}{3}$ \cr
                  &  $\{ 0, 0, 0, 0, 1, 0 \}$  & $\{ 0, 0, 0, 0, -1, 0 \}$ &  \cr
                  &  $\{ 0, 0, 0, 0, 1, -1 \}$  & $\{ 0, 0, 0, 0, -1, 1 \}$ & \cr
      \hline
                  &  \multicolumn{2}{c|}
      {$\frac{1}{2}(3b_2^2 + 4b_3^2 + 3b_4^2 + 4b_2{b_3} + 2b_2{b_4} + 4b_3{b_4}) < 2$, or} &  \cr
      $E_7(a_1)$  &  \multicolumn{2}{c|}{$2(b_2 + b_3)^2 + 2(b_3 + b_4)^2 + (b_2 + b_4)^2 < 4$} & \cr
    \cline{2-3}
      & $\{ 0, 0, 0, 0, 1, -1, 0 \}$  &  $\{ 0, 0, 0, 0, -1, 1, 0 \}$ & $\frac{3}{2}$ \cr
      & $\{ 0, 0, 0, 0, 0, 1, -1 \}$ &  $\{ 0, 0, 0, 0, 0, -1, 1 \}$ &  \cr
      & $\{ 0, 0, 0, 0, 0, 0, 1 \}$  &  $\{ 0, 0, 0, 0, 0, 0, -1 \}$ & \cr
      & $\{ 0, 0, 0, 0, 1, 0, 0 \}$  &  $\{ 0, 0, 0, 0, -1, 0, 0 \}$ & \cr
      \hline
                  &  \multicolumn{2}{c|}
      {$\frac{1}{2}(3b_2^2 + 3b_3^2 + 3b_4^2 + 2b_2{b_3} - 2b_2{b_4} + 2b_3{b_4}) < 2$, or} &  \cr
      $E_7(a_2)$  &  \multicolumn{2}{c|}
      {$(b_2 + b_3)^2 + (b_2 - b_4)^2 + (b_3 + b_4)^2 + b_2^2  + b_3^2 + b_4^2 < 4$} & \cr
    \cline{2-3}
      & $\{ 0, 0, 0, 0, 1, 0, 0 \}$  &  $\{ 0, 0, 0, 0, -1, 0, 0 \}$ & $\frac{3}{2}$  \cr
      & $\{ 0, 0, 0, 0, 0, 1, 0 \}$  &  $\{ 0, 0, 0, 0, 0, -1, 0 \}$ & \cr
      & $\{ 0, 0, 0, 0, 0, 0, 1 \}$  &  $\{ 0, 0, 0, 0, 0, 0, -1 \}$ & \cr
      & $\{ 0, 0, 0, 0, 1, -1, 1 \}$  & $\{ 0, 0, 0, 0, -1, 1, -1 \}$ & \cr
      \hline
                 &  \multicolumn{2}{c|}
      {$\frac{1}{2}(4b_2^2 + 3b_3^2 + 3b_4^2 + 4b_2{b_3} - 4b_2{b_4} - 2b_3{b_4}) < 2$, or} &  \cr
      $E_7(a_3)$  &  \multicolumn{2}{c|}
      {$2(b_2 + b_3)^2 + 2(b_2 - b_4)^2 + (b_3 - b_4)^2  < 4$} & \cr
    \cline{2-3}
      & $\{ 0, 0, 0, 0, 0, 1, 0 \}$  & $\{ 0, 0, 0, 0, 0, -1, 0 \}$ & $\frac{3}{2}$ \cr
      & $\{ 0, 0, 0, 0, 0, 0, 1 \}$  & $\{ 0, 0, 0, 0, 0, 0, -1 \}$ & \cr
      & $\{ 0, 0, 0, 0, 1, 0, 1 \}$  & $\{ 0, 0, 0, 0, -1, 0, -1 \}$ & \cr
      & $\{ 0, 0, 0, 0, 1, -1, 0 \}$  & $\{ 0, 0, 0, 0, -1, 1, 0 \}$ & \cr
      \hline
        &  \multicolumn{2}{c|}
      {$\frac{1}{2}(3b_2^2 + 3b_3^2 + 3b_4^2 -2b_2{b_3} - 2b_2{b_4} - 2b_3{b_4}) < 2$, or} &  \cr
      $E_7(a_4)$  &  \multicolumn{2}{c|}
      {$(b_2 - b_3)^2 + (b_2 - b_4)^2 + (b_3 - b_4)^2 + b_2^2 + b_3^3 + b_4^4  < 4$} & \cr
    \cline{2-3}
      & $\{ 0, 0, 0, 0, 1, 0, 0 \}$  & $\{ 0, 0, 0, 0, -1, 0, 0 \}$ & $\frac{3}{2}$ \cr
      & $\{ 0, 0, 0, 0, 0, 1, 0 \}$  & $\{ 0, 0, 0, 0, 0, -1, 0 \}$ & \cr
      & $\{ 0, 0, 0, 0, 0, 0, 1 \}$  & $\{ 0, 0, 0, 0, 0, 0, -1 \}$ & \cr
      & $\{ 0, 0, 0, 0, 1, 1, 1 \}$  & $\{ 0, 0, 0, 0, -1, -1, -1 \}$ & \cr
      \hline\hline
    \end{tabular}
  \vspace{2mm}
  \caption{\hspace{3mm}(cont.) $\beta$-unicolored linkage diagrams obtained as solutions of inequality
    $\mathscr{B}^{\vee}_{\Gamma}(\gamma^{\nabla}) < 2$}
  \label{homog_inequal_1}
  \end{table}

\begin{table}[H]
  \centering
  \renewcommand{\arraystretch}{1.7}
  \small
  \begin{tabular} {||c|l|l|l||}
  \hline \hline
      Diagram & \multicolumn{2}{c|}{$\mathscr{B}^{\vee}_{\Gamma}(\gamma^{\nabla}) = p < 2$,}  & $p$ \cr
              & \multicolumn{2}{c|}{$\gamma^{\nabla}$ -  $\beta$-unicolored linkage diagram} &  \\
      \hline\hline
       $D_5(a_1)$ &  \multicolumn{2}{c|}
      {$\frac{1}{4}(4b_2^2) < 2 ,\quad $ or $\quad b_2^2  < 2$} &  \cr
    \cline{2-3}
      & $\{ 0, 0, 0, 0, 1\}$  \quad\quad\quad\quad & $\{ 0, 0, 0, 0, -1\}$ & $1$ \cr
      \hline
       $D_6(a_1)$ &  \multicolumn{2}{c|}
      {$\frac{1}{2}(2b_2^2 + 2b_3^2) < 2 ,\quad $ or $\quad b_2^2 + b_3^2  < 2$} &  \cr
    \cline{2-3}
      & $\{ 0, 0, 0, 0, 1, 0 \}$   & $\{ 0, 0, 0, 0, -1, 0 \}$ & $1$ \cr
      & $\{ 0, 0, 0, 0, 0, 1 \}$  & $\{ 0, 0, 0, 0, 0, -1 \}$ & \cr
      \hline
        &  \multicolumn{2}{c|}
      {$\frac{1}{2}(2a_4^2 + 4b_1^2 + 4b_2^2 + 4a_4{b_2}) < 2,$ or} &  \cr
       $D_6(a_2)$  &  \multicolumn{2}{c|}{$\quad 2b_1^2 + b_2^2 + (a_4 + b_2)^2 < 2$} & \cr
    \cline{2-3}
      & $\{ 0, 0, 0, 1, 0, 0 \}$   & $\{ 0, 0, 0, -1, 0, 0 \}$ & $1$ \cr
      & $\{ 0, 0, 0, 1, 0, -1 \}$  & $\{ 0, 0, 0, -1, 0, 1 \}$ & \cr
      \hline
        &  \multicolumn{2}{c|}
      {$\frac{1}{4}(4a_4^2 + 4b_2^2 + 8b_3^2 + 8a_4{b_3}) < 2,$ or} &  \cr
       $D_7(a_1)$  &  \multicolumn{2}{c|}
       {$\quad b_2^2 + b_3^2 + (a_4 + b_3)^2  < 2$} & \cr
    \cline{2-3}
      & $\{ 0, 0, 0, 1, 0, 0, 0 \}$  & $\{ 0, 0, 0, -1, 0, 0, 0 \}$ & $1$ \cr
      & $\{ 0, 0, 0, 0, 0, 1, 0 \}$  & $\{ 0, 0, 0,  0, 0, -1, 0 \}$ & \cr
      & $\{ 0, 0, 0, 1, 0, 0, -1 \}$  & $\{ 0, 0, 0, -1, 0, 0, 1 \}$ & \cr
      \hline
        &  \multicolumn{2}{c|}
      {$\frac{1}{4}(4a_4^2 + 8b_2^2 + 4b_3^2 + 8a_4{b_2}) < 2,$ or} &  \cr
       $D_7(a_2)$  &  \multicolumn{2}{c|}
       {$\quad b_2^2 + b_3^2 + (a_4 + b_2)^2  < 2$} & \cr
    \cline{2-3}
      & $\{ 0, 0, 0, 1, 0, 0, 0 \}$  & $\{ 0, 0, 0, -1, 0, 0, 0 \}$ & $1$ \cr
      & $\{ 0, 0, 0, 0, 0, 0, 1 \}$  & $\{ 0, 0, 0,  0, 0, 0, -1 \}$ & \cr
      & $\{ 0, 0, 0, 1, 0, -1, 0 \}$  & $\{ 0, 0, 0, -1, 0, 1, 0 \}$ & \cr
      \hline\hline
         \end{tabular}
  \vspace{2mm}
  \caption{\hspace{3mm}(cont.) $\beta$-unicolored linkage diagrams obtained as solutions of inequality
    $\mathscr{B}^{\vee}_{\Gamma}(\gamma^{\nabla}) < 2$}
  \label{homog_inequal_2}
  \end{table}

  \begin{table}[H]
  \centering
  \renewcommand{\arraystretch}{1.6}
  \small
  \begin{tabular} {||c|l|l|l||}
  \hline \hline
      Diagram & \multicolumn{2}{c|}{$\mathscr{B}^{\vee}_{\Gamma}(\gamma^{\nabla}) = p < 2$,}  & $p$ \cr
              & \multicolumn{2}{c|}{$\gamma^{\nabla}$ -  $\beta$-unicolored linkage diagram} &  \\
    \hline\hline  %%% E_6
                  &  \multicolumn{2}{c|}{$\frac{1}{3}(4b_2^2 + 4b_3^2 + 4b_2{b_3}) < 2$, or} &  \cr
    $E_6$
                  &  \multicolumn{2}{c|}{$(b_2 + b_3)^2 + b_2^2  + b_3^2 < 3$} & \cr
    \cline{2-3}
    & \multicolumn{2}{c|}{(coincide with $E_6(a_1)$, $E_6(a_2)$)} & \cr
                  &  $\{ 0, 0, 0, 0, 0, 1 \}$ \quad\quad\quad\quad\quad\quad & $\{ 0, 0, 0, 0, 0, -1 \}$ & $\frac{4}{3}$ \cr
                  &  $\{ 0, 0, 0, 0, 1, 0 \}$  & $\{ 0, 0, 0, 0, -1, 0 \}$ &  \cr
                  &  $\{ 0, 0, 0, 0, 1, -1 \}$  & $\{ 0, 0, 0, 0, -1, 1 \}$ & \cr
      \hline
   %%%% E7
                  &  \multicolumn{2}{c|}
      {$\frac{1}{2}(3a_4^2 + 4b_2^2 + 8b_3^2 + 8b_2{b_3} + 4a_4{b_2} + 8a_4{b_3}) < 2$, or} &  \cr
      $E_7$  &  \multicolumn{2}{c|}{$(a_4 + 2b_2 + 2b_3)^2 + (a_4 + 2b_3)^2 + a_4^2 < 4$} & \cr
    \cline{2-3}
      & $\{ 0, 0, 0, 1, 0, 0, 0 \}$  &  $\{ 0, 0, 0, -1, 0, 0, 0 \}$ & $\frac{3}{2}$ \cr
      & $\{ 0, 0, 0, 1, 0, -1, 0 \}$ &  $\{ 0, 0, 0, -1, 0, 1, 0 \}$ &  \cr
      & $\{ 0, 0, 0, 1, 0, 1, -1 \}$  &  $\{ 0, 0, 0, -1, 0, -1, 1 \}$ & \cr
      & $\{ 0, 0, 0, 1, 0, 0, -1 \}$  &  $\{ 0, 0, 0, -1, 0, 0, 1 \}$ & \cr
      \hline
   %%%% D5
     $D_5$   &  \multicolumn{2}{c|}
      {$\frac{1}{4}(4b_2^2) < 2$, or $b_2^2 < 2$} &  \cr
    \cline{2-3}
      & $\{ 0, 0, 0, 0, 1 \}$  &  $\{ 0, 0, 0, 0, -1 \}$ & $1$  \cr
   %%%% D6
      \hline
                 &  \multicolumn{2}{c|}
      {$\frac{1}{2}(2a_4^2 + 4b_2^2 + 4a_4{b}_2) < 2$, or} &  \cr
      $D_6$  &  \multicolumn{2}{c|}
      {$(2b_2 + a_4)^2 + a_4^2  < 2$} & \cr
    \cline{2-3}
      & $\{ 0, 0, 0, 1, 0, 0 \}$  & $\{ 0, 0, 0, 1, 0, -1 \}$ & $1$ \cr
      & $\{ 0, 0, 0, -1, 0, 0 \}$  & $\{ 0, 0, 0, -1, 0, 1 \}$ & \cr
      \hline
                 &  \multicolumn{2}{c|}
      {$\frac{1}{4}(8a_4^2 + 16a_4{b}_2 + 8b_3{a}_4 + 12b_2^2 + 8b_2{b}_3 + 4b_3^2) < 2$, or} &  \cr
      $D_7$  &  \multicolumn{2}{c|}
      {$(a_4 + b_2 + b_3)^2 + (a_4 + b_2)^2 + b_2^2  < 2$} & \cr
    \cline{2-3}
       & $\{ 0, 0, 0, 0, 0, 0, 1 \}$  & $\{ 0, 0, 0, 0, 0, 0, -1 \}$ & $1$ \cr
       & $\{ 0, 0, 0, -1, 0, 0, 1 \}$  & $\{ 0, 0, 0, 1, 0, 0, -1 \}$ & \cr
       & $\{ 0, 0, 0, -1, 0, 1, 0 \}$  & $\{ 0, 0, 0, 1, 0, -1, 0 \}$ & \cr
      \hline\hline
    \end{tabular}
  \vspace{2mm}
  \caption{\hspace{3mm}(cont.) $\beta$-unicolored linkage diagrams obtained as solutions of inequality
    $\mathscr{B}^{\vee}_{\Gamma}(\gamma^{\nabla}) < 2$}
  \label{homog_inequal_3}
  \end{table}

\newpage
\subsection{Linkage diagrams $\gamma^{\nabla}_{ij}(6)$ per loctets and components}
  \label{sec_diagr_per_comp}
  \index{$\alpha$-set, subset of roots corresponding to $w_1$}
~\\

 \index{loctet ! - of type $L_{ij}$}
 \index{$L_{ij}$, loctet of type $(ij) \in \{(12), (13), (23)\}$}
\begin{table}[H]
  \small
  \centering
  \renewcommand{\arraystretch}{1}
  \begin{tabular} {|| c | c || l | l | l ||}
  \hline \hline
    Diagram & Comp. & \multicolumn{3}{c||}{Linkage diagrams $\gamma^{\nabla}_{ij}(6)$ for the loctet of type $L_{ij}$}    \\
            \cline{3-5}
            &           & \multicolumn{1}{c|}{Type $L_{12}: (1, 1, 0, -1, ...)$}
                        & \multicolumn{1}{c|}{Type $L_{13}: (1, 0, 1, -1, ...)$}
                        & \multicolumn{1}{c||}{Type $L_{23}: (0, 1, 1, -1, ...)$} \\
    \hline \hline
    $D_5(a_1)$ & $1$  & $(1, 1, 0, -1, 0) \quad L^a_{12}$
                      & $(1, 0, 1, -1, 0) \quad  L^a_{13}$  & -  \cr
               & $2$  & $(1, 1, 0, -1, 1) \quad  L^b_{12}$
                      & $(1, 0, 1, -1, -1) \quad  L^b_{13} $  & -  \cr
               & $3$  & - & - &  $(0, 1, 1, -1, 0) \quad  L^a_{23}$  \\
   \hline
    $D_6(a_1)$ & $1$  & $(1, 1, 0, -1, 1, -1) \quad L^d_{12}$
                      & $(1, 0, 1, -1, -1, -1) \quad L^d_{13}$  & -  \cr
               &      & $(1, 1, 0, -1, 0, 0) \quad L^a_{12}$
                      & $(1, 0, 1, -1, 0, 0) \quad  L^a_{13}$    & -  \cr
               & $2$  & $(1, 1, 0, -1, 0, -1) \quad L^b_{12}$
                      & $(1, 0, 1, -1, 0, -1) \quad L^b_{13}$   & -  \cr
               &      & $(1, 1, 0, -1, 1, 0) \quad L^c_{12}$
                      & $(1, 0, 1, -1, -1, 0) \quad L^c_{13}$   & -  \cr
               & $3$  & - & - &  $(0, 1, 1, -1, 0, 0) \quad L^a_{23}$  \\
   \hline
    $D_6(a_2)$ & $1$  & $(1, 1, 0, 0, -1, 0) \quad L^a_{12}$
                      & $(1, 0, 1, -1, -1, 0) \quad L^b_{13}$  & -  \cr
               &      & $(1, 1, 0, 0, -1, 0) \quad L^c_{12}$
                      & $(1, 0, 1, 1, -1, -1) \quad L^d_{13}$   & -  \cr
               & $2$  & $(1, 1, 0, 1, -1, 0) \quad L^b_{12}$
                      & $(1, 0, 1, 0, -1, 0) \quad L^a_{13}$  & -  \cr
  $l(\alpha) = 4$     &      & $(1, 1, 0, -1, -1, 1) \quad L^d_{12}$
                      & $(1, 0, 1, 0, -1, -1) \quad L^c_{13}$  & -  \cr
               & $3$  & - & - &  $(0, 1, 1, 0, -1, 0) \quad L^a_{23}$  \\
   \hline
    $E_6(a_1)$ & $1$  & $(1, 1, 0, -1, 0, 0) \quad L^a_{12}$
                      & $(1, 0, 1, -1, 0, -1) \quad L^a_{13}$
                      & $(0, 1, 1, -1, 0, 0) \quad L^a_{23}$  \cr
               & $2$  & $(1, 1, 0, -1, 1, 0) \quad L^b_{12}$
                      & $(1, 0, 1, -1, -1, 0) \quad L^b_{13}$
                      & $(0, 1, 1, -1, 0, -1) \quad L^b_{23}$  \\
   \hline
    $E_6(a_2)$ & $1$  & $(1, 1, 0, -1, 0,  1) \quad L^a_{12}$
                      & $(1, 0, 1, -1, 0, 0) \quad L^a_{13}$
                      & $(0, 1, 1, -1, 0, -1) \quad L^a_{23}$  \cr
               & $2$  & $(1, 1, 0, -1, 1, 0) \quad L^b_{12}$
                      & $(1, 0, 1, -1, -1, 0) \quad L^b_{13}$
                      & $(0, 1, 1, -1, 0, 0) \quad L^b_{23}$  \\
    \hline
    $E_7(a_1)$ & $1$  & $(1, 1, 0, -1, 1, 0, -1) \quad L^b_{12}$
                      & $(1, 0, 1, -1, -1, 0, -1) \quad L^b_{13}$
                      & $(0, 1, 1, -1, 0, -1, 0) \quad L^b_{23}$  \cr
               &      & $(1, 1, 0, -1, 0, 0, 0) \quad L^a_{12}$
                      & $(1, 0, 1, -1, 0, -1, 0) \quad L^a_{13}$
                      & $(0, 1, 1, -1, 0, 0, 0) \quad L^a_{23}$  \\
   \hline
    $E_7(a_2)$ & $1$  & $(1, 1, 0, -1, 0, 0, -1) \quad L^a_{12}$
                      & $(1, 0, 1, -1, 0, -1, 0) \quad L^a_{13}$
                      & $(0, 1, 1, -1, 0, 0, -1) \quad L^a_{23}$  \cr
               &      & $(1, 1, 0, -1, 1, 0, 0) \quad L^b_{12}$
                      & $(1, 0, 1, -1, -1, 0, 0) \quad L^b_{13}$
                      & $(0, 1, 1, -1, 0, -1, 0) \quad L^b_{23}$  \\
   \hline
    $E_7(a_3)$ & $1$  & $(1, 1, 0, -1, 0, 1, -1) \quad L^a_{12}$
                      & $(1, 0, 1, -1, 0, 0, 0) \quad L^a_{13}$
                      & $(0, 1, 1, -1, 0, -1, -1) \quad L^a_{23}$  \cr
               &      & $(1, 1, 0, -1, 1, 0, 0) \quad L^b_{12}$
                      & $(1, 0, 1, -1, -1, 0, 0) \quad L^b_{13}$
                      & $(0, 1, 1, -1, 0, 0, 0) \quad L^b_{23}$  \\
   \hline
    $E_7(a_4)$ & $1$  & $(1, 1, 0, -1, 0, -1, 0) \quad L^a_{12}$
                      & $(1, 0, 1, -1, 0, 0, 1) \quad L^a_{13}$
                      & $(0, 1, 1, -1, 0, 1, 0) \quad L^a_{23}$  \cr
               &      & $(1, 1, 0, -1, 1, 0, 0) \quad L^b_{12}$
                      & $(1, 0, 1, -1, -1, 0, 0) \quad L^b_{13}$
                      & $(0, 1, 1, -1, 0, 0, -1) \quad L^b_{23}$  \\
   \hline
   $D_7(a_1)$  & $1$  & $(1, 1, 0,  1, -1, 0, -1) \quad L^c_{12}$
                      & $(1, 0, 1,  1, -1, 0, -1) \quad L^c_{13}$  & -  \cr
                &     & $(1, 1, 0, -1, -1, 0, 0) \quad L^d_{12}$
                      & $(1, 0, 1, -1, -1, 0, 0) \quad L^d_{13}$  & -  \cr
                &     & $(1, 1, 0,  0, -1, 1, 0) \quad L^e_{12}$
                      & $(1, 0, 1,  0, -1, -1, 0) \quad L^e_{13}$  & -  \cr
                &     & $(1, 1, 0,  0, -1, 1, -1) \quad L^f_{12}$
                      & $(1, 0, 1,  0, -1, -1, -1) \quad L^f_{13}$  & -  \cr
              & $2$  & $(1, 1, 0, 0, -1, 0, 0) \quad L^a_{12}$
                      & $(1, 0, 1, 0, -1, 0, 0) \quad L^a_{13}$  & -  \cr
    $l(\alpha) = 4$ & & $(1, 1, 0, 1, -1, 1 -1) \quad L^g_{12}$
                      & $(1, 0, 1, 1, -1, -1, -1) \quad L^g_{13}$   & -  \cr
                &     & $(1, 1, 0, -1, -1, 1, 0) \quad L^h_{12}$
                      & $(1, 0, 1, -1, -1, -1, 0) \quad L^h_{13}$  & -  \cr
                &     & $(1, 1, 0, 0, -1, 0, -1) \quad L^b_{12}$
                      & $(1, 0, 1, 0, -1, 0, -1) \quad L^b_{13}$  & -  \cr
               & $3$  & - & - &  $(0, 1, 1, 0, -1, 0, 0) \quad L^a_{23}$  \\
   \hline
    $D_7(a_2)$ & $1$  & $(1, 1, 0,  1, -1, 0, -1) \quad L^c_{12}$
                      & $(1, 0, 1,  1, -1, 0, -1) \quad L^b_{13}$  & -  \cr
                &     & $(1, 1, 0,  0, -1, 0, 0) \quad L^a_{12}$
                      & $(1, 0, 1, -1, -1, 0, 0) \quad L^d_{13}$  & -  \cr
                &     & $(1, 1, 0,  0, -1, -1, 0) \quad L^e_{12}$
                      & $(1, 0, 1,  1, -1, -1, 0) \quad L^h_{13}$  & -  \cr
                &     & $(1, 1, 0,  -1, -1, 1, -1) \quad L^g_{12}$
                      & $(1, 0, 1,  0, -1, -1, -1) \quad L^f_{13}$  & -  \cr
               & $2$  & $(1, 1, 0, 0, -1, 0, -1) \quad L^b_{12}$
                      & $(1, 0, 1, -1, -1, 0, -1) \quad L^c_{13}$  & -  \cr
    $l(\alpha) = 4$ & & $(1, 1, 0, 1, -1, 0, 0) \quad L^d_{12}$
                      & $(1, 0, 1, 0, -1, 0, 0) \quad L^a_{13}$   & -  \cr
                &     & $(1, 1, 0, -1, -1, 1, 0) \quad L^h_{12}$
                      & $(1, 0, 1, 0, -1, -1, 0) \quad L^e_{13}$  & -  \cr
                &     & $(1, 1, 0, 0, -1, 1, -1) \quad L^f_{12}$
                      & $(1, 0, 1, 1, -1, -1, -1) \quad L^g_{13}$  & -  \cr
               & $3$  & - & - &  $(0, 1, 1, 0, -1, 0, 0)$   \\
   \hline  \hline  %1
\end{tabular}
  \vspace{2mm}
  \caption{\hspace{3mm}Linkage diagrams $\gamma^{\nabla}_{ij}(6)$ for the Carter diagrams
  from the class $\mathsf{C4}$, $n < 8$.
   For $D_6(a_2)$, $D_7(a_1)$, $D_7(a_2)$, the length of the $\alpha$-set is $4$}
  \label{tab_seed_linkages_6}
  \end{table}

%% file: D_linkageSyst.tex
\section{\sc\bf The linkage systems for the Carter diagrams of rank $l < 8$}
 \label{sect_linkage_diagr}
  \index{linkage system  ! - $\mathscr{L}(D_l)$ and $\mathscr{L}(D_l(a_k))$ for $l < 8$}

~\\
  The linkage systems are similar to the weight systems (= weight diagrams)
  of the irreducible representations of the semisimple Lie algebras.
~\\
~\\
\subsection{The linkage systems
 $\mathscr{L}(D_4(a_1))$, $\mathscr{L}(D_5(a_1))$, $\mathscr{L}(D_6(a_1))$, $\mathscr{L}(D_6(a_2))$}
 \index{linkage system  ! - $\mathscr{L}(D_5)$ and $\mathscr{L}(D_5(a_1))$}
 \index{linkage system  ! - $\mathscr{L}(D_6)$ and $\mathscr{L}(D_6(a_k))$}
~\\

\begin{figure}[H]
\centering
\includegraphics[scale=0.8]{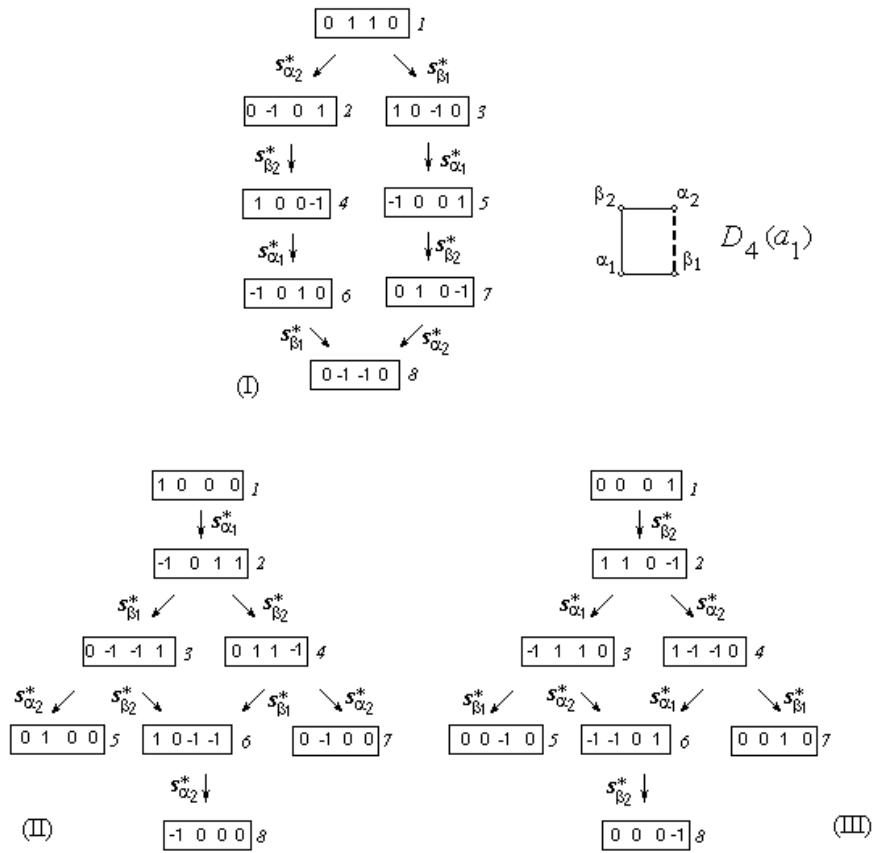}
\caption[\hspace{3mm}Three components of the linkage system $D_4(a_1)$, 3 components]{\hspace{3mm}Three components of the linkage system $D_4(a_1)$. There are $24$ linkage diagrams in the case $D_4(a_1)$}
%%%%%% The label must come after caption
\label{D4a1_linkages}
\end{figure}

\begin{figure}[H]
\centering
\includegraphics[scale=0.62]{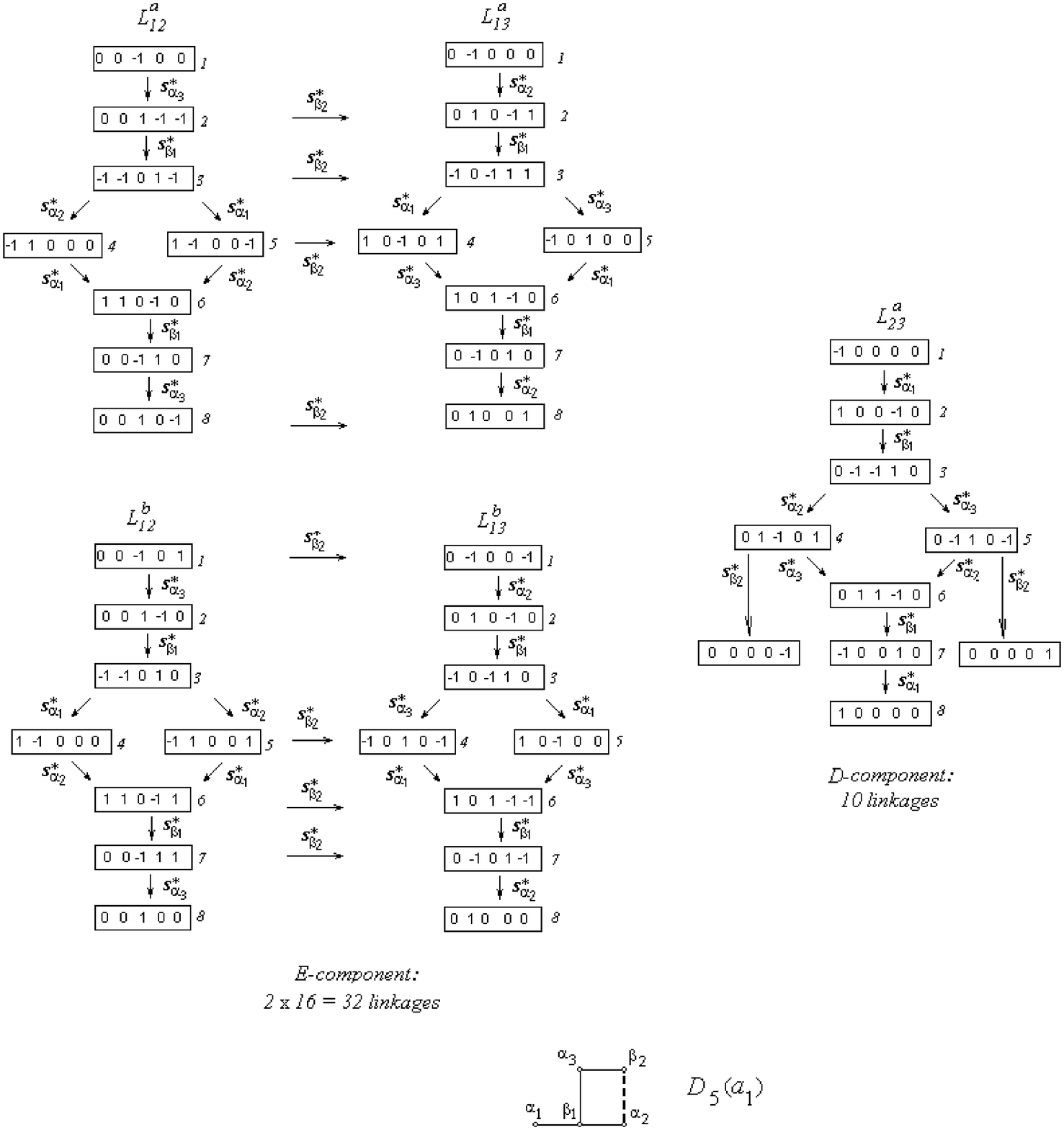}
\vspace{3mm} \caption[\hspace{3mm}The linkage system
 $\mathscr{L}(D_5(a_1))$, $3$ components, $5$ loctets]
 {\hspace{3mm}The linkage system $\mathscr{L}(D_5(a_1))$. There are one part of the $D$-component
 containing $10$ linkage diagrams, and two parts of the $E$-component
 containing $2\times16 = 32$ elements}
%%%%%% The label must come after caption
\label{D5a1_linkages}
\end{figure}

\begin{figure}[H]
\centering
\includegraphics[scale=0.5]{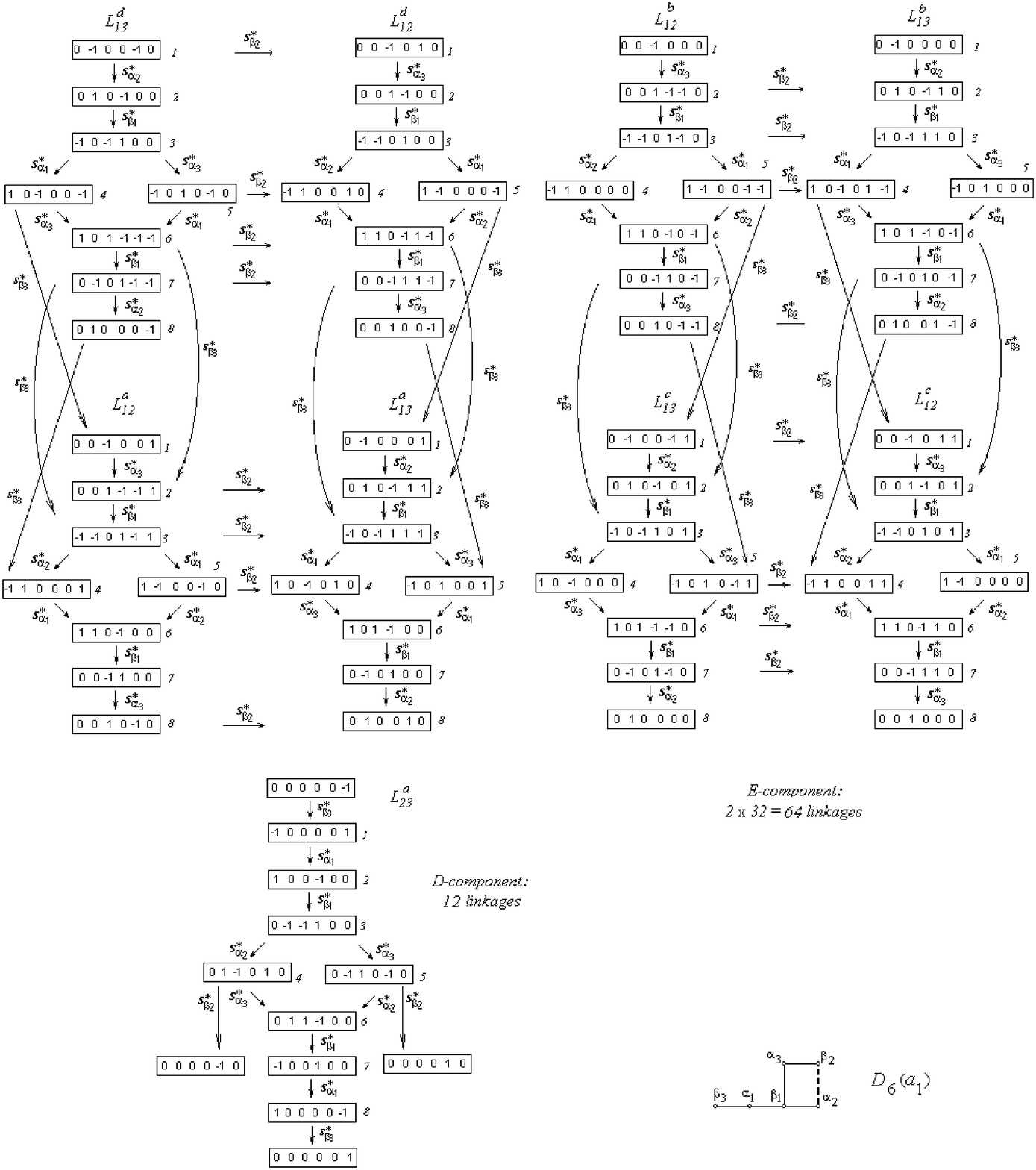}
\vspace{3mm}
 \caption[\hspace{3mm}The linkage system $\mathscr{L}(D_6(a_1))$, $3$ components, $9$ loctets]
 {\hspace{3mm}The linkage system  $\mathscr{L}(D_6(a_1))$.
 There are $12$ linkage diagrams, $1$ loctet in the single part of the $D$-component,
 and $2\times32 = 64$ linkage diagrams, $8$ loctets in two parts of the $E$-component}
%%%%%% The label must come after caption
\label{D6a1_linkages}
\end{figure}

\begin{figure}[H]
\centering
\includegraphics[scale=0.5]{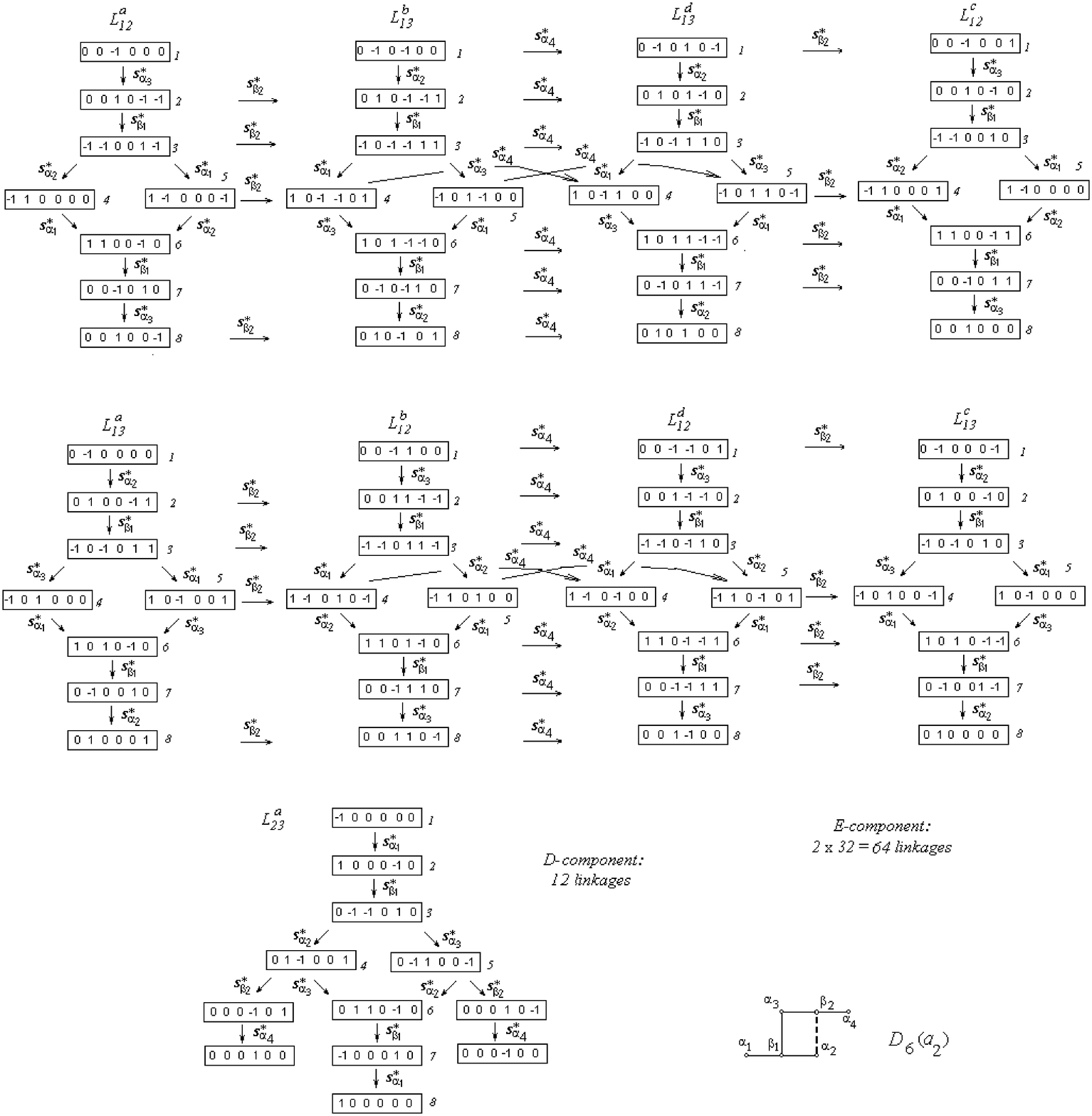}
\vspace{3mm} \caption[\hspace{3mm}The linkage system
 $\mathscr{L}(D_6(a_2))$, $3$ components, $9$ loctets]
 {\hspace{3mm}The linkage system $\mathscr{L}(D_6(a_2))$.
 There are $12$ linkages, $1$ loctet in the single part of the $D$-component,
 and $2\times32 =64$ linkages, $8$ loctets in two parts of the $E$-component}
%%%%%% The label must come after caption
\label{D6a2_linkages}
\end{figure}

\newpage
\subsection{The linkage systems $\mathscr{L}(E_6(a_i))$, $\mathscr{L}(E_6)$,
 $\mathscr{L}(E_7(a_i))$, $\mathscr{L}(D_5)$, $\mathscr{L}(D_6)$}
 \index{linkage system  ! - $\mathscr{L}(D_5)$ and $\mathscr{L}(D_5(a_1))$}
 \index{linkage system  ! - $\mathscr{L}(D_6)$ and $\mathscr{L}(D_6(a_k))$}
 \index{linkage system  ! - $\mathscr{L}(E_l)$ and $\mathscr{L}(E_l(a_k))$}
~\\

\begin{figure}[H]
\centering
\includegraphics[scale=1.4]{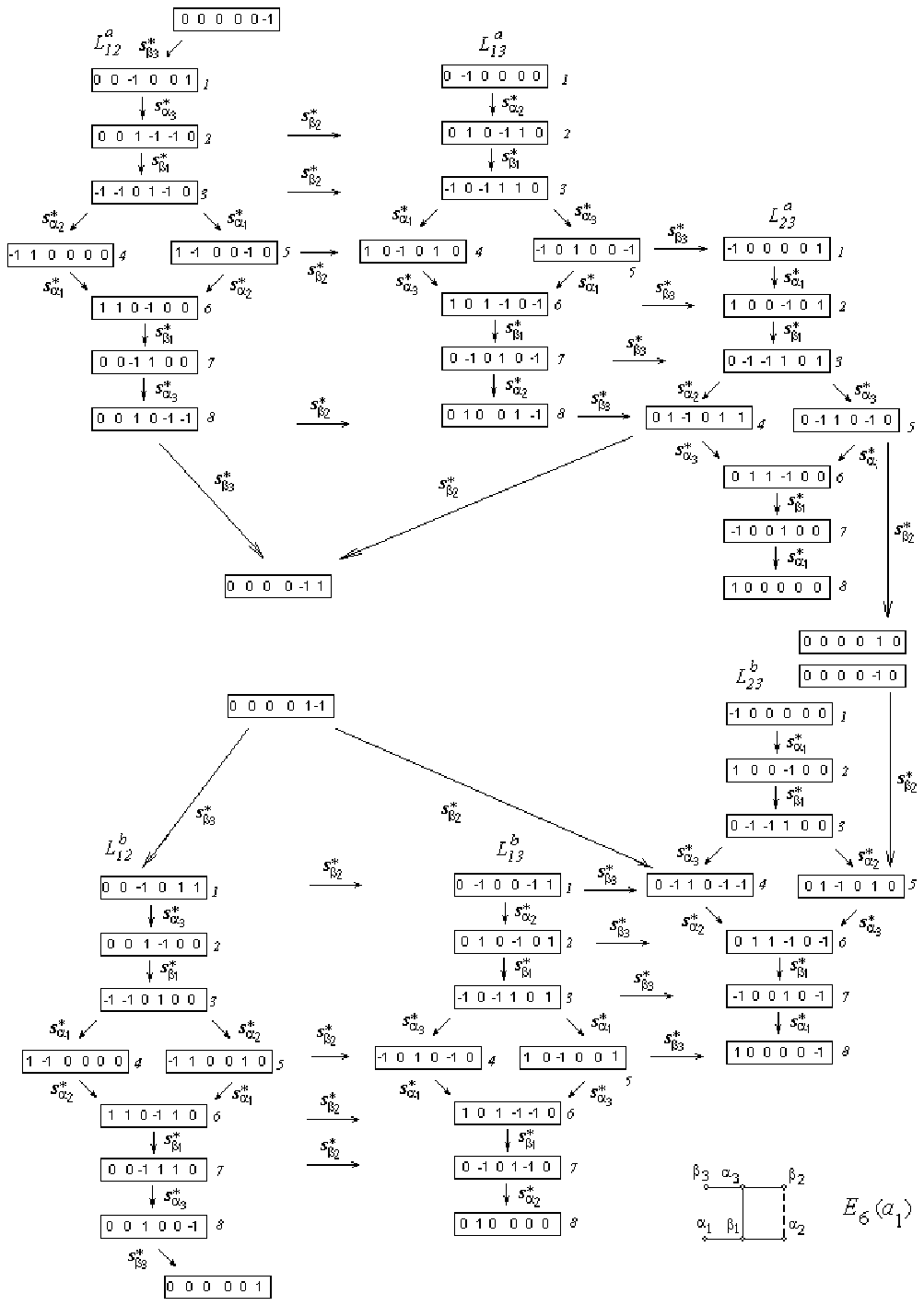}
 \vspace{3mm} \caption{\hspace{3mm}The linkage system $E_6(a_1)$. The only component is the $E$-component
  containing two parts, $54$ linkage diagrams, $6$ loctets}
%%%%%% The label must come after caption
\label{E6a1_linkages}
\end{figure}

\begin{figure}[H]
\centering
\includegraphics[scale=0.62]{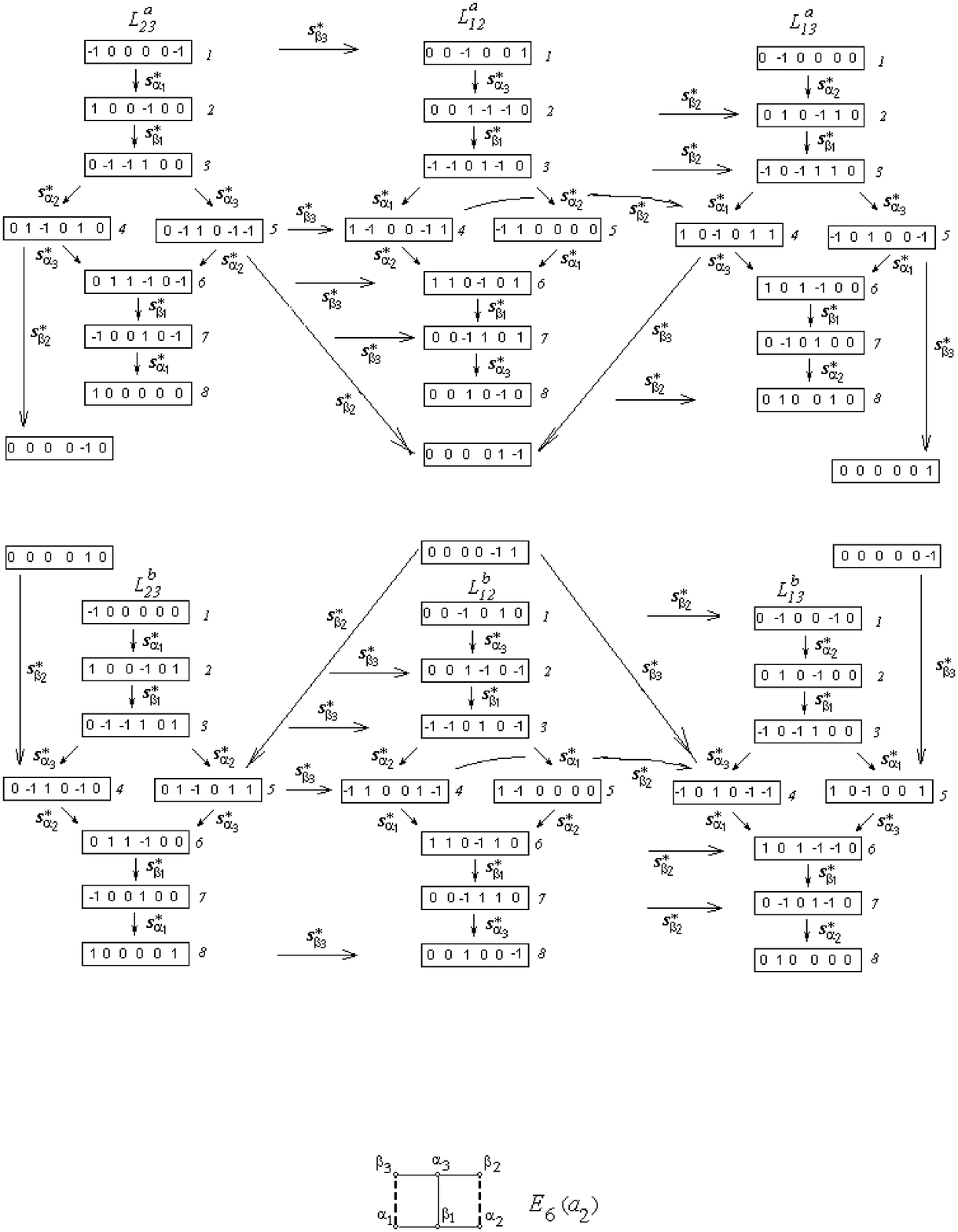}
 \vspace{3mm} \caption{\hspace{3mm}The linkage system $E_6(a_2)$. The only component is the $E$-component
  containing two parts, $54$ linkage diagrams, $6$ loctets}
%%%%%% The label must come after caption
\label{E6a2_linkages}
\end{figure}

\begin{figure}[H]
\centering
\includegraphics[scale=1.4]{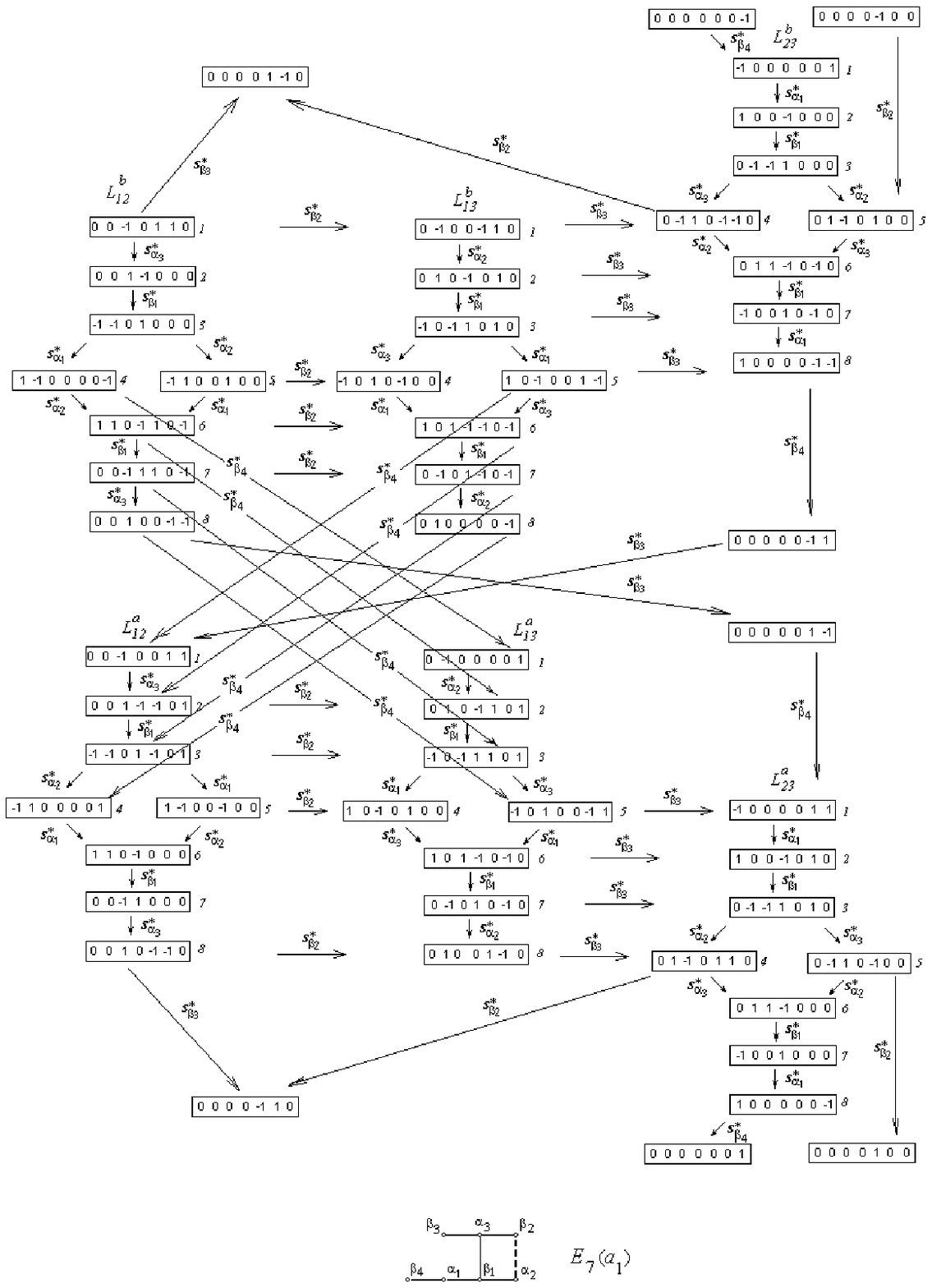}
\caption{\hspace{3mm}The linkage system $E_7(a_1)$. The only component is the $E$-component
 containing $56$ linkage diagrams, $6$ loctets}
%%%%%% The label must come after caption
\label{E7a1_linkages}
\end{figure}

\begin{figure}[H]
\centering
\includegraphics[scale=1.4]{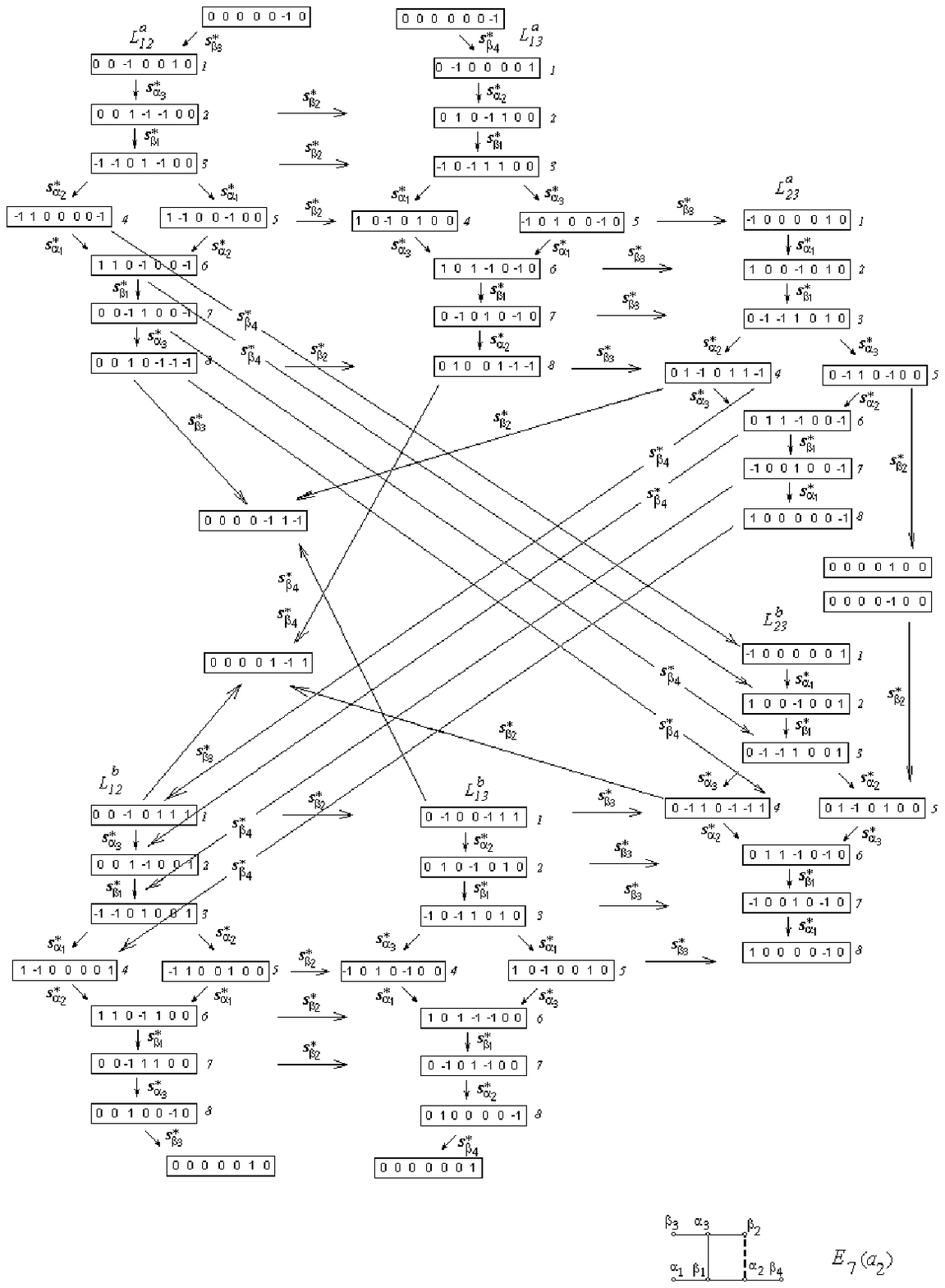}
\caption{\hspace{3mm}The linkage system $E_7(a_2)$. The only component is the $E$-component
 containing $56$ linkage diagrams, $6$ loctets}
%%%%%% The label must come after caption
\label{E7a2_linkages}
\end{figure}

\begin{figure}[H]
\centering
\includegraphics[scale=1.4]{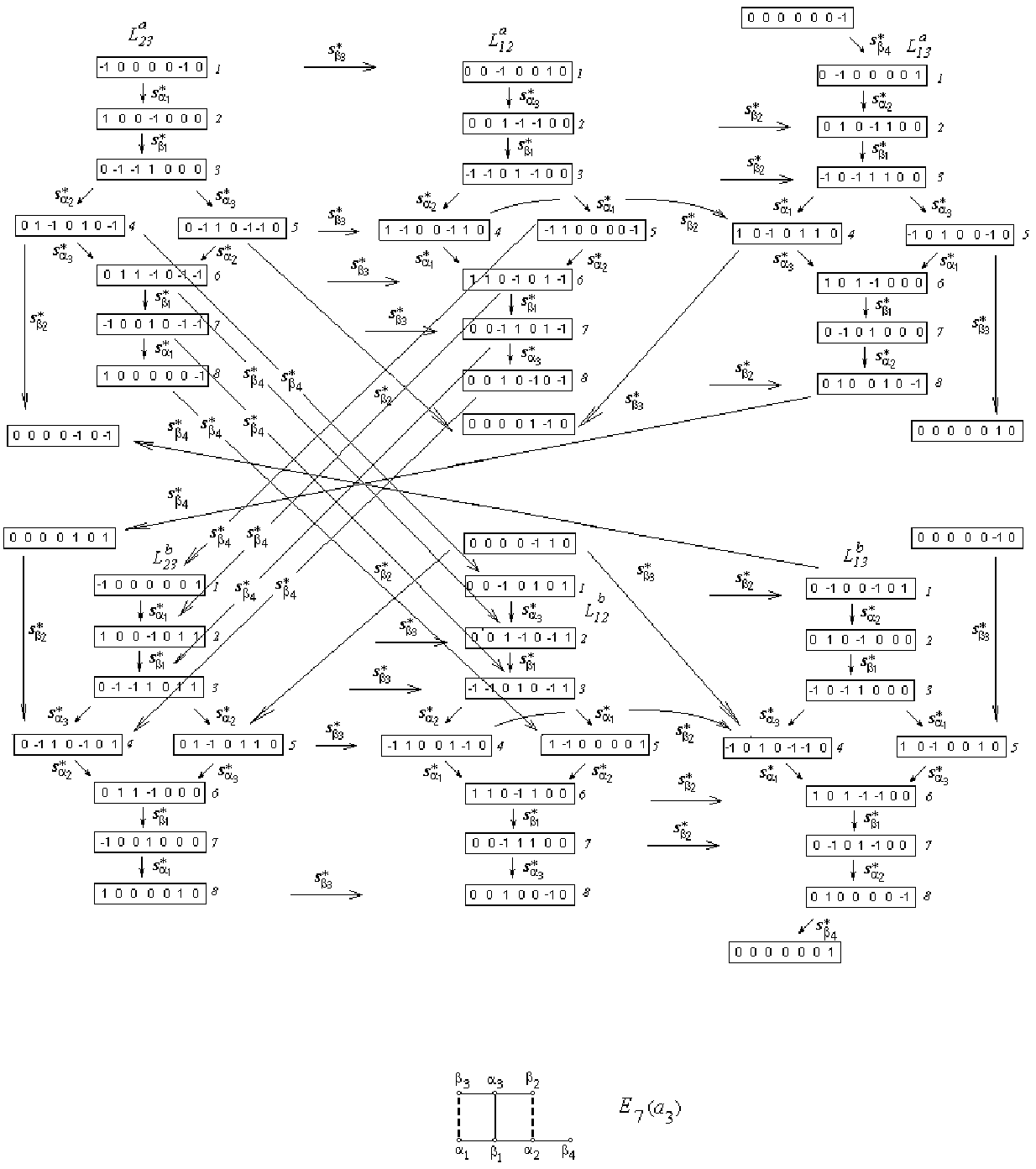}
\caption{\hspace{3mm}The linkage system $E_7(a_3)$, one component, $56$ linkage diagrams, $6$ loctets}
%%%%%% The label must come after caption
\label{E7a3_linkages}
\end{figure}

\begin{figure}[H]
\centering
\includegraphics[scale=1.4]{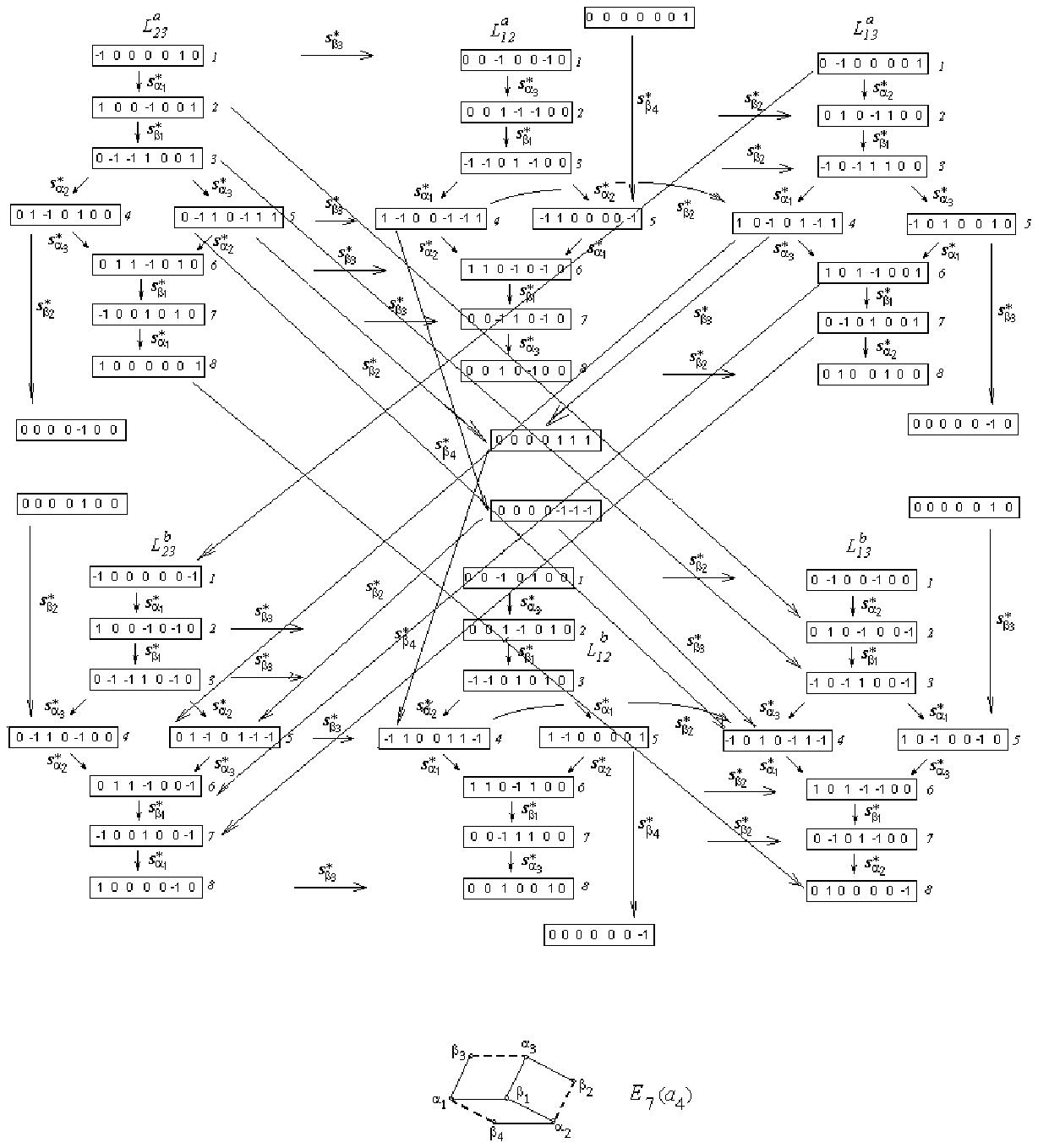}
\caption{\hspace{3mm}The linkage system $E_7(a_4)$,  one component, $56$ linkage diagrams,
 $6$ loctets}
%%%%%% The label must come after caption
\label{E7a4_linkages}
\end{figure}

\begin{figure}[H]
\centering
\includegraphics[scale=0.45]{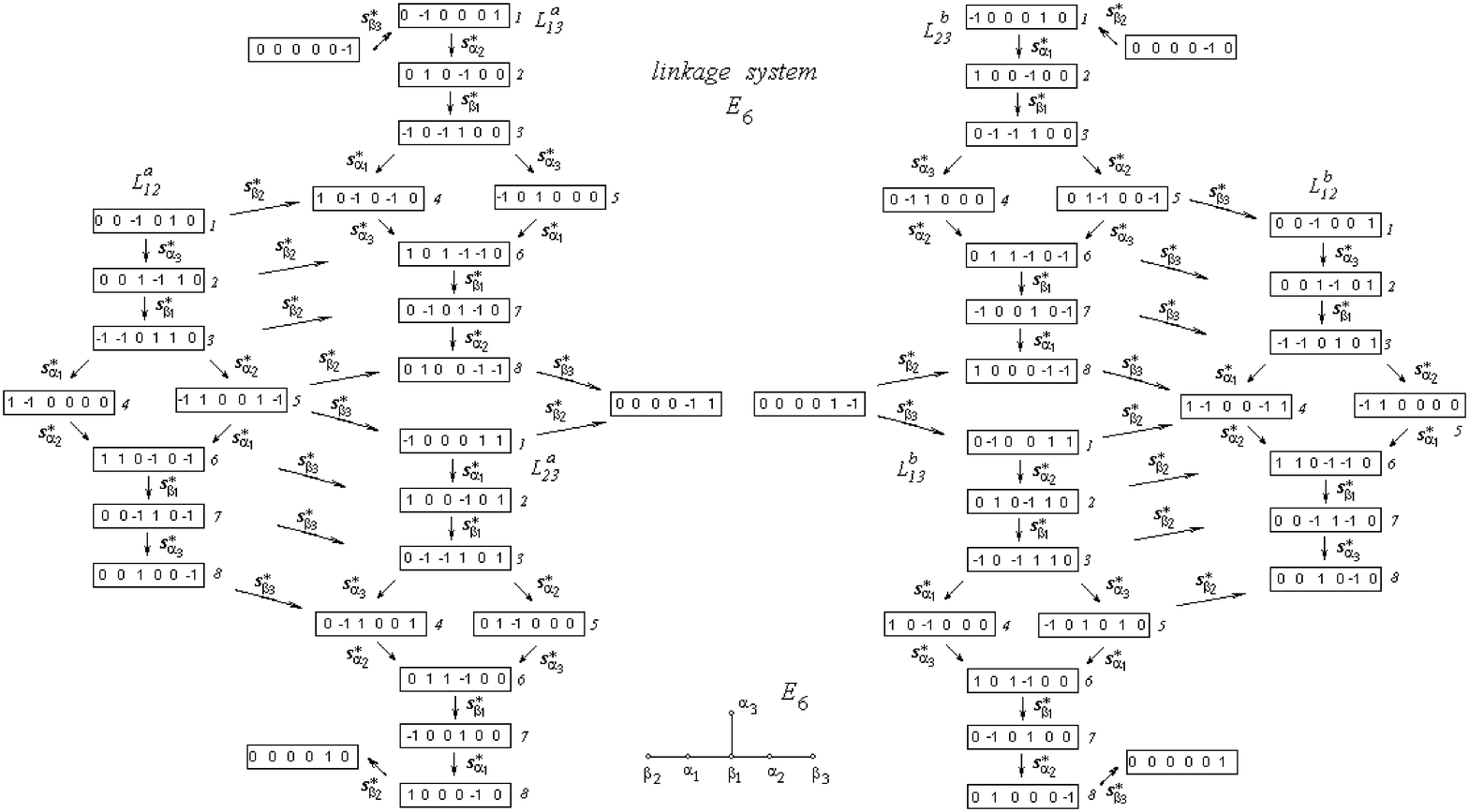}
~\\
\vspace{0.7cm}
\includegraphics[scale=0.45]{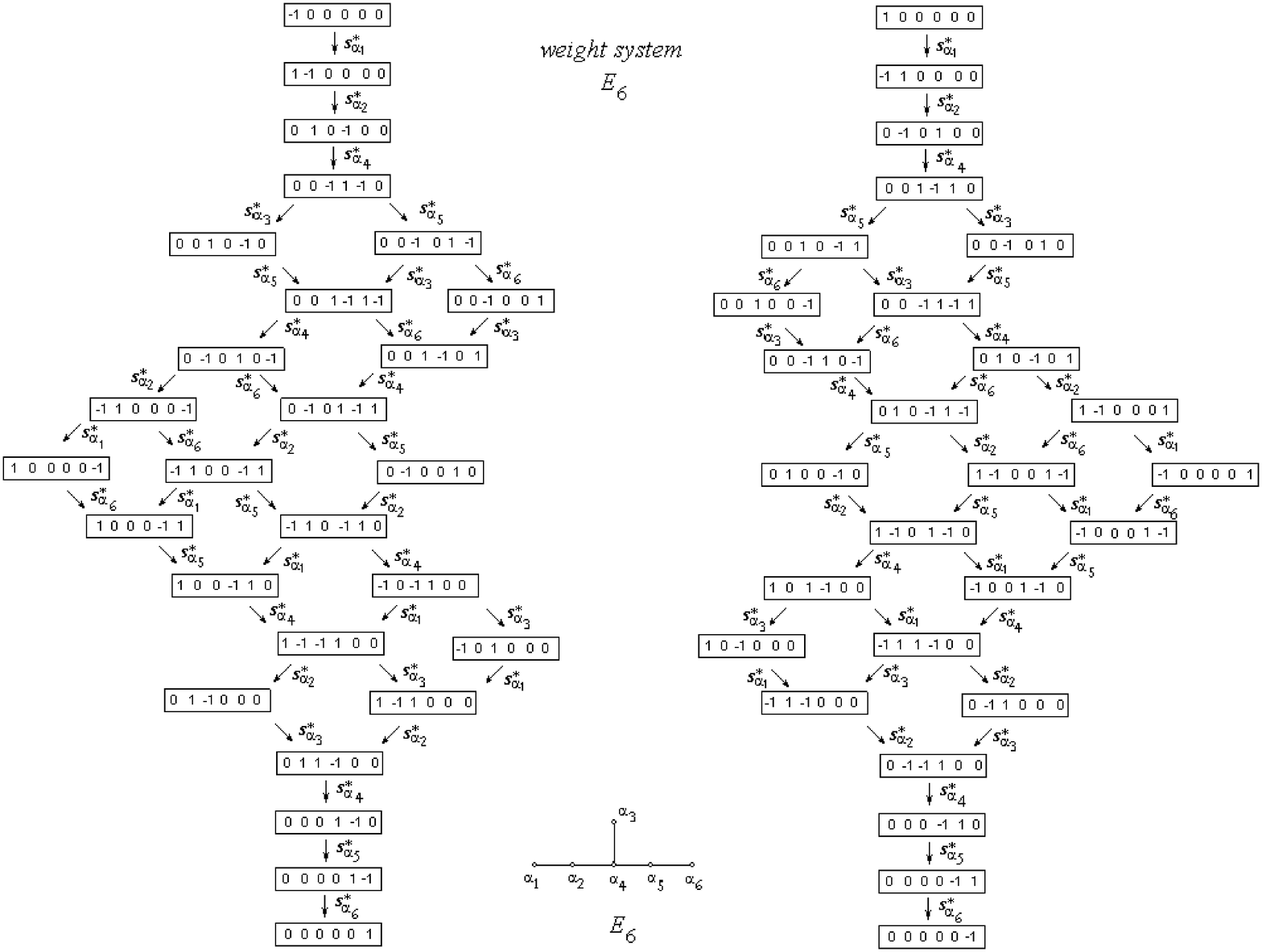}
\caption{\footnotesize \hspace{3mm}$\mathscr{L}(E_6)$: The $E$-component contains $2$ parts,
 each part contains $27$ linkages (top).
 %%The weight systems of ${\bf 27}$ and $\overline{\bf 27}$  of the semisimple Lie algebra $E_6$]
 Weight systems of representations ${\bf 27}$ and $\overline{\bf 27}$  of the semisimple Lie algebra $E_6$ (bottom)}
%%%%%% The label must come after caption
\label{27_weight_diagr_E6__2comp}
\end{figure}

\begin{figure}[H]
\centering
\includegraphics[scale=1.2]{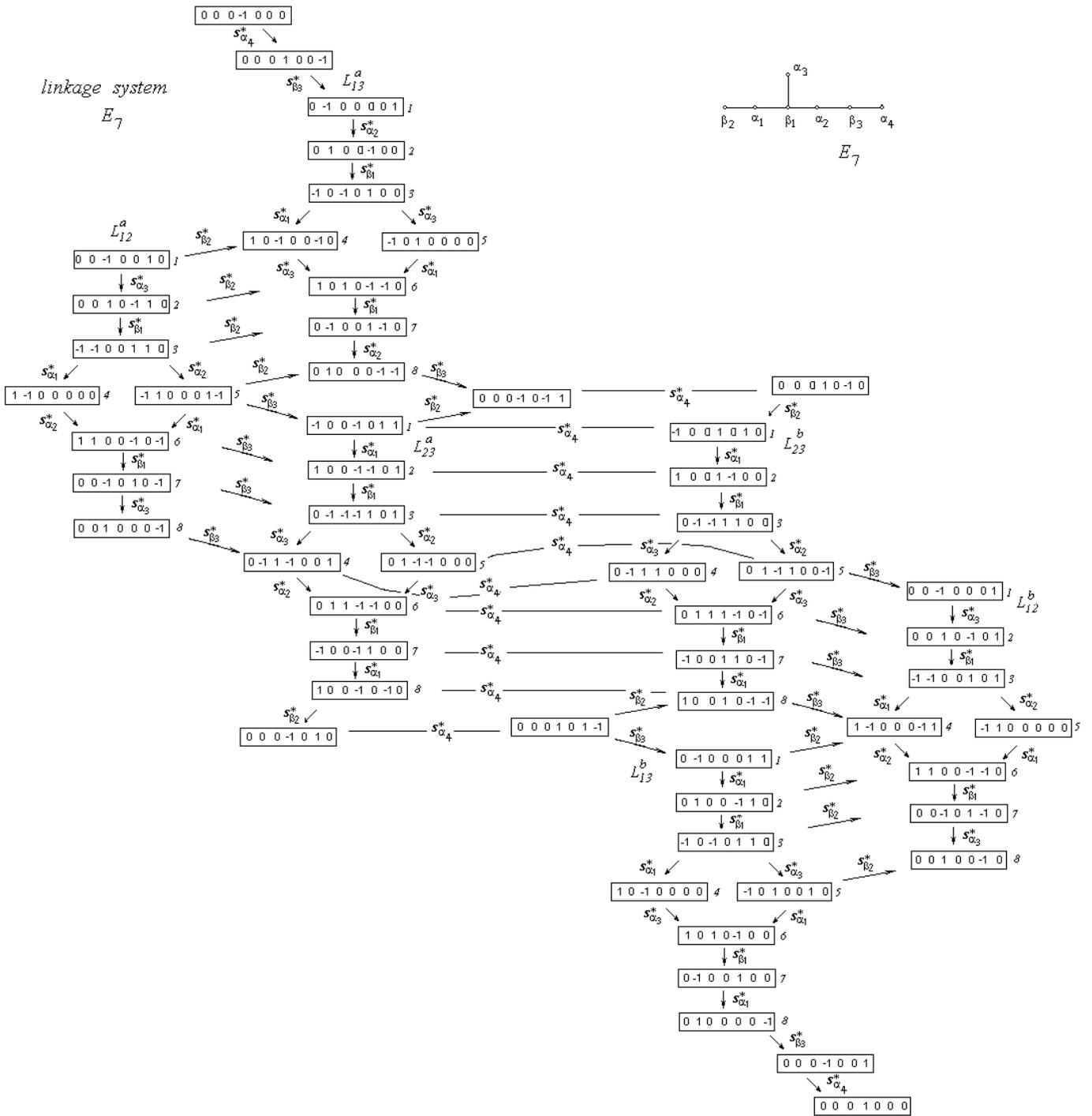}
\caption{\hspace{3mm}
The linkage system $E_7$, one component, $56$ elements  }
%%%%%% The label must come after caption
\label{E7pure_linkage_system}
\end{figure}

 \index{weight system ! - ${\bf 56}$  of  $E_7$}
\begin{figure}[H]
\centering
\includegraphics[scale=1.2]{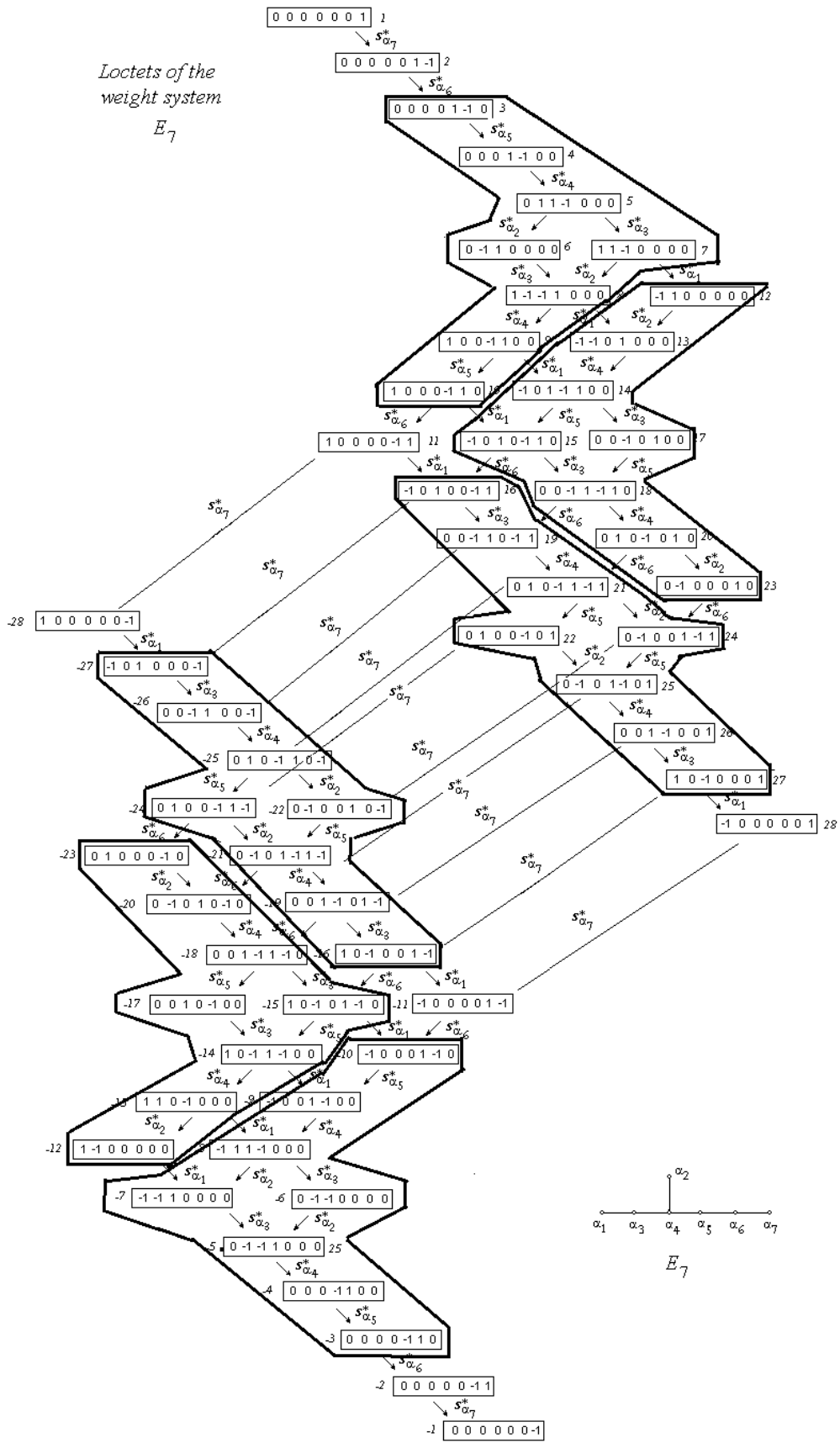}
\vspace{3mm}
\caption{\hspace{3mm}Loctets in the weight system of the fundamental representations ${\bf 56}$  of  $E_7$}
%%%%%% The label must come after caption
\label{56_weight_diagr_E7}
\end{figure}

\begin{figure}[H]
\centering
\includegraphics[scale=0.48]{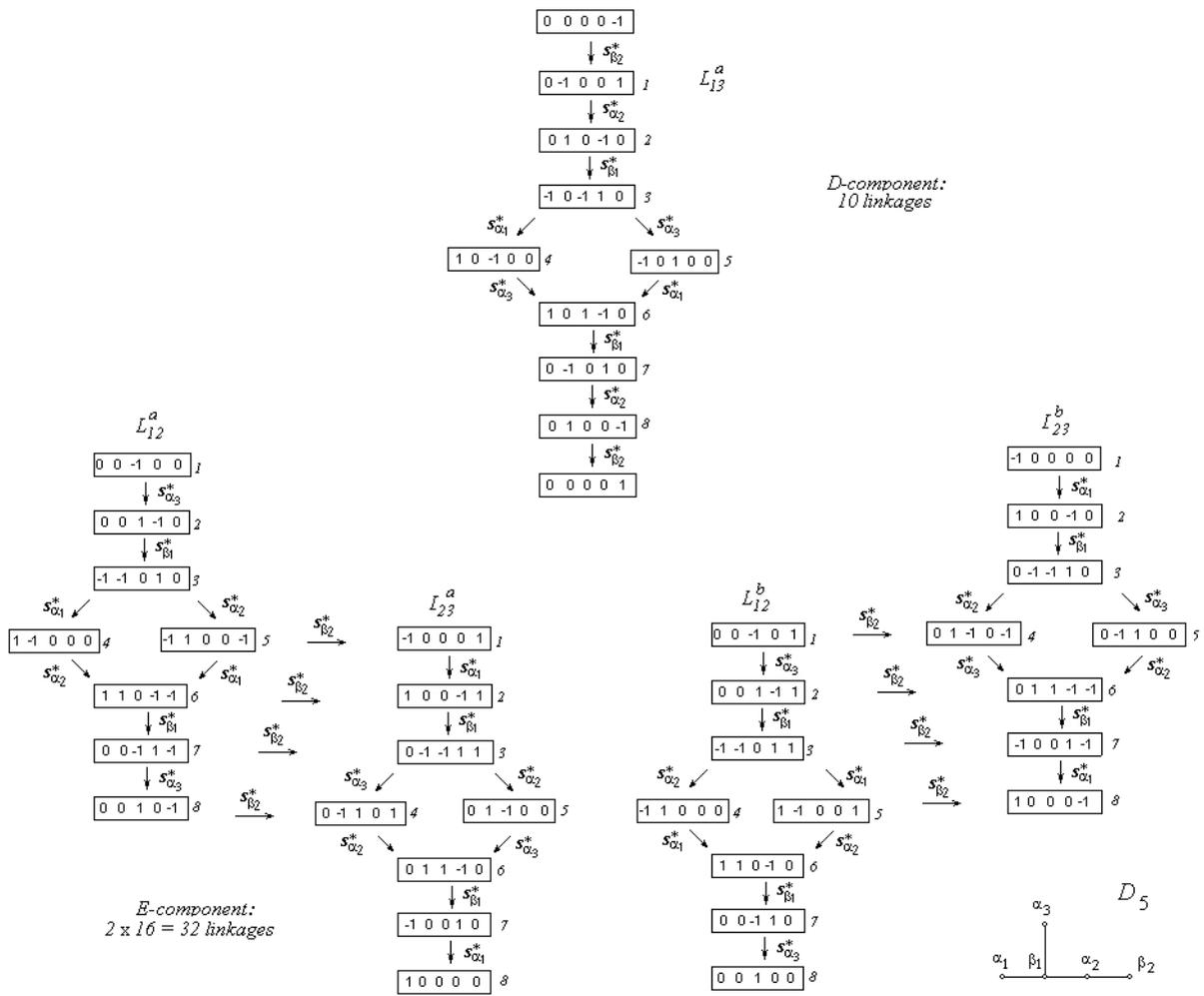}
\vspace{3mm}
 \caption[\hspace{3mm}Linkage system for $\mathscr{L}(D_5)$, 3 components, 10 + 2x16 linkages]
 {\hspace{3mm}The linkage system $\mathscr{L}(D_5)$: One part of the $D$-component
 containing $10$ linkage diagrams and two $2\times16 = 32$ element parts in the $E$-component}
%%%%%% The label must come after caption
\label{D5pure_loctets}
\end{figure}

\begin{figure}[H]
\centering
\includegraphics[scale=0.37]{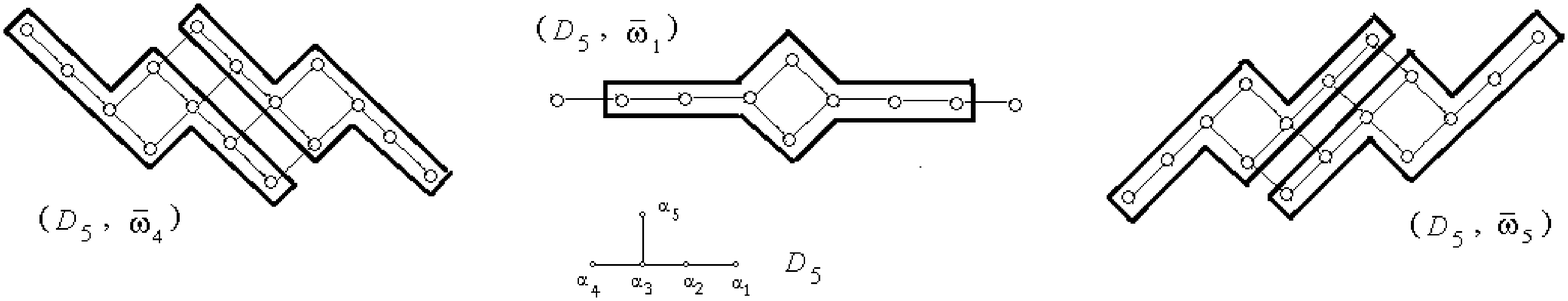}
\vspace{3mm}
\caption[\hspace{3mm}Loctets in the weight system for representations $(D_5, \overline{\omega}_1)$,
 $(D_5, \overline{\omega}_4)$ and $(D_5, \overline{\omega}_5)$]
 {\hspace{3mm}Loctets in the weight system for $3$ fundamental representations of $D_5$:
  $(D_5, \overline{\omega}_1$), $(D_5, \overline{\omega}_4$) and $(D_5, \overline{\omega}_5$)  }
%%%%%% The label must come after caption
\label{16x2_10_weight_diagr_D5}
\end{figure}

\begin{figure}[H]
\centering
\includegraphics[scale=0.36]{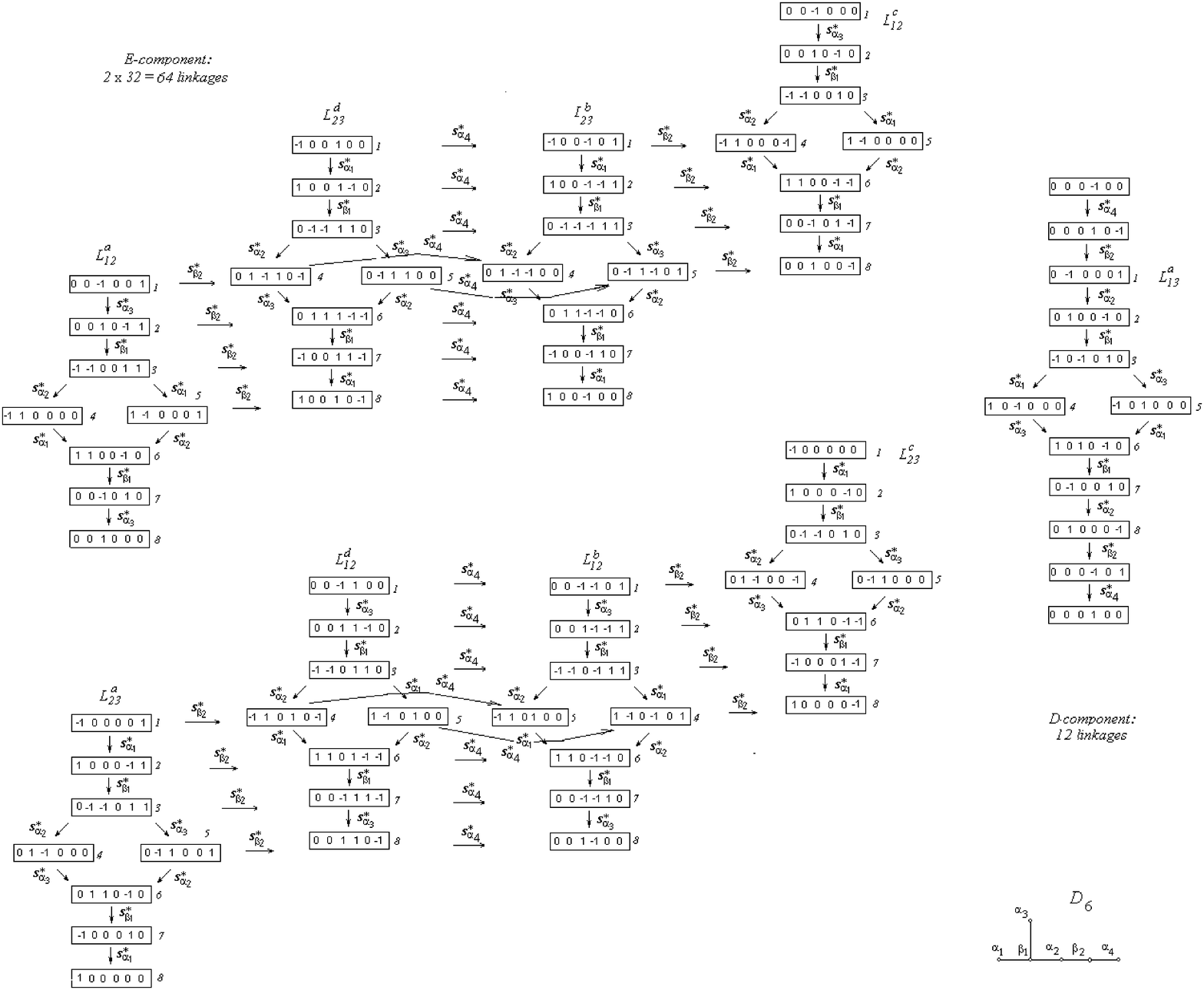}
\vspace{3mm}
 \caption[\hspace{3mm}Linkage system for $\mathscr{L}(D_6)$, 3 components, 12 + 2x32 linkages]
 {\hspace{3mm}The linkage system $\mathscr{L}(D_6)$,
 one part of the $D$-component containing $12$ elements and two parts of the $E$-component,
 each of which contains $2\times32 = 64$ elements}
%%%%%% The label must come after caption
\label{D6pure_loctets}
\end{figure}

\begin{figure}[H]
\centering
\includegraphics[scale=0.85]{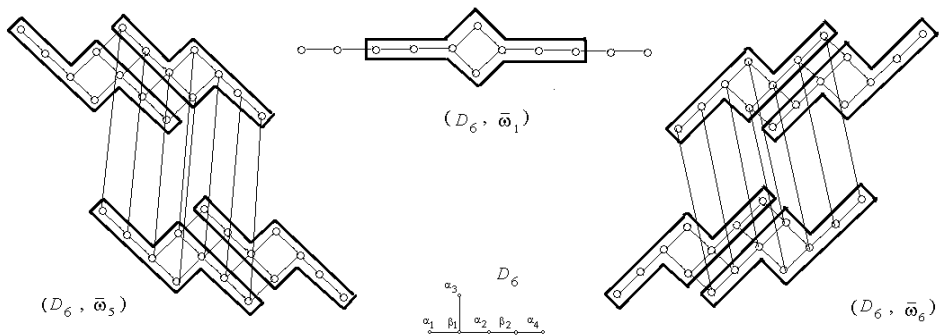}
\vspace{3mm}
\caption[\hspace{3mm}Loctets in the weight system for representations $(D_6, \overline{\omega}_i)$,
$(D_6, \overline{\omega}_5$) and $(D_6, \overline{\omega}_6$)]{\hspace{3mm}Loctets in the weight system for $3$ fundamental representations of $D_6$:
  $(D_6, \overline{\omega}_1$), $(D_6, \overline{\omega}_5$) and $(D_6, \overline{\omega}_6$)  }
%%%%%% The label must come after caption
\label{32x2_12_weight_diagr_D6}
\end{figure}

%%=========================
\newpage
\subsection{The linkage systems $\mathscr{L}(D_7(a_1))$, $\mathscr{L}(D_7(a_2))$, $\mathscr{L}(D_7)$}
 \index{linkage system  ! - $\mathscr{L}(D_7)$ and $\mathscr{L}(D_7(a_k))$}
~\\

\begin{figure}[H]
\centering
\includegraphics[scale=0.42,  angle=270]{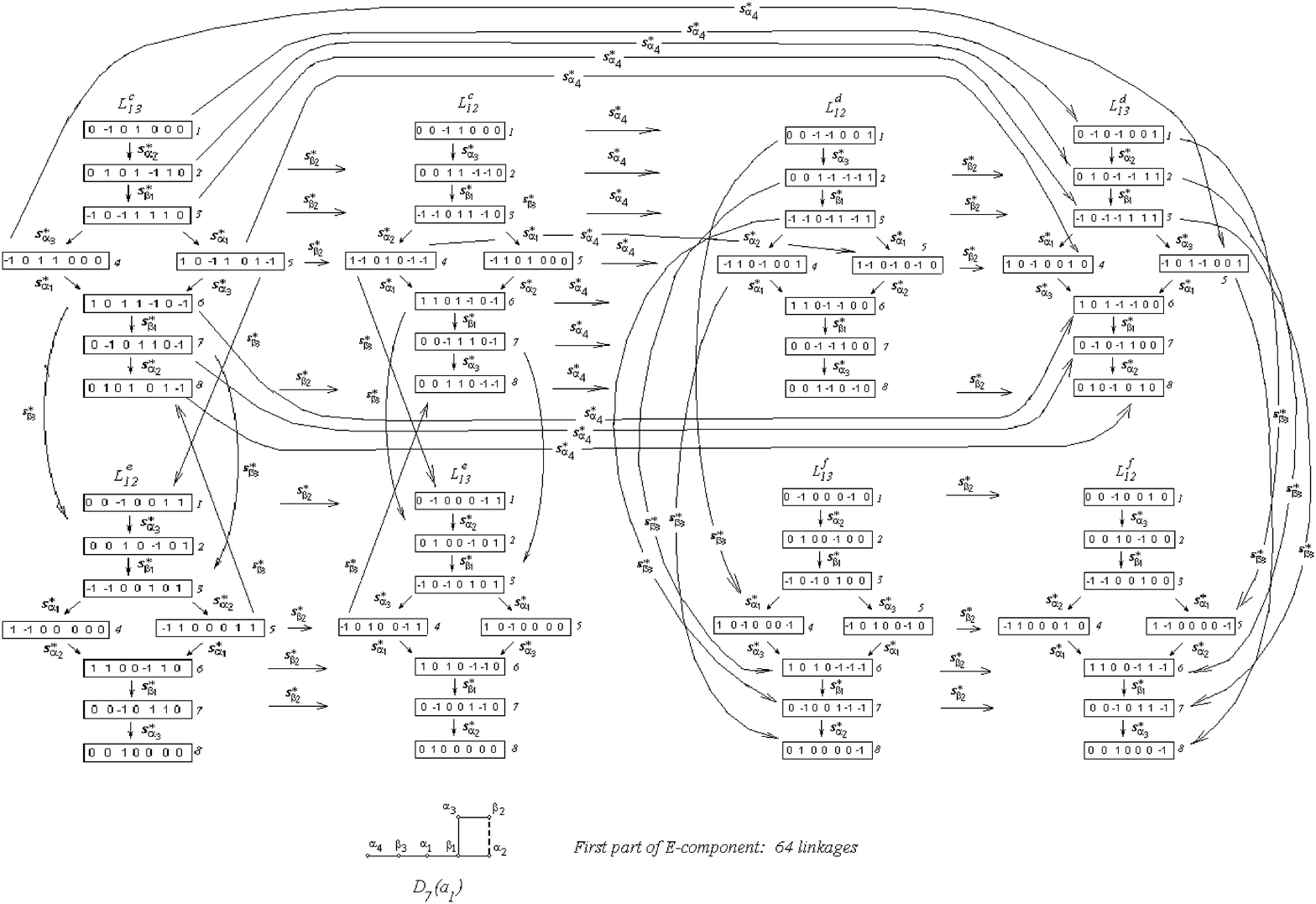}
 \caption[\hspace{3mm}The linkage system $\mathscr{L}(D_7(a_1))$, first part of the $E$-component]
 {\hspace{3mm}The linkage system $\mathscr{L}(D_7(a_1))$, $1$st part of the $E$-component, $64$ linkage diagrams, $8$ loctets}
%%%%%% The label must come after caption
\label{D7a1_linkages_cdef}
\end{figure}

\begin{figure}[H]
\centering
\includegraphics[scale=0.47,  angle=270]{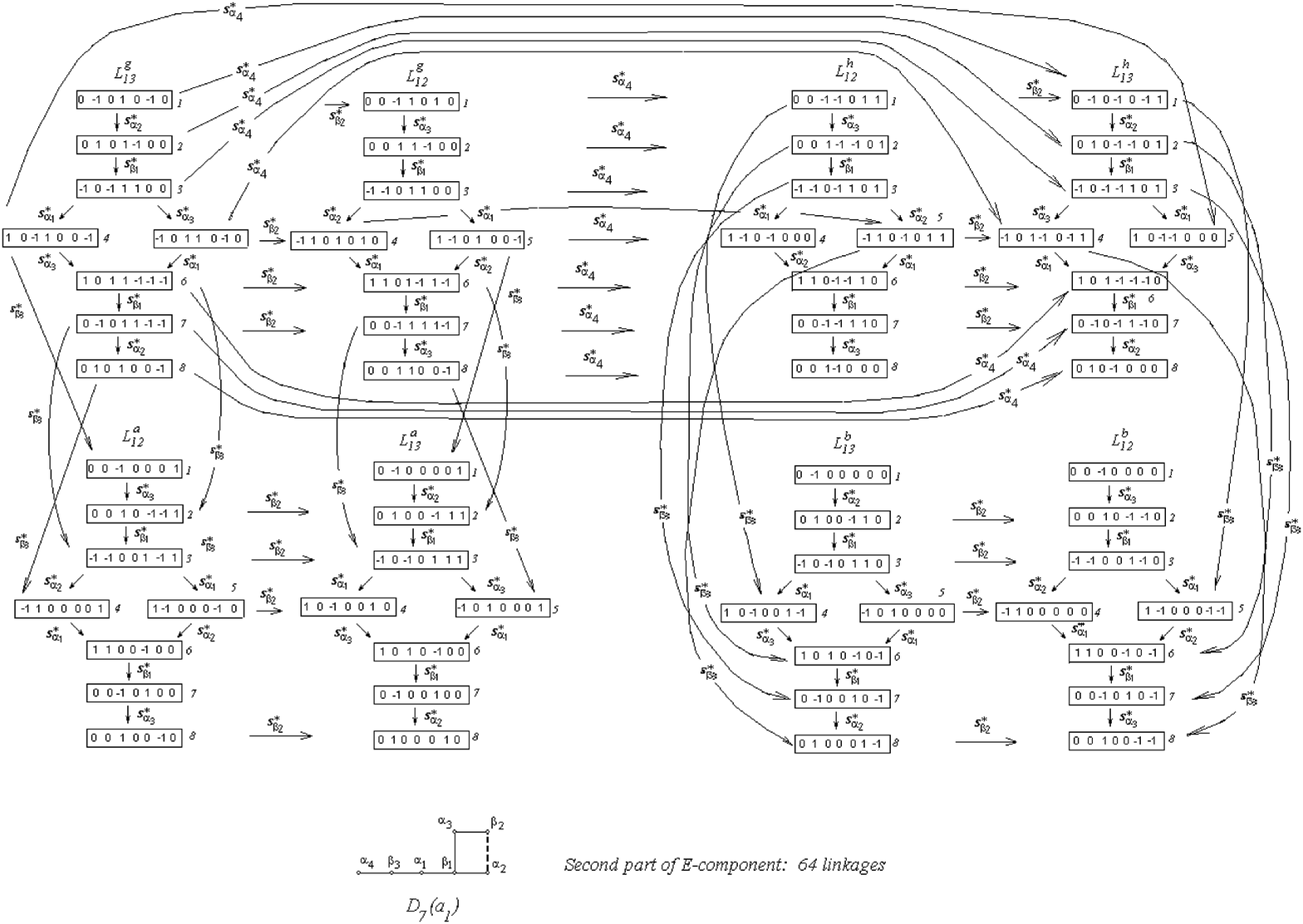}
\caption{\hspace{3mm}$D_7(a_1)$, second part of the $E$-component:
 $64$ linkages, $8$ loctets}
%%%%%% The label must come after caption
\label{D7a1_linkages_abgh}
\end{figure}

\begin{figure}[H]
\centering
\includegraphics[scale=0.47, angle=270]{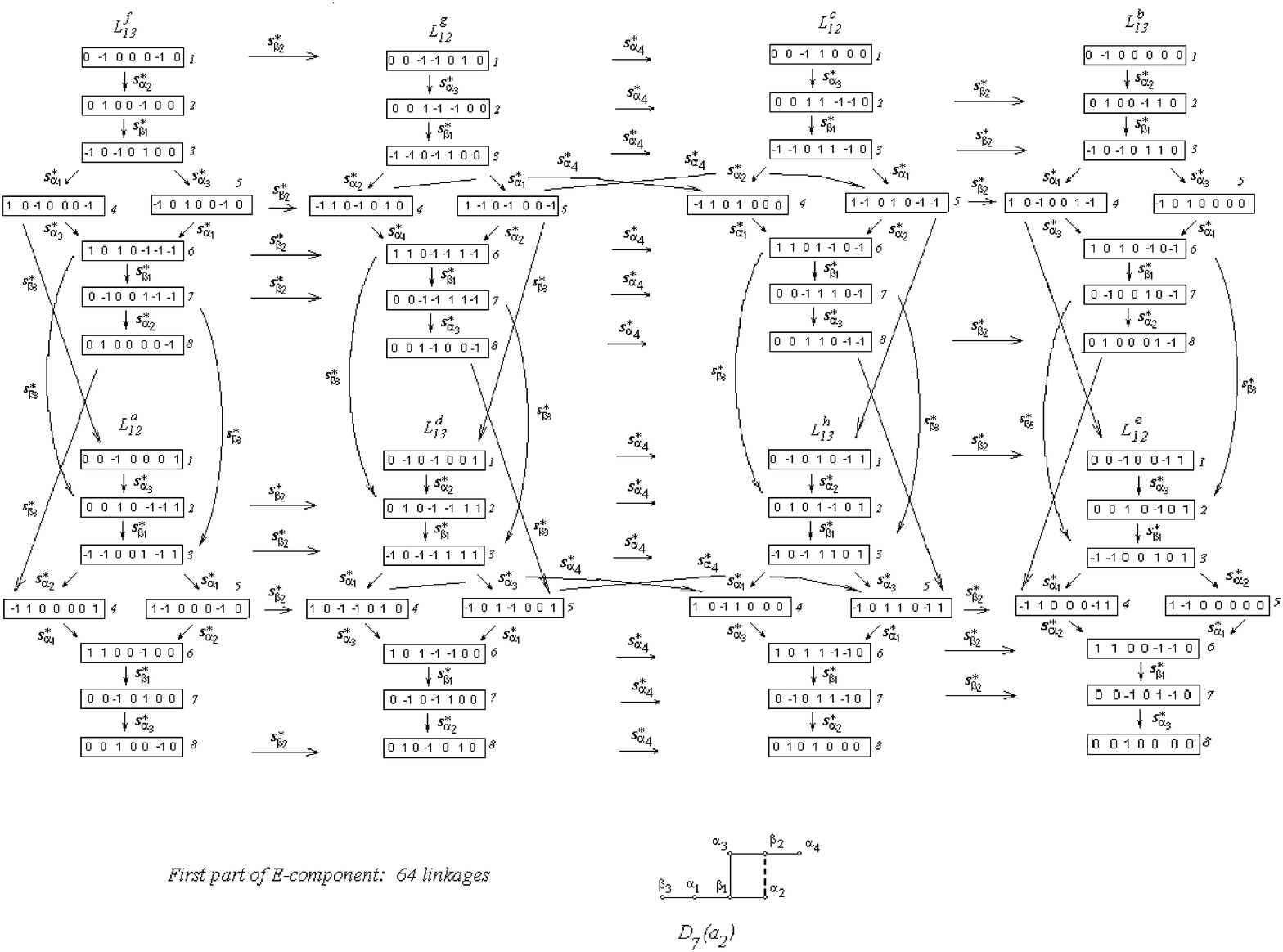}
\caption{\hspace{3mm}The linkage system $D_7(a_2)$, $1$st part of the $E$-component , $64$ linkage diagrams, $8$ loctets}
%%%%%% The label must come after caption
\label{D7a2_linkages_comp1}
\end{figure}

\begin{figure}[H]
\centering
\includegraphics[scale=0.55, angle=270]{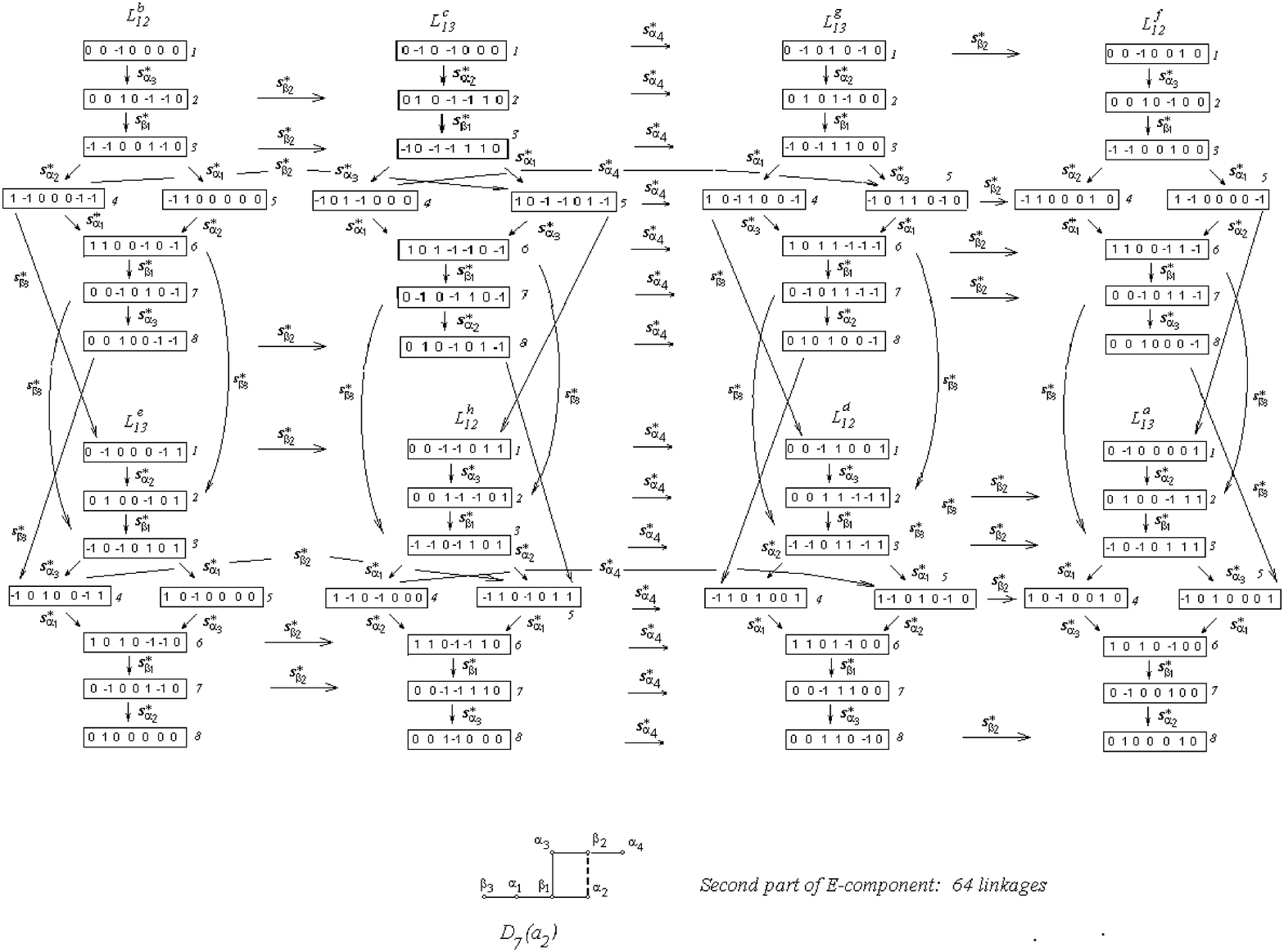}
\caption{\hspace{3mm} The linkage system $D_7(a_2)$, $2$nd part of the $E$-component, $64$ linkage diagrams, $8$ loctets}
%%%%%% The label must come after caption
\label{D7a2_linkages_comp2}
\end{figure}

\begin{figure}[H]
\centering
\includegraphics[scale=0.45, angle=270]{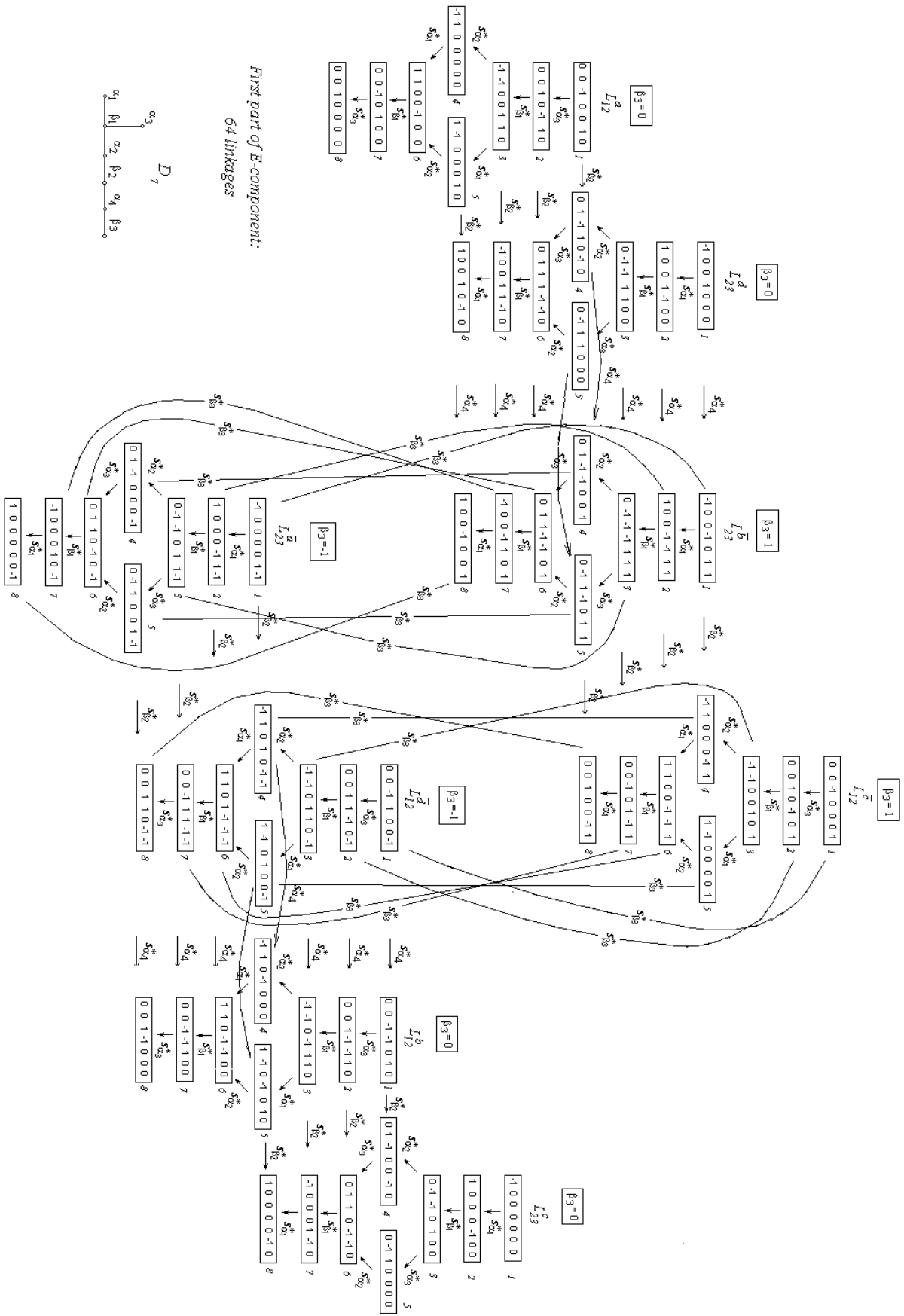}
\vspace{3mm}
\caption{\hspace{3mm}The first component of the linkage system $D_7$, 64 linkages}
%%%%%% The label must come after caption
\label{D7pu_loctets_comp1}
\end{figure}

\begin{figure}[H]
\centering
\includegraphics[scale=0.47, angle=270]{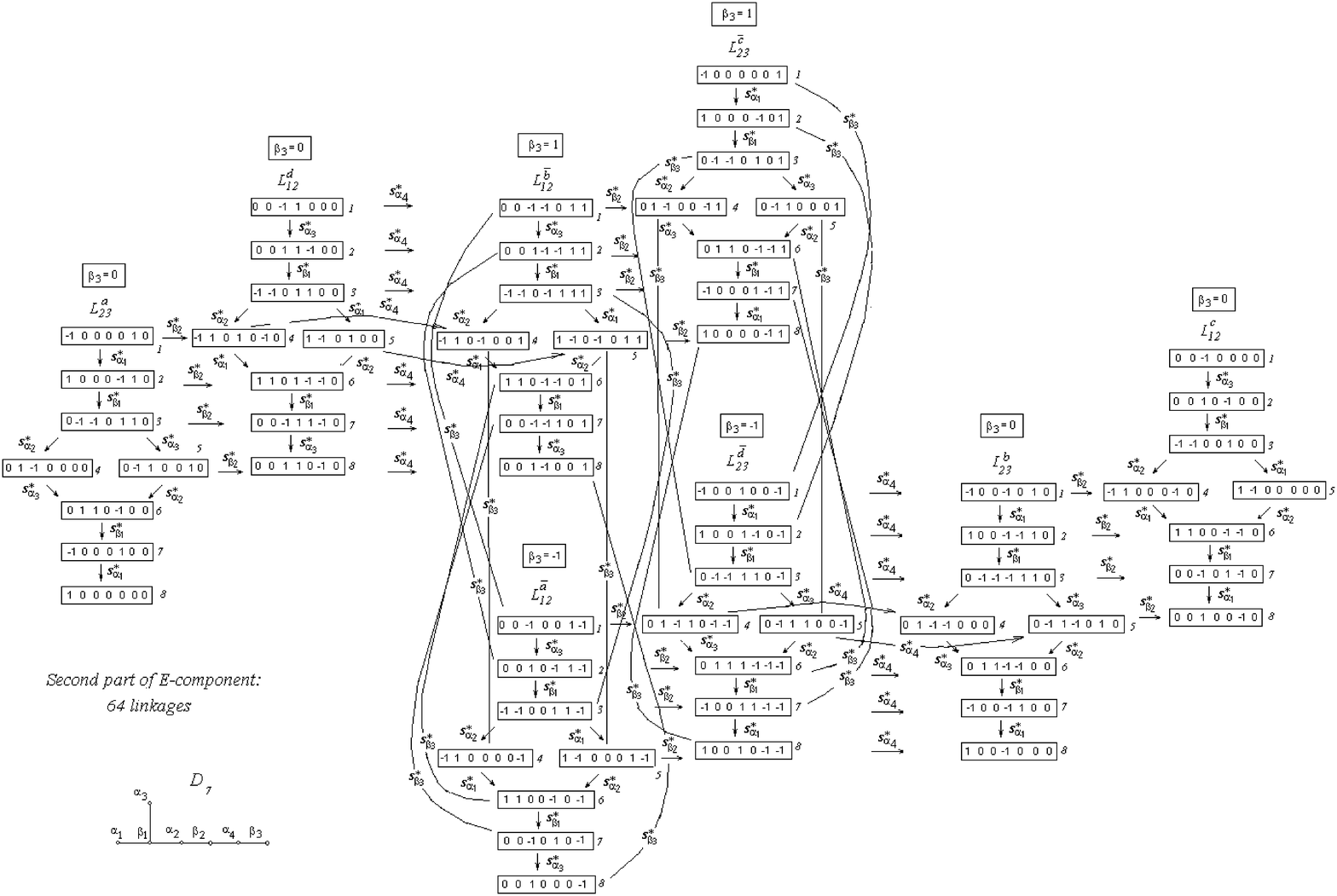}
\vspace{3mm}
\caption{\hspace{3mm}The second component of the linkage system $D_7$, 64 linkages}
%%%%%% The label must come after caption
\label{D7pu_loctets_comp2}
\end{figure}

\begin{figure}[H]
\centering
\includegraphics[scale=0.47, angle=270]{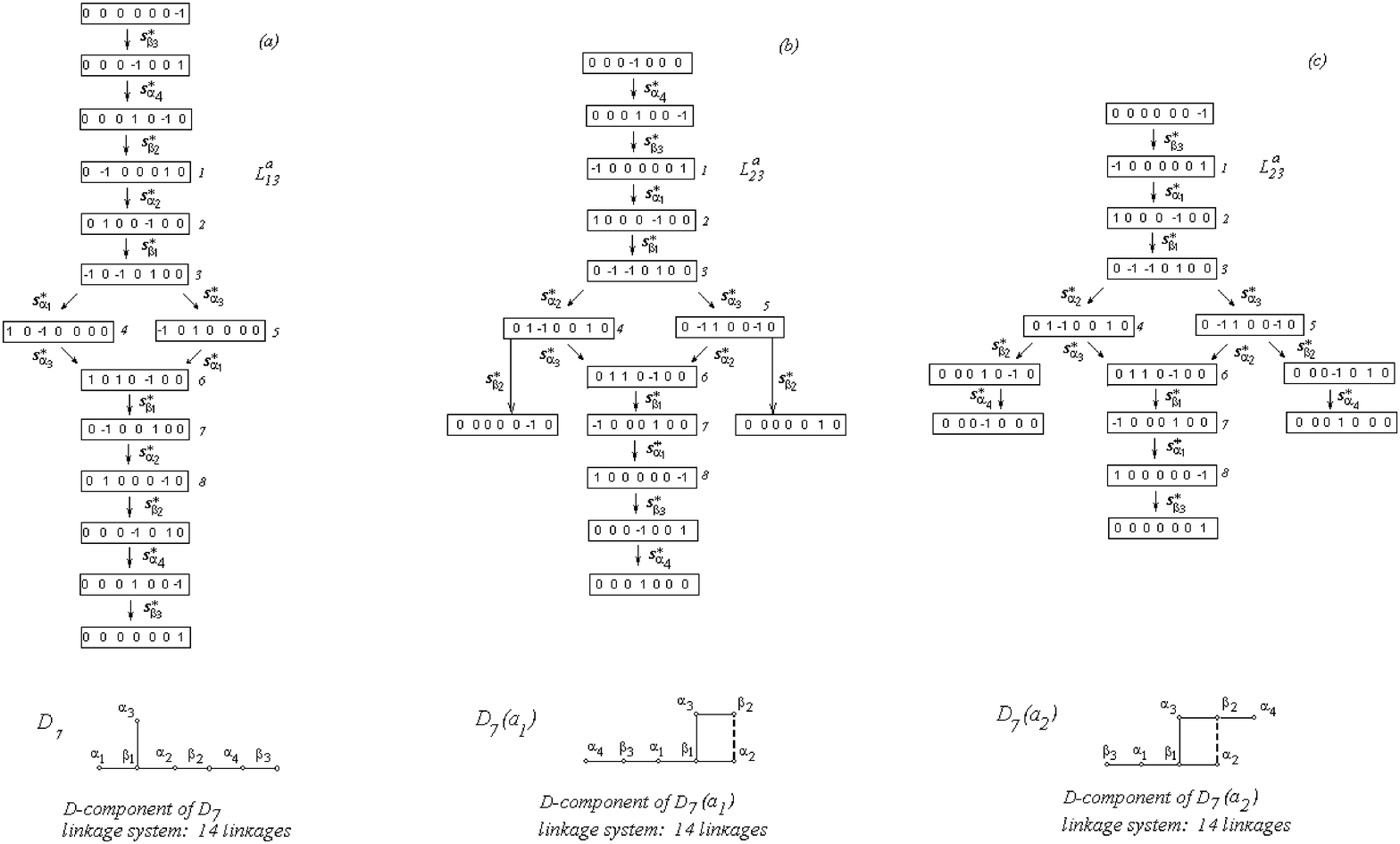}
\vspace{3mm}
 \caption[\hspace{3mm}The $D$-components for linkage systems
 $\mathscr{L}(D_7)$, $\mathscr{L}(D_7(a_1))$, $\mathscr{L}(D_7(a_2))$]
 {The $D$-components for linkage systems $\mathscr{L}(D_7)$, $\mathscr{L}(D_7(a_1))$, $\mathscr{L}(D_7(a_2))$}
%%%%%% The label must come after caption
\label{D7a1_D7a2_D7pu_loctets_comp3}
\end{figure}

%%=========================
\clearpage
\subsection{The linkage systems $\mathscr{L}(D_4)$}
~\\
\begin{figure}[H]
\centering
\includegraphics[scale=1.8]{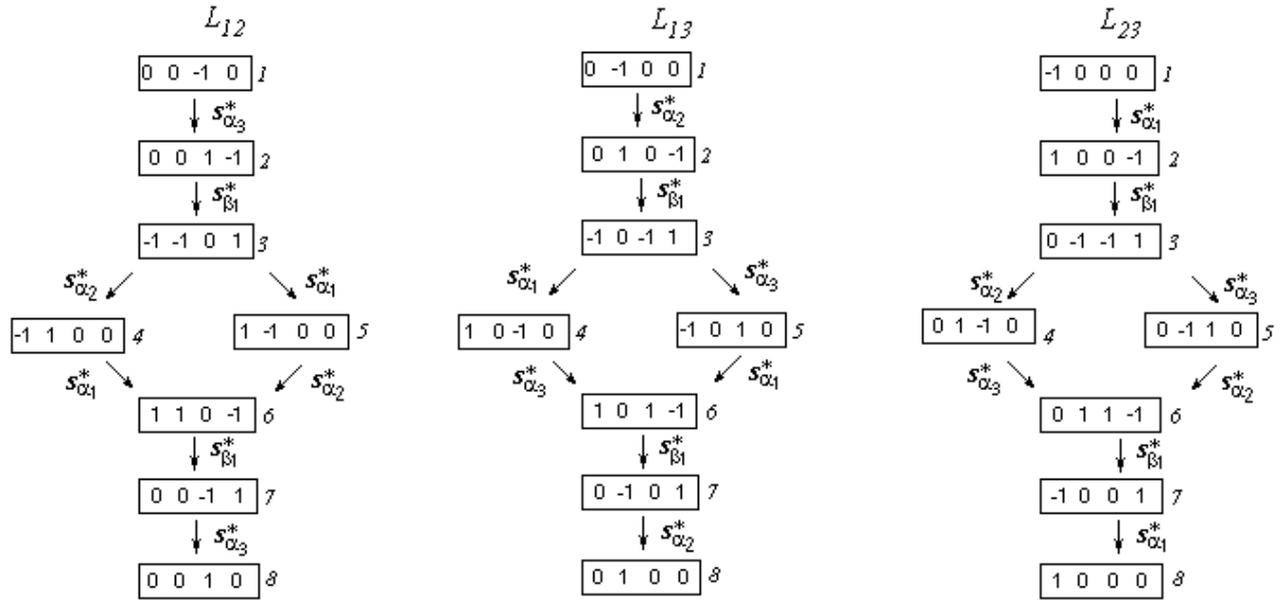}
\caption{\hspace{3mm}The linkage system $D_4$, $24$ linkages, $3$ loctets}
%%%%%% The label must come after caption
\label{D4_loctets}
\end{figure}

%% file: biblio-2.tex
%% The Carter theorem - Biblio
\pagebreak[4]